\documentclass[preprint,11pt]{elsarticle} 

\usepackage[bookmarks=true,colorlinks=true,linkcolor=blue]{hyperref}
\usepackage[usenames,dvipsnames,svgnames,table]{xcolor}
\usepackage{amsthm}
\usepackage[ruled]{algorithm2e}
\biboptions{sort&compress}

\usepackage{amssymb,amsmath}
\usepackage{calc}
\usepackage[margin=1.in]{geometry}
\usepackage{tikz}

\pdfoutput=1

\renewcommand{\url}[1]{}
\newcommand{\citeCount}[1]{}

\usepackage{soul}

\usepackage{ifthen}
\newcommand{\REV}[2]{%
  \ifthenelse{\equal{#1}{0}}{{{\color{brown}#2}}}{}%
  \ifthenelse{\equal{#1}{1}}{{{\color{green!50!black}#2}}}{}%
  \ifthenelse{\equal{#1}{2}}{{{\color{red!50!black}#2}}}{}%
  \ifthenelse{\equal{#1}{3}}{{{\color{blue!50!black}#2}}}{}%
}

\newcommand{\dx}{\Delta x}
\newcommand{\dy}{\Delta y}
\newcommand{\dt}{\Delta t}

\newcommand{\Dpx}{D_{+x}}
\newcommand{\Dmx}{D_{-x}}

\newcommand{\pd}[2]{\frac{\partial #1}{\partial #2}}
\newcommand{\pdn}[3]{\frac{\partial^#3 #1}{\partial #2^#3}}
\newcommand{\dd}[2]{\frac{d #1}{d #2}}
\newcommand{\ddn}[3]{\frac{d^#3 #1}{d #2^#3}}

\newcommand{\gv}{\mathbf{ g}}

\newcommand{\iv}{\mathbf{ i}}

\newcommand{\nv}{\mathbf{ n}}

\newcommand{\rv}{\mathbf{ r}}

\newcommand{\tv}{\mathbf{ t}}
\newcommand{\uv}{\mathbf{ u}}

\newcommand{\xv}{\mathbf{ x}}

\newcommand{\Uv}{\mathbf{ U}}
\newcommand{\Vv}{\mathbf{ V}}
\newcommand{\Wv}{\mathbf{ W}}

\newcommand{\zerov}{\mathbf{0}}

\newcommand{\Dc}{{\mathcal D}}

\newcommand{\Gc}{{\mathcal G}}

\newcommand{\Ic}{{\mathcal I}}

\newcommand{\Lc}{{\mathcal L}}
\newcommand{\Mc}{{\mathcal M}}

\newcommand{\Oc}{{\mathcal O}}

\newcommand{\tableFont}{\footnotesize}

\newcommand{\num}[2]{#1e#2} 

\newcommand{\trimfigcb}[3]{\includegraphics[width=#2, height=#3, clip, trim=17cm 2.35cm 1.65cm 2.35cm]{#1}}

\newcommand{\xcb}{0}
\newcommand{\ycb}{0}
\newcommand{\xccb}{0}
\newcommand{\xrcb}{0}
\newcommand\drawColorBarH[7]{ 
\renewcommand{\xccb}{\xcb+0.5*#4}
\renewcommand{\xrcb}{\xcb+#4}
\begin{scope}[#2] 
  \draw (\xcb,\ycb)   node[anchor=south,xshift=0cm,yshift=0cm,rotate=-90] {\trimfigcb{#1}{#3 cm}{#4 cm}};
  \draw (\xcb,\ycb)   node[anchor=north, xshift=0.15cm,yshift=0cm] {\scriptsize $#5$};
  \draw (\xccb,\ycb) node[anchor=north, xshift=0cm,yshift=0cm] {\scriptsize $#7$};
  \draw (\xrcb,\ycb)  node[anchor=north, xshift=0cm,yshift=0cm] {\scriptsize $#6$};
\end{scope}
}

\newlength{\tfwidth}
\newlength{\tfheight}
\newlength{\tfxa}
\newlength{\tfxb}
\newlength{\tfya}
\newlength{\tfyb}
\newcommand{\trimFigWithBox}[6]{%
\setlength\fboxsep{0pt}%
\setlength\fboxrule{1.0pt}
\fbox{\includegraphics[width=#2, clip, trim=#3 #4 #5 #6]{#1}}%
}
\newcommand{\trimFigNoBox}[6]{%
\setlength\fboxsep{1pt}
\setlength\fboxrule{0.0pt}
\fbox{\includegraphics[width=#2, clip, trim=#3 #4 #5 #6]{#1}}%
}
\newcommand{\trimFigHeightWithBox}[6]{%
\setlength\fboxsep{0pt}%
\setlength\fboxrule{1.0pt}
\fbox{\includegraphics[height=#2, clip, trim=#3 #4 #5 #6]{#1}}%
}
\newcommand{\trimFigHeightNoBox}[6]{%
\setlength\fboxsep{1pt}
\setlength\fboxrule{0.0pt}
\fbox{\includegraphics[height=#2, clip, trim=#3 #4 #5 #6]{#1}}%
}

\newsavebox\figBox

\newcommand{\trimw}[6]{%
\sbox\figBox{\includegraphics{#1}}
\setlength{\tfwidth}{\the\wd\figBox}
\setlength{\tfheight}{\the\ht\figBox}
\setlength{\tfxa}{\tfwidth*\real{#3}}%
\setlength{\tfxb}{\tfwidth*\real{#4}}%
\setlength{\tfya}{\tfheight*\real{#5}}%
\setlength{\tfyb}{\tfheight*\real{#6}}%
\trimFigNoBox{#1}{#2}{\tfxa}{\tfya}{\tfxb}{\tfyb}%
}

\newcommand{\trimwb}[6]{%

\sbox\figBox{\includegraphics{#1}}
\setlength{\tfwidth}{\the\wd\figBox}
\setlength{\tfheight}{\the\ht\figBox}
\setlength{\tfxa}{\tfwidth*\real{#3}}%
\setlength{\tfxb}{\tfwidth*\real{#4}}%
\setlength{\tfya}{\tfheight*\real{#5}}%
\setlength{\tfyb}{\tfheight*\real{#6}}%
\trimFigWithBox{#1}{#2}{\tfxa}{\tfya}{\tfxb}{\tfyb}%
}

\newcommand{\trimh}[6]{%
\sbox\figBox{\includegraphics{#1}}
\setlength{\tfwidth}{\the\wd\figBox}
\setlength{\tfheight}{\the\ht\figBox}
\setlength{\tfxa}{\tfwidth*\real{#3}}%
\setlength{\tfxb}{\tfwidth*\real{#4}}%
\setlength{\tfya}{\tfheight*\real{#5}}%
\setlength{\tfyb}{\tfheight*\real{#6}}%
\trimFigHeightNoBox{#1}{#2}{\tfxa}{\tfya}{\tfxb}{\tfyb}%
}

\newcommand{\trimhb}[6]{%

\sbox\figBox{\includegraphics{#1}}
\setlength{\tfwidth}{\the\wd\figBox}
\setlength{\tfheight}{\the\ht\figBox}
\setlength{\tfxa}{\tfwidth*\real{#3}}%
\setlength{\tfxb}{\tfwidth*\real{#4}}%
\setlength{\tfya}{\tfheight*\real{#5}}%
\setlength{\tfyb}{\tfheight*\real{#6}}%
\trimFigHeightWithBox{#1}{#2}{\tfxa}{\tfya}{\tfxb}{\tfyb}%
}
\usepackage{mathrsfs}
\usepackage{multirow}

\usepackage{caption}
\usepackage{subcaption}

\newcommand{\KL}{Kirchhoff-Love }

\newcommand{\csf}{C_\text{sf} } 
\newcommand{\df}{C_\text{df} } 

\newcommand{\ad}{\nu_\text{ad} } 

\newcommand{\Dad}{\Dc_\text{ad} } 

\newcommand{\Dpr}{D_{+r}}
\newcommand{\Dmr}{D_{-r}}
\newcommand{\Dps}{D_{+s}}
\newcommand{\Dms}{D_{-s}}
\newcommand{\Dor}{D_{0r}}
\newcommand{\Dos}{D_{0s}}

\newcommand{\dr}{h_r}
\newcommand{\ds}{h_s}

\begin{document}

\small

\begin{frontmatter}

  \title{
 Stable finite difference methods for Kirchhoff-Love plates on \\ overlapping grids
  }

\author[ul]{Longfei~Li\corref{cor1}\fnref{RCSThanks}}
\ead{longfei.li@louisiana.edu}

\author[ucla]{Hangjie~Ji}
\ead{hangjie@math.ucla.edu}

\author[lanl]{Qi Tang}
\ead{qtang@lanl.gov}

\address[ul]{Department of Mathematics, University of Louisiana at Lafayette, Lafayette, LA  70504, USA.}
\address[ucla]{Department of Mathematics, University of California Los Angeles, Los Angeles, CA 90095, USA.}
\address[lanl]{Los Alamos National Laboratory, Los Alamos, NM 87545, USA.}
\cortext[cor1]{Corresponding author.}

\fntext[RCSThanks]{Research supported by the  Louisiana Board of Regents Support Fund under contract No. LEQSF(2018-21)-RD-A-23.}

\begin{abstract}
In this work, we propose and develop  efficient  and accurate  numerical methods for solving the Kirchhoff-Love plate model in  domains with complex geometries. The algorithms proposed here employ  curvilinear finite-difference methods for spatial discretization of the governing PDEs   on general composite overlapping grids.
The coupling of different components of the composite overlapping  grid is through numerical interpolations.
However,  interpolations  introduce perturbation to the finite-difference discretization, which causes numerical instability for  time-stepping schemes used to advance the  resulted semi-discrete system.  
To address the instability, motivated by an upwind scheme for solving the wave equation, we propose to add a fourth-order hyper-dissipation  to the spatially discretized system to stabilize its  time integration; this additional dissipation term  captures  the essential upwinding effect of the original upwind scheme. 
The investigation of  strategies for incorporating the   upwind dissipation term  into several time-stepping schemes (both explicit and implicit) leads to the development of   four novel  algorithms. For each algorithm,
formulas for determining  a stable time step and a sufficient dissipation coefficient on curvilinear grids are derived  by performing  a local Fourier analysis. 
Quadratic eigenvalue problems for  a simplified model plate in 1D domain   are considered to reveal  the weak instability due to the presence of interpolating equations in the spatial discretization.
This  model problem is further   investigated for
the stabilization effects of the proposed  algorithms.  Carefully designed numerical experiments are carried out to  validate the accuracy and stability of the proposed  algorithms, followed by two benchmark problems to  demonstrate the capability and efficiency of our approach for solving realistic applications. Results that concern  the performance of the proposed algorithms are also presented. 

\end{abstract}


\end{frontmatter}

\clearpage
\tableofcontents

\clearpage
\section{Introduction}\label{sec:introduction}
\KL theory~\cite{reddy2006theory}  concerns the small deflection of  thin plates  subject to external loadings,
and it  is widely used in structural engineering for determining the stresses and deformations of thin-walled structures.
The theory simplifies the   solid mechanics    by assuming that   a 3D plate  can be  represented by its 2D mid-surface.  The dimension reduction of \KL  theory offers great convenience for studying   plates both analytically and numerically. This paper   aims for the development of accurate and efficient finite-difference schemes for solving the \KL plate model in domains with complex geometries.


The presence of the biharmonic operator in the governing PDE of the  \KL model poses a significant  challenge to its numerical methods.  To avoid the complication from the biharmonic  operator, numerical methods have been developed based on various  reformulations of  the system \cite{Ehrlich1971,ChenEtal2008,BecacheEtal2004} or mixed finite elements \cite{Ciarlet2002}; however, for some reformulations, 
it is challenging to support the  general boundary conditions (such as supported or free boundary conditions) beyond the essential  (i.e., clamped)  boundary conditions \cite{Ehrlich1971,ChenEtal2008}.
Solving the   \KL model directly using conforming finite element method (FEM) requires the use of  finite elements of class ${C}^1$. However,  ${C}^1$ elements are difficult to construct in multi-dimensions \cite{Ciarlet2002,Belytschko2005};  hence they  are  rarely used in practice except for beams (1D plates). To circumvent this issue, one could use non-conforming elements  such as Morley elements \cite{Morley1968,brenner1996,BrennerSung1999,MingXu2006,LiEtal2014} or a discontinuous Galerkin method~\cite{noels2008new} where high-order continuities along the element edges are only enforced weakly.

An alternative approach to handle the biharmonic operator is through finite differences. Based on  direct finite-difference approximations of the  biharmonic operator and all the other lower-order derivatives,   we have recently developed efficient and accurate finite-difference methods   for the \KL plate equation~\cite{NguyenEtal2020,JiLi2019}. These methods support  all the  boundary conditions (i.e., clamped, supported and free) and  are straightforward to implement, but are limited to simple domains due to  constraints with respect to meshing using a single structured grid. 

Supporting a general boundary condition is a key step towards fluid-structure interaction (FSI)  solvers involving plates.
For instance, in our previous work of~\cite{LiHenshaw2016}, 
we derived
a generalized interface coupling condition for the beam/shell coupled with an incompressible flow,
which can be viewed as a generalization of a non-standard boundary condition. 
Such an interface condition has been extended to other regimes such as the rigid body~\cite{rbinsmp2017, rbins2017, rbins3d2018}
or elastic solids~\cite{fibrmp2019, fibr2019} coupled with an incompressible flow.
In this work, as a first step towards an FSI solver involving plates in complex moving and deforming domain, we focus on the stable and robust schemes for solving the \KL plate with complex geometry  along with different boundary conditions. 

Many numerical methods have been developed based on solving  alternative models that are closely related to the \KL plate.  For example, 
the \KL plate model can be  regarded  as the  thin plate limit of the  Reissner–Mindlin  theory~\cite{reddy2006theory}, which concerns the rotation in addition to the    deflection of plates.  The Reissner-Mindlin  theory is  more accurate for thicker plates and is  easier to solve numerically because  only   ${C}^0$ elements  are  needed to approximate its unknowns (deflection and  rotation).  Therefore, good numerical methods developed for the  Reissner-Mindlin  plate can be exploited  to solve  the \KL limit. For this type of methods, 
one needs to pay special attention to 
address the shear and membrane locking phenomena \cite{daVeiga-2007,daVeiga-2008,daVeigaEtal2015} that are caused by the  inconsistency between the \KL plate and    Reissner-Mindlin   plate at zero thickness. Continuum based (CB) element method is another widely used method in commercial software and research for solving  plate/shell models   \cite{Stanley1985,KlinkelEtal1994,Belytschko2005}. The essence of the CB  methodology is to derive a simplified  model (referred to as CB shell elements)  for the thin structure   at   the discrete level by imposing   kinematic  assumptions to the discretization   of the  entire  3D solid.  

Isogeometric analysis (IGA) that uses NURBS basis functions for    more precise geometric representations   was proposed as a generalization of classical finite element analysis \cite{IGAoriginal2005,BazilevsEtal2006}. Since its introduction, IGA has gained increasing attention in engineering  and applied sciences communities for  simulating  challenging  PDEs in domains with complex geometries. 
Recently, IGA
 and its variations (e.g.,  isogeometric collocation or  Galerkin) have been well developed and widely used for problems with higher order derivatives such as various plate  and shell models; see for example  \cite{KiendlEtal2009,kiendl2010bending, KiendlEtal2015, nguyen2015extended, ZouEtal2021,daVeigaEtal2015,BensonEtal2010} and the references therein. 

 In this paper, we propose to develop finite-difference based  numerical methods for solving \KL plates on composite overlapping grids. A composite overlapping grid refers to the  collection of  logically rectangular curvilinear component grids that cover  the entire domain and overlap where they meet.  Overlapping grids, also known as overset or Chimera grids, are   often used for the  efficient and accurate solution  of  PDEs on regions of complex geometry \cite{CGNS}.  The novel algorithms presented here   are based on the common spatial discretization of the PDE  on composite overlapping grid, which involves curvilinear finite-difference approximations for spatial derivatives on each component grid and interpolating formulas for coupling  solutions on overlapping  component grids. Four time-stepping methods, both explicit and implicit, are considered for the temporal integration  of the spatially discrete system. One numerical challenge of our approach lies in the weak instability caused by the presence of interpolating equations in the discrete system,  which  breaks the nice  symmetric property of the finite-difference discretization on structured component grids.  Motivated by an upwind scheme that was developed for solving the wave equation on overlapping grids \cite{sosup2012,mxsosup2018}, we propose to add a fourth-order hyper-dissipation to the spatially discretized system to  stabilize its time integration; this additional dissipation term captures the essential upwinding effect of the original upwind scheme. Analysis and  numerical tests  are carefully carried out to validate the numerical properties of the methods.  Finally, two  benchmark  problems involving plates with complex geometries  are presented to illustrate that our  finite-difference based approach  is  well-suited for solving plate models arising in  realistic applications.

 The remainder of the paper is organized as follows. In Section~\ref{sec:governingEqns}, we describe  the governing equation and boundary conditions for the problem considered, followed by a detailed  presentation of the proposed numerical algorithms in Section~\ref{sec:algorithm}.  In Section~\ref{sec:stabilityAnalysis}, we analyze the stability of our algorithms and derive formulas for determining stable time steps and sufficient dissipation that can be used in actual computations.  In Section~\ref{sec:model1d},  a simplied 1D model problem is considered to illustrate the unstable modes caused by the interpolating equations, and to investigate the effects of our stabilization strategies. Numerical results for convergence studies and benchmark simulations are discussed in Section~\ref{sec:results}. Finally, some concluding remarks are given in  Section~\ref{sec:conclusions}.  

\section{Governing Equations}\label{sec:governingEqns}
We consider the \KL plate model for  isotropic and homogeneous material  that incorporates various physics including  bending, tension and linear restoration (i.e., the Winkler foundation known in engineering community  \cite{vlasov1966beams,flugge2013viscoelasticity}). Let  $w(x,y,t)$ with $(x,y)\in\Omega\subset\mathbb{R}^2$ denote  the transverse displacement of the plate, then the  \KL  model can be  described by an initial-boundary value problem (IBVP).
Specifically, the  governing  PDE for the displacement    is given by 
\begin{equation} \label{eq:KLShell}
{\rho}{H}\pdn{w}{t}{2}=-{K}w+{T}\nabla^2w-{D} \nabla^4w+f, 
\end{equation}
where $f=f(x,y,t)$ is a given external forcing,  $\rho$ is density, $H$ is thickness, $K$ the linear stiffness coefficient that acts as a linear restoring force, $T$ is the tension coefficient, and $D = EH^3/(12(1 - \nu^2))$ represents the flexural rigidity with $\nu$ and $E$ being the Poisson's ratio and Young's modulus, respectively.

On any  boundary of the domain, $ \forall (x,y)\in\partial\Omega$, we may impose one  of the following boundary conditions,
\begin{alignat}{3}
  \text{clamped:} &\quad w=0, &&\quad \pd{ w}{\nv}=0; \label{eq:clampedBC}\\
  \text{supported:}  & \quad  w=0, &&\quad \pdn{ w}{\nv}{2}+\nu\pdn{ w}{\tv}{2}=0; \label{eq:supportedBC}\\
  \text{free:}  & \quad   \pdn{ w}{\nv}{2}+\nu\pdn{ w}{\tv}{2}=0, &&\quad   \pd{}{\nv}\left[\pdn{ w}{\nv}{2}+\left(2-\nu\right)\pdn{ w}{\tv}{2}\right]=0, \label{eq:freeBC}
\end{alignat}
where   $\partial/\partial{\nv}$ and  $\partial/\partial{\tv}$ are the normal and tangential derivatives defined on  the boundary of the domain. It is important to point out that, at corners between two free boundaries \eqref{eq:freeBC},  a corner condition that imposes zero forcing, $\partial^2w/\partial x\partial y=0$, must be included  \cite{bilbao2008family,JiLi2019}.

The initial state of the plate is defined by 
\begin{equation}\label{eq:IC}
w(x,y,0)=w_0(x,y) \quad\text{and}\quad \pd{w}{t}(x,y,0)=v_0(x,y),
\end{equation}
where $w_0$ and $v_0$  prescribe  the  initial displacement and velocity, respectively.

\section{Numerical Methods}\label{sec:algorithm}

We aim to develop efficient and accurate numerical methods for solving  \KL plates   with general geometries using composite overlapping grids.  First, 
we discuss composite overlapping grids and the associated discretization approach,
followed by the  presentation of four time-stepping schemes for the stable integration of the spatially discretized system.
We are also interested  in the velocity  $v(x,y,t)=\partial_t{w}(x,y,t)$ and acceleration   $a(x,y,t)=\partial^2_t{w}(x,y,t)$ of   plates; therefore,   all of our numerical methods are designed to   solve  for  $v$ and $a$ to the same accuracy  as  the displacement solution $w$. Note that $v$ and $a$  are crucial information for multi-physics problems such as FSI applications; the  accurate computation of these quantities  is essential for any future development of FSI solvers involving \KL plates.

\subsection{Composite overlapping grids}
Composite overlapping grids are efficient and  powerful techniques that are often used  for
the solution of PDEs on domains with complicated shapes \cite{CGNS}. 
In Figure~\ref{fig:cgExample}, we show an example of composite overlapping grid  for a square plate with a circular hole cut at the center. In general, a composite overlapping grid, $\Gc$, consists of a set of structured component grids, $\{\Gc^{(n)}\}$,
$g=1,\ldots,{\mathcal N}$, that cover the entire plate domain $\Omega$;  the component grids overlap where they
meet. Solutions on the different component grids in the overlapped region are coupled by interpolation.
Typically, boundary-fitted curvilinear grids are used near the
boundaries to resolve the shape of a plate, while one or more
background Cartesian grids are used to handle the bulk of the  domain for efficiency. 
Each component grid $\Gc^{(n)}$ is a logically rectangular, curvilinear grid in 2D, and is defined by a smooth 
mapping that maps the unit square into the subdomain covered by $\Gc^{(n)}$; i.e.,
$
\xv = \gv(\rv)$
with
 $\rv=(r,s)\in[0,1]^2$ and $ \xv=(x,y)\in\mathbb{R}^2$.  The  original $(x,y)$-space is referred to as physical space, and the unit square $(r,s)$-space is called the reference space.

{
\newcommand{\figWidth}{6cm}
\def\xa{6.}
\def\ya{5.5}
\newcommand{\trimfig}[2]{\trimw{#1}{#2}{0.12}{0.12}{0.12}{0.12}}
\begin{figure}[h]
\begin{center}
\begin{tikzpicture}[scale=1]
  \useasboundingbox (0.0,0.0) rectangle (\xa,\ya);  

\draw(-0.5,-0.5) node[anchor=south west,xshift=0pt,yshift=0pt] {\trimfig{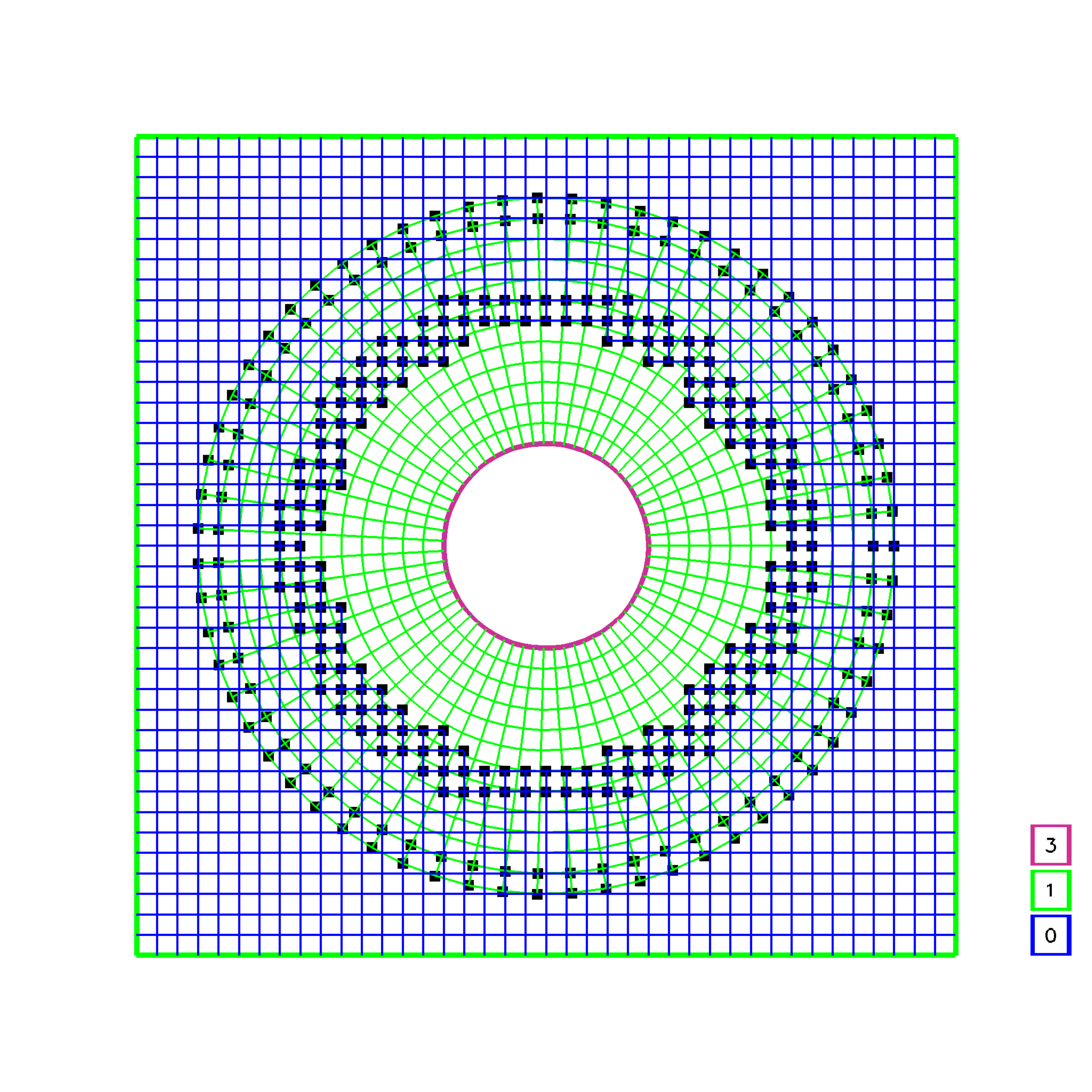}{\figWidth}};

%
\end{tikzpicture}

\end{center}
\caption{ Composite overlapping grid $\Gc_1$ for a square plate with a circular hole cut at the center. The side of the square is 4 and the radius of the circle is 1/2. Interpolation points on both component grids are highlighted.} \label{fig:cgExample}
\end{figure}
}

Grid points on a composite overlapping grid  are classified as discretization points, interpolation points or unused points. On the  discretization points of each component grid,  curvilinear finite-difference methods are used to discretize the spatial derivatives of the PDE and boundary conditions;    ghost points are used to aid the discretization of physical boundary  conditions.   Numerical  solutions  between different component grids are coupled together with interpolating equations.  The interpolation points of the example composite grid are highlighted in Figure~\ref{fig:cgExample}.

We  build a sequence of refined meshes denoted as $\{\Gc_n\}$ for each problem  in practice, where $\Gc_1$  represents the base grid that has a target grid spacing $h=1/10$. For curved grids (non-Cartesian), the grid spacings are not constant. So the target grid spacing is used as a guidance, and  we try to maintain the cell volumes as uniform as possible.
The  grid shown in Figure~\ref{fig:cgExample} is the base grid  $\Gc_1$ for this  example grid.  Finer grid $\Gc_n$ is  obtained from refining the base grid by a factor of $n$; that is, the target grid spacing for $\Gc_n$ is $1/(10n)$. For accuracy studies,  we typically perform convergence studies on the sequence of grids,  $\{\Gc_1,\Gc_2,\Gc_4,\Gc_8,\Gc_{16}\}$.

\subsection{Spatial discretization on curvilinear grid}
To discretize the governing equation \eqref{eq:KLShell} on a component grid, we first transform the problem  in  physical  space to  the reference space, and then approximate the partial derivatives  in the $(r, s)$-space using   standard central difference formulas. The transformed derivatives are obtained using the  derivatives of the  mapping $\xv = \gv(\rv)$ and  chain rules.

Specifically, the \KL plate equation \eqref{eq:KLShell} on the  reference domain can be written as
\newcommand{\cc}[1]{{\color{black} #1}}
\def\wo{W}
\def\wr{\cc{W_{r}}}
\def\ws{\cc{W_{s}}}
\def\wrr{\cc{W_{rr}}}
\def\wrs{\cc{W_{rs}}}
\def\wss{\cc{W_{ss}}}
\def\wrrr{\cc{W_{rrr}}}
\def\wrrs{\cc{W_{rrs}}}
\def\wrss{\cc{W_{rss}}}
\def\wsss{\cc{W_{sss}}}
\def\wrrrr{\cc{W_{rrrr}}}
\def\wrrrs{\cc{W_{rrrs}}}
\def\wrrss{\cc{W_{rrss}}}
\def\wrsss{\cc{W_{rsss}}}
\def\wssss{\cc{W_{ssss}}}
\begin{align}
  \rho H \pdn{W}{t}{2}= \Lc W+F, \label{eq:KLShellReference}
\end{align}
where the plate operator $\Lc$ is defined by 
\begin{align*}
 \Lc W = &-KW+a_1(\rv) \wr+a_2(\rv) \ws+b_{11}(\rv)\wrr+b_{12}(\rv)\wrs+b_{22}(\rv)\wss \nonumber \\
&+c_{111}(\rv)\wrrr+c_{112}(\rv)\wrrs+c_{122}(\rv)\wrss+c_{222}(\rv)\wsss\\
&+d_{1111}(\rv)\wrrrr+d_{1112}(\rv)\wrrrs+d_{1122}(\rv)\wrrss+d_{1222}(\rv)\wrsss+d_{2222}(\rv)\wssss.
\end{align*}
Here $ W(\rv) \equiv w(\xv) =w(\gv(\rv))$,  $ F(\rv,t) \equiv f(\xv,t) =f(\gv(\rv),t)$  and the  coefficients $a_i(\rv), b_i(\rv),c_i(\rv),d_i(\rv)$  are   derivatives of  the mapping  $\gv$ with respect to  $r$ and $ s$.  For a concise presentation of the paper, the complete  expression  of \eqref{eq:KLShellReference} as well as the coefficients 
 is presented  in \ref{sec:KLShellReferenceFormula}.

Derivatives of $W$ in $(r,s)$-space  are approximated using standard finite difference formulas on the unit square grid. To facilitate the discussion of discretizing \eqref{eq:KLShellReference}, we introduce the following basic finite-difference operators,
$$
\Dpr W_{ij}=\frac{W_{i+1,j}-W_{i,j}}{\dr},~\Dmr W_{ij}=\frac{W_{i,j}-W_{i-1,j}}{\dr},~\Dor W_{ij}=\frac{W_{i+1,j}-W_{i-1,j}}{2\dr},
$$
where $W_{ij}$ denotes the  solution on grid $\rv_{i,j}$, and $\dr$ is the grid spacing in the $r$-direction of the unit square grid. Similar operators, $\Dps$, $\Dms$, and $\Dos$, can be defined in the $s$-direction.

\begin{figure}[ht]
  \centering
  \includegraphics[width=6cm]{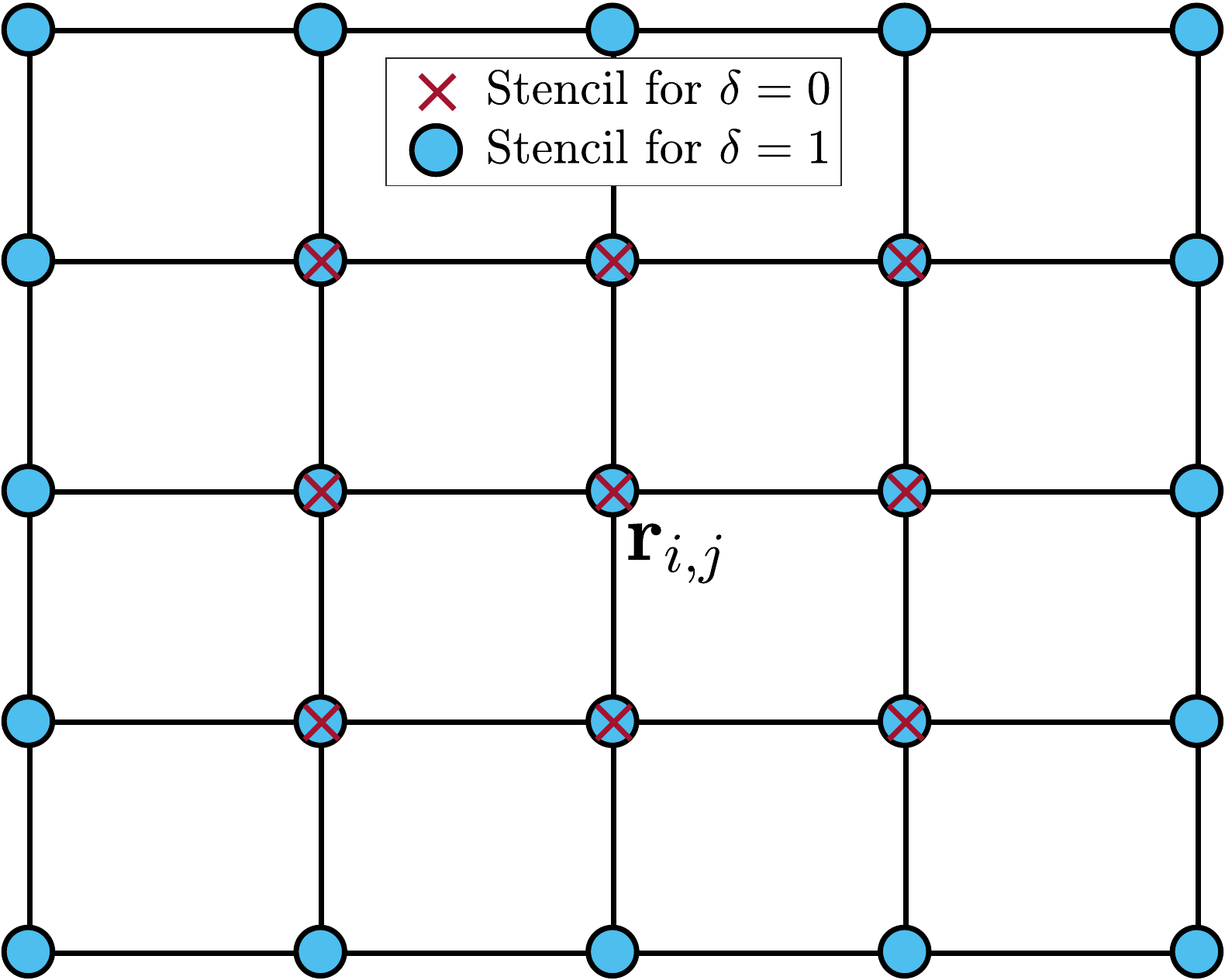}
  \caption{Stencil of the finite-difference operators.\label{fig:stencil} }
\end{figure}

We approximate all the third- and fourth-order derivatives using 2nd-order accurate  central finite difference formulas. Thus, the discrete operators for all the third-order derivatives are
\begin{alignat}{4}\label{eq:3rdDerivFD}
   D_{rrr}= \Dpr\Dmr\Dor, &\quad D_{rrs} =\Dpr\Dmr\Dos, &\quad D_{rss}=\Dor\Dps\Dms,   &\quad D_{sss}=\Dps\Dms\Dos,
\end{alignat}
and those for all the fourth-order derivatives are
\begin{equation}\label{eq:4thDerivFD}
  \begin{alignedat}{3}
     &D_{rrrr}=(\Dpr\Dmr)^2, &&\quad D_{rrrs}=\Dpr\Dmr\Dor\Dos, &&\quad D_{rrss}=\Dpr\Dmr\Dps\Dms,\\
     & D_{rsss}=\Dor\Dos\Dps\Dms,&&\quad D_{ssss}=(\Dps\Dms)^2. &&
     \end{alignedat}
\end{equation}
Note that the difference operators defined in \eqref{eq:3rdDerivFD} and \eqref{eq:4thDerivFD} utilize a stencil that  takes  5 points in each direction of  the $(r,s)$-space 
as is shown in Figure~\ref{fig:stencil}.

For the lower order (1st or 2nd) derivatives in \eqref{eq:KLShellReference}, we have the freedom to approximate the derivatives   either maintaining  the same   2nd-order accuracy by taking  a smaller stencil, or maintaining  the same stencil size as all the other difference operators in  \eqref{eq:3rdDerivFD} and \eqref{eq:4thDerivFD}  by implementing  4th-order accurate schemes.  Specifically, the finite difference operators for the first and second order derivatives can be defined by
\begin{alignat}{3}
  &D_r= \Dor- \delta \frac{\dr^2}{6} \Dor\Dpr\Dmr,  && \quad D_s= \Dos- \delta \frac{\ds^2}{6} \Dos\Dps\Dms, && \label{eq:1stDerivFD}\\
  & D_{rr}= \Dpr\Dmr -\delta \frac{\dr^2}{12}(\Dpr\Dmr)^2, && \quad D_{ss}= \Dps\Dms -\delta \frac{\ds^2}{12}(\Dps\Dms)^2,  &&  \quad D_{rs}=D_rD_s. \label{eq:2ndDerivFD}
\end{alignat}
These formulas are 2nd-order accurate if $\delta=0$ (same order), and  are 4nd-order accurate if  $\delta=1$ (same stencil); see  Figure~\ref{fig:stencil} for the stencils for $\delta=0$ and $\delta=1$.

 Using the previously defined finite-difference operators, we are now ready  to present the discrete  equations  for approximating   \eqref{eq:KLShellReference} on a component grid, given by:
\def\Dr{D_{r}W_{ij}}
\def\Ds{D_{s}W_{ij}}
\def\Drr{D_{rr}W_{ij}}
\def\Drs{D_{rs}W_{ij}}
\def\Dss{D_{ss}W_{ij}}
\def\Drrr{D_{rrr}W_{ij}}
\def\Drrs{D_{rrs}W_{ij}}
\def\Drss{D_{rss}W_{ij}}
\def\Dsss{D_{sss}W_{ij}}
\def\Drrrr{D_{rrrr}W_{ij}}
\def\Drrrs{D_{rrrs}W_{ij}}
\def\Drrss{D_{rrss}W_{ij}}
\def\Drsss{D_{rsss}W_{ij}}
\def\Dssss{D_{ssss}W_{ij}}
\begin{equation}
  \rho H\ddn{W_{ij}}{t}{2}=\Lc_hW_{ij}+F_{ij},\label{eq:KLShellDiscrete}
\end{equation}
where the discrete plate operator $\Lc_h$ is given by 
  \begin{align*}
\Lc_hW_{ij}=&-Kw+a_1(\rv) \Dr+a_2(\rv) \Ds+b_{11}(\rv)\Drr+b_{12}(\rv)\Drs+b_{22}(\rv)\Dss\\
&+c_{111}(\rv)\Drrr+c_{112}(\rv)\Drrs+c_{122}(\rv)\Drs+c_{222}(\rv)\Dsss+d_{1111}(\rv)\Drrrr\\
&+d_{1112}(\rv)\Drrrs+d_{1122}(\rv)\Drrss+d_{1222}(\rv)\Drsss+d_{2222}(\rv)\Dssss. 
    \end{align*}
Given the order of accuracy for all the difference operators defined in \eqref{eq:3rdDerivFD} -- \eqref{eq:2ndDerivFD}, we expect the truncation error of \eqref{eq:KLShellDiscrete} to be 2nd-order.

If a physical boundary is  present on  this component grid,  the discretized boundary condition can be readily derived by replacing the derivatives with the corresponding finite-difference approximations.

\subsection{Time-stepping schemes}\label{sec:tsSchemes}
For numerical purposes, we propose to modify \eqref{eq:KLShellDiscrete} by adding  an artificial hyper-dissipation term. The modified plate equation is given by 
\begin{equation}
\ddn{W_{ij}}{t}{2}=\frac{1}{  \rho H}\left(\Lc_hW_{ij}+F_{ij}\right)+\ad\Dad\dd{W_{ij}}{t},\label{eq:KLShellReferenceAD}
\end{equation}
where $\ad$ is the dissipation parameter and $\Dad=-\dr^4D_{rrrr}-\ds^4D_{ssss}$ is the difference operator for the artificial hyper-dissipation.
The dissipation term, inspired by the  upwind schemes for  the wave equation on overlapping grids \cite{sosup2012,mxsosup2018,ssmx2018},   incorporates the essential   upwinding effect into the problem  and    is indispensable   for stabilizing the plate simulations on composite overlapping grids. It is important to point out that with $\ad=\Oc(1/h^2)$ the dissipation term is formally $\Oc(h^2)$ so that its addition does not affect the overall accuracy of the spatial  discretization;  here  $h=\min(\dr,\ds)$.

In this paper, we propose the following time-stepping schemes to advance the  modified discrete \KL plate equations \eqref{eq:KLShellReferenceAD} in time. Let $W^n_{\iv}$, $V^n_{\iv}$, and $A^n_{\iv}$ denote the numerical solutions of the plate displacement, velocity and acceleration at time $t_n$, respectively. Here $\iv=(i,j)$ is a multi-index for the grid point, and $t_n=n\dt$ with  a fixed  time-step $
\dt$. The aim  of the following  time-stepping algorithms is to determine the numerical  solutions  at a new time given solutions at previous time levels. If  boundary conditions for velocity or acceleration need to be applied by a particular scheme,  they can be derived  by  taking  appropriate time derivatives of the displacement boundary conditions given in \eqref{eq:clampedBC} -- \eqref{eq:freeBC}.

Solving \eqref{eq:KLShellReferenceAD} in the second-order form directly  leads to the algorithms referred to by us  as  the C2 and UPC2 schemes, which  are listed in Algorithms~\ref{alg:c2} and~\ref{alg:upc2}, respectively. Both schemes advance the system in time by approximating  the second-order time derivative with the  second-order centered finite-difference formula (C2). The difference between C2 and UPC2  lies in the treatment of  $\partial_tW_{\iv}$ in the artificial dissipation term. In particular, we approximate $\partial_tW_{\iv}$ with  a backward time difference for the  C2 scheme,  while we  include  the upwind dissipation term using a predictor-corrector scheme for the UPC2 scheme, so that both schemes remain  explicit with or without the artificial dissipation. But how we deal with the dissipation  leads to different time step restrictions that is to be discussed in Section~\ref{sec:timeStep}. 

\begin{algorithm}[h!]
 {\bf Input:} {displacement  at two previous time levels; i.e.,   $W_\iv^n$ and $W_\iv^{n-1}$ }\\
 {\bf Output:} {new displacement  $W_\iv^{n+1}$, and current velocity and acceleration $V_\iv^{n}$,$A_\iv^{n}$  }\\
 {\bf Procedures:} \\
 {\em Advance the equation using centered time difference scheme }
 \begin{align}
   &A_\iv^{n}=\frac{1}{  \rho H}\left(\Lc_hW^{n}_{\iv}+F^{n}_{\iv}\right)\nonumber\\
   & W_\iv^{n+1}= 2W_\iv^{n}-W_{\iv}^{n-1}+\dt^2A_\iv^{n}+\dt\ad\left(\Dad W^{n}_\iv  - \Dad W^{n-1}_\iv \right)\label{eq:C2updating}\\
   &V_\iv^{n}=\frac{W_\iv^{n+1}-W_\iv^{n-1}}{2\dt}\nonumber
 \end{align}   
     {\bf Remark:} {\em $V_\iv^{n}$ and  $A_\iv^{n}$ are obtained from post-processing the displacement solutions, which are not involved in the main updating formula \eqref{eq:C2updating}; therefore, boundary conditions are enforced  for $W_\iv^{n+1}$ only to fill in the values at ghost and/or boundary grid points.  Interpolation routines  are called at the end of each step to fill in solutions for  the interpolation points of composite overlapping grids.}
  \caption{C2 time-stepping scheme}\label{alg:c2}
\end{algorithm}

\begin{algorithm}[h!]
 {\bf Input:} {displacement  at two previous time levels; i.e.,   $W_\iv^n$ and $W_\iv^{n-1}$ }\\
 {\bf Output:} {new displacement  $W_\iv^{n+1}$, and current velocity and acceleration $V_\iv^{n}$,$A_\iv^{n}$  }\\
 {\bf Procedures:} \\
 {\em Stage I:  predict the solution using C2 scheme  without artificial dissipation }
 \begin{align}
   &A_\iv^{n}=\frac{1}{  \rho H}\left(\Lc_hW^{n}_{\iv}+F^{n}_{\iv}\right)\nonumber\\
   & W_\iv^{p}= 2W_\iv^{n}-W_{\iv}^{n-1}+\dt^2A_\iv^{n}\label{eq:UPC2updatingP}
 \end{align}
     {\em Stage II:  correct  the solution with the artificial dissipation }
     \begin{align}
       &W_\iv^{n+1}= W_\iv^{p}+\frac{1 }{2}\dt\ad\left(\Dad W^{p}_\iv-\Dad W^{n-1}_\iv \right)\label{eq:UPC2updatingC}\\
       &V_\iv^{n}=\frac{W_\iv^{n+1}-W_\iv^{n-1}}{2\dt}\nonumber
     \end{align}
{\bf Remark:} {\em $V_\iv^{n}$ and  $A_\iv^{n}$ are obtained from post-processing the displacement solutions, which are not involved in the main updating formulas \eqref{eq:UPC2updatingP} and  \eqref{eq:UPC2updatingC}; therefore, boundary conditions are enforced  for $W_\iv^{n+1}$ only to fill in the values at ghost and/or boundary grid points.  Interpolation routines  are called at the end of both stages to fill in solutions for  the interpolation points of composite overlapping grids.}
 \caption{UPC2 time-stepping scheme}\label{alg:upc2}
\end{algorithm}

We have previously developed two numerical methods referred to as the PC22 and NB2 schemes to solve a generalized \KL model written in  first order form  \cite{NguyenEtal2020}. However, both schemes become unstable if applied directly for solving the plate equation on composite overlapping grids due to the presence of interpolating equations in the spatial discretization. In this paper,  we propose to stabilize  the schemes by solving the  first order form of the modified  plate equation \eqref{eq:KLShellReferenceAD},
\begin{equation}\label{eq:KLShellReferenceAD1stOrderForm}
  \left\{
  \begin{aligned}
&\dd{W_{ij}}{t}=V_{ij}\\  
&\dd{V_{ij}}{t}=A_{ij}+\ad\Dad V_{ij}
  \end{aligned}
  \right.
  \quad\text{with} \quad A_{ij}=\frac{1}{  \rho H}\left(\Lc_hW_{ij}+F_{ij}\right).
\end{equation}

For each step of the  PC22 time-stepping scheme, the algorithm advances \eqref{eq:KLShellReferenceAD1stOrderForm} by taking a second-order accurate Adams-Bashforth  predictor, followed by  a second-order Adams-Moulton corrector. Specifically, the scheme is summarized in Algorithm~\ref{alg:pc22}. 

\begin{algorithm}[h!]
 {\bf Input:} {solutions at two previous time levels; i.e.,   $(W_\iv^n,V_\iv^n,A_\iv^n)$ and $(W_\iv^{n-1},V_\iv^{n-1},A_\iv^{n-1})$ }\\
 {\bf Output:} {solutions at the new time level; i.e., $(W_\iv^{n+1},V_\iv^{n+1},A_\iv^{n+1})$  }\\
  {\bf Procedures:} \\
 {\em Stage I:  predict solutions using a second-order Adams-Bashforth (AB2) predictor}
 \begin{equation*}
   \begin{aligned}
     &W_\iv^p= W_\iv^n+\dt\left(\frac{3}{2}V_\iv^{n}-\frac{1}{2}V_\iv^{n-1}\right)\\
     &V_\iv^p= V_\iv^n+\dt\left[\frac{3}{2}(A_\iv^{n}+\ad\Dad V^n_{ij})-\frac{1}{2}(A_\iv^{n-1}+\ad\Dad V^{n-1}_{ij})\right] \\
     &A^{p}_\iv=\frac{1}{  \rho H}\left(\Lc_hW^p_{ij}+F^{n+1}_{ij}\right).
   \end{aligned}
 \end{equation*}
 {\em Stage II:   correct solutions using a  second-order Adams-Moulton (AM2) corrector}
  \begin{equation*}
   \begin{aligned}
     &W_\iv^{n+1}= W_\iv^n+\dt\left(\frac{1}{2}V_\iv^{n}+\frac{1}{2}V_\iv^{p}\right)\\
     &V_\iv^{n+1}= V_\iv^n+\dt\left[\frac{1}{2}(A_\iv^{n}+\ad\Dad V^n_{ij})+\frac{1}{2}(A_\iv^{p}+\ad\Dad V^{p}_{ij})\right] \\
     &A^{n+1}_\iv=\frac{1}{  \rho H}\left(\Lc_hW^{n+1}_{ij}+F^{n+1}_{ij}\right).
   \end{aligned}
  \end{equation*}
  
  {\bf Remark}:{\em Boundary conditions are applied after both the predictor and corrector stages to fill in the solutions of $W$ and $V$ at ghost and/or boundary grid points. We do not apply boundary conditions for  $A$ since it is derived from post-processing $W$ solutions and is not  part of the main updating formulas.      Interpolation routines  are called at the end of both stages to fill in solutions for  the interpolation points of composite overlapping grids. }
 \caption{PC22 time-stepping scheme}\label{alg:pc22}
\end{algorithm}

The NB2 scheme takes advantage of a second-order accurate version of the well-known  Newmark-Beta scheme \cite{Newmark59} and implements it in a predictor-corrector fashion. The complete algorithm for the NB2 scheme is given in Algorithm~\ref{alg:nb2}.

\begin{algorithm}[h!]
{\normalsize
 {\bf Input:} {solutions at the previous time level; i.e.,   $(W_\iv^n,V_\iv^n,A_\iv^n)$  }\\
 {\bf Output:} {solutions at the new time level; i.e., $(W_\iv^{n+1},V_\iv^{n+1},A_\iv^{n+1})$  }\\
 {\bf Procedures:} \\
 {\em Stage I.  compute a first-order prediction for displacement and velocity}
 \begin{equation*} 
\begin{aligned}
  &W_{\iv}^{p}= W_{\iv}^{n}+\dt V_{\iv}^n +\frac{\dt^2}{2} (1-2\beta)A_{\iv}^{n}  \\
  &V_{\iv}^{p} = V_{\iv}^{n}+\dt  (1-\gamma)A_{\iv}^{n}
\end{aligned}
 \end{equation*}
{\em Stage II.   solve a system of equations for acceleration at $t_{n+1}$}
\begin{equation*}
 \left( {\rho}{H}-\beta \dt^2\Lc_h-\gamma\dt\ad\Dad\right)A^{n+1}_\iv={\Lc_h} W_{\iv}^{p}+F_{\iv}^{n+1}+\ad\Dad V_{\iv}^{p}
\end{equation*}

{\em Stage III.  correct the  displacement and velocity using the  new acceleration solution}\\
\begin{equation*} 
\begin{aligned}
  &W_{\iv}^{n+1}= W_{\iv}^{n}+\dt V_{\iv}^n +\frac{\dt^2}{2}\left[ (1-2\beta)A_{\iv}^{n}+2\beta A_{\iv}^{n+1} \right]  \\
  &V_{\iv}^{n+1} = W_{\iv}^{n}+\dt \left[ (1-\gamma)a_{\iv}^{n}+\gamma A_{\iv}^{n+1} \right]
\end{aligned}
\end{equation*}

{\bf Remark}: {\em A second-order accurate Newmark-Beta scheme with $\beta = 1/4$ and $\gamma = 1/2$ is used here.  Boundary conditions are applied after stages I and III to fill in the solutions of $W$ and $V$ at ghost and/or boundary grid points; for composite overlapping grids, interpolation routines are called  as well.  For stage II,  equations for acceleration at ghost and boundary nodes are given by  boundary  conditions, while equations  at the interpolation points are given by a 4th order polynomial interpolation formula.} 
}
 \caption{NB2 time-stepping scheme}\label{alg:nb2}
\end{algorithm}

\clearpage
\section{Stability Analysis}\label{sec:stabilityAnalysis}
Strategies for determining the dissipation coefficient $\ad$ in \eqref{eq:KLShellReferenceAD} and a stable time step $\dt$ for each of the algorithms proposed in Section~\ref{sec:tsSchemes} are discussed here. 
Formulas for curvilinear grids are derived   in the usual way by freezing coefficients and using a local Fourier analysis.

To be specific, the discrete Fourier transformation of the homogeneous version of  \eqref{eq:KLShellReferenceAD} for  a given  wave number pair $\omega=(\omega_r,\omega_s)$ is found to be
\begin{equation}
\ddn{\hat{W}_{\omega}}{t}{2}=\frac{1}{\rho H}\hat{Q}(\xi_r,\xi_s;\rv) \hat{W}_{\omega}-\ad\left(16\sin^4(\xi_r/2)+16\sin^4(\xi_s/2)\right) \dd{\hat{W}_{\omega}}{t}, \label{eq:plateEqnDFT}
\end{equation}
where $\hat{Q}(\xi_r,\xi_s;\rv) $ denotes the  Fourier transform (symbol) of the discrete operator $\Lc_h$, and $\hat{W}_{\omega}$ is the Fourier coefficient of $W_{jk}$ that is defined by $W_{jk}=\hat{W}_{\omega}e^{i\xi_r j}e^{i\xi_s k}$ with  $\xi_r=2\pi\omega_r\dr$ and  $\xi_s=2\pi\omega_s\ds$.  The specific expression of  $\hat{Q}(\xi_r,\xi_s;\rv) $ is given in \ref{sec:formulaDFT} to save space. 
A particular  time-stepping scheme is stable for solving \eqref{eq:KLShellReferenceAD} provided it is stable for  \eqref{eq:plateEqnDFT}  $\forall (\xi_r,\xi_s)\in[-\pi,\pi]^2$ and $\forall \rv\in[0,1]^2$.  A sufficient condition for the stability of the scheme can be derived by considering the worst case scenario of \eqref{eq:plateEqnDFT}.

With a slight abuse of notations,   we consider the following test problem to study  the stability issues for solving the plate  problem,
\begin{equation}\label{eq:testProblem}
\ddn{w}{t}{2}=Q w -\mu \pd{w}{t} ~~\text{with}~~\mu\geq0,
\end{equation}
where $Q=\hat{Q}_M/(\rho H)$ such that $|\hat{Q}_M|=\max\{|\hat{Q}(\xi_r,\xi_s;\rv) |\}$ and $\mu=32\ad$. Note that $\hat{Q}_M$  and $\mu$ are attained at  $\xi_r=\pm\pi$ and  $\xi_s=\pm\pi$ which correspond to the plus-minus modes in each direction of the $(r,s)$-space.

\subsection{Dissipation coefficient determination}
First, we look at the test problem \eqref{eq:testProblem} analytically, whose  
 general solution   is given by 
\begin{equation}\label{eq:generalSolutionTP}
w(t)=c_1e^{\xi_+t}+c_2e^{\xi_-t}~~\text{with}~~
\xi_{\pm}=\frac{-\mu\pm\sqrt{\mu^2+4Q}}{2}.
\end{equation}
Here  $c_1$ and $c_2$ are constants to be determined by the initial conditions. Since the plate operator $\Lc$ in \eqref{eq:KLShellReference} is self-adjoint, the Fourier transform (or symbol) of the corresponding discrete operator $\Lc_h$ on a single grid should be a real negative number (c.f., \ref{sec:formulaDFT}). In this case, the test problem is well-posed in the sense that the solution \eqref{eq:generalSolutionTP} is bounded over time.

However, on composite overlapping grids, the symmetry of the differentiation matrix associated with $\Lc_h$ is spoiled by the interpolating equations.  As a result, 
   for some modes, the symbol of $\Lc_h$, as well as the eigenvalue of the differential matrix, is perturbed  off the real axis that can  cause the solution \eqref{eq:generalSolutionTP} to grow exponentially if  there is no sufficient dissipation.  Thus, we investigate imposing an artificial hyper-dissipation 
to stabilize the perturbation induced by interpolation. For this purpose,   we consider the test problem \eqref{eq:testProblem} with perturbed $Q$. In particular,  we  assume  $Q=m+i n$, where $m=\hat{Q}_M<0$ and $n\neq0$  is induced by perturbation of interpolation. In this case, $\xi_{\pm}$ in the general solution  \eqref{eq:generalSolutionTP} can be written as 
\begin{equation}
\xi_{\pm}=\frac{-\mu\pm\sqrt{r}e^{i\theta/2}}{2}=\frac{-\mu\pm\left(\sqrt{r}\cos(\theta/2)+i\sqrt{r}\sin(\theta/2)\right)}{2},
\end{equation}
where
$$
r=\sqrt{(\mu^2+4m)^2+16n^2},~~\theta=\arctan2\left({4n},{\mu^2+4m}\right).
$$
An artificial  dissipation is sufficient to stabilize the problem if 
$
\Re(\xi_\pm) \leq 0
$; therefore, a lower bound for $\mu$ can be derived from this restriction that is given by
\begin{equation}
  \mu\geq \sqrt{r}|\cos(\theta/2)|. \label{eq:muLowerBound}
\end{equation}
It is not straightforward to use \eqref{eq:muLowerBound} to determine the strength of  dissipation  to be added  in the algorithms  since it is hard to get an estimate for $n$ and $\theta$  without numerically computing the eigenvalues, which can be very expensive computationally.

On the other hand, we do not want the dissipation to be too strong to  damp out the oscillatory nature of the solution \eqref{eq:generalSolutionTP}, so we restrict
$\mu^2+4m<0$.  Given that $m=\hat{Q}_M<0$ and $\mu=32\ad$, we may specify  the dissipation coefficient as 
\begin{equation}\label{eq:dissipationCoefficient}
\ad = \df\frac{\sqrt{-\hat{Q}_M}}{16},
\end{equation}
where $\df\in[0,1]$ is a dissipation factor that can be used  to control the strength of the artificial   dissipation in practice.

\subsection{Time step determination}\label{sec:timeStep}

We now analyze the stability for the various time-stepping  methods using the test problem \eqref{eq:testProblem} with perturbation ($Q=m+in$) and dissipation taken into consideration.    An estimate of a stable time step $\dt$ is obtained  for each time-stepping scheme following the stability analysis.

\medskip
\noindent  {\bf C2 scheme}.

The test problem if advanced  with the C2 scheme leads to the following difference equation,
    $$
   \frac{w^{n+1}-2w^{n}+w^{n-1}}{\dt^2}= Qw^{n}-\mu \frac{w^{n}-w^{n-1}}{\dt},
   $$
whose characteristic equation is  given by
   $$
      \zeta^2-(2+\dt^2Q-\mu\dt)\zeta+(1-\mu\dt)=0.
   $$
      The stability of the scheme depends on the root condition $| \zeta_\pm|\leq1$.
      
Without dissipation, we find that the scheme is nondissipative ($|\xi_{\pm}|=1$) and  works only for the unperturbed case ($Q=\hat{Q}_M$) subject to  a time step restriction given by 
 $$
 \dt \leq \frac{2}{\sqrt{-\hat{Q}_M}}.
 $$

 With dissipation, the time step constraint derived from solving the root condition $| \zeta_\pm|\leq1$ is 
 $$
   \dt\leq\min\left\{\frac{4\sqrt{2}}{\mu+\sqrt{\mu^2-8\sqrt{2}m}}, \frac{2}{\mu}\right\},
$$
   where the assumption that the perturbation induced imaginary part is much smaller in magnitude (i.e., $|n|\ll |m|$) has been used.

   Noting that $m=\hat{Q}_M$, we propose the following formula for determining  $\dt$ in   Algorithm~\ref{alg:c2},
   \begin{equation}   \label{eq:c2dt}
     \dt=\left\{
     \begin{aligned}
        &\csf\frac{2}{\sqrt{-\hat{Q}_M}}, &&\quad \text{if}\quad \mu=0, \\
        &\csf\min\left\{\frac{4\sqrt{2}}{\mu+\sqrt{\mu^2-8\sqrt{2}\hat{Q}_M}}, \frac{2}{\mu}\right\},  &&\quad \text{if}\quad \mu\neq 0. 
     \end{aligned}
     \right.
   \end{equation}
   Here $\csf \in (0, 1]$ is a stability factor (sf) that multiplies an estimate of the largest stable time step based on the above  analysis, and $\mu=32\ad$ with the  dissipation coefficient   defined in \eqref{eq:dissipationCoefficient}.

\medskip
\noindent  {\bf UPC2 scheme}.

The UPC2 scheme  applied to  the test problem leads to 
\begin{align*}
   &\frac{w^{p}-2w^{n}+w^{n-1}}{\dt^2}= Qw^{n}, \\
   &\frac{w^{n+1}-w^{p}}{\dt^2}= -\mu\frac{w^{p}-w^{n-1}}{2\dt}.
\end{align*}
Combining the predictor and corrector steps by eliminating $w^p$,  the UPC2 scheme is essentially a 3-step scheme given by
$$
w^{n+1}-(1-\frac{1}{2}\mu\dt)(2+\dt^2Q)w^{n}+(1-\mu\dt)w^{n-1}=0.
$$
Similarly, the root condition $|\zeta_{\pm}|\leq1$ for the characteristic equation
implies the following time step restriction
$$
\dt\leq\min\left\{\frac{2}{\sqrt{-m}} ,\frac{2}{\mu}\right\}.
$$
Therefore,   we use  the following time step for Algorithm~\ref{alg:upc2},
\begin{equation}
  \label{eq:upc2dt}
\dt =\csf\min\left\{\frac{2}{\sqrt{-\hat{Q}_M}} ,\frac{2}{\mu}\right\},
\end{equation}
where  $\csf \in (0, 1]$ and $\mu=32\ad$ with $\ad$ given in \eqref{eq:dissipationCoefficient}.

\medskip
\noindent  {\bf PC22 scheme}.

The PC22 scheme solves the test problem in the first-order form,
$$
\dd{\uv}{t}=A \uv\quad\text{with}\quad
A =
\begin{bmatrix}
0 & 1\\
 Q & -\mu
\end{bmatrix},\quad\text{and}\quad \uv=
\begin{bmatrix}w\\v\end{bmatrix}. 
  $$
  The  time-stepping scheme is stable  for  the first order system provided the scheme is stable for all the eigenvalues of $A$. The 3-step updating formula for the PC22 scheme is found to be
  $$
\uv^{n+1}-\left(1+z+\frac{3}{4}z^2\right)\uv^n+\frac{1}{4}z^2\uv^{n-1}=0,
$$
where $z=\dt\lambda_A$ and $\lambda_A=\left(-\mu\pm\sqrt{\mu^2+4Q}\right)/2$ are the eigenvalues of $A$.

Following \cite{NguyenEtal2020}, by  requiring $z$ to be inside a super-ellipse that approximates the stability region of the scheme, a stable  time step for Algorithm~\ref{alg:pc22} can be derived;
that is, 
\begin{equation}
  \label{eq:pc22dt}
\dt = \csf\left( \left|\frac{\mu}{2r_a}\right|^\frac{3}{2}+\left|\frac{\sqrt{-\mu^2-4\hat{Q}_M}}{2r_b}\right|^\frac{3}{2}\right)^{-\frac{2}{3}}~~\text{with} ~~r_a=1.75, ~r_b=1.2.
\end{equation}
Similarly, we may take $\csf \in (0, 1]$ and $\mu=32\ad$ with $\ad$ given in \eqref{eq:dissipationCoefficient}.

\medskip
\noindent  {\bf NB2 scheme}.

The NB2 scheme for the test problem \eqref{eq:testProblem} can be written as
  $$
  \begin{bmatrix}
    w^{n+1}\\v^{n+1}
  \end{bmatrix}
  =A^{-1}B \begin{bmatrix}
    w^{n}\\v^{n}
  \end{bmatrix}
  $$
  where
  $$
  A=\begin{bmatrix} 1-Q\,\beta\,{\mathrm{\dt}}^2 & \beta\,{\mathrm{\dt}}^2\,\mu \\ -Q\,\mathrm{\dt}\,\gamma & \mathrm{\dt}\,\gamma\,\mu +1 \end{bmatrix}
~~\text{and}~~
B=\begin{bmatrix} 1-Q\,{\mathrm{\dt}}^2\,\left(2\,\beta-1\right) & \frac{\mu \,\left(2\,\beta-1\right)\,{\mathrm{\dt}}^2}{2}+\mathrm{\dt}\\ -Q\,\mathrm{\dt}\,\left(\gamma-1\right) & \mathrm{\dt}\,\mu \,\left(\gamma-1\right)+1 \end{bmatrix}.
$$
The scheme is stable if $|\mathrm{eig}(A^{-1}B)|\leq1$. Noting that $\beta=1/4$ and $\gamma=1/2$ for the NB2 scheme, the eigenvalues of $A^{-1}B$ if $\mu=0$ are found to be
  \begin{align*}
    \lambda_1=-\frac{\frac{3\,Q\,{\mathrm{\dt}}^2}{2}-2\,\mathrm{\dt}\,\sqrt{\frac{9\,Q^2\,{\mathrm{\dt}}^2}{16}+4\,Q}+4}{4\,\left(\frac{Q\,{\mathrm{\dt}}^2}{4}-1\right)}
    \quad\text{and}\quad
\lambda_2=-\frac{2\,\mathrm{\dt}\,\sqrt{\frac{9\,Q^2\,{\mathrm{\dt}}^2}{16}+4\,Q}+\frac{3\,Q\,{\mathrm{\dt}}^2}{2}+4}{4\,\left(\frac{Q\,{\mathrm{\dt}}^2}{4}-1\right)}.
  \end{align*}
  We can see  that  the NB2 scheme is nondissipative since   $|\lambda_1|=|\lambda_2|=1$  if there is no perturbation and no artificial  dissipation added, i.e., $Q=\hat{Q}_M<0$ and $\mu=0$.
As a result,   any perturbation caused by interpolation  introduces weak instability to    the scheme (see the next section for more details on weak instability). Given  sufficient artificial dissipation to suppress the weak instability,  the NB2 scheme is implicit in time and stable for any time step.  For accuracy reasons, we choose its time step based on the condition for the explicit PC22  scheme \eqref{eq:pc22dt}, but with a much larger stability factor. Typically, we choose $\csf\in[1,50]$  in applications.

\medskip
For convenience, we summarize  the time step and dissipation  formulas in Table~\ref{tab:timeStep}. Except for the implicit NB2 scheme, all the other schemes are explicit. We normally set $\csf=0.9$ for the explicit schemes.   From our experience,  it is generally sufficient to set $\df = 0.1$  to damp out the weak instability caused  by the inclusion of interpolating equations in the discretization on composite overlapping grids.

\begin{table}[h]\tableFont 
\caption{Summary of the  time step and dissipation  results.} \label{tab:timeStep}
\vspace{-0.2in}
\begin{center}
\rowcolors{1}{white}{gray!20}
\newcolumntype{M}[1]{>{\raggedright\arraybackslash}p{#1}}
\newcolumntype{P}[1]{>{\centering\arraybackslash}p{#1}}
\begin{tabular}{P{2cm}M{8cm}m{5cm}}
\hline
Name & \multicolumn{1}{c}{Formula} &  \multicolumn{1}{c}{Comments}\\\hline
{C2} & $\dt=\left\{
     \begin{aligned}
        &\csf\frac{2}{\sqrt{-\hat{Q}_M}}~~ \text{if}~~ \mu=0 \\
        &\csf\min\left\{\frac{4\sqrt{2}}{\mu+\sqrt{\mu^2-8\sqrt{2}\hat{Q}_M}}, \frac{2}{\mu}\right\}~~ \text{if}~~ \mu\neq 0 
     \end{aligned}
     \right.$ & Without dissipation, the scheme is nondissipative and works only for the single grid. With dissipation, the scheme has a reduced time-step. We choose  $\csf\in(0,1]$.\\
UPC2&$\displaystyle\dt =\csf\min\left\{\frac{2}{\sqrt{-\hat{Q}_M}} ,\frac{2}{\mu}\right\}$ & The scheme has no reduction in time-step  with artificial dissipation. We choose  $\csf\in(0,1]$ \\
PC22&$\displaystyle \dt = \csf\left( \left|\frac{\mu}{2r_a}\right|^\frac{3}{2}+\left|\frac{\sqrt{-\mu^2-4\hat{Q}_M}}{2r_b}\right|^\frac{3}{2}\right)^{-\frac{2}{3}}$ & The radii of the super ellipse are $r_a=1.75, r_b=1.2.$. We choose  $\csf\in(0,1]$.\\
NB2&$\displaystyle \dt = \csf\left( \left|\frac{\mu}{2r_a}\right|^\frac{3}{2}+\left|\frac{\sqrt{-\mu^2-4\hat{Q}_M}}{2r_b}\right|^\frac{3}{2}\right)^{-\frac{2}{3}}$   & The time step is chosen to be the same as the PC22 scheme but with a larger $\csf$; typically, we choose   $\csf\in[1,50]$ for accuracy reasons. \\\hline
Dissipation  & $\displaystyle\mu=32\ad~~\text{with}~~\ad = \df\frac{\sqrt{-\hat{Q}_M}}{16}$ & $\df\in[0,1]$ is the dissipation factor used to control the strength of the artificial dissipation. \\
\hline
\end{tabular}
\end{center}
\end{table}

\section{Model Problem in 1D}\label{sec:model1d}
A simplified model plate that consists of bending only  in 1D  is considered here  to demonstrate the numerical properties of the proposed schemes. On a simple composite overlapping grid,  we  illustrate the~\emph{weak instability} caused by the presence of interpolating equations, and investigate the   stabilizing effects of the various strategies of incorporating   artificial dissipation to the system.

\subsection{Problem setup}
Specifically, the   model equation reads 
\begin{equation}\label{eq:modelEqn}
  \pdn{w}{t}{2}= - \pdn{w}{x}{4}, \quad x\in[0,1].
\end{equation}
and the corresponding boundary conditions at the end points ($x=0$, $x=1$) are given by
\begin{alignat}{3}
  \text{clamped}: &\quad w=0, && \quad\pd{w}{x}=0,\label{eq:modelClampedBC}\\ 
  \text{supported}: &\quad w=0, && \quad \pdn{w}{x}{2}=0,\label{eq:modelSupportedBC}\\
      \text{free}: &\quad \pdn{w}{x}{2}=0, && \quad \pdn{w}{x}{3}=0.\label{eq:modelFreeBC}
\end{alignat}
Note that in 1D domain, the \KL  plate is the same as the well-known Euler-Bernoulli beam model.




To demonstrate the~\emph{weak instability} due to the presence of interpolating equations  in the discretized system, we consider the model equation \eqref{eq:modelEqn} subject to  the supported boundary condition \eqref{eq:modelSupportedBC} on both ends of the domain  as an example, and discretize the IBVP  using  both  a single  grid ($\Gc_s$) and  a composite overlapping grid  ($\Gc_c=\Gc^{(1)}\cup\Gc^{(2)}$) for comparison;   a diagram for the mesh setup is shown  in Figure~\ref{fig:MPgrids}. For the composite grid,   interpolation points of each  component  grid  are identified using  solid-square markers, and are connected with their donor points from  the other component grid. Note that, for  consistency  with the finite-difference  discretization,  a five-point stencil is used for  interpolation (i.e., 4th-order polynomial interpolation). We also investigate the stabilizing effects in all the  algorithms introduced in Section~\ref{sec:algorithm}  using this model problem. Thus,  the artificial dissipation term is  incorporated into the simplified model problem  \eqref{eq:modelEqn}, and are spatially discretized on the grids $\Gc_s$ and $\Gc_c$ as well.

{ 
  \newcommand{\figWidth}{16cm}
\def\xa{16}
\def\ya{7}
\newcommand{\trimfig}[2]{\trimw{#1}{#2}{0.}{0.}{0.}{0.}}
\begin{figure}[h]
  \begin{center}
\resizebox{12cm}{!}
{
\begin{tikzpicture}[scale=1]
  \useasboundingbox (0.0,0.0) rectangle (\xa,\ya);  

\draw(-0.5,-0.2) node[anchor=south west,xshift=0pt,yshift=0pt] {\trimfig{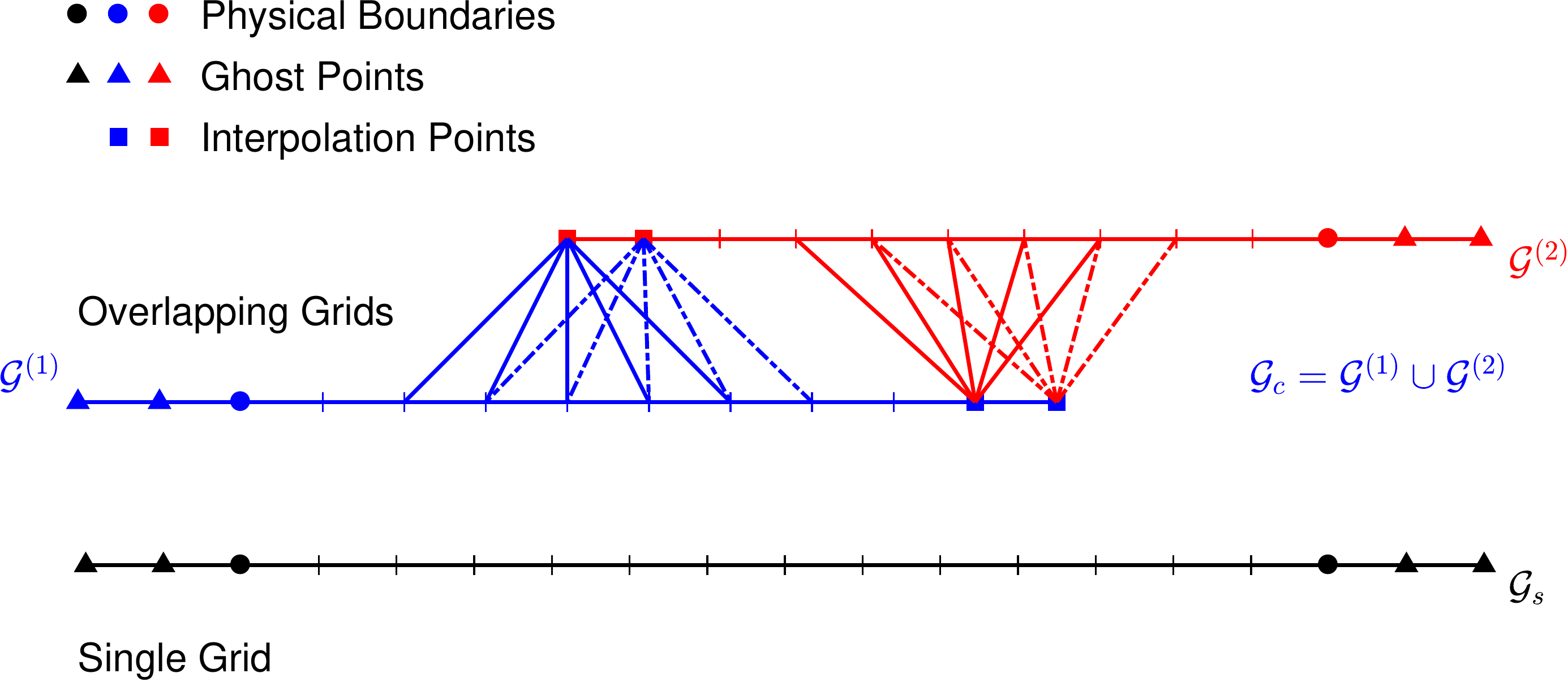}{\figWidth}};

\def\xstart{.9}
\def\ystart{1.1} 
\def\dy{1}
\def\dx{.79}

\draw({\xstart+2*\dx},\ystart)  node[anchor=north,xshift=0pt,yshift=-3pt] { $0$};
\draw({\xstart+16*\dx},\ystart)  node[anchor=north,xshift=2pt,yshift=-3pt] { $1$};

\draw(\xstart,\ystart)  node[anchor=south , xshift=3pt,yshift=3pt] {\footnotesize $x_{-2}$};
\draw({\xstart+\dx},\ystart)  node[anchor=south,xshift=3pt,yshift=3pt] {\footnotesize $x_{-1}$};
\draw({\xstart+2*\dx},\ystart)  node[anchor=south,xshift=0pt,yshift=3pt] {\footnotesize $x_{0}$};
\draw({\xstart+3*\dx},\ystart)  node[anchor=south,xshift=0pt,yshift=3pt] {\footnotesize $x_{1}$};
\draw({\xstart+4*\dx},\ystart)  node[anchor=south,xshift=0pt,yshift=3pt] {\footnotesize $x_{2}$};
\draw({\xstart+6*\dx},\ystart)  node[anchor=south,xshift=0pt,yshift=3pt] {\footnotesize $\cdots\cdots$};

\draw({\xstart+8*\dx},\ystart)  node[anchor=south,xshift=0pt,yshift=3pt] {\footnotesize $x_{i-1}$};
\draw({\xstart+9*\dx},\ystart)  node[anchor=south,xshift=0pt,yshift=3pt] {\footnotesize $x_i$};
\draw({\xstart+10*\dx},\ystart)  node[anchor=south,xshift=0pt,yshift=3pt] {\footnotesize $x_{i+1}$};

\draw({\xstart+12*\dx},\ystart)  node[anchor=south,xshift=0pt,yshift=3pt] {\footnotesize $\cdots\cdots$};
\draw({\xstart+14*\dx},\ystart)  node[anchor=south,xshift=6pt,yshift=3pt] {\footnotesize $x_{N-2}$};
\draw({\xstart+15*\dx},\ystart)  node[anchor=south,xshift=6pt,yshift=3pt] {\footnotesize $x_{N-1}$};
\draw({\xstart+16*\dx},\ystart)  node[anchor=south,xshift=2pt,yshift=3pt] {\footnotesize $x_{N}$};
\draw({\xstart+17*\dx},\ystart)  node[anchor=south,xshift=4pt,yshift=3pt] {\footnotesize $x_{N+1}$};
\draw({\xstart+18*\dx},\ystart)  node[anchor=south,xshift=7pt,yshift=3pt] {\footnotesize $x_{N+2}$};

\draw(\xstart,{\ystart+\dy})  node[anchor=south , xshift=0pt,yshift=0pt,blue] {\footnotesize $x^{(1)}_{-2}$};
\draw({\xstart+\dx},{\ystart+\dy})  node[anchor=south,xshift=2pt,yshift=0pt,blue] {\footnotesize $x^{(1)}_{-1}$};
\draw({\xstart+2*\dx},{\ystart+\dy})  node[anchor=south,xshift=0pt,yshift=0pt,blue] {\footnotesize $x^{(1)}_{0}$};
\draw({\xstart+3*\dx},{\ystart+\dy})  node[anchor=south,xshift=2pt,yshift=0pt,blue] {\footnotesize $x^{(1)}_{1}$};
\draw({\xstart+4*\dx},{\ystart+\dy})  node[anchor=south,xshift=2pt,yshift=0pt,blue] {\footnotesize $x^{(1)}_{2}$};
\draw({\xstart+6*\dx},{\ystart+\dy})  node[anchor=south,xshift=0pt,yshift=0pt,blue] {\footnotesize $\cdots\cdots$};
\draw({\xstart+8*\dx},{\ystart+\dy})  node[anchor=south,xshift=0pt,yshift=0pt,blue] {\footnotesize $\cdots\cdots$};
\draw({\xstart+10*\dx},{\ystart+\dy})  node[anchor=south,xshift=0.45cm,yshift=0pt,blue] {\footnotesize $x^{(1)}_{N_1}$};
\draw({\xstart+11*\dx},{\ystart+\dy})  node[anchor=south,xshift=0.5cm,yshift=0pt,blue] {\footnotesize $x^{(1)}_{N_1+1}$};
\draw({\xstart+12*\dx},{\ystart+\dy})  node[anchor=south,xshift=0.75cm,yshift=0pt,blue] {\footnotesize $x^{(1)}_{N_1+2}$};

\def\xstart{5.9}
\def\ystart{4.5} 
\def\dy{0.05}
\def\dx{.79}

\draw(\xstart,{\ystart+\dy})  node[anchor=south , xshift=0pt,yshift=0pt,red] {\footnotesize $x^{(2)}_{-2}$};
\draw({\xstart+\dx},{\ystart+\dy})  node[anchor=south,xshift=0pt,yshift=0pt,red] {\footnotesize $x^{(2)}_{-1}$};
\draw({\xstart+2*\dx},{\ystart+\dy})  node[anchor=south,xshift=0pt,yshift=0pt,red] {\footnotesize $x^{(2)}_{0}$};
\draw({\xstart+3*\dx},{\ystart+\dy})  node[anchor=south,xshift=0pt,yshift=0pt,red] {\footnotesize $x^{(2)}_{1}$};
\draw({\xstart+4*\dx},{\ystart+\dy})  node[anchor=south,xshift=0pt,yshift=0pt,red] {\footnotesize $x^{(2)}_{2}$};
\draw({\xstart+6*\dx},{\ystart+\dy})  node[anchor=south,xshift=-0.2cm,yshift=0pt,red] {\footnotesize $\cdots\cdots$};
\draw({\xstart+8*\dx},{\ystart+\dy})  node[anchor=south,xshift=-0.2cm,yshift=0pt,red] {\footnotesize $\cdots\cdots$};
\draw({\xstart+10*\dx},{\ystart+\dy})  node[anchor=south,xshift=-0.1cm,yshift=0pt,red] {\footnotesize $x^{(2)}_{N_2}$};
\draw({\xstart+11*\dx},{\ystart+\dy})  node[anchor=south,xshift=0.cm,yshift=0pt,red] {\footnotesize $x^{(2)}_{N_2+1}$};
\draw({\xstart+12*\dx},{\ystart+\dy})  node[anchor=south,xshift=0.15cm,yshift=0pt,red] {\footnotesize $x^{(2)}_{N_2+2}$};

%
\end{tikzpicture}
} 
  \end{center}
  \vspace{-0.2in}
  \caption{The single and composite overlapping grids used for discretizing  the model problem and its boundary conditions.}\label{fig:MPgrids} 
\end{figure}
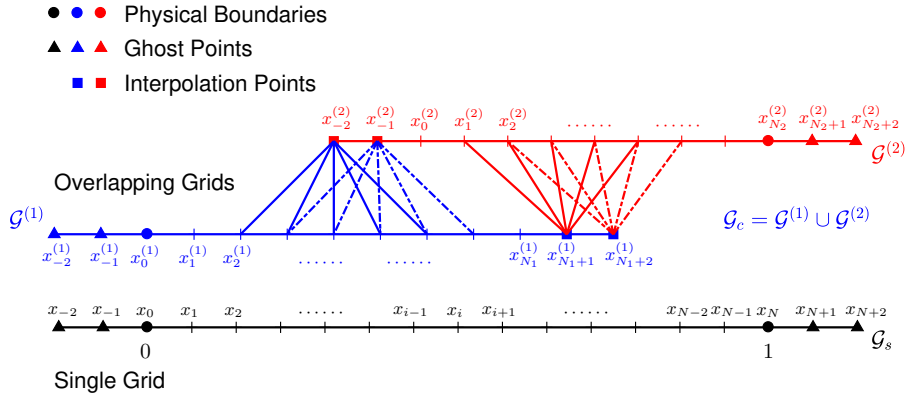
}

Let $W_i(t)$ be the approximation of $w(x_i,t)$, then   the discrete system on $\Gc_s$ is given by
\begin{align*}
  \ddn{W_i}{t}{2}= - (\Dpx\Dmx)^2{W_i}-\ad h_x^4(\Dpx\Dmx)^2\dd{W_i}{t},  &\quad  i = 1,2,\dots, N-1; ~ j =1,2, \\
  W_{i_b}=\Dpx\Dmx W_{i_b}=0,   &\quad i_b=0, N,\\
  D^5_{\pm x}W_{i_g}=0,    &\quad i_g=-2, N+2. 
\end{align*}
The last equation corresponds to the fifth-order extrapolation for the second ghost line, where the extrapolation operator is either  $\Dpx^5$ or $\Dmx^5$  so that it extrapolates into the interior of the grid.

On the composite overlapping grid $\Gc_c$ and denoting   $W_i^{(j)}\approx w\left(x^{(j)}_i,t\right)$,  the model problem is discretized by 
\begin{align*}
  \ddn{W^{(j)}_i}{t}{2}= - (\Dpx\Dmx)^2{W^{(j)}_i}-\ad h_x^4(\Dpx\Dmx)^2\dd{W^{(j)}_i}{t},  &\quad  i = 1,2,\dots, N_j-1; j=1,2, \\
  W^{(j)}_{i_b}=\Dpx\Dmx W^{(j)}_{i_b}=0,   &\quad i_b=0 ~ \text{if}~ j=1;~i_b=N_2~ \text{if}~ j=2,\\
  D^5_{\pm x}W^{(j)}_{i_g}=0,    &\quad i_g=-2 ~ \text{if}~ j=1;~i_g=N_2+2~ \text{if}~ j=2,\\
  W^{(1)}_{i_p}- \sum_{l=1}^5c_{2,l}(i_p)  W^{(2)}_{d_{2,l}(i_p) } = 0,  & \quad i_p=N_1+1, N_1+2, \\
  W^{(2)}_{i_p}- \sum_{l=1}^5c_{1,l}(i_p)   W^{(1)}_{d_{1,l}(i_p) } = 0,  & \quad i_p=-1,-2.
\end{align*}
 Here the last  two conditions are the interpolating equations that couple the numerical  solutions on the two component grids  together. In the interpolating equations, $d_{j,l}(i_p)$ for $l=1,\dots,5$ represent the indices of the donor points on $\Gc^{(j)}$ for interpolating the grid function  at  $x_{i_p}^{(k)}$ on the other component grid $\Gc^{(k)}$;  the corresponding interpolating coefficients are denoted by  $c_{j,l}(i_p) $. The relationship between the interpolation points and their donor points are schematically depicted in Figure~\ref{fig:MPgrids}.

 If we denote  $\Wv$  as  the solution vector, where $\Wv=\left[W_{-2},\dots, W_{N+2} \right]^T$ for the single grid case and $\Wv=\left[W^{(1)}_{-2},\dots, W^{(1)}_{N_1+2},W^{(2)}_{-2},\dots, W^{(2)}_{N_2+2}  \right]^T$ for the composite grid case, then  the above spatially discretized systems, including the boundary and interpolation conditions,  can be concisely written into the following  matrix form,
 \begin{equation}
   M_I \ddn{\Wv}{t}{2}= M_{L}\Wv +\ad M_{\text{ad}}\dd{\Wv}{t}.
 \label{eq:modelEqnMatrixForm}
 \end{equation}
 Here  $M_I$ is identity matrix except for the rows corresponding to the ghost, boundary and interpolation points, which are set to be zero  to accommodate  the boundary and interpolation conditions. Following \eqref{eq:dissipationCoefficient}, the dissipation coefficient $\ad$ for this 1D  model problem is specifically given by
 \begin{equation}
  \ad=\df\frac{1}{2h_x^2}. 
 \label{eq:adMP}
 \end{equation}

\subsection{Stability issues}
 For the purpose of demonstrating the stability issues of the various  time-stepping methods introduced in Section~\ref{sec:tsSchemes}, we consider the spatially discretized system \eqref{eq:modelEqnMatrixForm}  obtained from  the   one-dimensional mesh setup  shown in Figure~\ref{fig:MPgrids} with $N_1=15$, $N_2=9$ and $N=20$.

 \medskip
\noindent {\bf C2 and UPC2 schemes.}

The matrix equation \eqref{eq:modelEqnMatrixForm}, if solved with the C2 or UPC2 scheme, leads to  the following  general updating formula,
\begin{equation}
Q_2\Wv^{n+1}+Q_1\Wv^{n}+Q_0\Wv^{n-1}=\zerov,  \label{eq:updatingFormulaMP}
\end{equation}
where $\Wv^{n}$ is the  numerical solution at time level $t_n$.
The coefficient matrices in the updating formula \eqref{eq:updatingFormulaMP} for the C2 scheme are  given by 
\begin{equation*}
Q_2=M_I,\quad Q_1=-2M_I-\dt^2M_L-\dt \ad M_{\text{ad}}, \quad Q_0=M_I+\dt \ad M_{\text{ad}},  \label{eq:QsForC2}
\end{equation*}
while  those for the UPC2 scheme are given by 
\begin{equation*}
Q_2=M_I,\quad Q_1=-\left(M_I+\frac{\dt}{2}\ad M_\text{ad}\right)\left(2M_I+\dt^2M_L\right), \quad Q_0=M_I+\dt \ad M_{\text{ad}}. \label{eq:QsForUPC2}
\end{equation*}
Note that   C2 and UPC2 schemes are the same scheme if no artificial dissipation is added  (i.e., $\ad=0$). 

The solution  of the difference  equation \eqref{eq:updatingFormulaMP} is related to the quadratic eigenvalue problem (also known as the characteristic equation),
\begin{equation}
Q_2\lambda^2+Q_1\lambda+Q_0=0,   \label{eq:quadraticEigMP}
\end{equation}
which is derived from seeking separable solutions of the form $\Wv^n=\lambda^n\Wv_0$.
The stability of the time-stepping scheme  requires that $|\lambda|\leq 1$ for all the eigenvalues of \eqref{eq:quadraticEigMP}.  Therefore,  in Figure~\ref{fig:polyEigsUPC2}, we plot 
the  numerically computed eigenvalues for the schemes with and without artificial dissipation on a complex plane together   with a unit circle  $|z|=1$ as reference; a scheme would be unstable if eigenvalues are observed  outside of the unit circle.

{ 
  \newcommand{\figWidth}{5cm}
  \newcommand{\figWidthLgd}{2.8cm}
\def\xa{16}
\def\ya{5.5}
\newcommand{\trimfig}[2]{\trimw{#1}{#2}{0.}{0.}{0.}{0.}}
\begin{figure}[h]
\begin{center}
\begin{tikzpicture}[scale=1]
  \useasboundingbox (0.0,0.0) rectangle (\xa,\ya);  

\draw(-0.5,-0.2) node[anchor=south west,xshift=0pt,yshift=0pt] {\trimfig{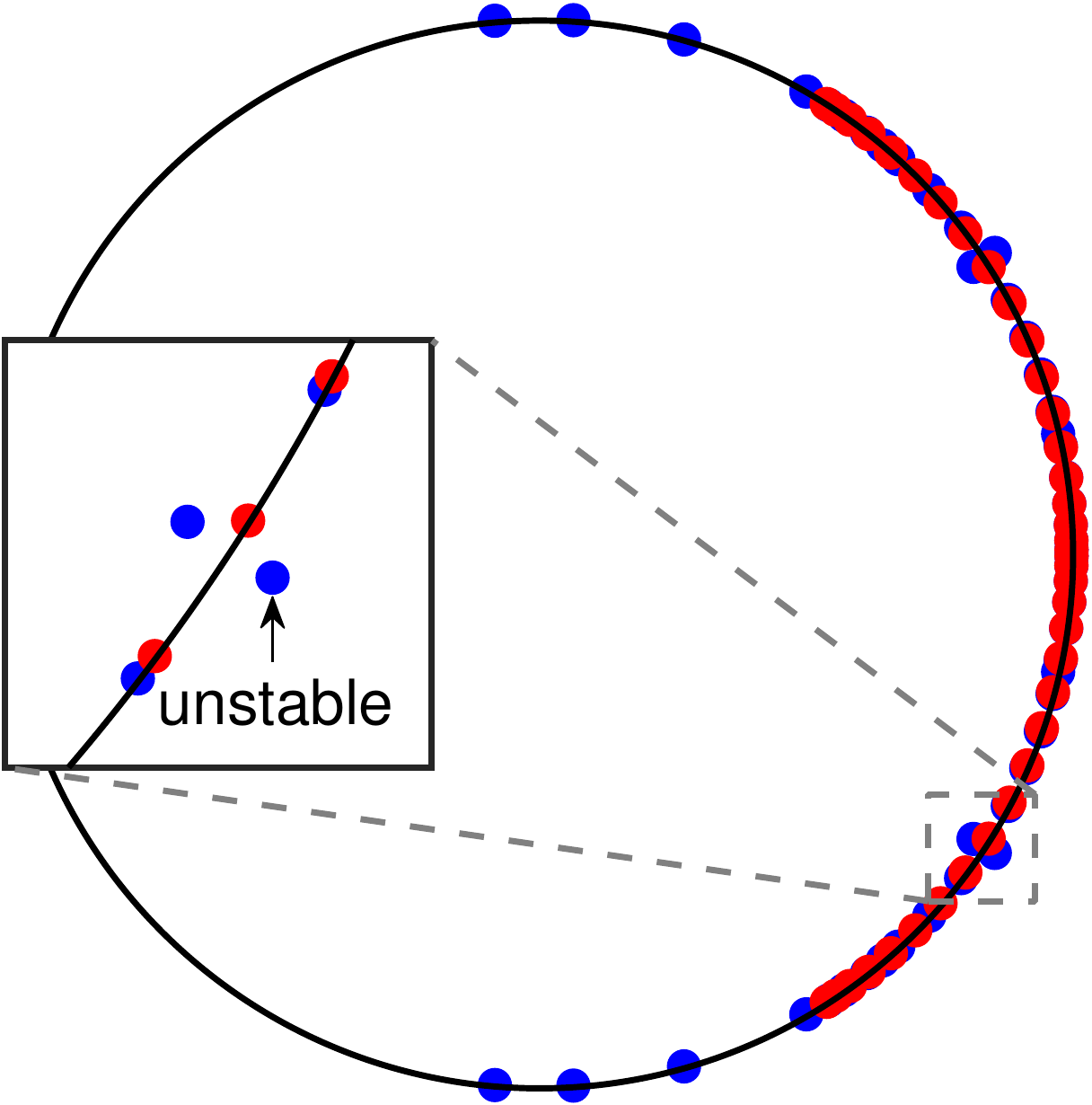}{\figWidth}};
\draw(5.,-0.2) node[anchor=south west,xshift=0pt,yshift=0pt] {\trimfig{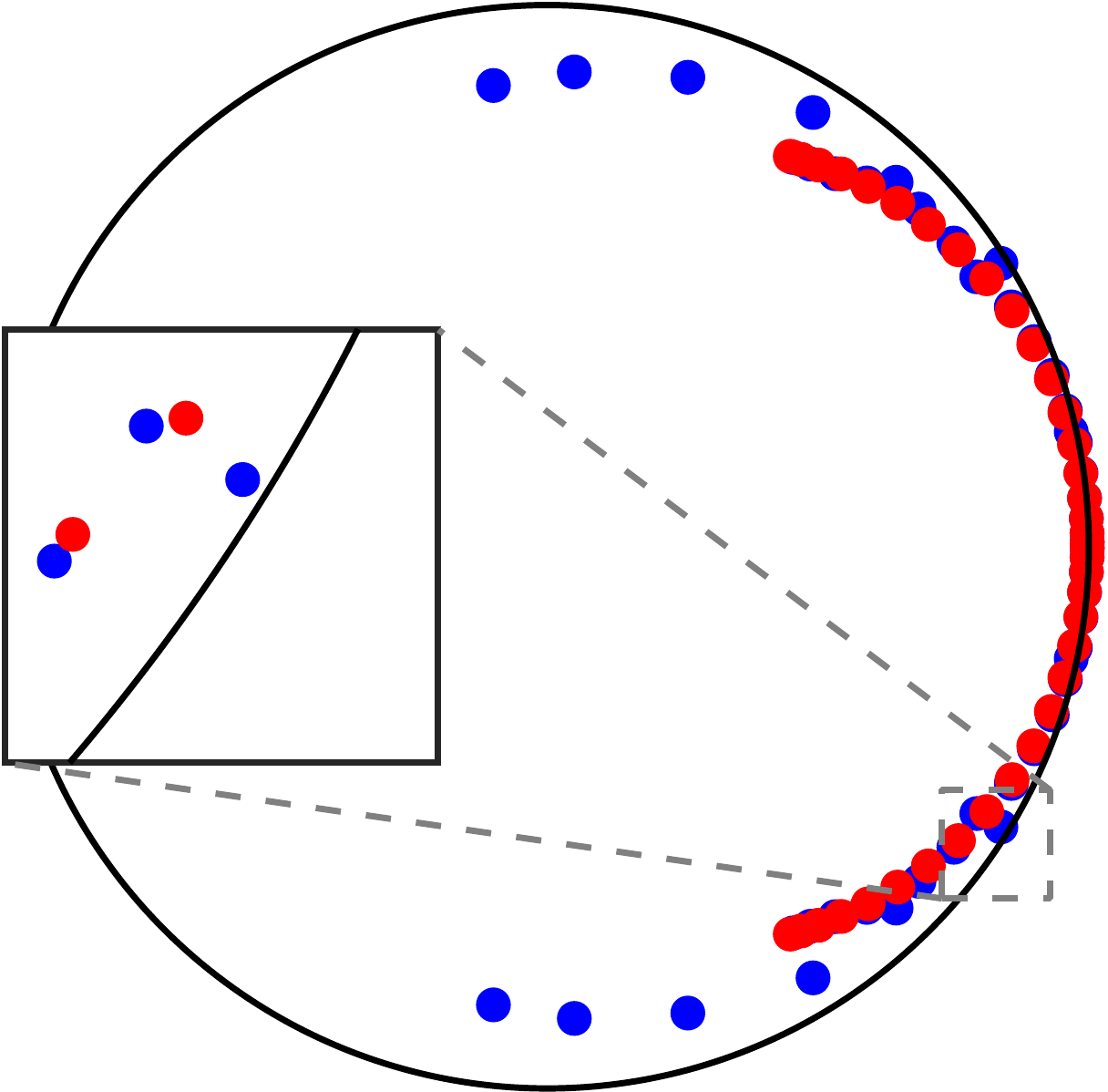}{\figWidth}};
\draw(10.5,-0.2) node[anchor=south west,xshift=0pt,yshift=0pt] {\trimfig{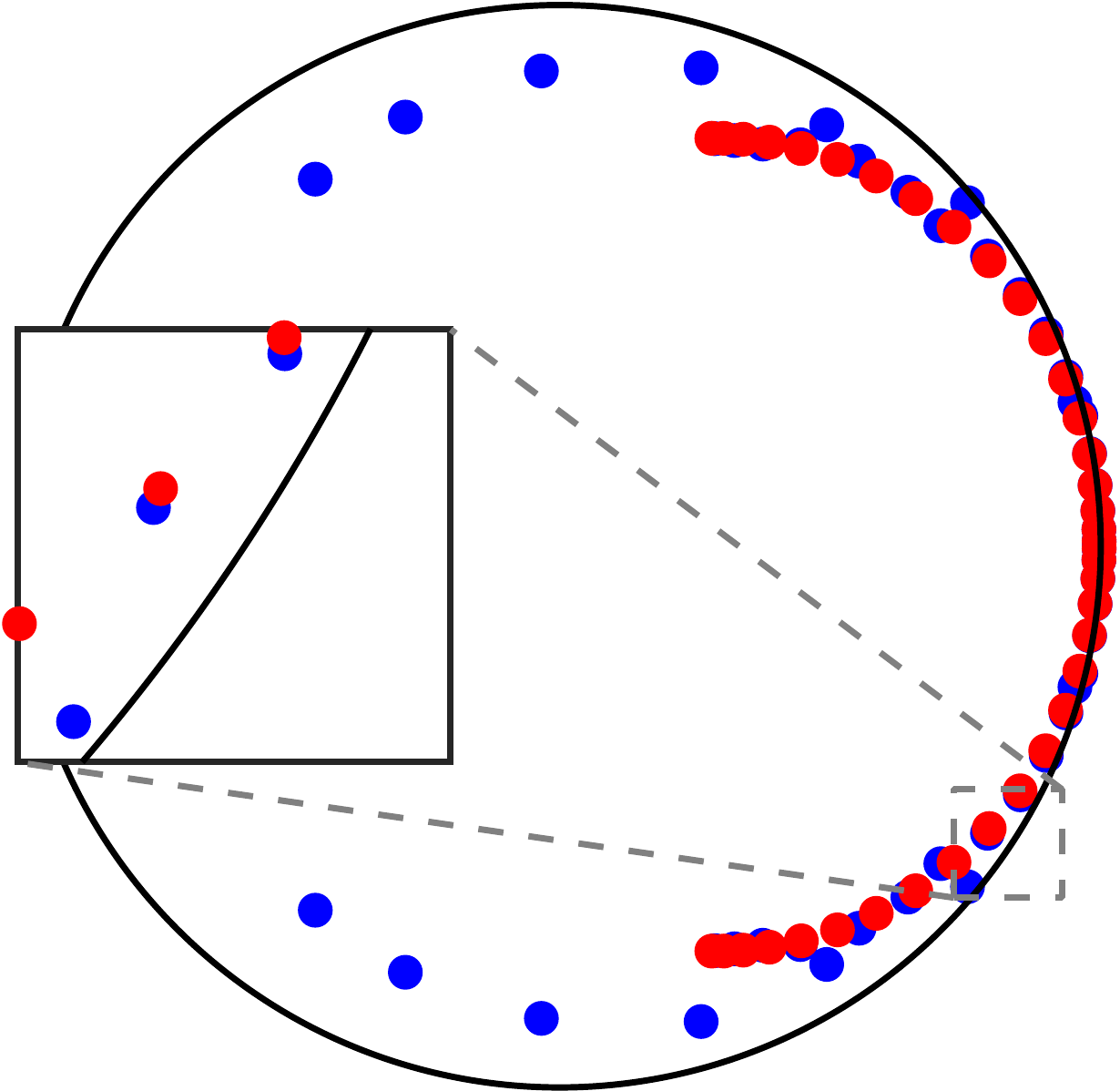}{\figWidth}};

   \draw(3.4,4.4) node[anchor=south west,xshift=0pt,yshift=0pt] {\trimfig{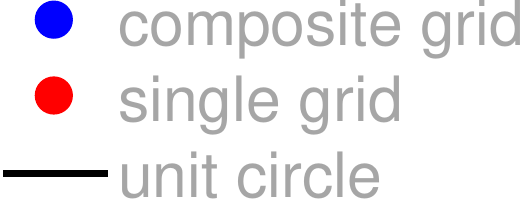}{\figWidthLgd}};

\draw(2.5,0)  node[anchor=north,xshift=0pt,yshift=0pt] {\footnotesize no dissipation };
\draw(8.5,0)  node[anchor=north,xshift=0pt,yshift=0pt] {\footnotesize  with dissipation};
\draw(13.5,0)  node[anchor=north,xshift=0pt,yshift=0pt] {\footnotesize  UPC2 with dissipation};

%
\end{tikzpicture}

\end{center}
\caption{Eigenvalues of    C2 scheme with no artificial dissipation (left), C2 scheme with artificial dissipation (middle), UPC2 scheme (right).} \label{fig:polyEigsUPC2}
\end{figure}
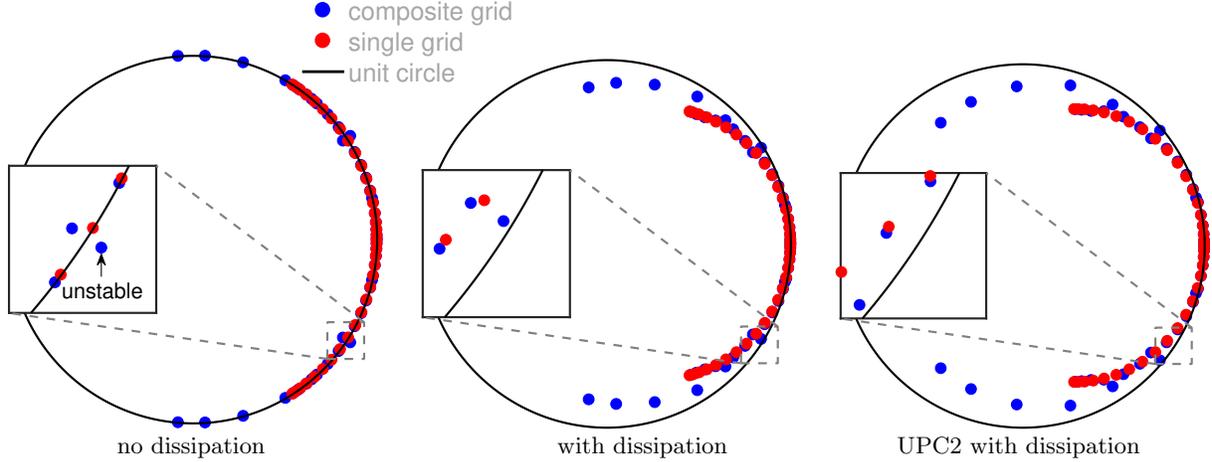
}

If no artificial dissipation is included in the scheme, the C2 scheme is non-dissipative and stable on the single grid $\Gc_s$; this is seen  in the left image of Figure~\ref{fig:polyEigsUPC2} that  all the eigenvalues for the single grid case remain on the unit circle.  Due to the perturbation caused by the interpolating equations, we can see in the same image  that, for the composite grid case,  there are complex conjugate pairs of eigenvalues  with one pair lying  inside the unit circle and the other outside that causes weak instability for the scheme.

The weak instability can be stabilized by including artificial dissipation in an explicit manner  as the C2 scheme described in Algorithm~\ref{alg:c2} or in an upwind predictor-corrector approach as the UPC2 scheme  described  in Algorithm~\ref{alg:upc2}. For both schemes, setting  the  dissipation $\ad$ in \eqref{eq:adMP}  with   $\df=0.1$ is sufficient for stabilization.  The eigenvalues for the C2 scheme with dissipation and the UPC2 scheme are shown in the middle and right images of  Figure~\ref{fig:polyEigsUPC2}, respectively. It is clear that  all the eigenvalues, especially the high frequency modes,  are  damped so that no more unstable modes are observed outside the unit circle.

\medskip
\noindent {\bf PC22 scheme.}

To solve \eqref{eq:modelEqnMatrixForm} with the PC22 time-stepping scheme, we need to convert the equation into the first-order form,
\begin{equation}\label{eq:modelEqn1stOrderForm}
\Ic\dd{\Uv}{t}=\Mc\Uv,
\end{equation}
where  
$$
\Uv=\begin{bmatrix}
\Wv\\
\Vv
 \end{bmatrix},\quad
  \Ic=\begin{bmatrix}
  M_I & \zerov\\
 \zerov  &  M_I 
  \end{bmatrix},\quad\text{and}\quad 
 \Mc= \begin{bmatrix}
    \zerov & M_v\\
    M_L & \ad M_\text{ad}
  \end{bmatrix}.
  $$
  Here $M_v$ equals the  identity  matrix except for the rows corresponding to the ghost, boundary and interpolation points, which are given by the boundary conditions of the velocity component. Note  that the boundary conditions for velocity are  derived from the  time derivative  of the displacement ones.

Integrating \eqref{eq:modelEqn1stOrderForm} using the PC22 time-stepping scheme leads to a three step updating formula,
$$
Q_2\Uv^{n+1}+Q_1\Uv^{n}+Q_0\Uv^{n-1}=\zerov,
$$
where
$$
Q_2=\Ic, \quad Q_1= -\Ic-\dt \Mc-\frac{3}{4}\dt^2 \Mc^2,\quad Q_0=\frac{1}{4}\dt^2 \Mc^2.
$$
The stability issues and the stabilizing effects of this scheme are then investigated by solving the corresponding  eigenvalue problem with and without artificial dissipation; the results are shown in  Figure~\ref{fig:polyEigsPC22}.

From the results for  $\ad=0$ that are shown in the  left image of   Figure~\ref{fig:polyEigsPC22},  we can see  that  PC22  is a dissipative scheme, and it is stable for solving the model problem on the single grid. However, the dissipation of the scheme itself is not enough  to stabilize the algorithm applied  on the  composite overlapping grid  $\Gc_c$;  extra artificial dissipation are essential for  stability. With $\df=0.1$, it is then observed in the right image of  Figure~\ref{fig:polyEigsPC22}  that the unstable modes are further damped by the artificial dissipation such that the corresponding eigenvalues fall inside of the unit circle.

{ 
\newcommand{\figWidth}{5cm}
\def\xa{11}
\def\ya{5.5}
\newcommand{\trimfig}[2]{\trimw{#1}{#2}{0.}{0.}{0.}{0.}}
  \newcommand{\figWidthLgd}{3cm}

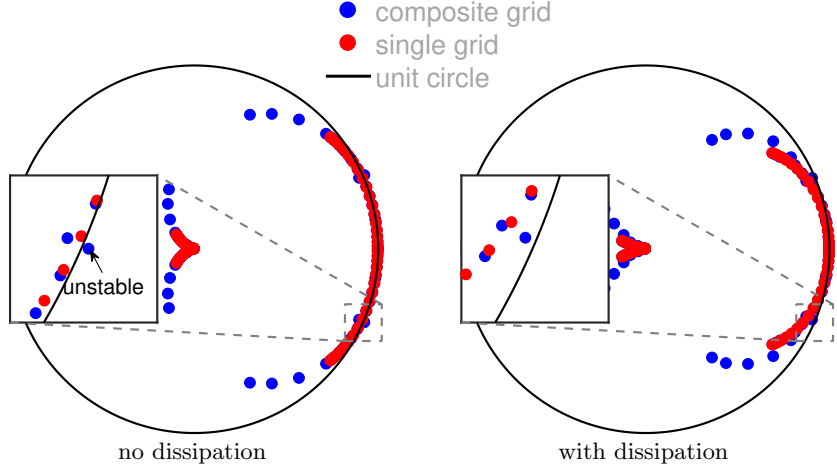
\begin{figure}[h]
\begin{center}
\begin{tikzpicture}[scale=1]
  \useasboundingbox (0.0,0.0) rectangle (\xa,\ya);  

\draw(-0.5,-0.2) node[anchor=south west,xshift=0pt,yshift=0pt] {\trimfig{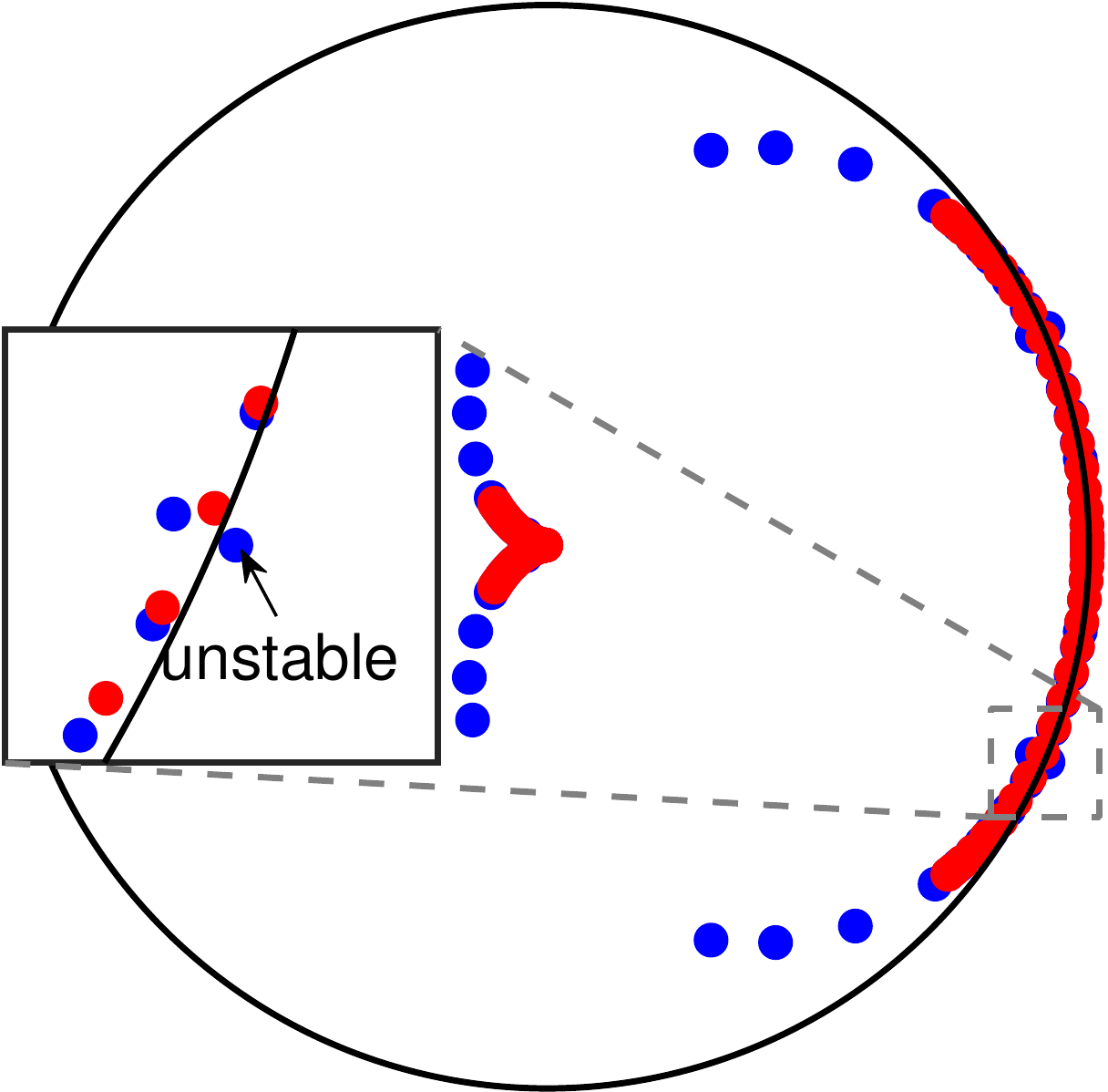}{\figWidth}};
\draw(5.5,-0.2) node[anchor=south west,xshift=0pt,yshift=0pt] {\trimfig{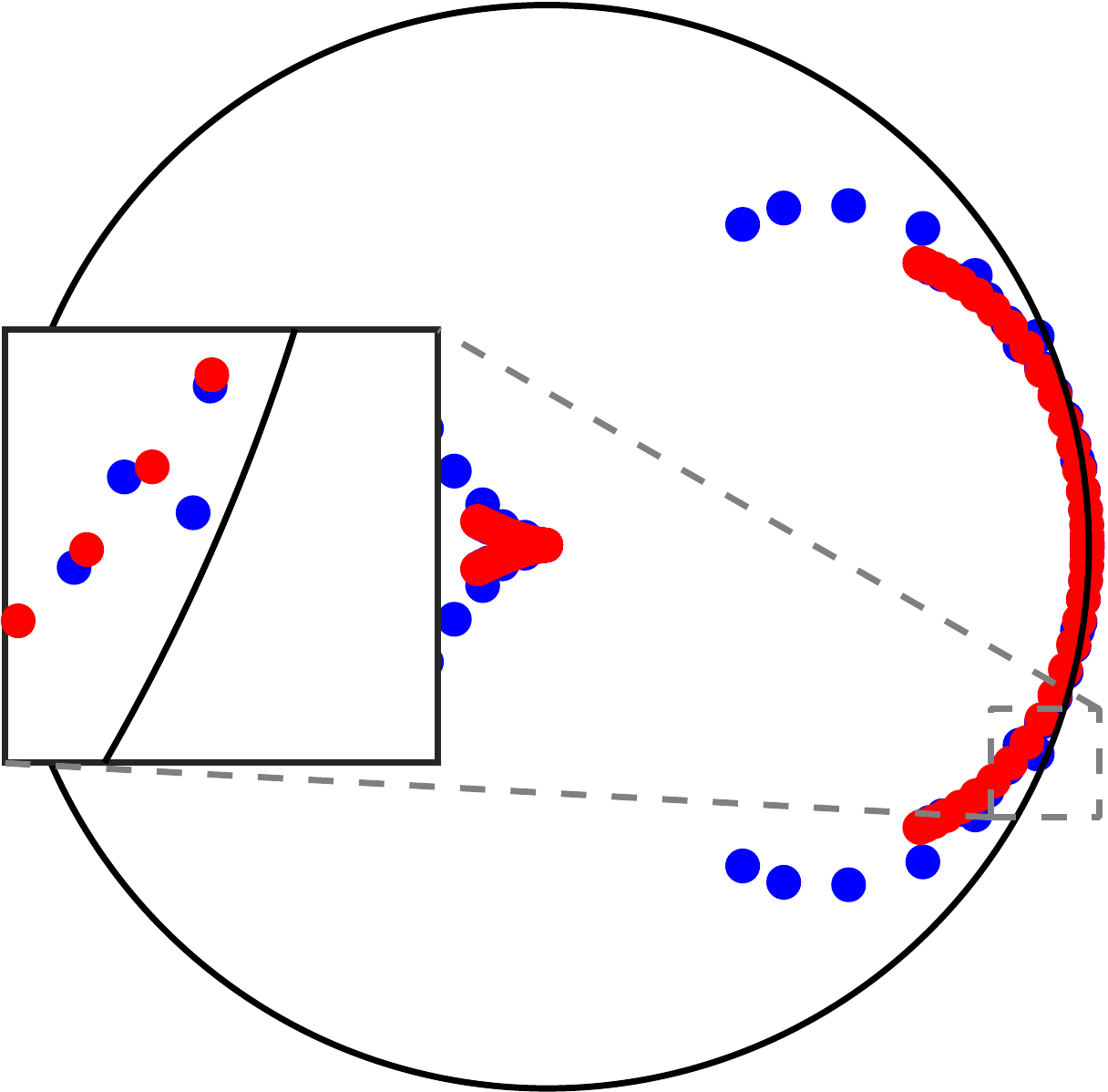}{\figWidth}};

   \draw(3.7,4.4) node[anchor=south west,xshift=0pt,yshift=0pt] {\trimfig{fig/PolyEigLegend}{\figWidthLgd}};

\draw(2.5,0)  node[anchor=north,xshift=0pt,yshift=0pt] {\footnotesize no dissipation };
\draw(8.5,0)  node[anchor=north,xshift=0pt,yshift=0pt] {\footnotesize  with dissipation};

%
\end{tikzpicture}

\end{center}
\caption{Eigenvalues of  PC22 scheme.} \label{fig:polyEigsPC22}
\end{figure}
}

\medskip
\noindent {\bf NB2 scheme. }

The NB2 scheme also solves the model problem in the first-order form \eqref{eq:modelEqn1stOrderForm}. It is a two-step scheme and has the following updating formula,
$$
Q_2\Uv^{n+1}+Q_1  \Uv^{n}=\zerov,
$$
where
$${\footnotesize
Q_2=
\begin{bmatrix}
M_I-\dt^2\beta M_L &  -\dt^2\beta \ad M_\text{ad}\\
 -\dt\gamma M_L& M_I-\dt\gamma\ad M_\text{ad}
\end{bmatrix}
\quad\text{and}\quad
Q_1=-\begin{bmatrix}
M_I+\frac{\dt^2}{2}(1-2\beta)M_L &  \dt M_I+\frac{\dt^2}{2}(1-2\beta)\ad M_\text{ad}\\
          \dt(1-\gamma)M_L &  M_I+\dt (1-\gamma)\ad M_\text{ad}
\end{bmatrix}.
}
$$

{ 
  \newcommand{\figWidth}{5cm}
  \newcommand{\figWidthLgd}{3cm}
\def\xa{11}
\def\ya{5.5}
\newcommand{\trimfig}[2]{\trimw{#1}{#2}{0.}{0.}{0.}{0.}}
\begin{figure}[h]
\begin{center}
\begin{tikzpicture}[scale=1]
  \useasboundingbox (0.0,0.0) rectangle (\xa,\ya);  

  \draw(-0.5,-0.4) node[anchor=south west,xshift=0pt,yshift=0pt] {\trimfig{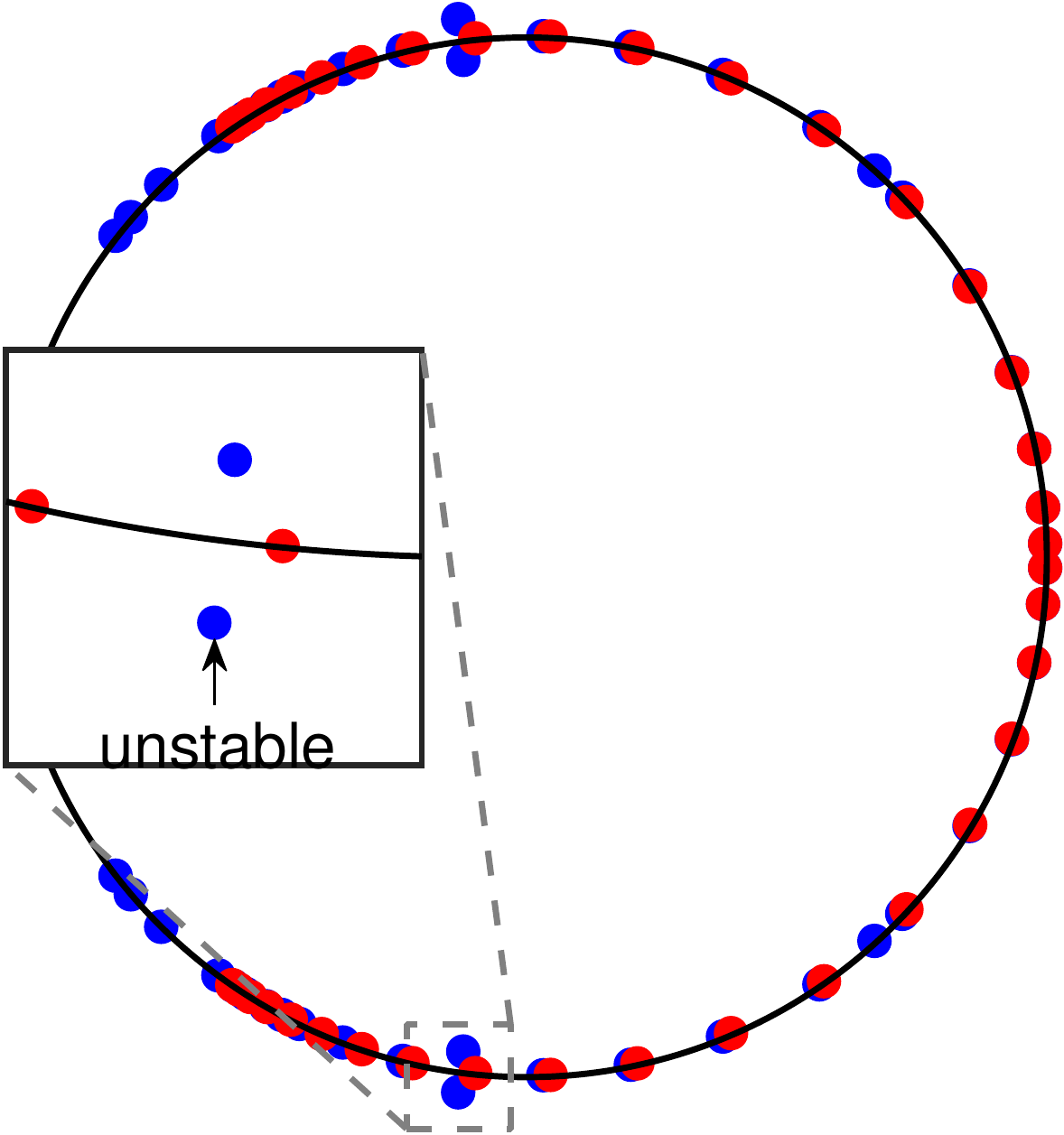}{\figWidth}};
  \draw(5.5,-0.4) node[anchor=south west,xshift=0pt,yshift=0pt] {\trimfig{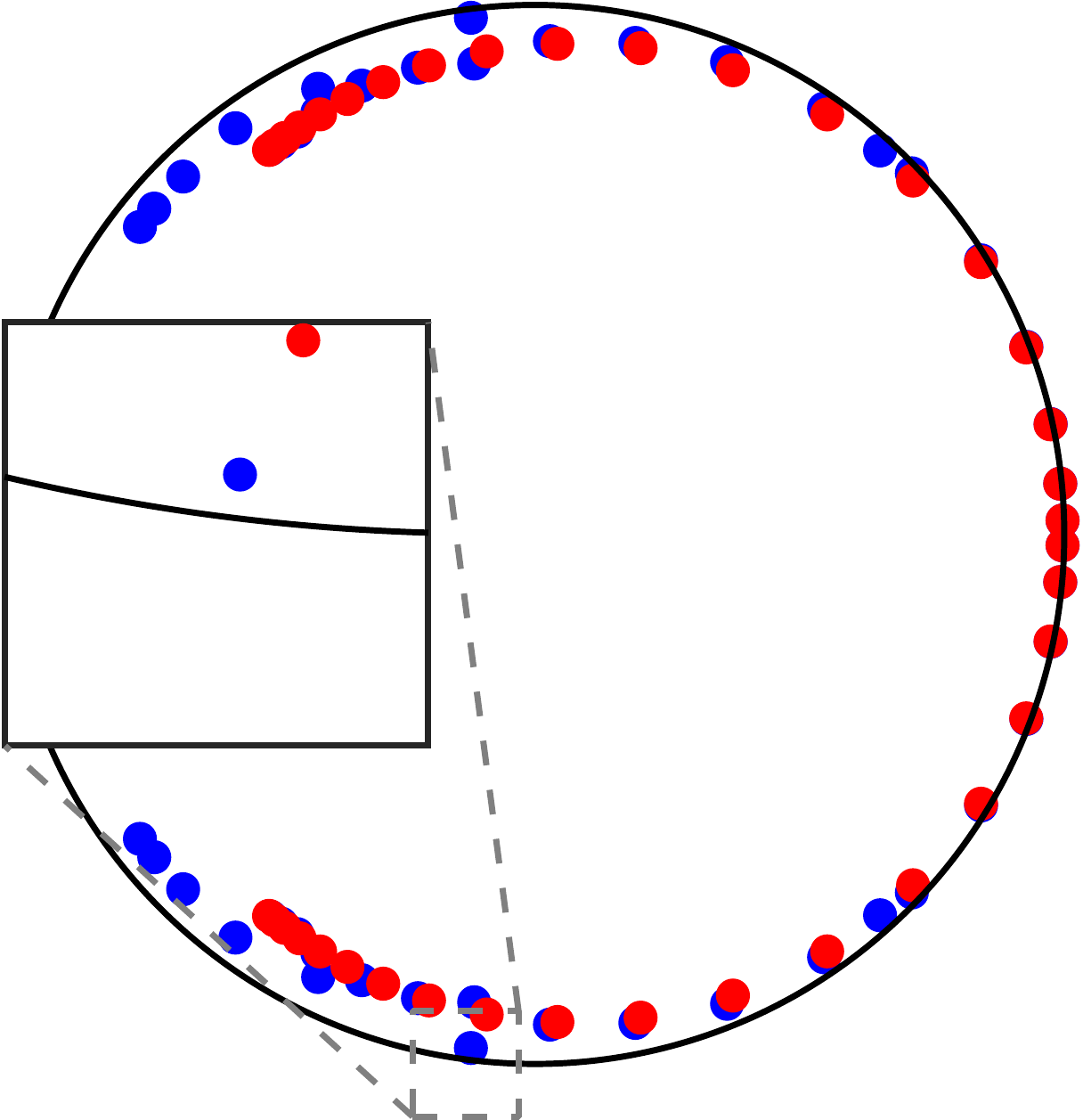}{\figWidth}};

   \draw(3.7,4.4) node[anchor=south west,xshift=0pt,yshift=0pt] {\trimfig{fig/PolyEigLegend}{\figWidthLgd}};

\draw(2.5,-0.2)  node[anchor=north,xshift=0pt,yshift=0pt] {\footnotesize no dissipation };
\draw(8.5,-0.2)  node[anchor=north,xshift=0pt,yshift=0pt] {\footnotesize  with dissipation};

%
\end{tikzpicture}

\end{center}
\caption{Eigenvalues of  NB2 scheme.} \label{fig:polyEigsNB2}
\end{figure}
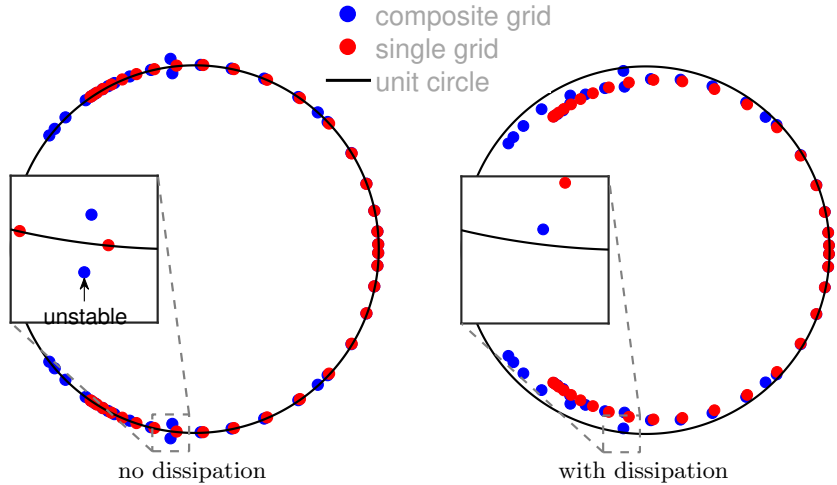
}

Computed eigenvalues for the NB2 scheme are plotted in Figure~\ref{fig:polyEigsNB2}. We can see from  the plot for $\ad=0$ (left image) that NB2 scheme is nondissipative and is stable for the single grid $\Gc_s$.  The presence of interpolation in the discretization formulas  on the composite  grid $\Gc_s$ perturbs the eigenvalues so that some modes become  unstable.  In the right image of Figure~\ref{fig:polyEigsNB2}, we demonstrate that incorporating the artificial  dissipation with $\df=0.1$ is sufficient to stabilize these unstable modes.

 
\section{Numerical results}\label{sec:results}
Numerical results for a series of test problems designed to validate the properties and performances of the proposed algorithms are now presented. Mesh refinement studies   using methods of manufactured solutions and  problems with  known analytical solutions are first conducted to  verify  the  accuracy and stability of all the algorithms. Using  a test case of   a plate  with multiple holes, we then  demonstrate the efficiency and capability of our approach  for  solving problems  with complicated geometrical configurations.   Finally,  to illustrate that the proposed algorithms can be  applied  in  realistic applications,  we explore the responding vibrations of  a guitar soundboard subject  to an initial pulse. Results for comparing the performance of the proposed algorithms for solving the guitar problem are also presented.

\subsection{Convergence study}
In this section, we conduct  convergence study  for a square plate using the method of manufactured solutions (MMS),  a technique often used  to construct exact solutions  for numerical validations  \cite{Roache1998}. In particular,  the forcing term in  \eqref{eq:KLShell} is specified so that a chosen (manufactured) function becomes an exact solution to the forced equations.  The accuracy of our proposed schemes are also carefully verified using  a simple circular plate model, whose   analytical solutions are available  for   certain boundary conditions.

For all the test problems with known exact solutions, errors  of  the numerical solutions  are measured in the maximum norm. Mathematically, we define  the error of the displacement and its maximum norm   at  time level $t_n$ as
$$
  E(w)=w_e(\xv_\iv,t_n)-W^n_{\iv} \quad\text{and} \quad ||E(w)||=\max_{\xv_{\iv}\in\Gc}|E(w)|, 
$$
  where $w_e$ is the exact  displacement. Similar definitions for  the  velocity and the  acceleration solutions are given  accordingly. 
  
\subsubsection{Same order  {\em vs.} same stencil}

In the first  test, we aim to  settle the   $\delta$  value to be used  in formulas  \eqref{eq:1stDerivFD} and \eqref{eq:2ndDerivFD} for approximating the 1st- and 2nd-order derivatives involved in the IBVP of the \KL model.   For this purpose, we compare the accuracy of the numerical results obtained using   the same order ($\delta=0$) formulas with  the  same stencil ($\delta=1$) formulas. 
We  focus on this finite-difference  discretization issue by keeping the setup of this numerical test as simple as possible.

Specifically, we consider solving \eqref{eq:KLShell} with parameters $ \rho=1$, $ H=1$, $ K=0$, $ T=0$ and $ D=1$ for the 
 following manufactured   solution on a  unit square domain,
\begin{equation}\label{eq:MMSOrderVsStencil}
u_e(x,y,t)=\cos(\pi x)\cos(\pi t).
\end{equation}
To avoid  distractions from other  factors  involved in our algorithms such as interpolation on overlapping grids and   artificial dissipation, this simple manufactured solution problem is solved using the C2  scheme without any  artificial dissipation  on   a  single grid as is shown in the left image of  Figure~\ref{fig:orderVsStencilGridsAContours}.  

{
\begin{table}[h]\tableFont 
  \begin{center}
    \begin{tabular}{c|cccccc||cccccc}
      \hline
                   \multicolumn{7}{c}{ same order ($\delta=0$) }   &  \multicolumn{6}{c}{same stencil ($\delta=1$)}\\\hline
grid             &  $||E(w)|| $ & ratio &  $ ||E(v)|| $ & ratio &  $||E( a)|| $ & ratio & $||E( w)|| $ & ratio &  $||E( v )||$ & ratio &  $||E( a )||$ & ratio \\\hline 
$\Gc_1$      & \num{4.7}{-3} &      & \num{1.7}{-1} &      & \num{17.2}{0} &       &  \num{9.7}{-3} &      & \num{1.9}{-1} &      & \num{1.7}{0} &       \\ 
 $\Gc_2$     & \num{1.3}{-3} &  3.58 & \num{3.9}{-2} &  4.34 & \num{11.8}{0} &  1.46    &  \num{2.4}{-3} &  3.98 & \num{4.8}{-2} &  4.00 & \num{4.1}{-1} &  4.26  \\ 
 $\Gc_4$     & \num{3.1}{-4} &  4.23 & \num{8.0}{-3} &  4.84 & \num{11.6}{0} &  1.02   &  \num{6.1}{-4} &  4.00 & \num{1.2}{-2} &  4.02 & \num{1.0}{-1} &  3.87  \\ 
 $\Gc_8$     & \num{7.9}{-5} &  3.93 & \num{2.5}{-3} &  3.22 & \num{8.6}{0} &  1.35    & \num{1.5}{-4} &  3.99 & \num{3.0}{-3} &  3.99 & \num{2.6}{-2} &  4.01  \\ 
 $\Gc_{16}$  & \num{1.9}{-5} &  4.07 & \num{6.0}{-4} &  4.12 & \num{5.0}{0} &  1.70  &   \num{3.8}{-5} &  4.00 & \num{7.5}{-4} &  4.00 & \num{6.5}{-3} &  4.01  \\  \hline
    rate                    &  $1.99$       &      &  $2.02$       &      &  $0.40$       &        & $2.00$   &            &  $2.00$       &      &  $2.00$       &       \\ \hline 
\end{tabular}
\caption{ Maximum-norm errors at $t=0.1$ and estimated convergence rates using the  manufactured solution \eqref{eq:MMSOrderVsStencil}. The numbers in the ratio columns provide the ratio of the errors at the current grid  to that on the next coarser grid. Results shown here are obtained using the C2 scheme for the plate subject to free boundary conditions. 
    }\label{tab:orderVsStencilFree}
\end{center}
\end{table}
}

{
\newcommand{\figWidth}{5cm}
\def\xa{16.}
\def\ya{5.}
\newcommand{\trimfig}[2]{\trimw{#1}{#2}{0.12}{0.12}{0.12}{0.12}}
\begin{figure}[h]
\begin{center}
\begin{tikzpicture}[scale=1]
  \useasboundingbox (0.0,0.0) rectangle (\xa,\ya);  

\draw(-0.5,0.0) node[anchor=south west,xshift=0pt,yshift=0pt] {\trimfig{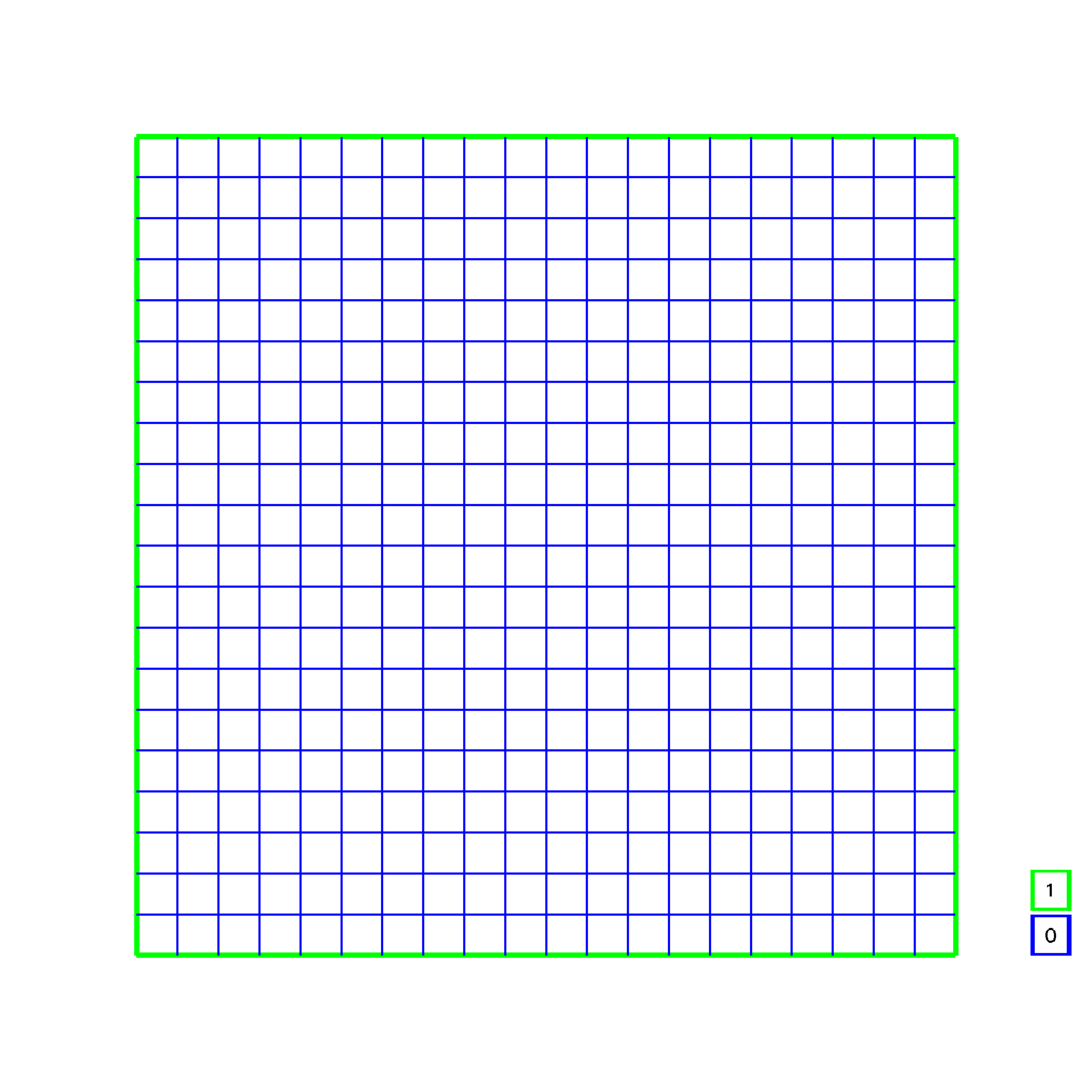}{\figWidth}};
\draw(5.,0.0) node[anchor=south west,xshift=0pt,yshift=0pt] {\trimfig{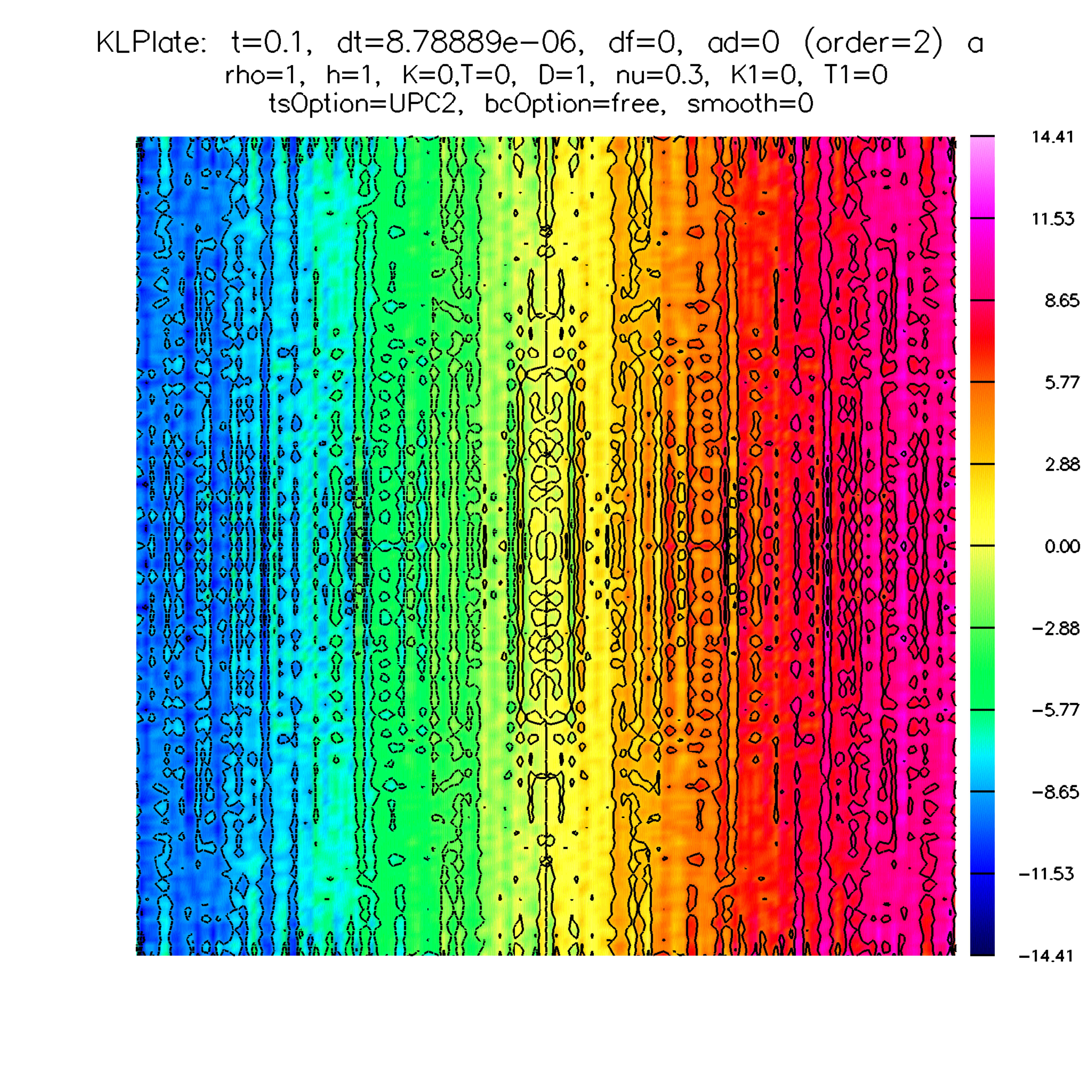}{\figWidth}};
\draw(10.5,0.0) node[anchor=south west,xshift=0pt,yshift=0pt] {\trimfig{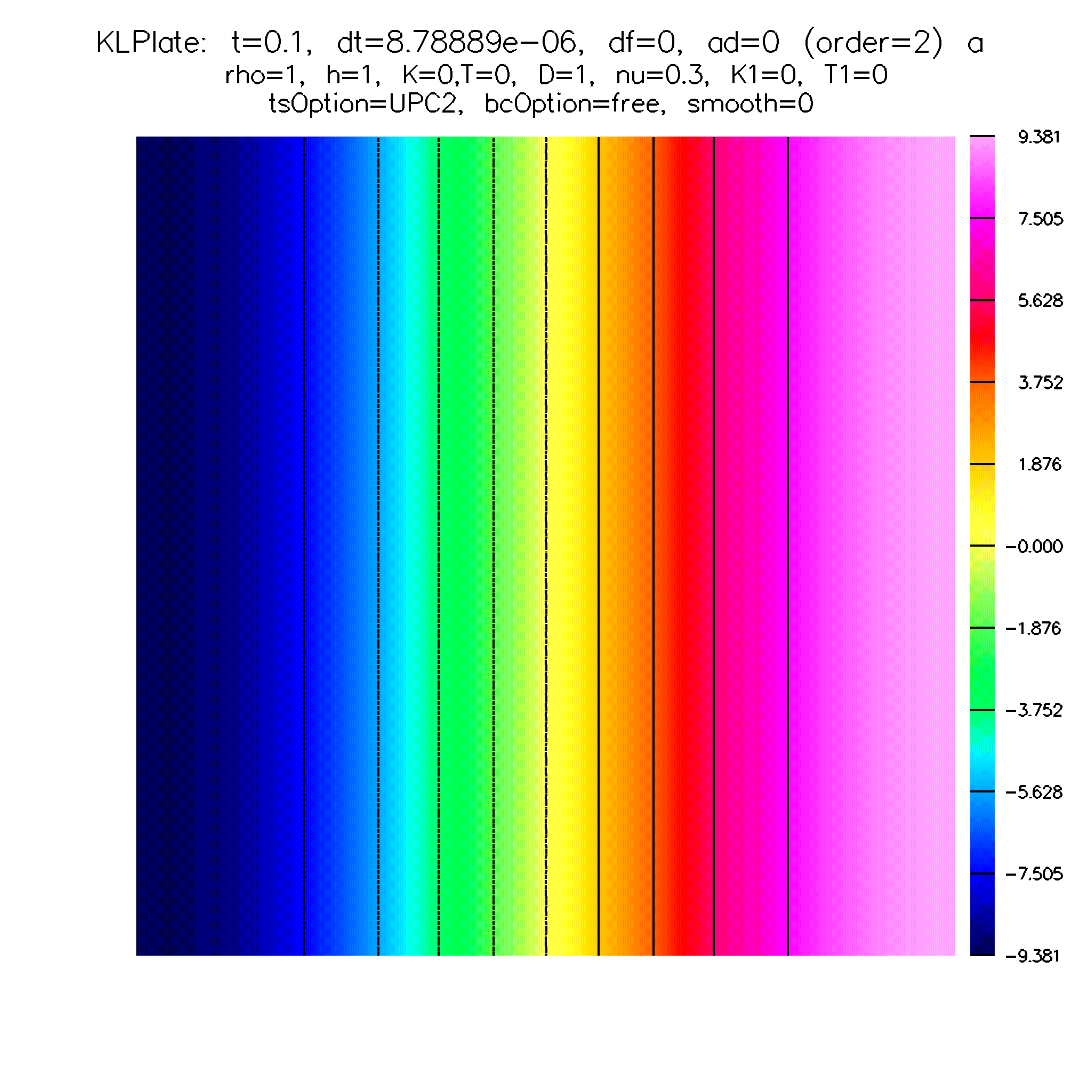}{\figWidth}};

\newcommand{\cbWidth}{.2}
\newcommand{\cbHeight}{3}

\drawColorBarH{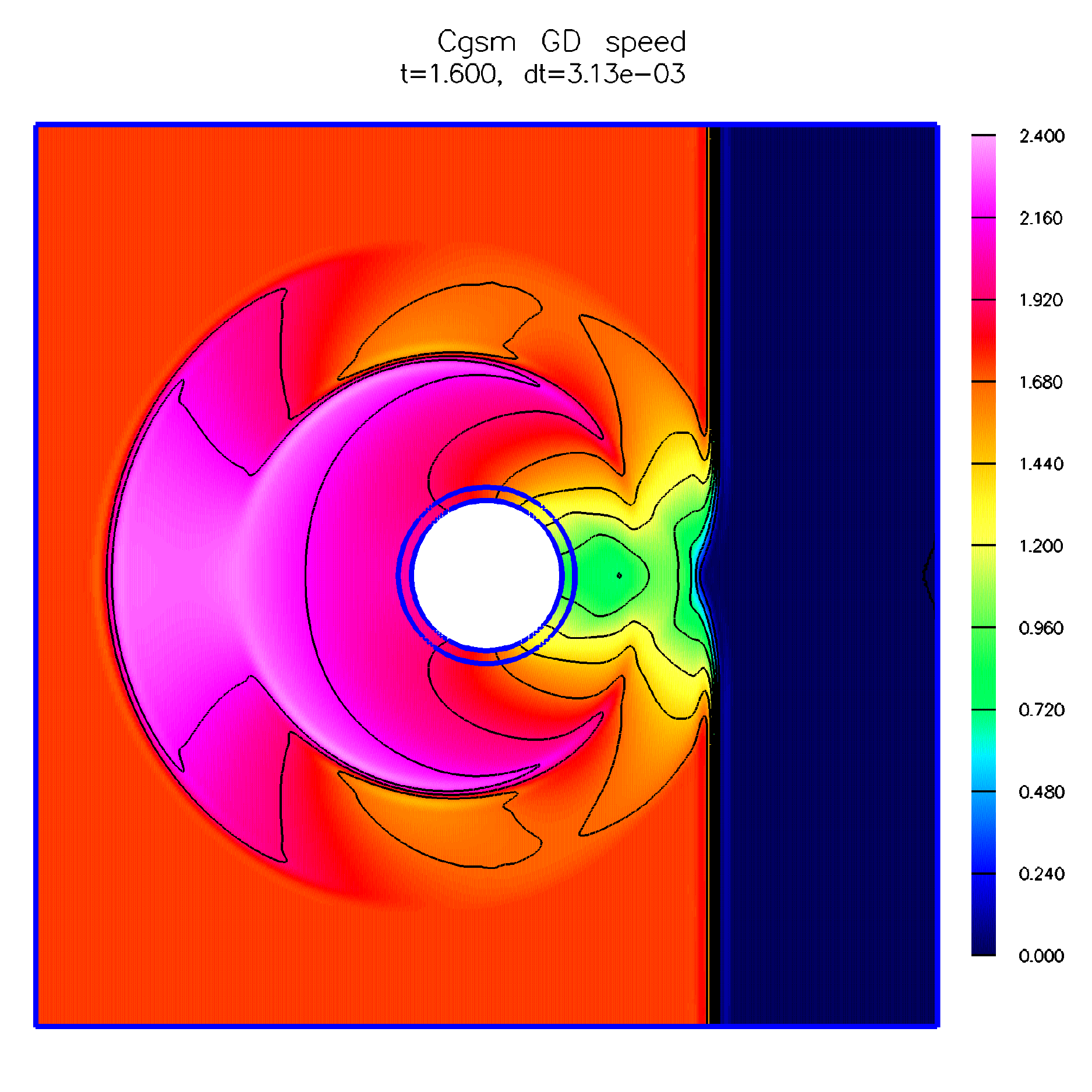}{xshift=6.6cm,yshift=-0.2cm}{\cbWidth}{\cbHeight}{$-14.41$}{$14.41$}{a(x,y,0.1)}
\drawColorBarH{fig/colourBarLines}{xshift=12.2cm,yshift=-0.2cm}{\cbWidth}{\cbHeight}{$-9.381$}{$9.381$}{a(x,y,0.1)}

%
\end{tikzpicture}

\end{center}
\caption{ Computational grids and numerical results obtained using UPC2 scheme for free boundary condition. Left: unit square grid $\Gc_2$. Middle: acceleration contour at $t=0.1$ for the same order  case ($\delta=0$).  Right: acceleration contour at $t=0.1$ for the same stencil  case ($\delta=0$).} \label{fig:orderVsStencilGridsAContours}
\end{figure}
}

Mesh refinement studies   for the plate with free boundaries are performed using both the   same order  and  the same stencil formulas;  maximum-norm errors for all the solution components ($w$,$v$ and $a$) at $t=0.1$ are   collected in Table~\ref{tab:orderVsStencilFree}.   For second order-accurate schemes, we expect the ratio to be around $4$, and the estimated rates to be  around $2$.   From this simple numerical  example, we observe  in Table~\ref{tab:orderVsStencilFree}  that the accuracies  for $w$ and  $v$ meet the  expectation regardless of the $\delta$ values; however,  the accuracy for  $a$ is greatly affected by the choice of $\delta$. To further examine  how $\delta$ influences  the acceleration solution, we plot the contours of $a(x,y,0.1)$ for $\delta=0$ and $\delta=1$ in the middle and right images of Figure~\ref{fig:orderVsStencilGridsAContours}, respectively. As is seen in those plots, the acceleration of the $\delta=0$ case  appears to be very noisy and inaccurate, while that of the $\delta=1$ case is smooth and accurate  as expected.

{
\newcommand{\figWidth}{5.cm}
\def\xa{16}
\def\ya{9}
\newcommand{\trimfig}[2]{\trimw{#1}{#2}{0.}{0.}{0.}{0.}}
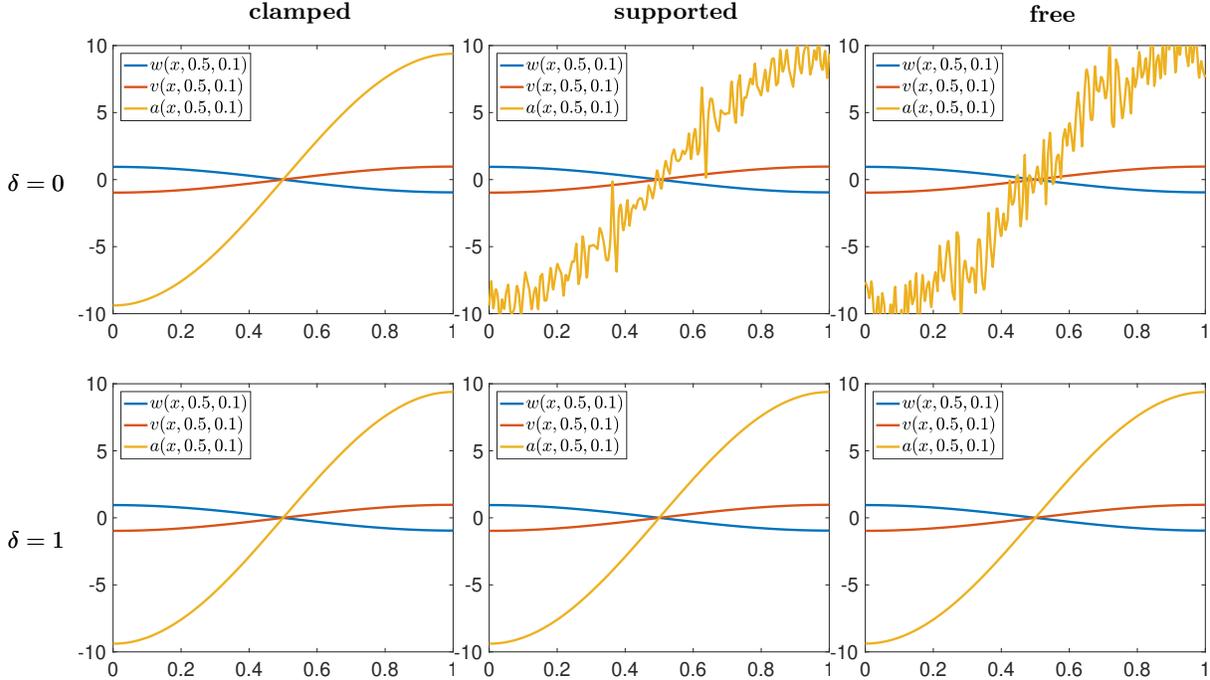
\begin{figure}[h]
\begin{center}
\begin{tikzpicture}[scale=1]
  \useasboundingbox (0.0,0.0) rectangle (\xa,\ya);  

\def\xs{0.}

\draw(0.5,4.5) node[anchor=south west,xshift=\xs cm,yshift=0pt] {\trimfig{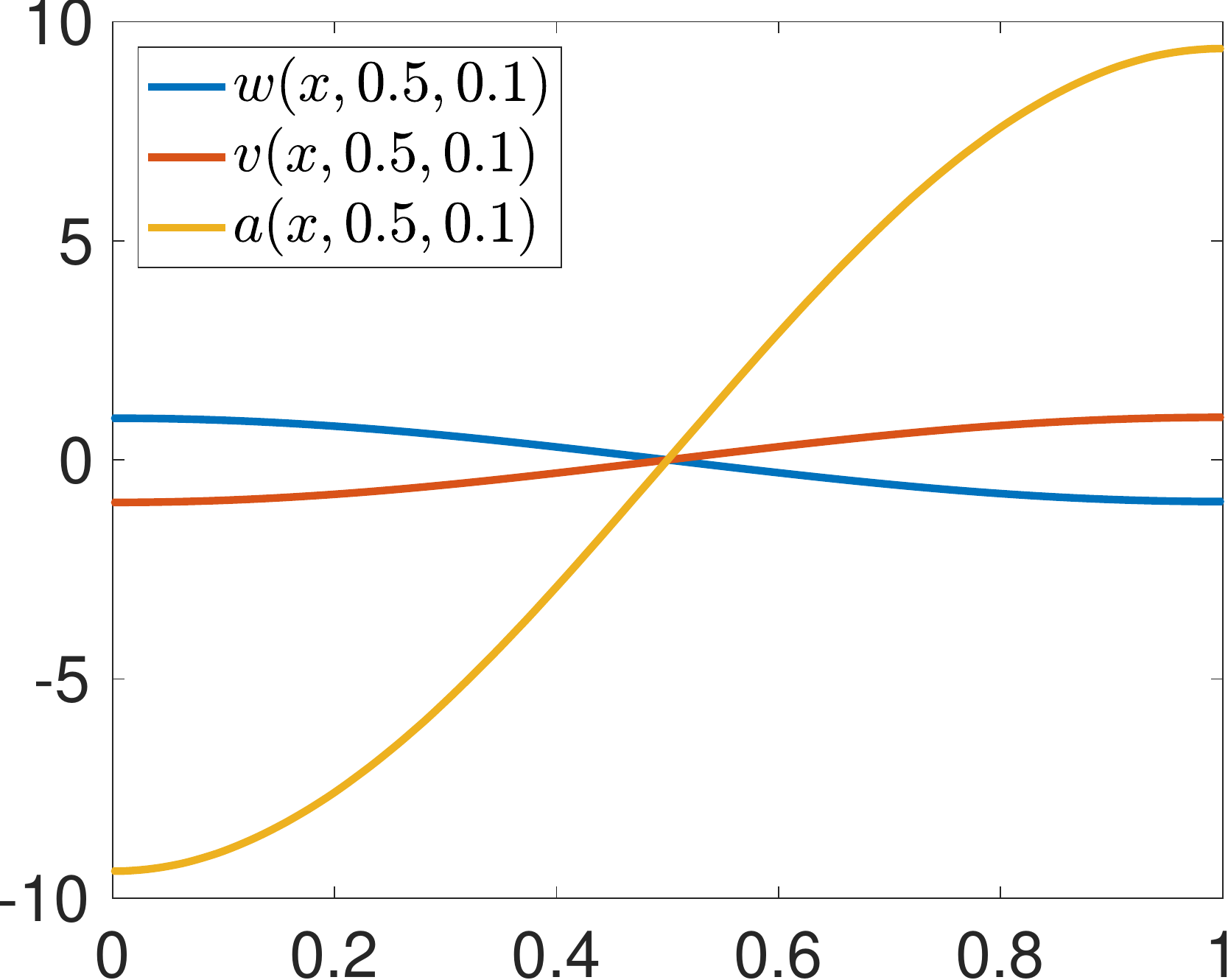}{\figWidth}};
\draw(5.5,4.5) node[anchor=south west,xshift=\xs cm,yshift=0pt] {\trimfig{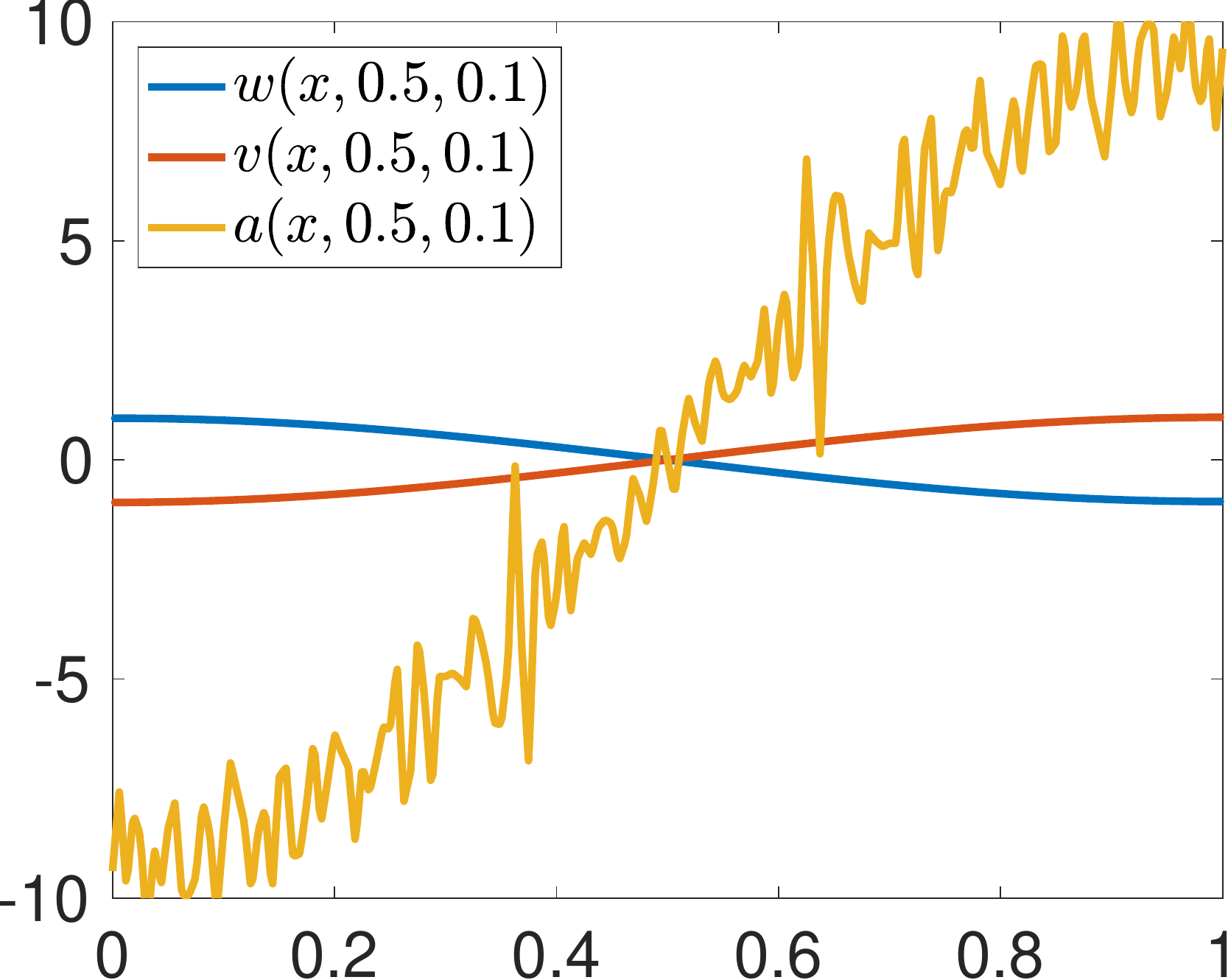}{\figWidth}};
\draw(10.5,4.5) node[anchor=south west,xshift=\xs cm,yshift=0pt] {\trimfig{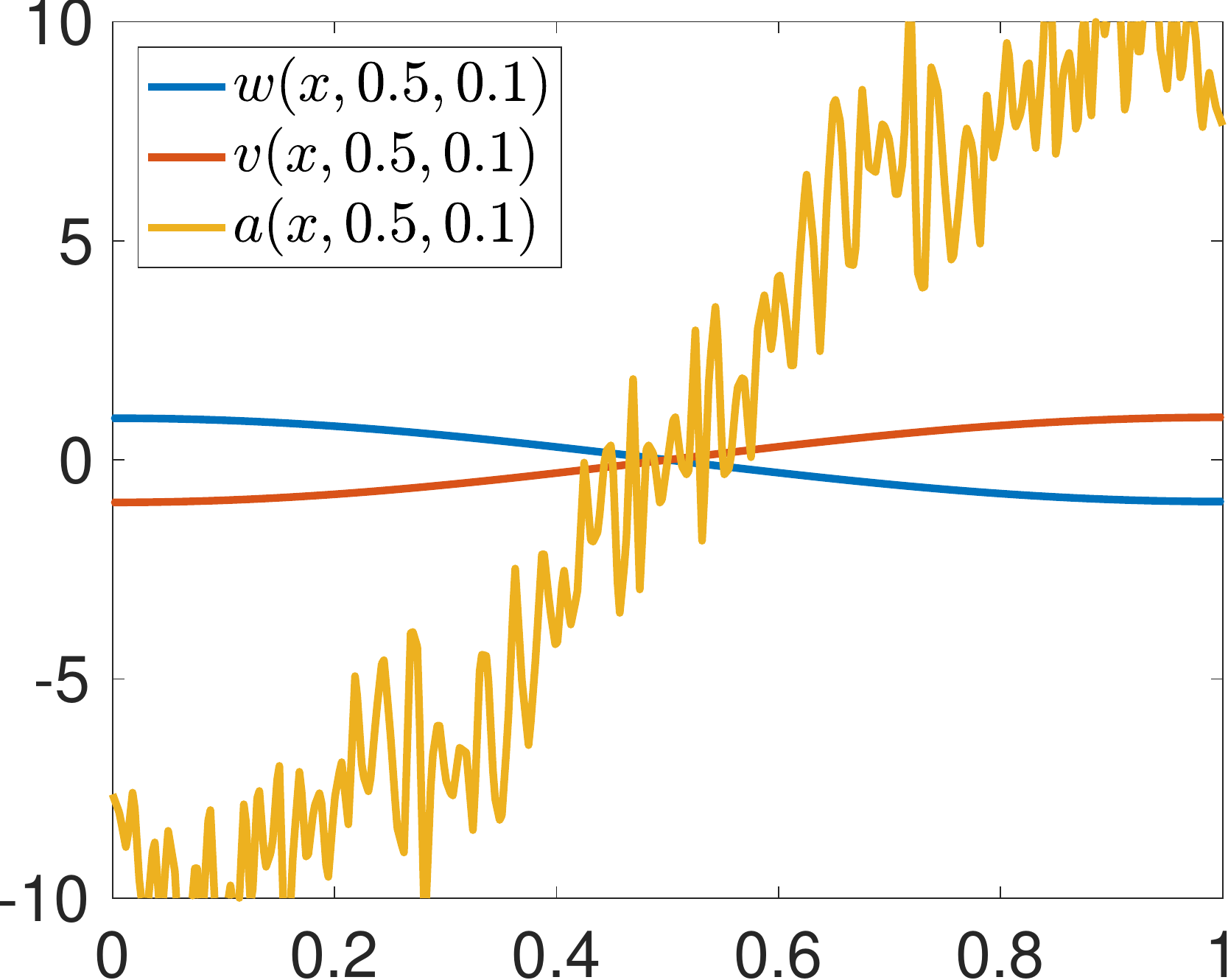}{\figWidth}};

\draw(0.5,0.0) node[anchor=south west,xshift=\xs cm,yshift=0pt] {\trimfig{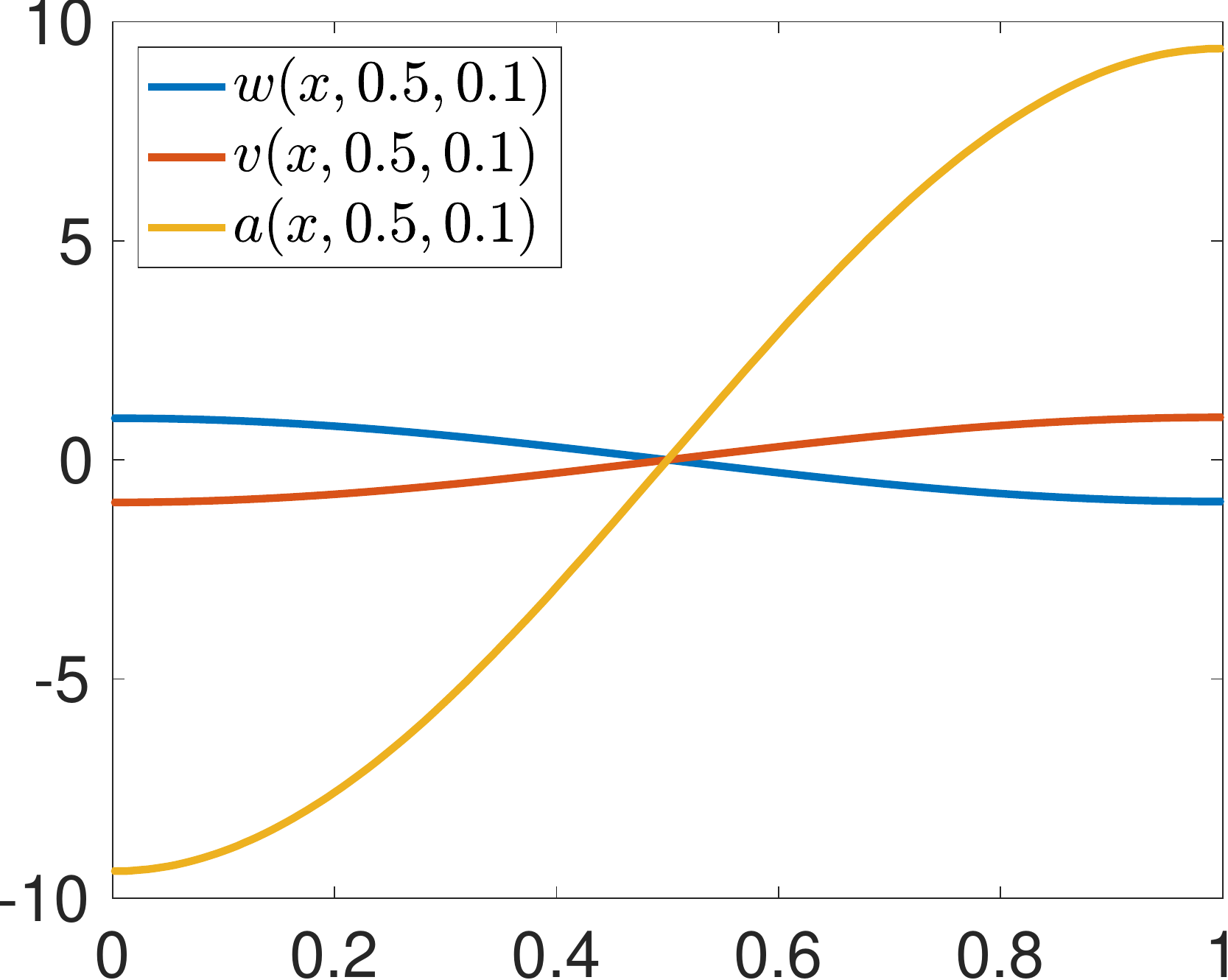}{\figWidth}};
\draw(5.5,0.0) node[anchor=south west,xshift=\xs cm,yshift=0pt] {\trimfig{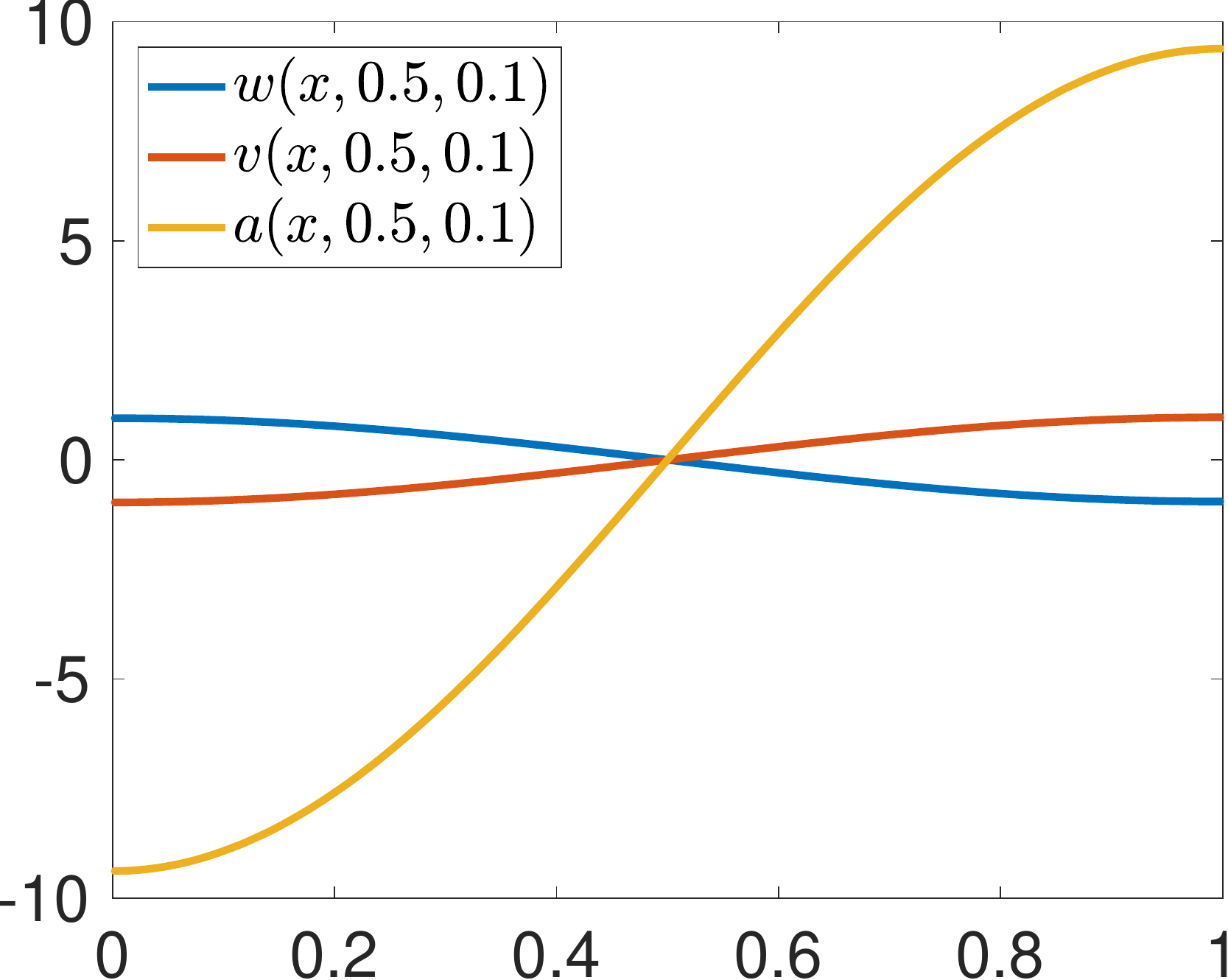}{\figWidth}};
\draw(10.5,0.0) node[anchor=south west,xshift=\xs cm,yshift=0pt] {\trimfig{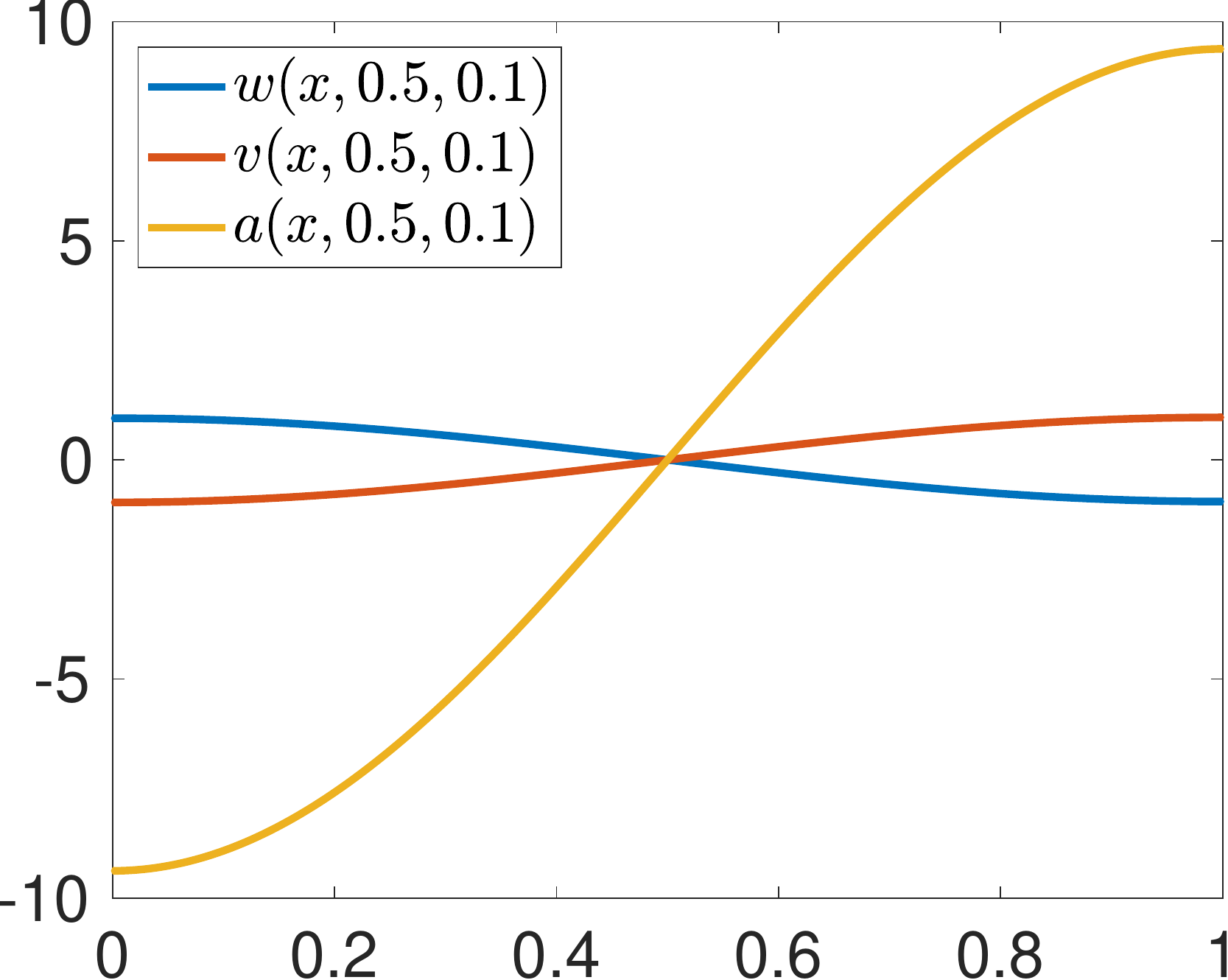}{\figWidth}};

\draw(1,6.75)  node[anchor=east,xshift=0cm,yshift=0pt] {\footnotesize\bf $\delta=0$};
\draw(1,2.)  node[anchor=east,xshift=0cm,yshift=0pt] {\footnotesize\bf $\delta=1$};

\draw(4,9)  node[anchor=center,xshift=0cm,yshift=0pt] {\footnotesize\bf clamped};
\draw(9,9)  node[anchor=center,xshift=0cm,yshift=0pt] {\footnotesize\bf supported};
\draw(14,9)  node[anchor=center,xshift=0cm,yshift=0pt] {\footnotesize\bf free};

\draw(1,6.75)  node[anchor=east,xshift=0cm,yshift=0pt] {\footnotesize\bf $\delta=0$};
\draw(1,2.)  node[anchor=east,xshift=0cm,yshift=0pt] {\footnotesize\bf $\delta=1$};

%
\end{tikzpicture}
\end{center}
\vspace{-0.2in}
\caption{  Numerical solutions along the horizontal center line of the unit square domain. Same order formulas ($\delta=0$) are used to generate the  solutions shown on the top row, while same stencil formulas ($\delta=1$) are used for those shown  on the bottom row. All results  are computed  on grid $\Gc_{16}$.} \label{fig:orderVsStencilLinePlots}
\end{figure}
}

The exact solution \eqref{eq:MMSOrderVsStencil} is intentionally specified to be independent of  $y$ variable, so that we can focus on understanding   the discretization issue in  $x-$direction only, and treat the problem as if it were in a 1D domain. Therefore, in Figure~\ref{fig:orderVsStencilLinePlots}, we plot  all the numerical solutions along  the horizontal center line of the  unit square domain; i.e., $w(x,0.5,0.1)$ ,     $v(x,0.5,0.1)$ and  $a(x,0.5,0.1)$. From these plots, we can see that  using the same order formulas for computing lower-order derivatives causes degradation of the acceleration accuracy for   plates with supported and free boundaries.  This is because, if $\delta=0$, the 
second-order derivative  $\partial^2_{x}w$ that appears   in both free and supported boundary conditions are approximated to 2nd-order accuracy with  a 3-point stencil  in the $x-$direction; this condition helps  determine   the $w$ solution on the first layer of ghost points. The  acceleration is related to the fourth-order derivative  $\partial^4_{x}w$ whose accuracy relies on the smoothness of $w$ solution across a 5-point stencil. Near the boundaries, values of $w$ on the ghost points are involved in the calculation of the acceleration. If derived using  the 3-point stencil,  the $w$ solution on the ghost points is less smooth than the other grid points that causes  a  boundary-layer error in the acceleration solution. The  boundary-layer error   then  propagates  through the entire domain that deteriorates  the overall  accuracy of the acceleration solution. However, if the same stencil formulas ($\delta=1$) are used instead, the smoothness of the solution is consistent across the ghost and interior points. Therefore, the acceleration accuracy remains 2nd-order that is consistent with the truncation error of the spatial discretization.





Based on the conclusion of this test, we  proceed with the same stencil formulas by setting  $\delta=1$  in \eqref{eq:1stDerivFD} and \eqref{eq:2ndDerivFD} for the rest of the numerical tests.

\subsubsection{Manufactured solutions}
Here we perform   an exhaustive convergence study for   all the numerical algorithms proposed in Section~\ref{sec:tsSchemes} subject to all the possible boundary conditions \eqref{eq:clampedBC} --\eqref{eq:freeBC} using the method of manufactured solutions.
To examine the stabilizing effects of these algorithms, we solve this problem using composite overlapping grid and artificial dissipation.

Specifically, we consider the following exact solution
$$
u_e(x,y,t)=\cos(\pi x)\cos(\pi y)\cos(\pi t),
$$
and solve the manufactured solution problem  for a unit square  plate that is discretized by a sequence of refined composite overlapping grids. The computational grid  $\Gc_2$ for this problem is shown in Figure~\ref{fig:tz_NB2_f_unitSquareCGe16Err}. The  physical parameters of the plate equation \eqref{eq:KLShell} are  $\rho=H=K=T=D=1$,  and the Poisson's ratio that is needed by the supported \eqref{eq:supportedBC} and free \eqref{eq:freeBC} boundary conditions is $\nu=0.3$.

In Figure~\ref{fig:tz_NB2_f_unitSquareCGe16Err},  results of the  plate with free boundaries are provided. Particularly, the  errors    shown in the figure are  computed at $t=0.1$ using the   numerical solutions of the NB2 algorithm~\ref{alg:nb2}. Errors of the other three algorithms are similar, so their plots are omitted here to save space.   We observe that the errors of all the solution components ($w$, $v$, $a$) are well behaved in that the magnitudes are small and they are smooth throughout the domain.

{
  \newcommand{\figWidth}{3.8cm}

\def\xa{16.}
\def\ya{4.}
\newcommand{\trimfig}[2]{\trimw{#1}{#2}{0.12}{0.12}{0.12}{0.12}}
\begin{figure}[h]
\begin{center}
\begin{tikzpicture}[scale=1]
  \useasboundingbox (0.0,0.0) rectangle (\xa,\ya);  

  \def\xs{-0.2}
   \begin{scope}[xshift=-0.25cm,yshift=0 cm]
   \draw(0,0) node[anchor=south west,xshift=\xs cm,yshift=0pt] {\trimfig{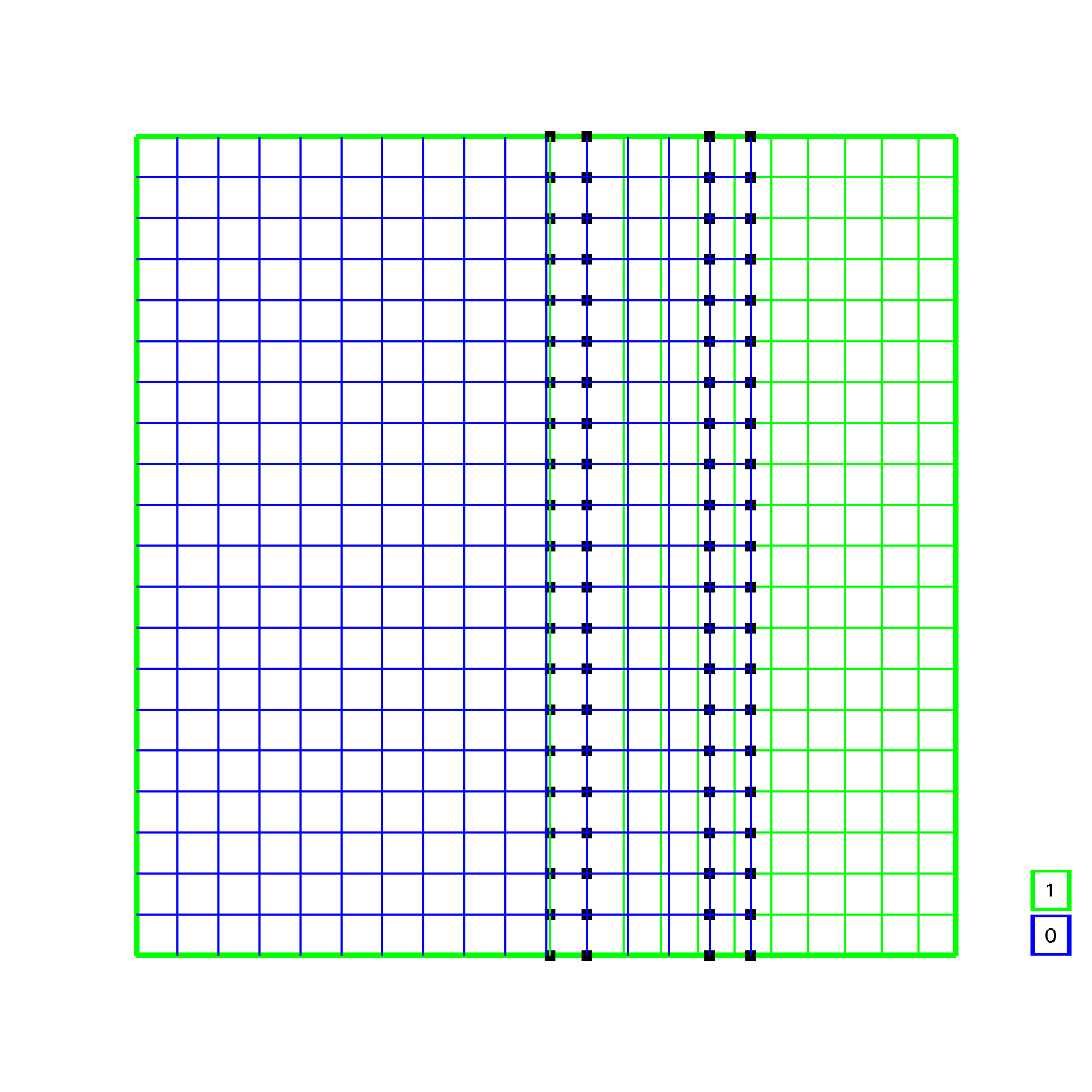}{\figWidth}};
   \draw[<-,black, very thick,xshift=0cm] (2.3,2.1)  -- (1.4,-0.)  node[black,anchor=north,xshift=0.8cm,yshift=0.1cm]{ \color{red!50!black}{interpolation  points}};
   \draw[<-,black, very thick,xshift=0cm] (3.,2.1)  -- (1.4,-0.) ;
   \end{scope}
\draw(3.75,0.0) node[anchor=south west,xshift=\xs cm,yshift=0pt] {\trimfig{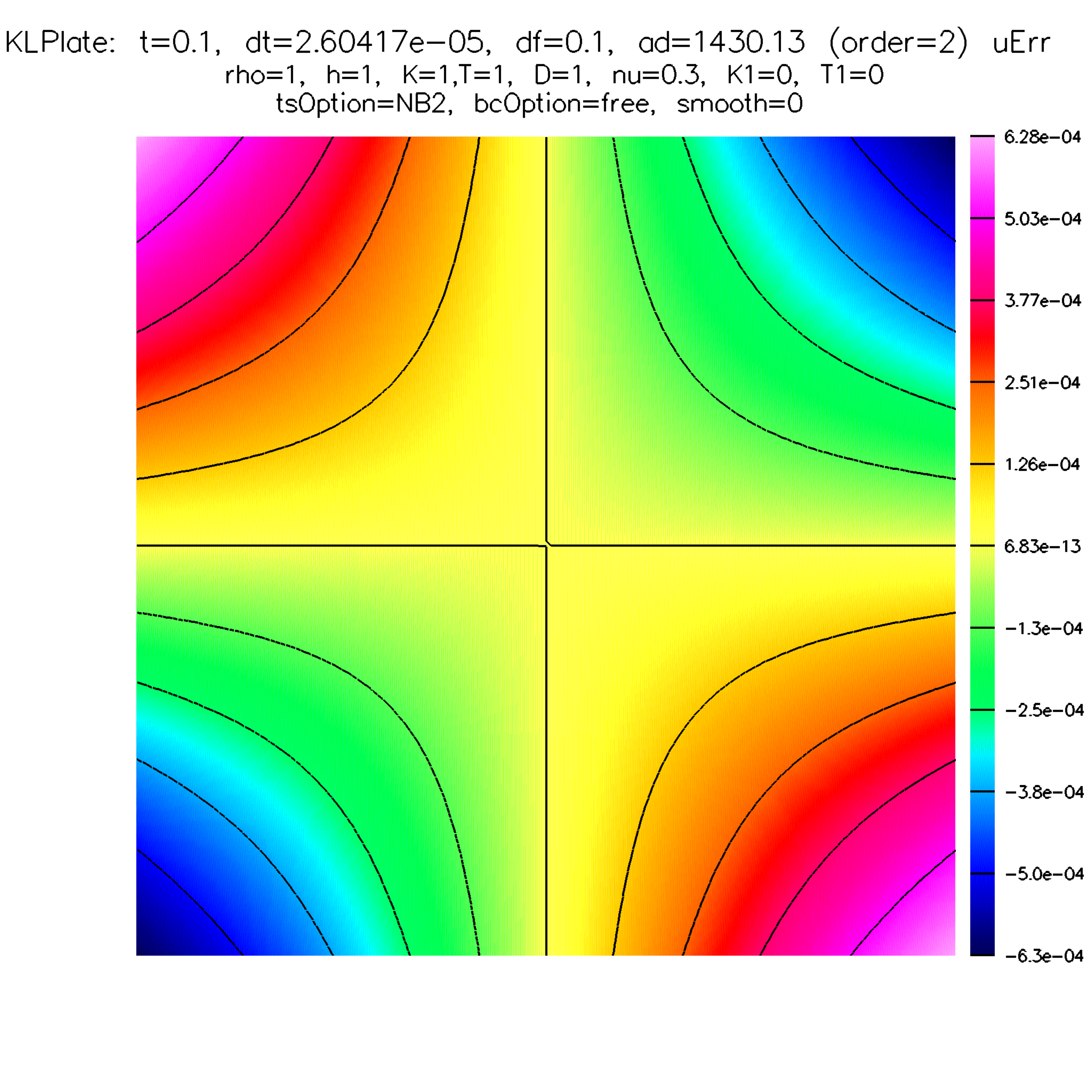}{\figWidth}};
\draw(7.75,0.0) node[anchor=south west,xshift=\xs cm,yshift=0pt] {\trimfig{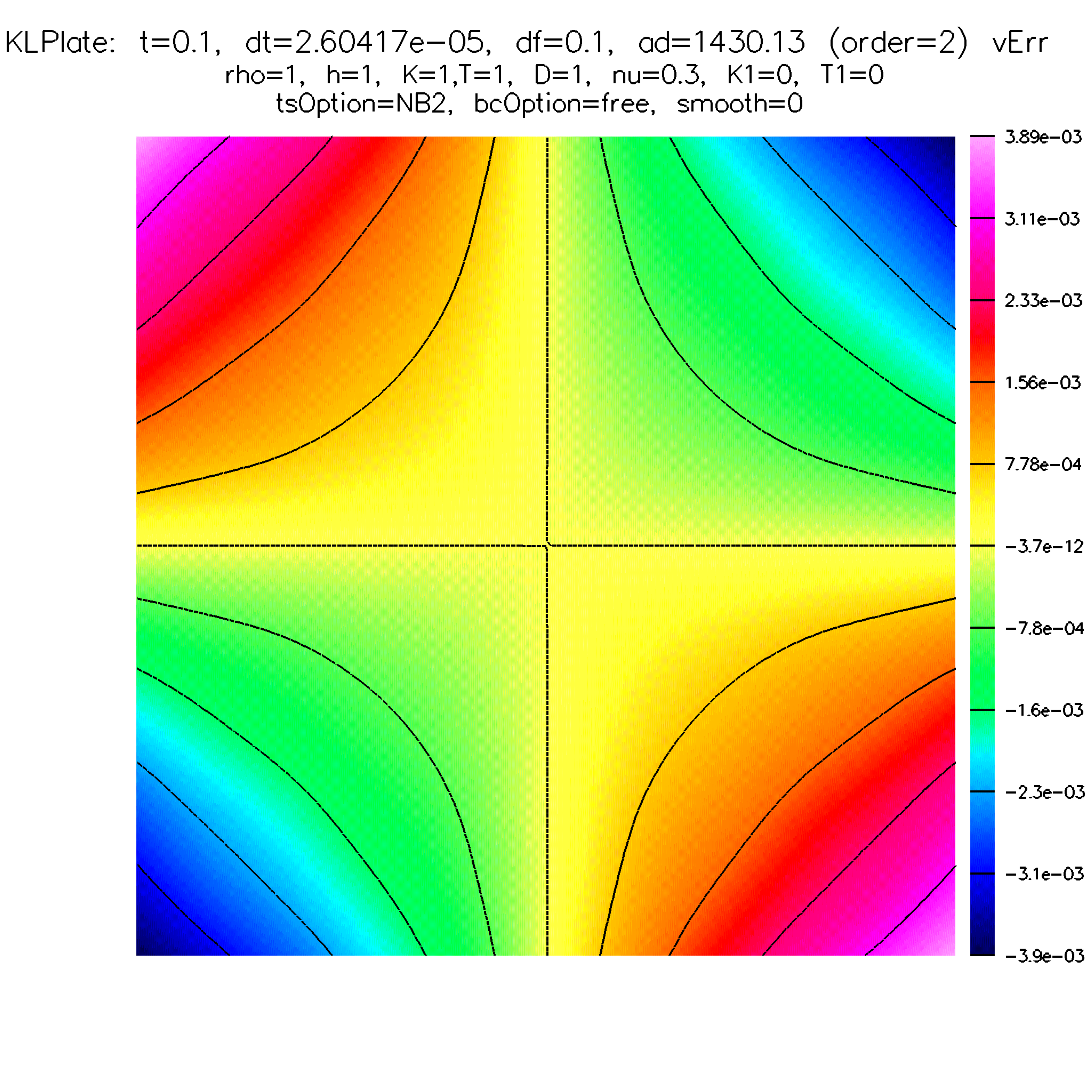}{\figWidth}};
\draw(11.75,0.0) node[anchor=south west,xshift=\xs cm,yshift=0pt] {\trimfig{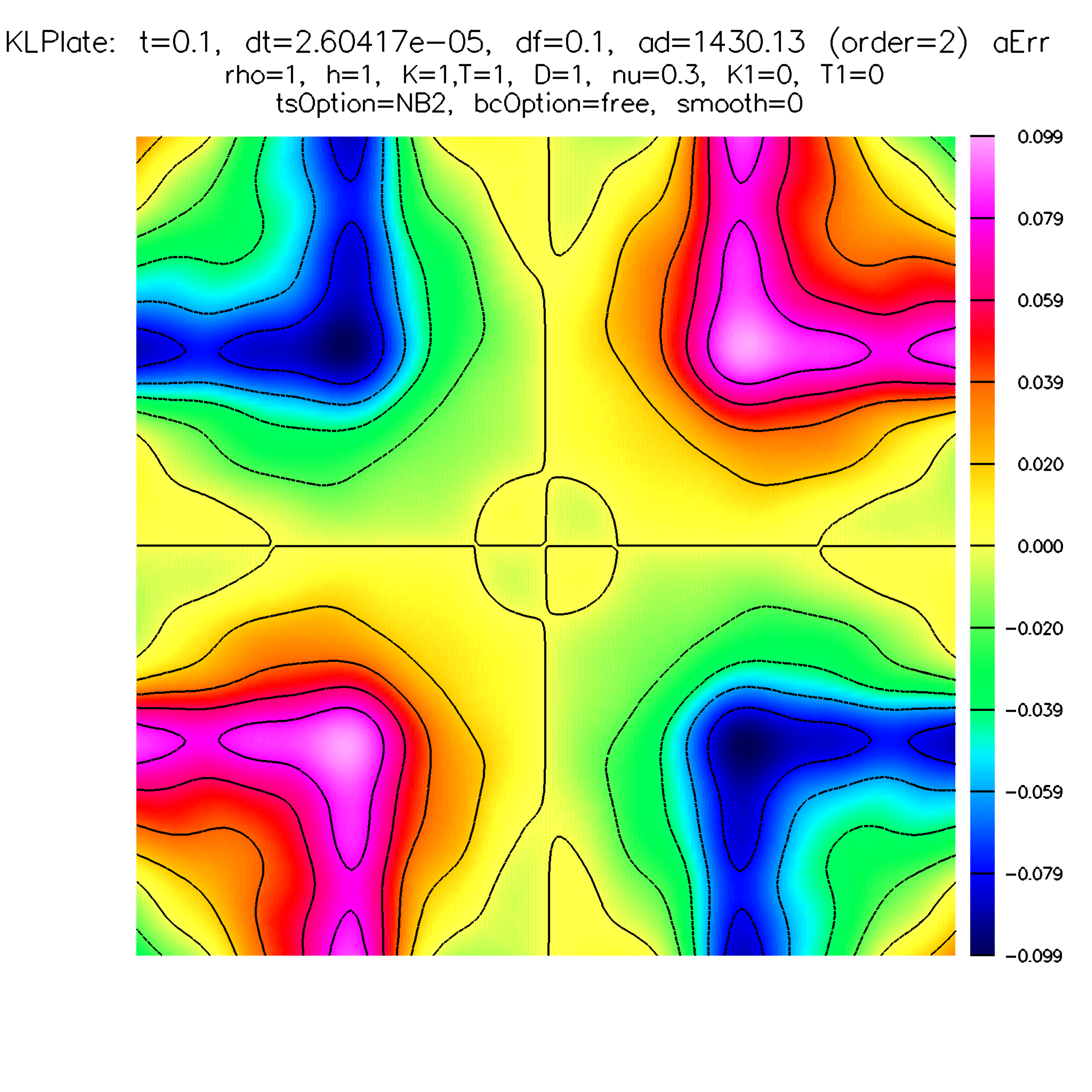}{\figWidth}};

\newcommand{\cbWidth}{.2}
\newcommand{\cbHeight}{3}

\drawColorBarH{fig/colourBarLines}{xshift=4.4cm,yshift=-.2cm}{\cbWidth}{\cbHeight}{$-6.3e-4$}{$6.3e-4$}{E(w)}
\drawColorBarH{fig/colourBarLines}{xshift=8.4cm,yshift=-.2cm}{\cbWidth}{\cbHeight}{$-3.9e-3$}{$3.9e-3$}{E(v)}
\drawColorBarH{fig/colourBarLines}{xshift=12.4cm,yshift=-.2cm}{\cbWidth}{\cbHeight}{$-0.099$}{$0.099$}{E(a)}


%
\end{tikzpicture}

\end{center}
\caption{Computational grids $\Gc_2$, and  the errors at  $t=0.1$   for the plate  with free boundary conditions.   Results shown here  are computed using   NB2 algorithm on the $\Gc_{16}$ grid.} \label{fig:tz_NB2_f_unitSquareCGe16Err}
\end{figure}
}

{
  \newcommand{\figWidth}{4cm}
  \newcommand{\figWidthLgd}{1.8cm}
\def\xa{16}
\def\ya{7}
\newcommand{\trimfig}[2]{\trimw{#1}{#2}{0.}{0.}{0.}{0.}}
\newcommand{\trimfigLgd}[2]{\trimw{#1}{#2}{0.675}{0.}{0.09}{0.475}}

\begin{figure}[h]
\begin{center}
\begin{tikzpicture}[scale=1]
  \useasboundingbox (0.0,0.0) rectangle (\xa,\ya);  

\draw(-0.5,0.0) node[anchor=south west,xshift=0pt,yshift=0pt] {\trimfig{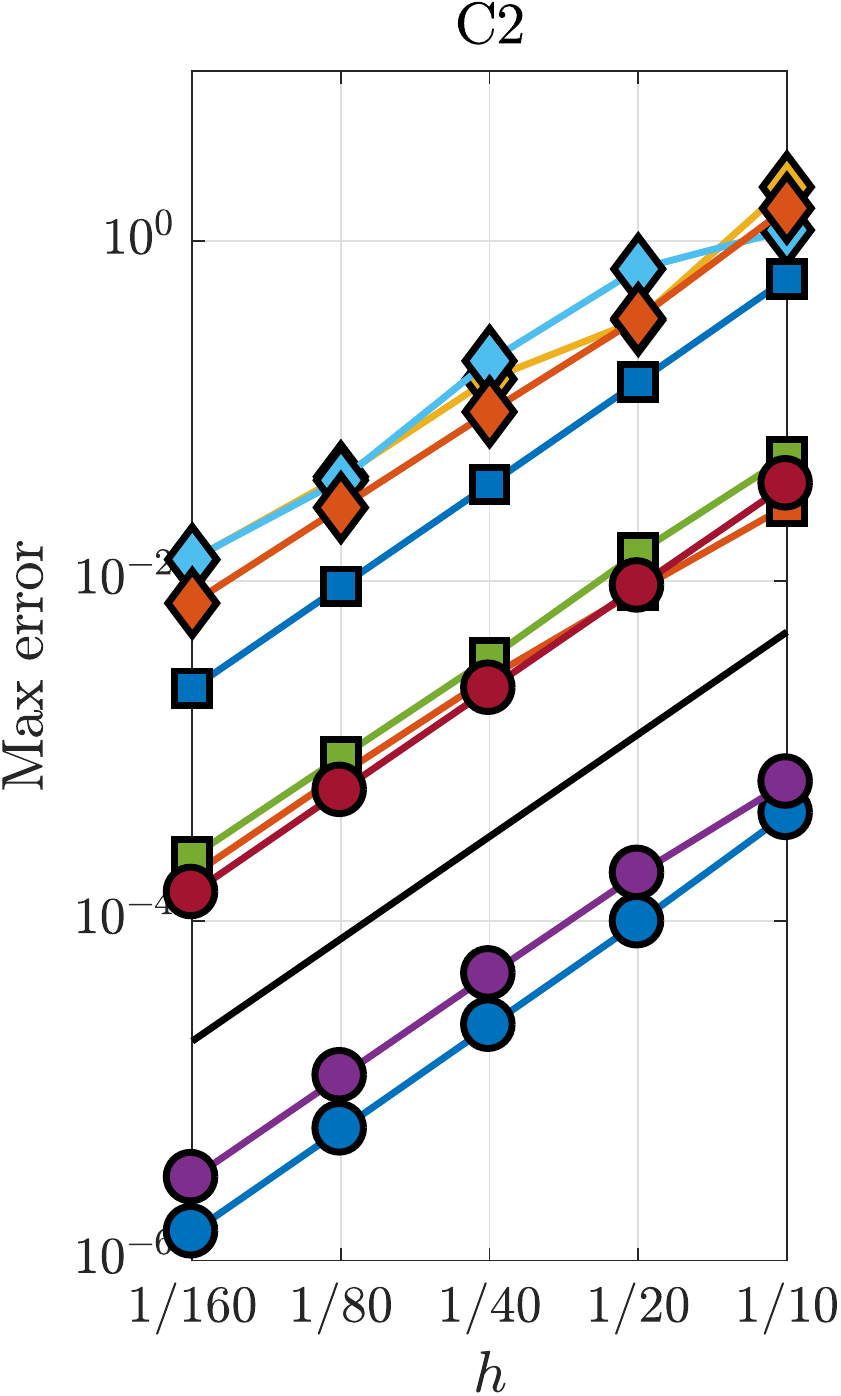}{\figWidth}};
\draw(3.5,0.0) node[anchor=south west,xshift=0pt,yshift=0pt] {\trimfig{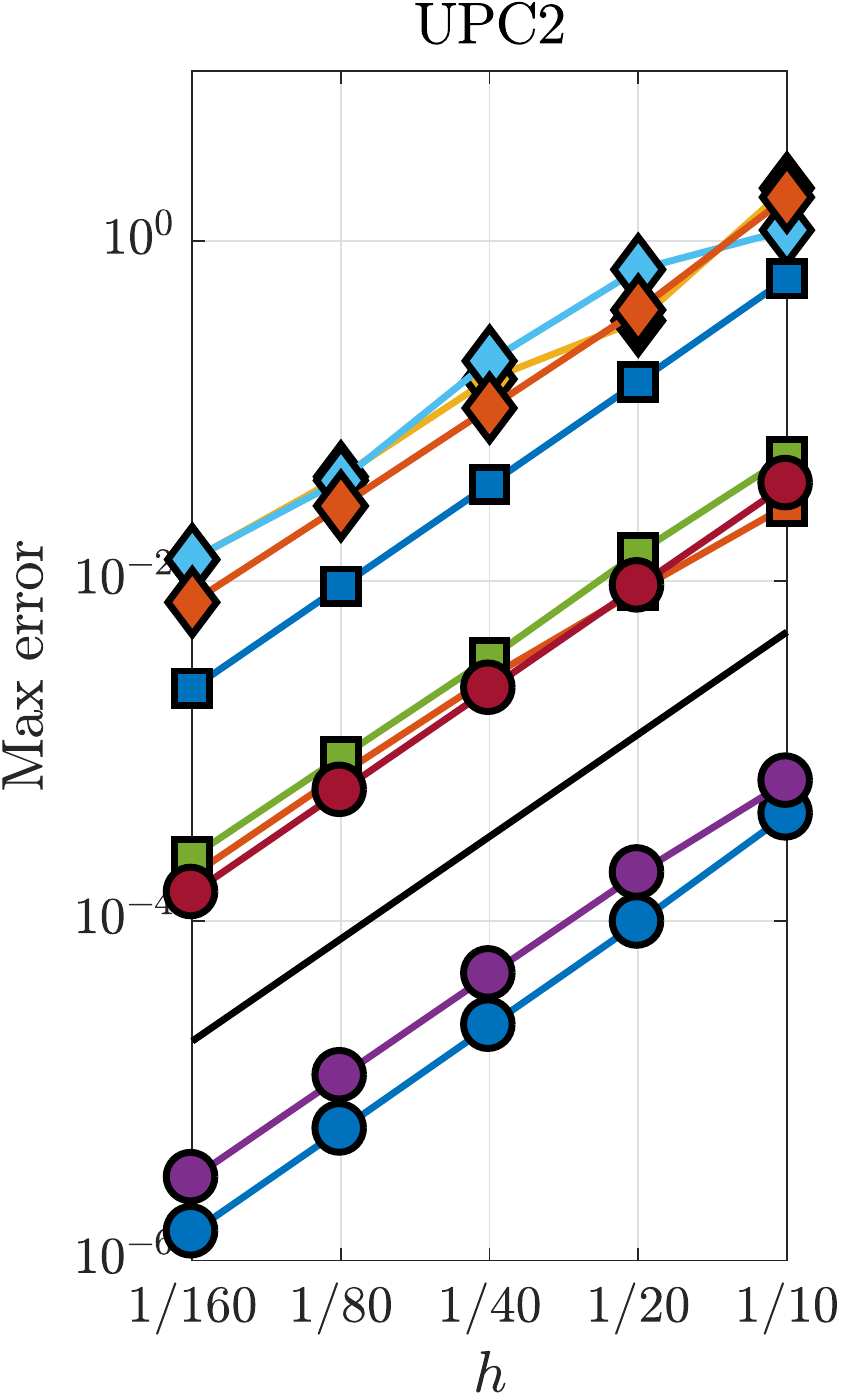}{\figWidth}};
\draw(7.5,0.0) node[anchor=south west,xshift=0pt,yshift=0pt] {\trimfig{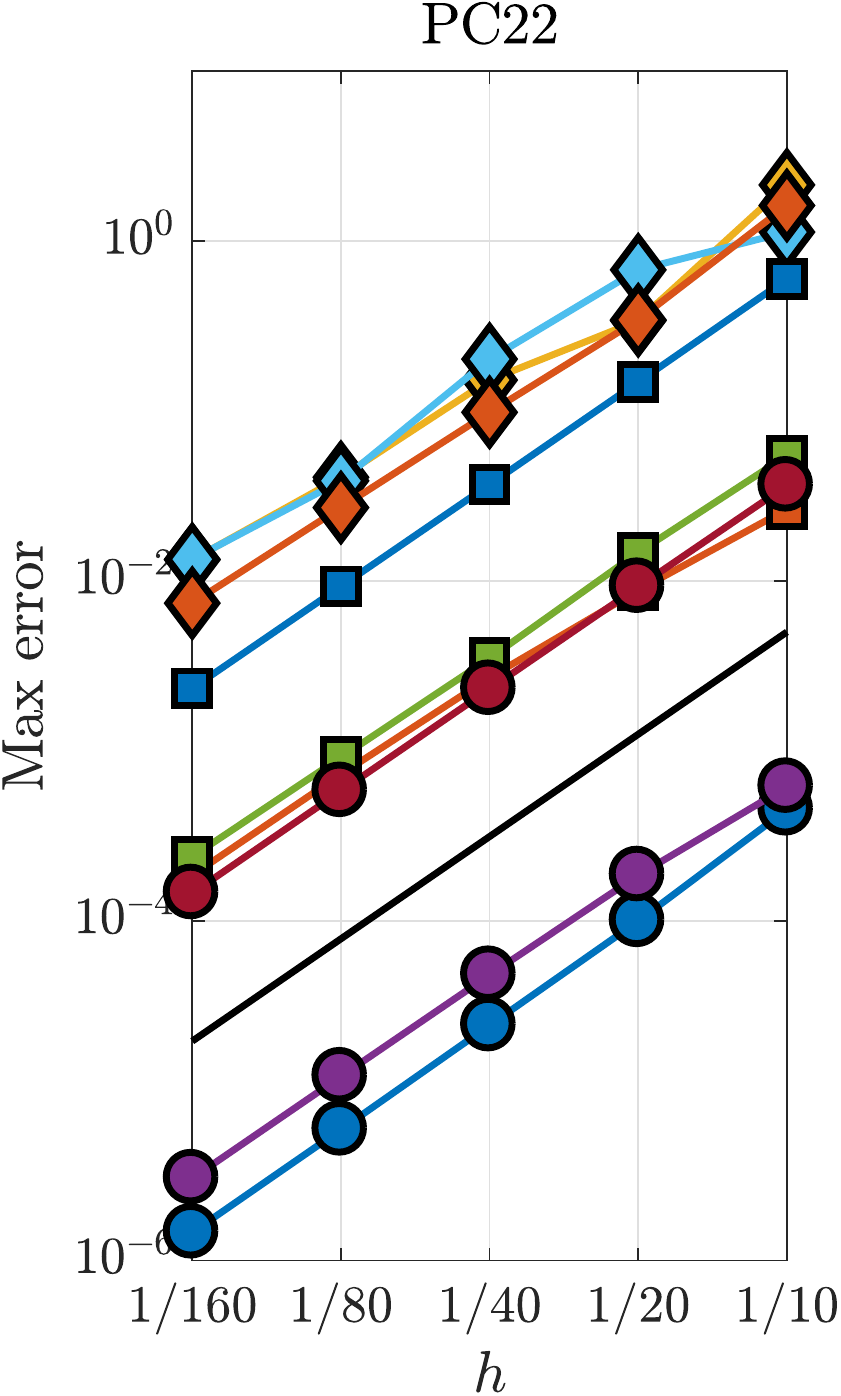}{\figWidth}};
\draw(11.5,0.0) node[anchor=south west,xshift=0pt,yshift=0pt] {\trimfig{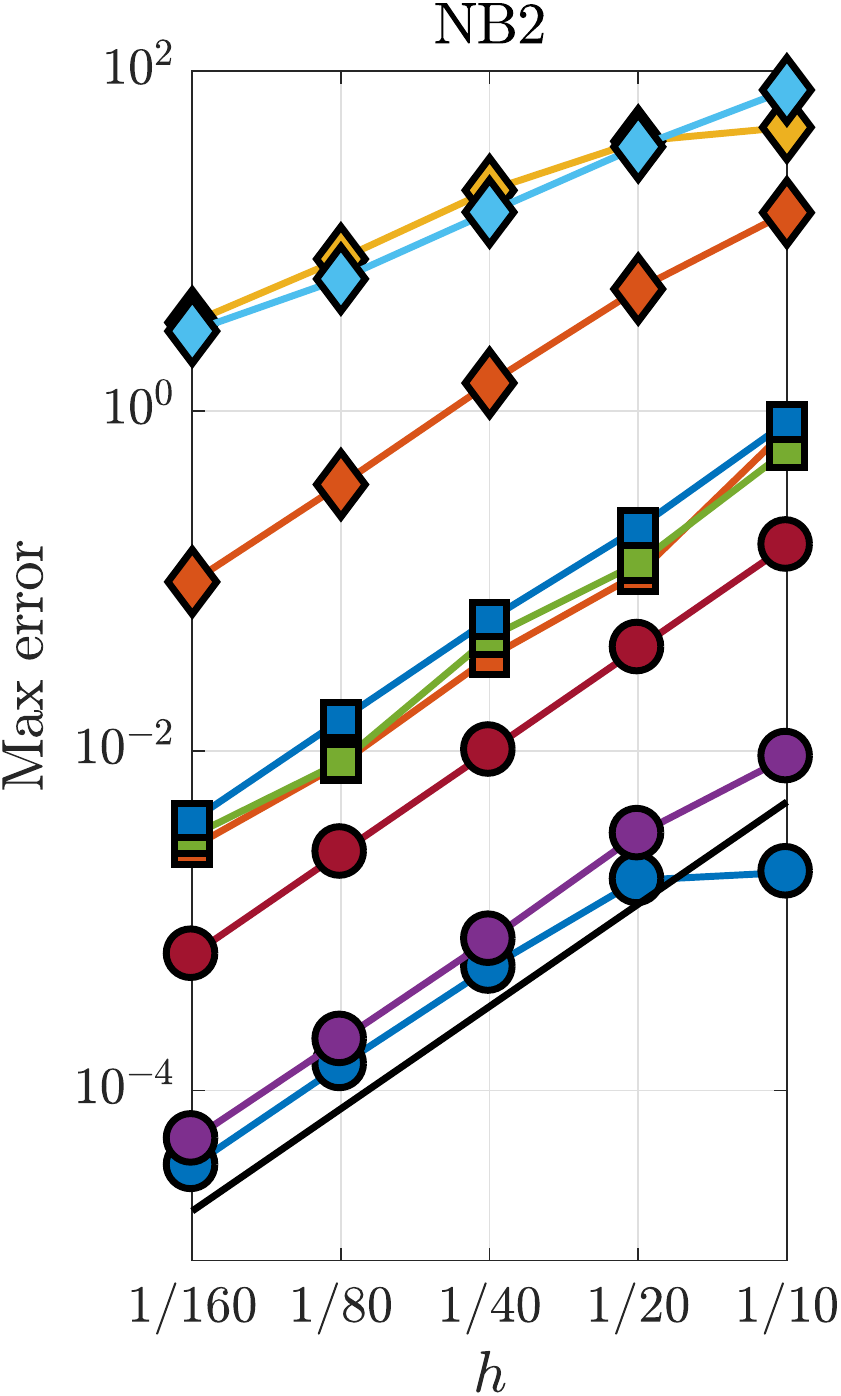}{\figWidth}};

\draw(-0.8,4.2) node[anchor=south west,xshift=0pt,yshift=0pt] {\trimfigLgd{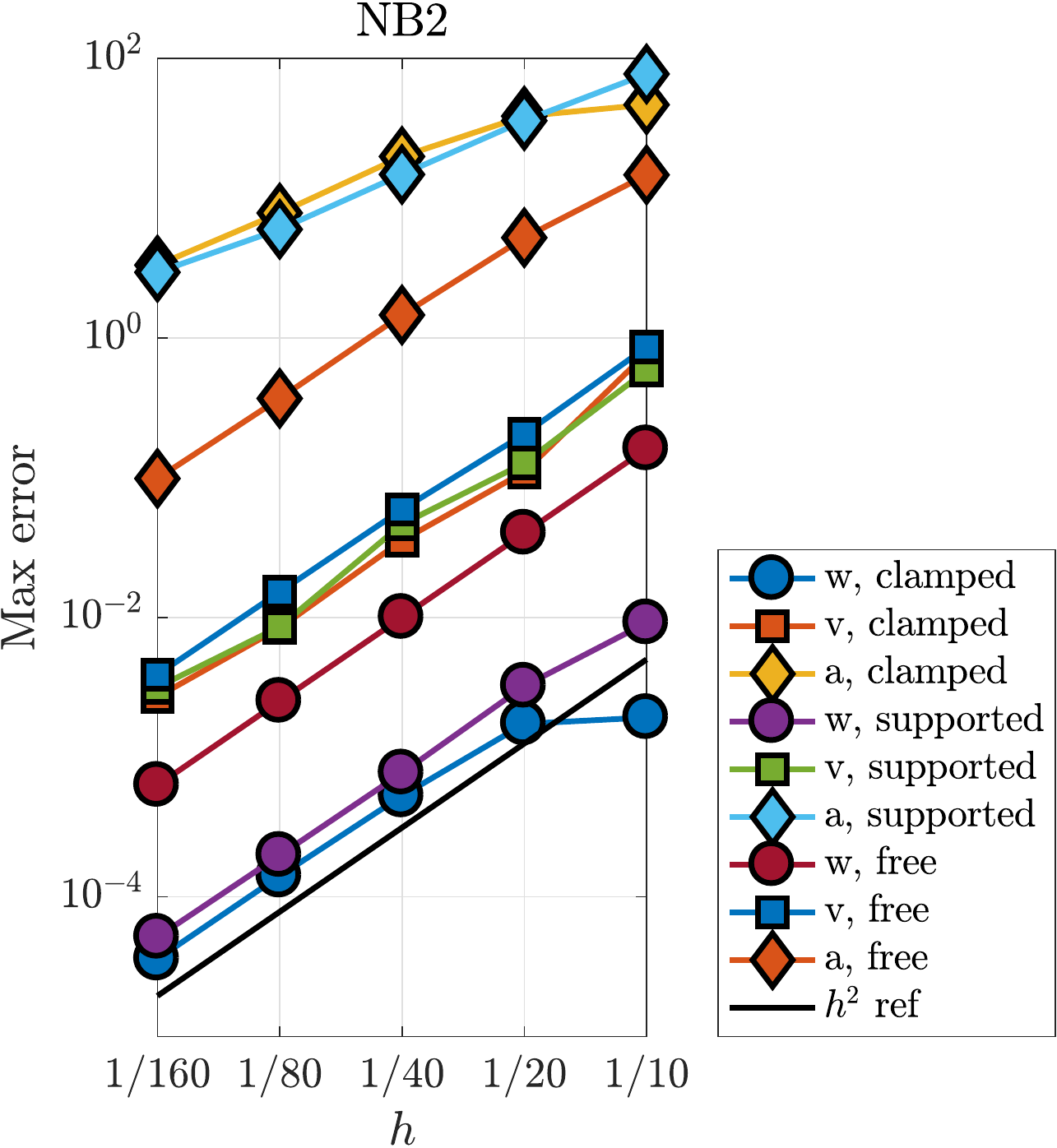}{\figWidthLgd}};


%
\end{tikzpicture}

\end{center}
\vspace{-0.3in}
\caption{Convergence rates of the manufactured solution problem  for all the numerical schemes and boundary conditions.  The stability factor is chosen to be $\csf=0.9$ for the explicit schemes (C2,UPC2,PC22), while $\csf=5$ for the implicit NB2 scheme. The dissipation factor is  $df=0.1$ for all the computations. } \label{fig:unitSquareCGeRate}
\end{figure}
}

Convergence studies of this test  are presented in Figure~\ref{fig:unitSquareCGeRate}. Here the manufactured solution problem is solved comprehensively with all of the proposed algorithms (C2,UPC2, PC22, NB2) subject to each of the boundary conditions (clamped, supported and free) on  a sequence of refined meshes. The time step for each calculation is determined according to the formulas summarized  in Table~\ref{tab:timeStep}, where the stability factor is chosen to be $\csf=0.9$ for the explicit schemes (C2,UPC2,PC22) and  $\csf=5$ for the implicit NB2 scheme. Since this problem is solved on composite overlapping grids, the artificial dissipation  \eqref{eq:dissipationCoefficient} with $\df=0.1$ is needed  in the algorithms for stabilizing the simulations.

In  Figure~\ref{fig:unitSquareCGeRate}, the maximum-norm errors of all the numerical solutions  are plotted against  the target grid spacings in  log-log scale, together with a reference line indicating 2nd-order accuracy.  As  are  observed  from these plots,  2nd-order accuracies are  achieved in  the solutions of  $w$,$v$ and $a$  by all the proposed algorithms subject to any of the considered boundary conditions.  The comprehensive convergence study  validates the  accuracy and stability properties of our approach for solving \KL plate model on composite overlapping grid.

\subsubsection{Analytical solutions for circular plates}

{
\begin{figure}[h]
     \centering
     \begin{subfigure}[b]{0.45\textwidth}
         \centering
         \includegraphics[width=\textwidth]{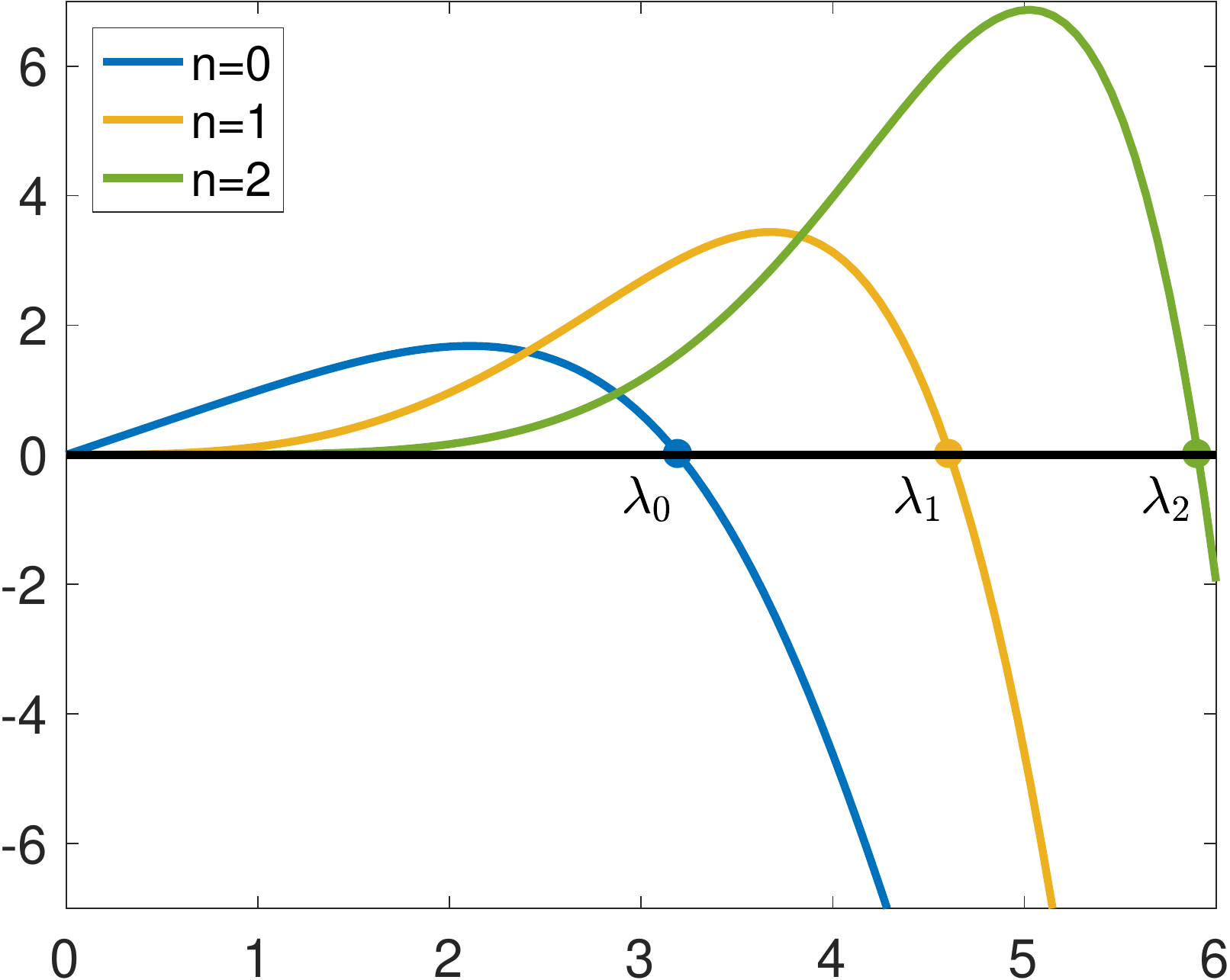}
         \caption{clamped $\phi_n^c(\lambda)$}
         \label{fig:transcendentalFunction_c}
     \end{subfigure}
     \hfill
     \begin{subfigure}[b]{0.45\textwidth}
         \centering
         \includegraphics[width=\textwidth]{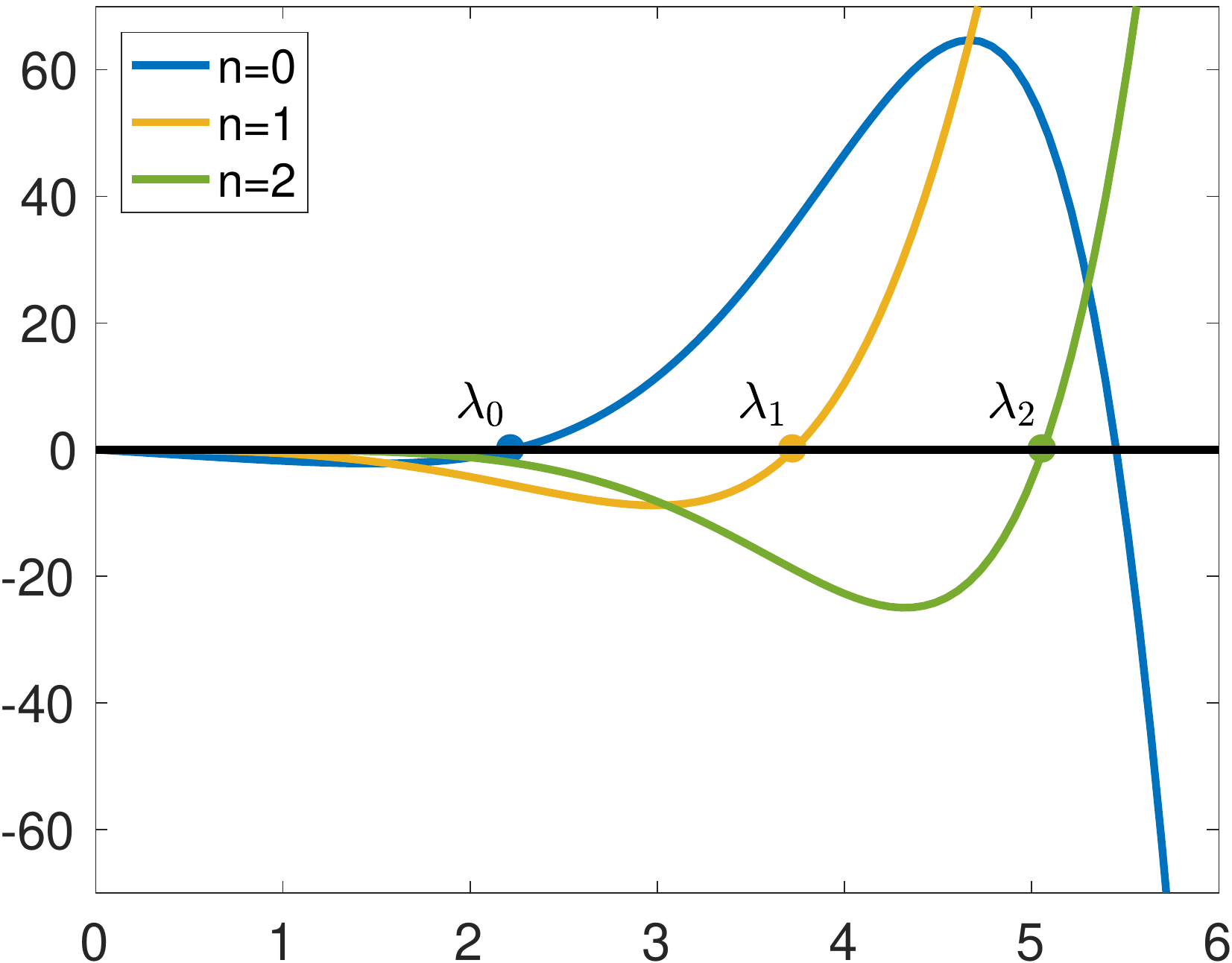}
         \caption{supported $\phi_n^s(\lambda)$}
         \label{fig:transcendentalFunction_s}
     \end{subfigure}
     \caption{Plots of the first three  transcendental functions  and their smallest roots $\lambda_n$ with $n=0,1,2$.} \label{fig:transcendentalFunction}
\end{figure}
}

Analytical solutions to  \eqref{eq:KLShell}  are available for  circular plates   with   clamped \eqref{eq:clampedBC} or  supported \eqref{eq:supportedBC} boundaries   \cite{Wah1962}.  We use these exact solutions to further verify the numerical properties of our algorithms for solving real plate problems  as opposed to the manufactured ones  in the previous tests.
Given  the  circular  domain $\Omega=\{\xv\in\mathbb{R}^2 : |\xv|\leq a\}$ and the parameters $K=T=f=0$,  analytical solutions  for  the \KL  model are  derived  in polar coordinates by separation of variables $w(r,\theta,t)=R(r)\Theta(\theta)T(t)$. The general solution satisfying the boundary condition $w=0$ reads,
\begin{equation}\label{eq:generalSolutionDiskPlate}
  w(r,\theta,t)=A_n\left[J_n\left(\frac{\lambda r}{a}\right)-\frac{J_n(\lambda)}{I_n(\lambda)}I_n\left(\frac{\lambda r}{a}\right)\right] \left(\cos(n\theta)+ \gamma_n \sin(n\theta)\right)(\sin \omega t +C_n\cos \omega t),
\end{equation}
  where  the natural frequency  is
  $
  \omega=\lambda^2\sqrt{D/(\rho H)},
  $
  and 
  $\lambda$ is related to the  eigenvalue of $ \Delta^2$, i.e.,
  $
    \Delta^2 w = \lambda^4 w
  $. 
    In the general solution \eqref{eq:generalSolutionDiskPlate},   $J_n$ is the Bessel's function of the first kind  and $I_n$ is the  modified Bessel's function of the first kind; here $n$ represents the order of the Bessel's functions.
The remaining coefficients, $A_n$, $\gamma_n$, and $C_n$, are to be determined by initial conditions.

{
  \newcommand{\figWidth}{3.8cm}

\def\xa{16.}
\def\ya{4.5}
\newcommand{\trimfig}[2]{\trimw{#1}{#2}{0.12}{0.12}{0.12}{0.12}}
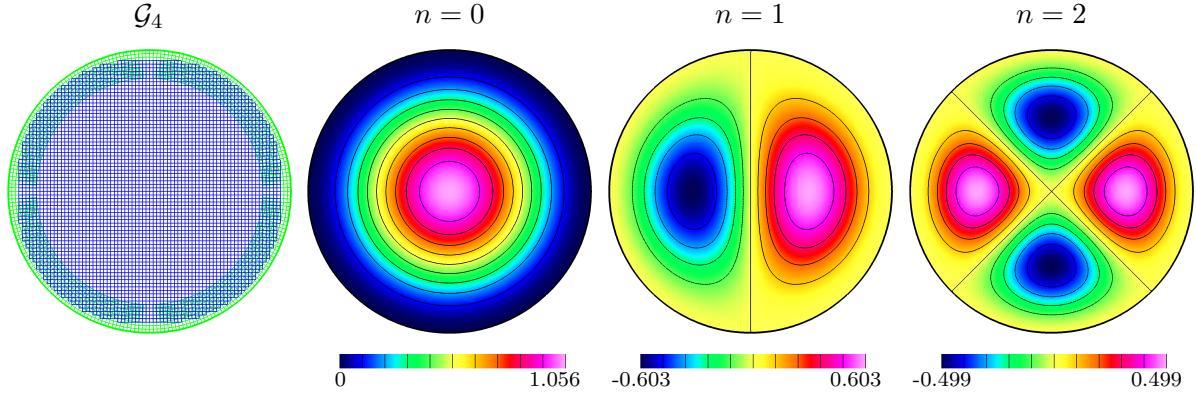
\begin{figure}[ht]
\begin{center}
\begin{tikzpicture}[scale=1]
  \useasboundingbox (0.0,0.0) rectangle (\xa,\ya);  

\def\xs{-0.2}
\draw(-0.25,0.0) node[anchor=south west,xshift=\xs cm,yshift=0pt] {\trimfig{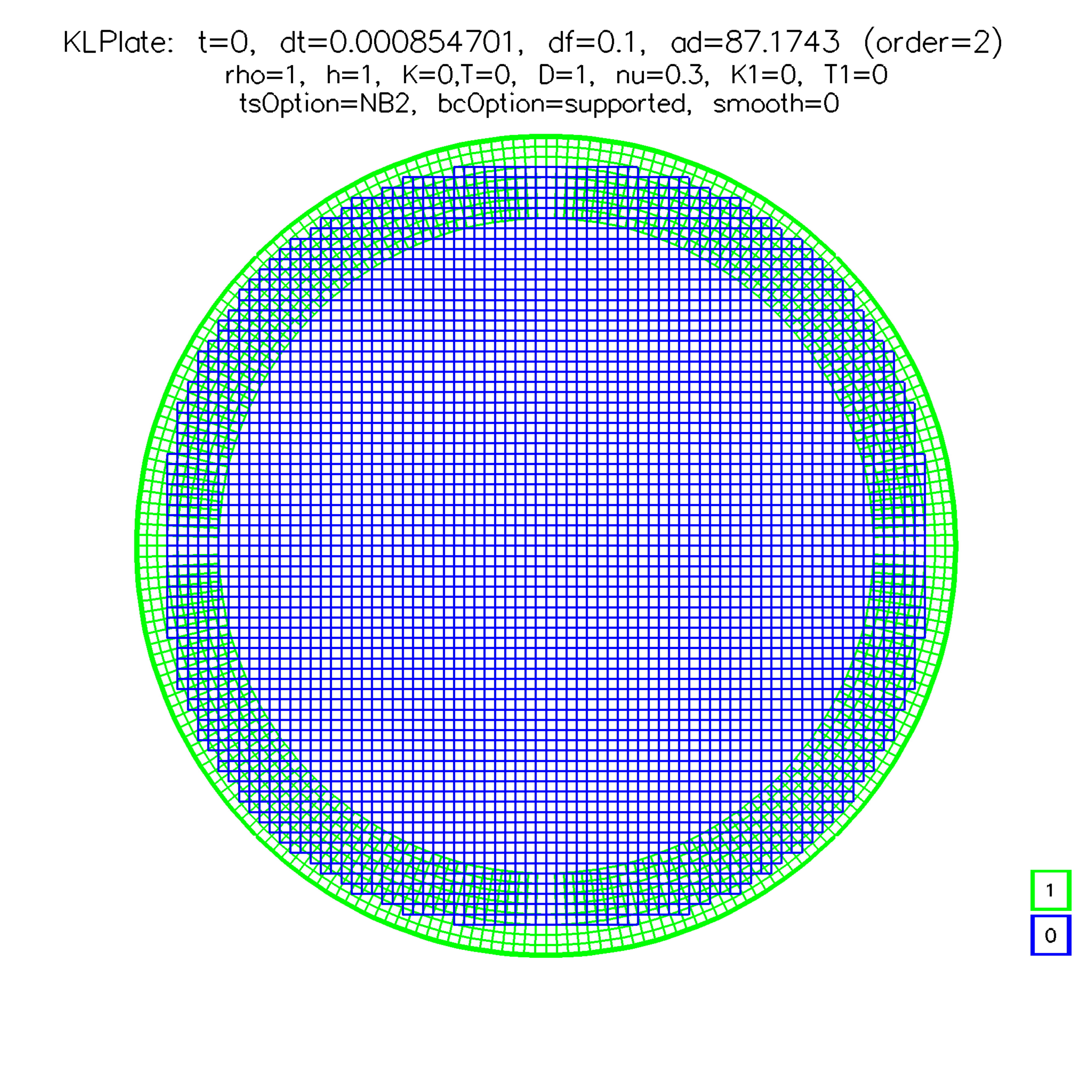}{\figWidth}};
\draw(3.75,0.0) node[anchor=south west,xshift=\xs cm,yshift=0pt] {\trimfig{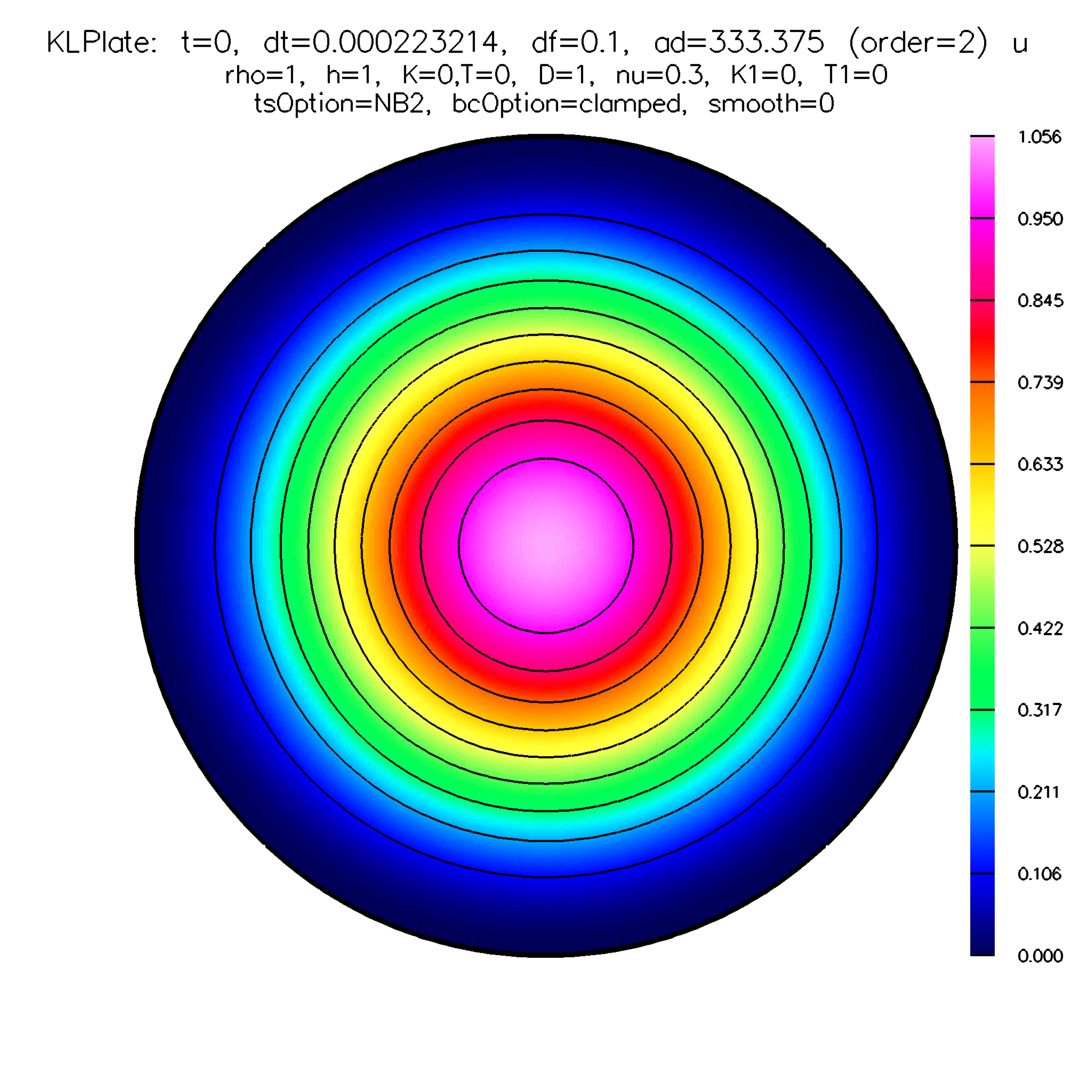}{\figWidth}};
\draw(7.75,0.0) node[anchor=south west,xshift=\xs cm,yshift=0pt] {\trimfig{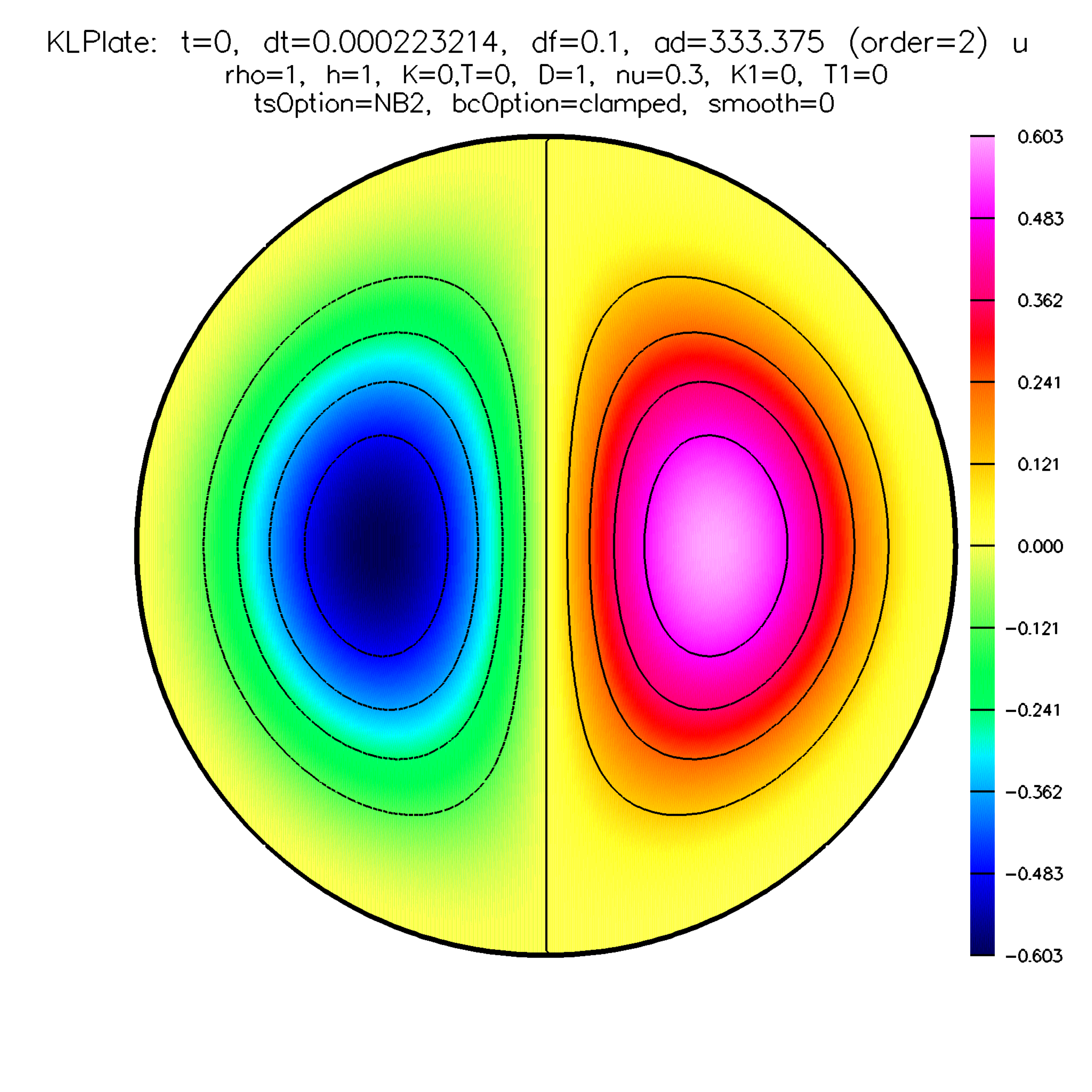}{\figWidth}};
\draw(11.75,0.0) node[anchor=south west,xshift=\xs cm,yshift=0pt] {\trimfig{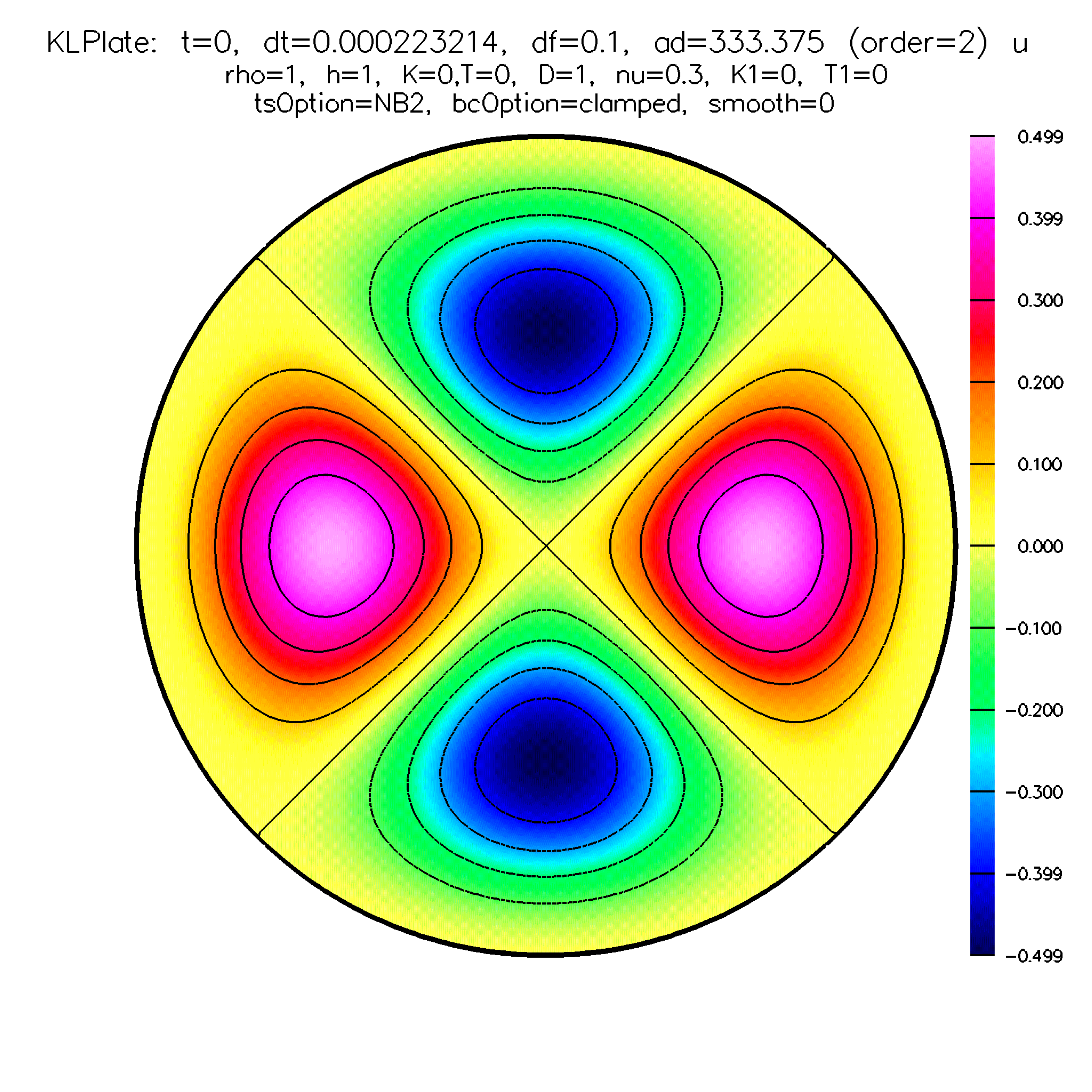}{\figWidth}};

\newcommand{\cbWidth}{.2}
\newcommand{\cbHeight}{3}

\drawColorBarH{fig/colourBarLines}{xshift=4.4cm,yshift=-.2cm}{\cbWidth}{\cbHeight}{$0$}{$1.056$}{}
\drawColorBarH{fig/colourBarLines}{xshift=8.4cm,yshift=-.2cm}{\cbWidth}{\cbHeight}{$-0.603$}{$0.603$}{}
\drawColorBarH{fig/colourBarLines}{xshift=12.4cm,yshift=-.2cm}{\cbWidth}{\cbHeight}{$-0.499$}{$0.499$}{}

\draw(2.,4.7)  node[anchor=north,xshift=0pt,yshift=0pt] { $\Gc_4$};
\draw(6.,4.7)  node[anchor=north,xshift=0pt,yshift=0pt] { $n=0$};
\draw(10,4.7)  node[anchor=north,xshift=0pt,yshift=0pt] { $n=1$ };
\draw(14.,4.7)  node[anchor=north,xshift=0pt,yshift=0pt] { $n=2$ };

%
\end{tikzpicture}

\end{center}
\caption{Computational grids $\Gc_4$ and  the initial displacements for $n=0,1,2$  subject to the clamped boundary conditions. Initial conditions for the supported boundary condition are similar with slightly different scales (omitted here  to save space) } \label{fig:dpnClampedICs}
\end{figure}
}

 The    values   of $\lambda$   depend on  the boundary conditions. 
 Now, we consider deriving  its values for the clamped  \eqref{eq:clampedBC}   and supported  \eqref{eq:supportedBC}   boundary conditions, respectively.

 \medskip
 \noindent {\bf  Clamped.}
 Clamped boundary  conditions in polar coordinates are  $$w=\pd{w}{r}=0.$$
 Enforcing them to the general solution  \eqref{eq:generalSolutionDiskPlate} leads to  
  $
  J_{n}(\lambda) I_{n+1}(\lambda)+ I_n(\lambda)J_{n+1}(\lambda) =0.
  $
The  values of $\lambda$ are determined by finding  the roots of  the following  transcendental function,
  \begin{equation}\label{eq:transcendentalFunctionClapmed}
  \phi^c_n(\lambda) =  J_{n}(\lambda) I_{n+1}(\lambda)+ I_n(\lambda)J_{n+1}(\lambda). 
 \end{equation}
  For each $n$,  the transcendental function \eqref{eq:transcendentalFunctionClapmed}   has  infinite number of roots. 
For the purpose of designing  numerical tests, we  take  the smallest root of \eqref{eq:transcendentalFunctionClapmed} for $n=0,1$ and $2$ (c.f.,  Figure~\ref{fig:transcendentalFunction_c}.), which  are 
  \begin{equation}\label{eq:lambdac}
  \lambda_0 \approx 3.196220616582554,\quad \lambda_1 \approx   4.610899879386510,\quad \lambda_2 \approx   5.905678237243653.
 \end{equation}

 \medskip
 \noindent {\bf  Supported.} The supported boundary conditions in polar coordinates are
 $$w=\pdn{w}{r}{2}+\nu\left(\frac{1}{r}\pd{w}{r}+\frac{1}{r^2}\pdn{w}{\theta}{2}\right)=0.$$  
    Similarly, enforcing them to the general solution leads to the definition of the following  transcendental function,
  $$
  \phi^s_n(\lambda) =  J_{n}(\lambda) I_{n+1}(\lambda)+ I_n(\lambda)J_{n+1}(\lambda)-\frac{2\lambda}{1-\nu} J_{n}(\lambda) I_{n}(\lambda).
  $$
  Assuming  $\nu=0.3$,  we plot the first three $ \phi^s_n(\lambda) $ in Figure~\ref{fig:transcendentalFunction_s}, whose   smallest roots  are 
  \begin{equation}\label{eq:lambdas}
  \lambda_0 \approx  2.221519534965056 ,\quad \lambda_1 \approx    3.728024285469852 ,\quad \lambda_2 \approx    5.060958083288190
  \end{equation}






{ 
  \newcommand{\figWidth}{4cm}
  \newcommand{\figWidthLgd}{1.8cm}
\def\xa{16}
\def\ya{21}
\newcommand{\trimfig}[2]{\trimw{#1}{#2}{0.}{0.}{0.}{0.}}
\newcommand{\trimfigLgd}[2]{\trimw{#1}{#2}{0.675}{0.}{0.09}{0.6}}

\begin{figure}[h]
\begin{center}
\begin{tikzpicture}[scale=1]
  \useasboundingbox (0.0,0.0) rectangle (\xa,\ya);  

\draw(-0.5,14) node[anchor=south west,xshift=0pt,yshift=0pt] {\trimfig{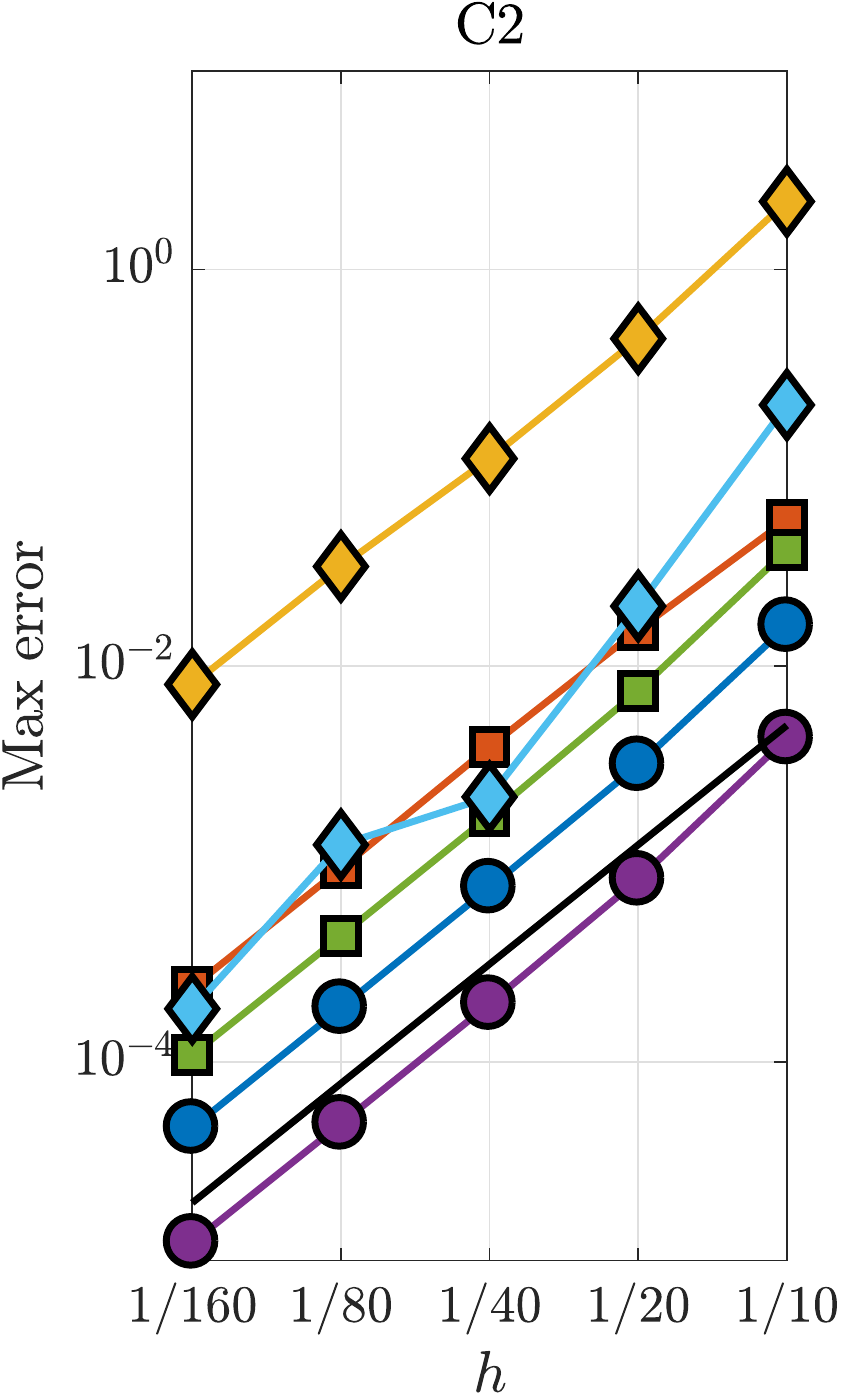}{\figWidth}};
\draw(3.5,14) node[anchor=south west,xshift=0pt,yshift=0pt] {\trimfig{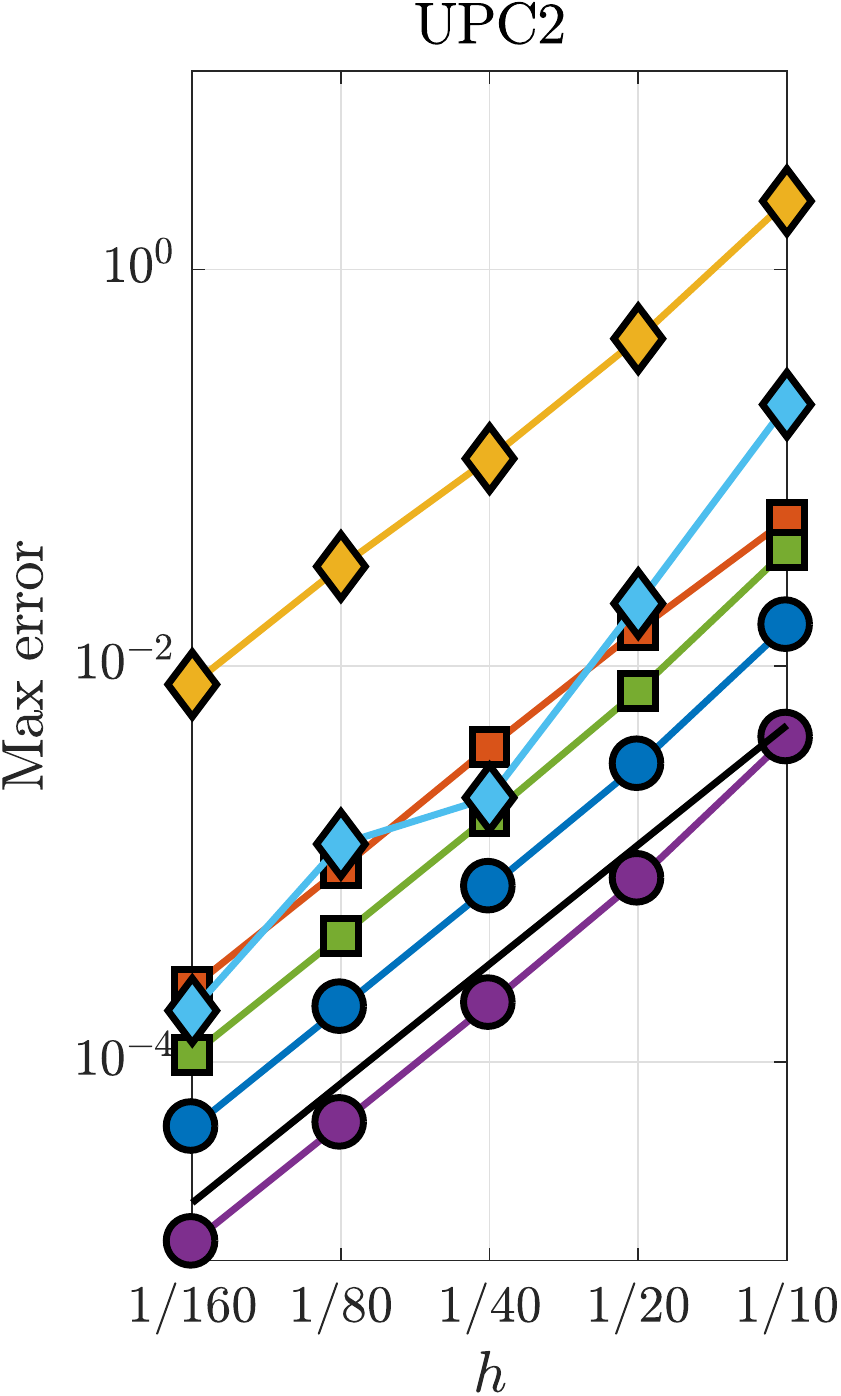}{\figWidth}};
\draw(7.5,14) node[anchor=south west,xshift=0pt,yshift=0pt] {\trimfig{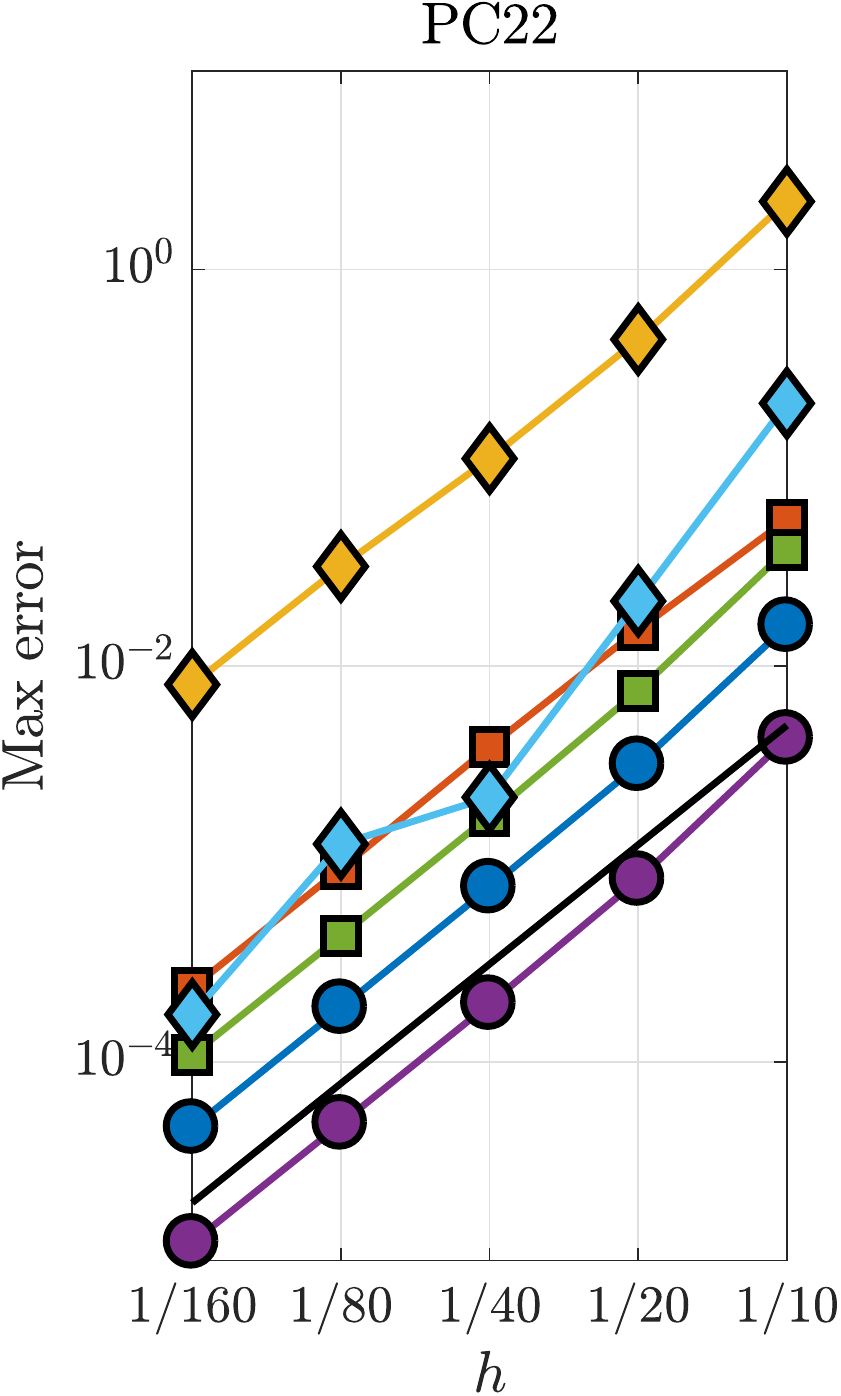}{\figWidth}};
\draw(11.5,14) node[anchor=south west,xshift=0pt,yshift=0pt] {\trimfig{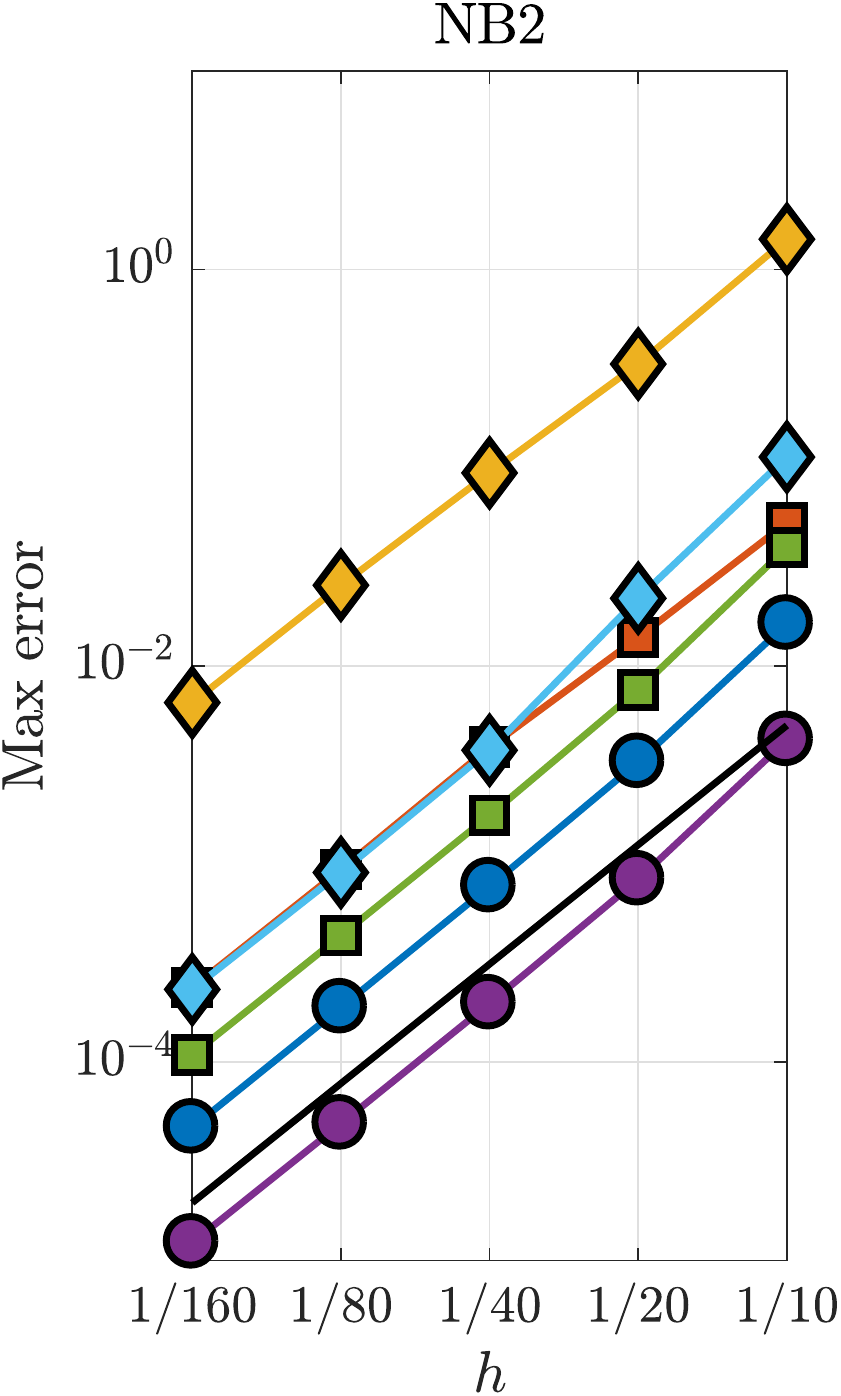}{\figWidth}};
\draw(8,20.75) node[anchor=south west,xshift=0pt,yshift=0pt]  {n=0};

\draw(-0.5,7) node[anchor=south west,xshift=0pt,yshift=0pt] {\trimfig{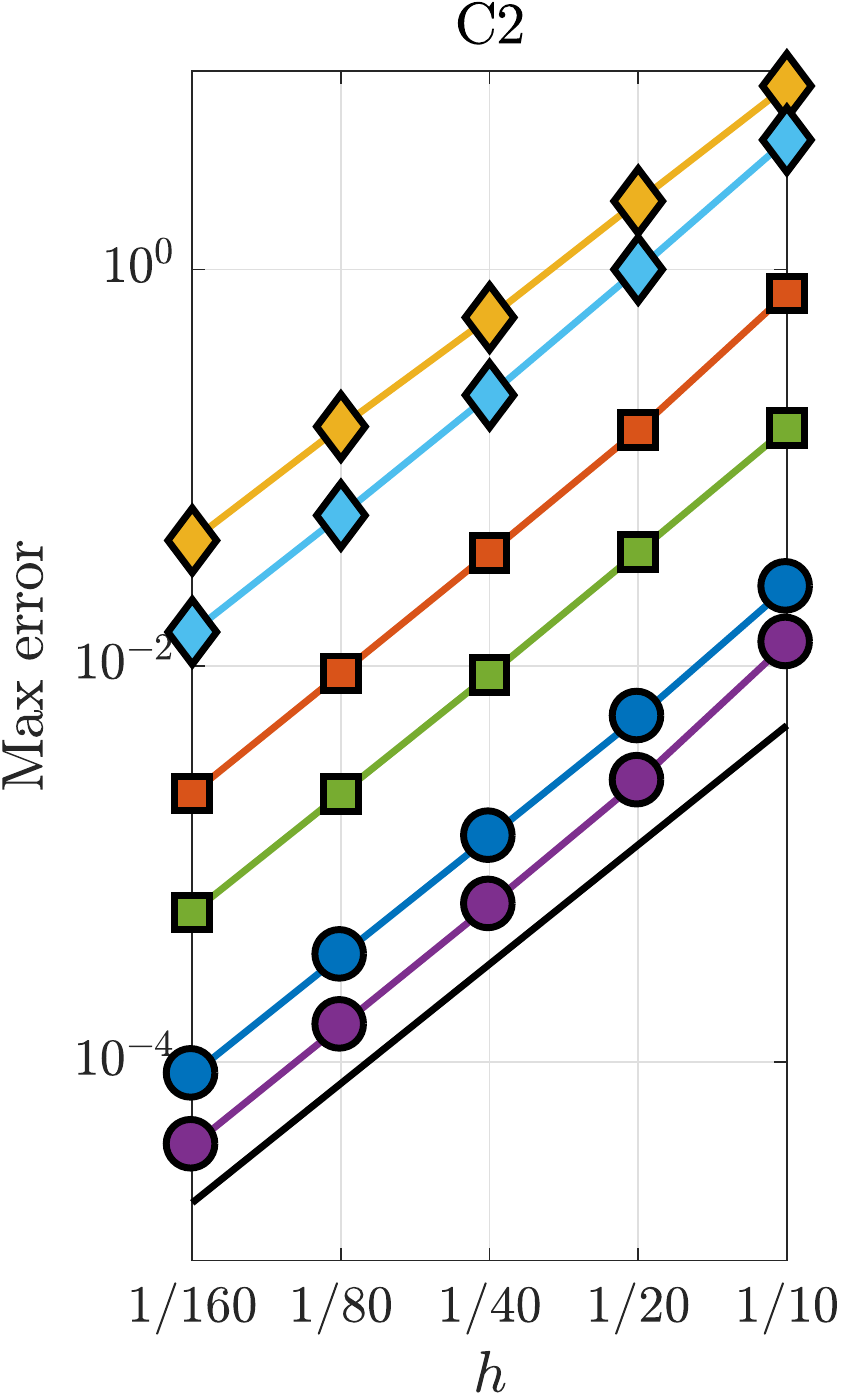}{\figWidth}};
\draw(3.5,7) node[anchor=south west,xshift=0pt,yshift=0pt] {\trimfig{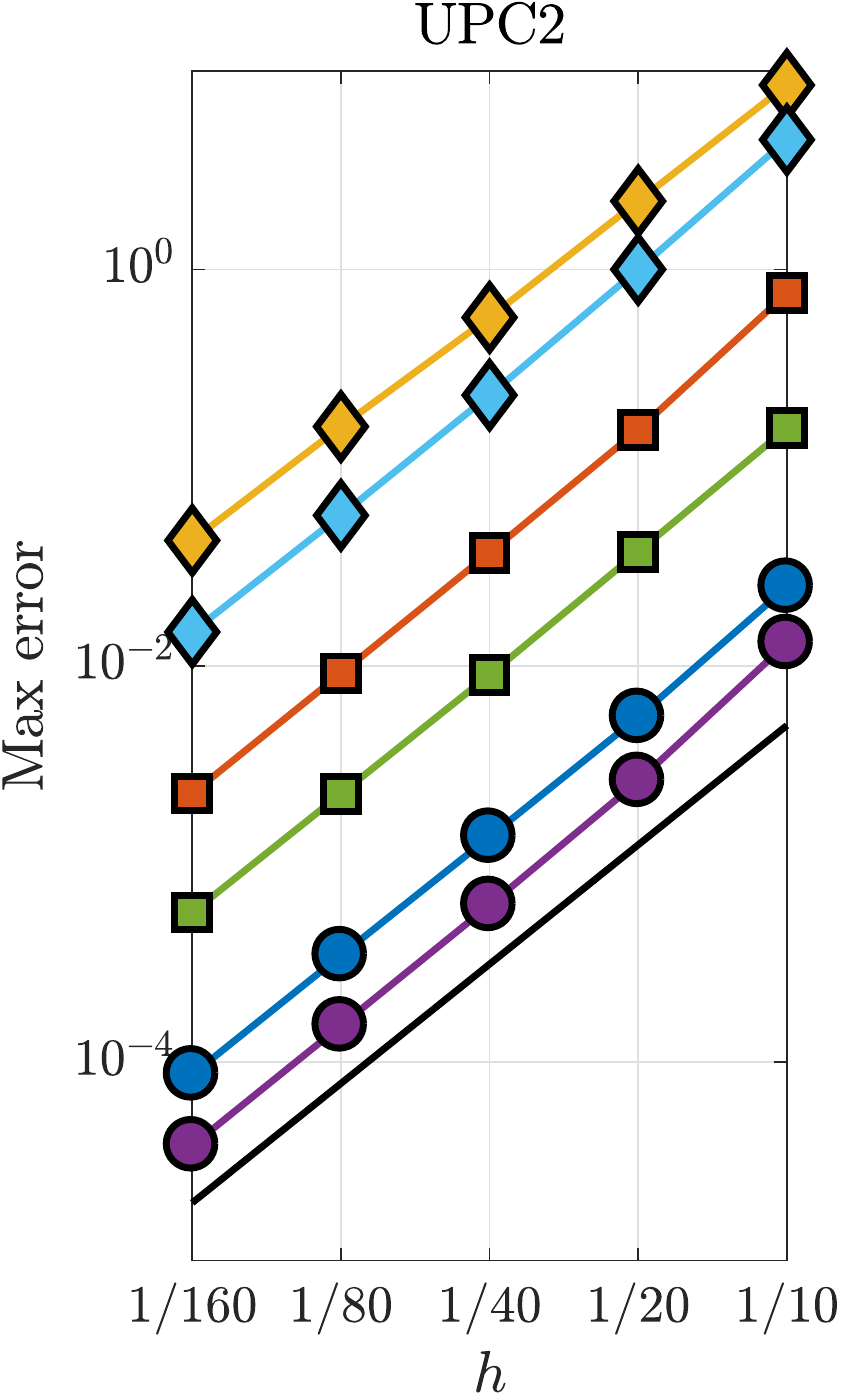}{\figWidth}};
\draw(7.5,7) node[anchor=south west,xshift=0pt,yshift=0pt] {\trimfig{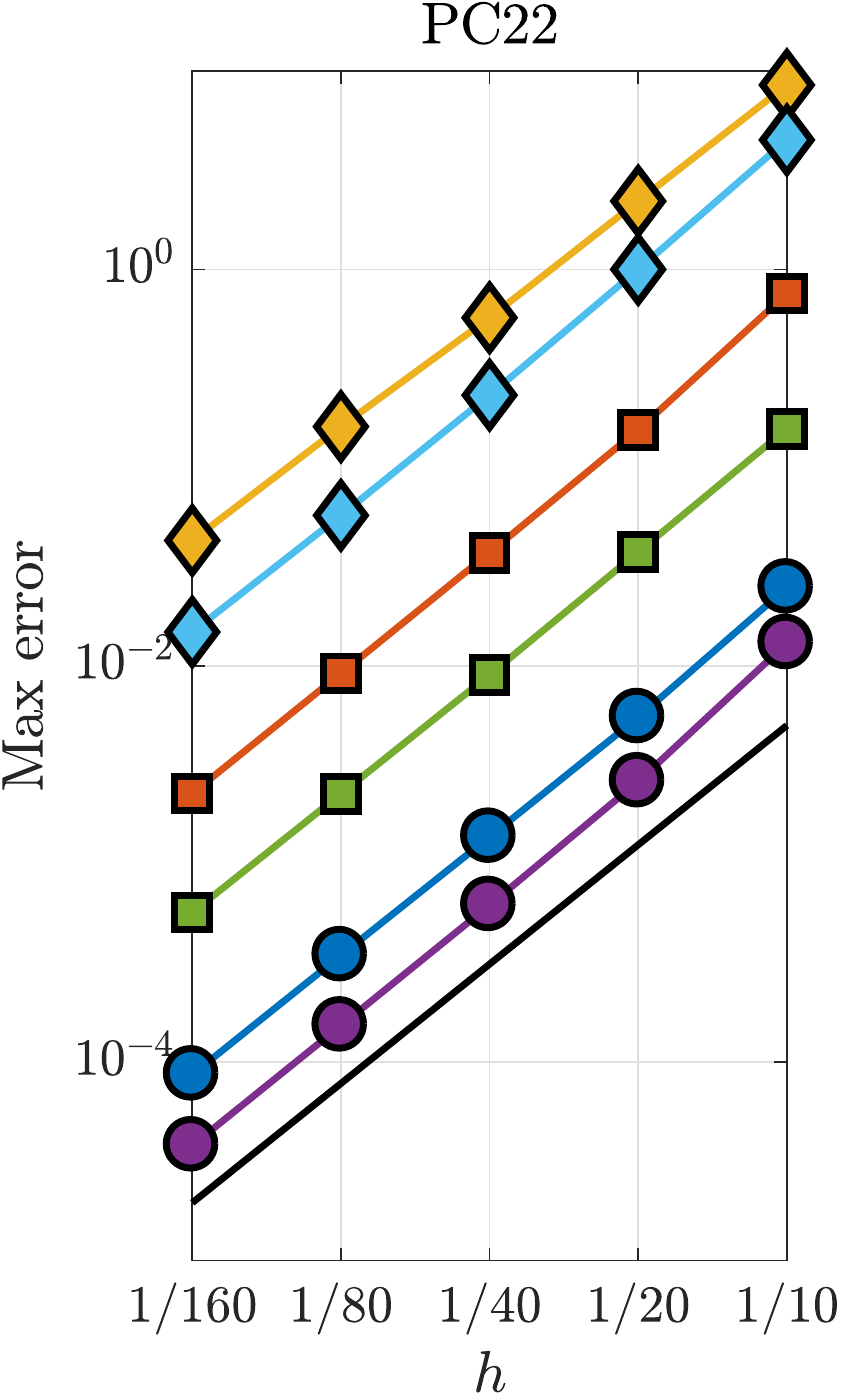}{\figWidth}};
\draw(11.5,7) node[anchor=south west,xshift=0pt,yshift=0pt] {\trimfig{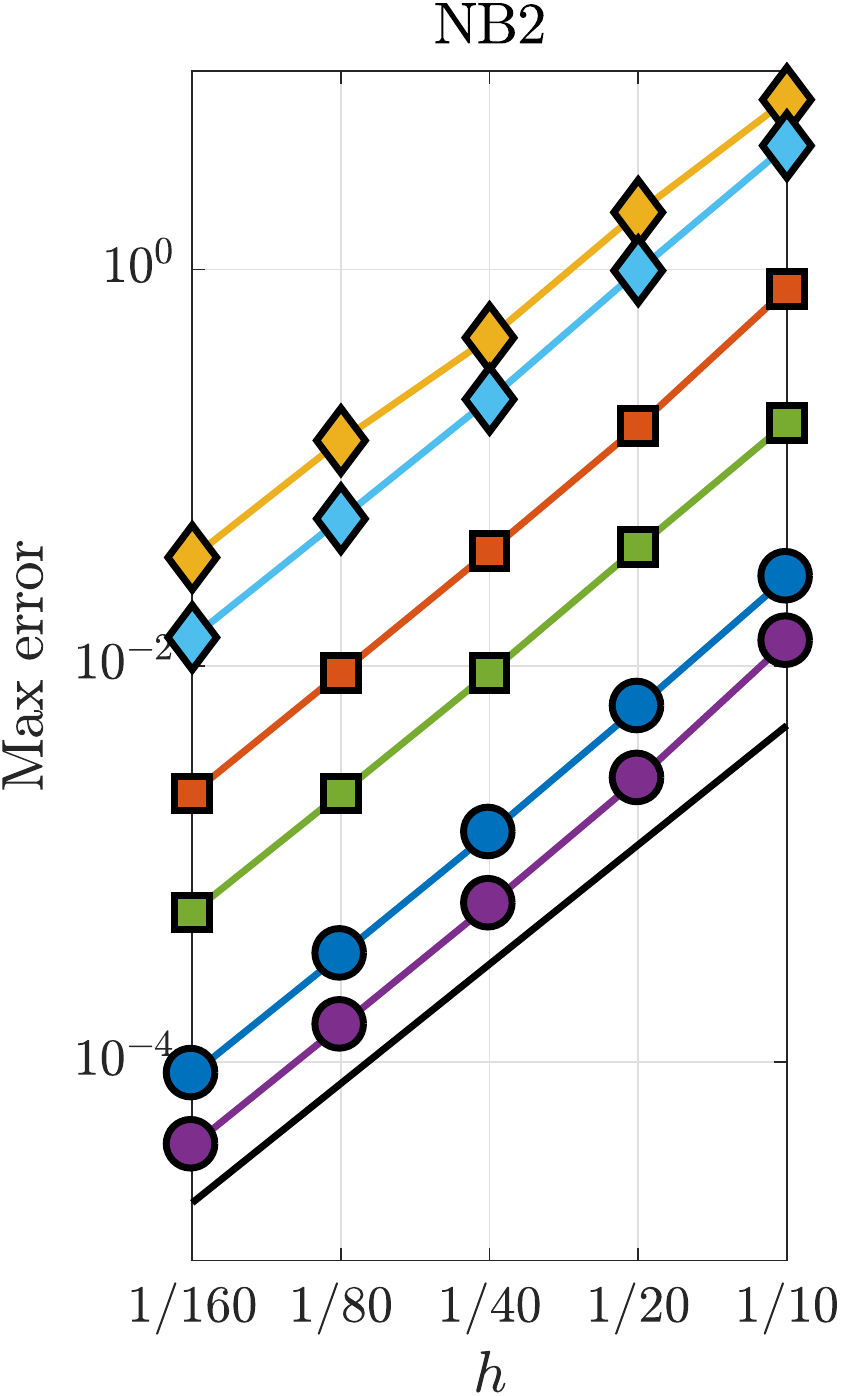}{\figWidth}};
\draw(8,13.75) node[anchor=south west,xshift=0pt,yshift=0pt]  {n=1};

\draw(-0.5,0.0) node[anchor=south west,xshift=0pt,yshift=0pt] {\trimfig{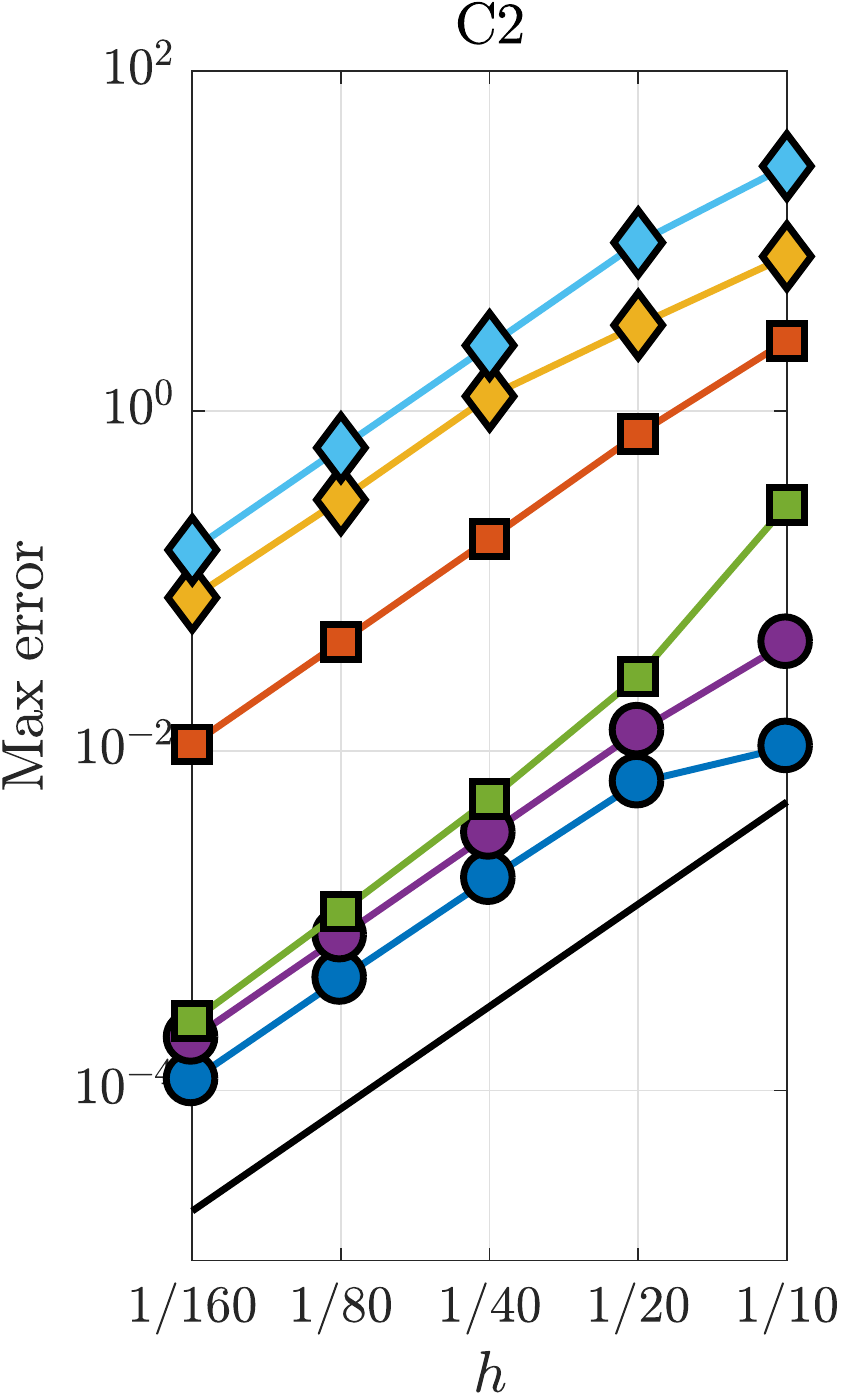}{\figWidth}};
\draw(3.5,0.0) node[anchor=south west,xshift=0pt,yshift=0pt] {\trimfig{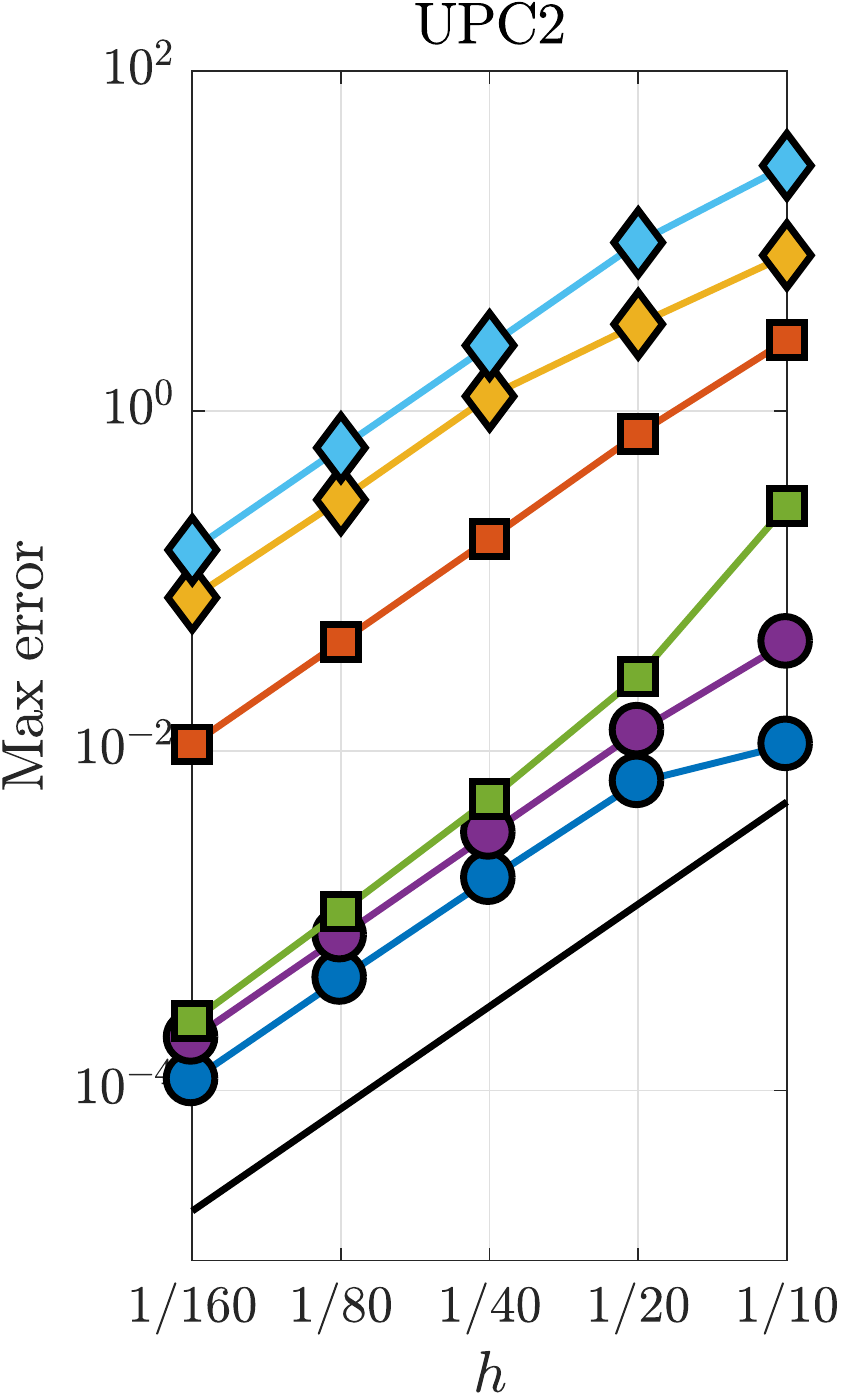}{\figWidth}};
\draw(7.5,0.0) node[anchor=south west,xshift=0pt,yshift=0pt] {\trimfig{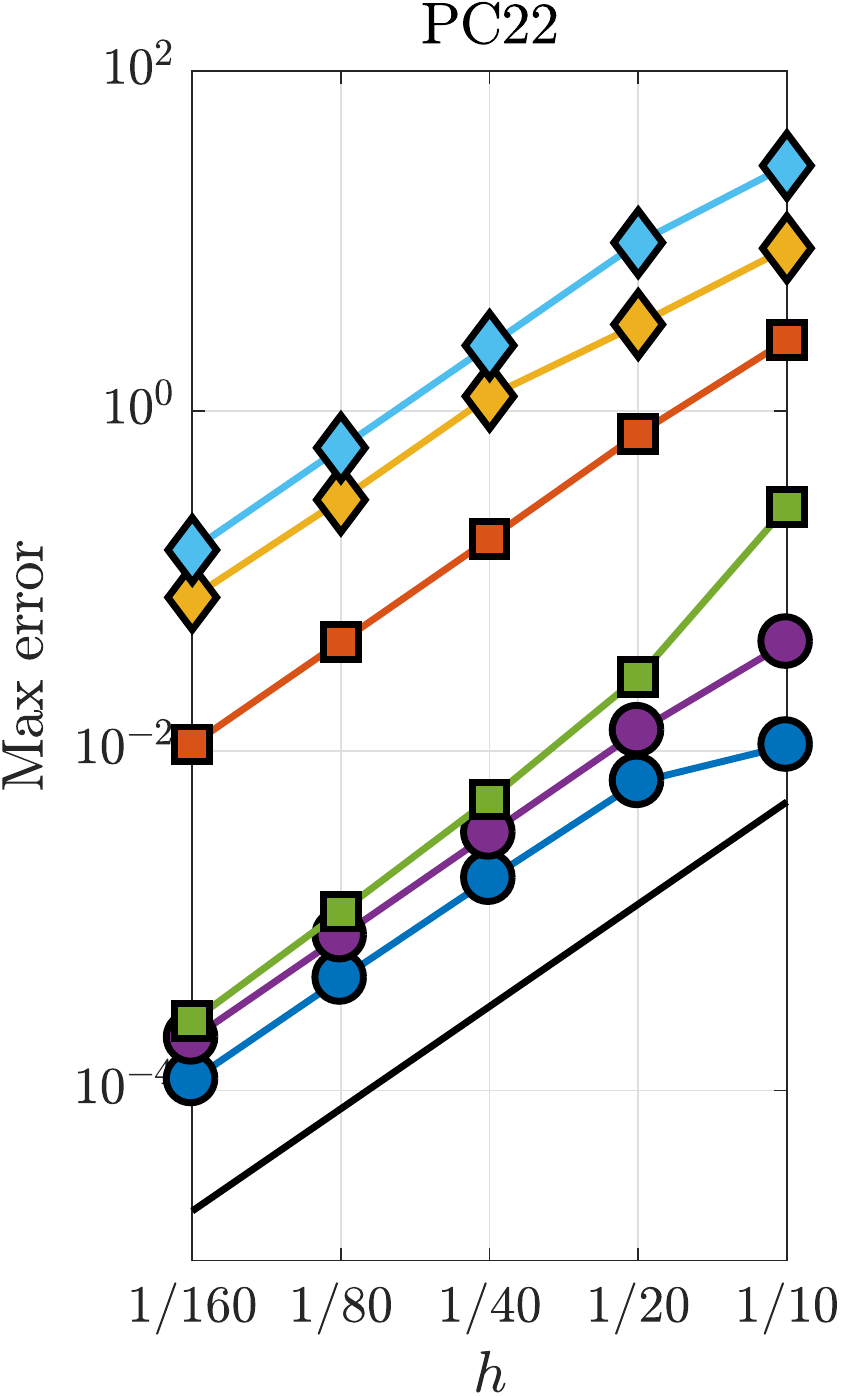}{\figWidth}};
\draw(11.5,0.0) node[anchor=south west,xshift=0pt,yshift=0pt] {\trimfig{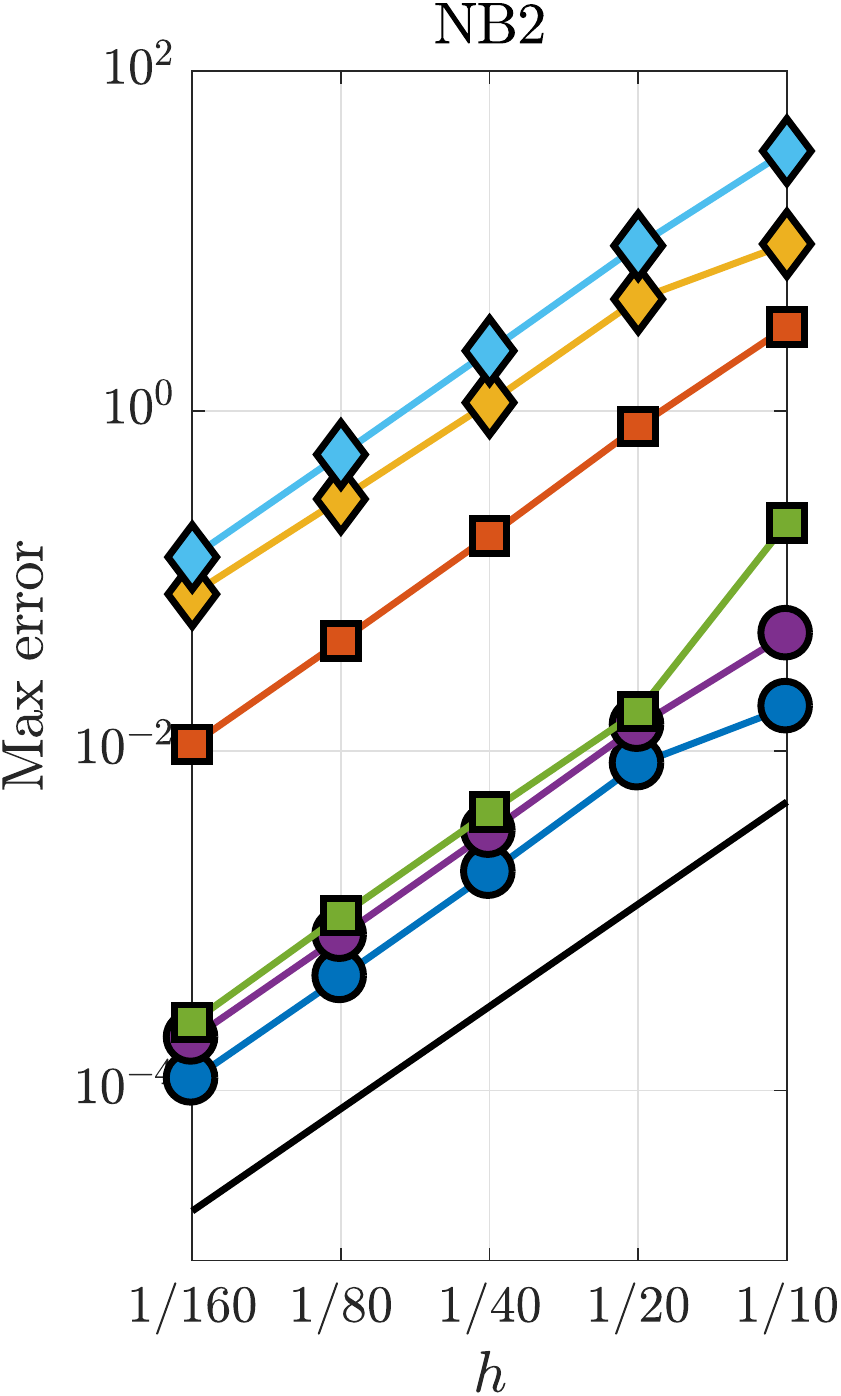}{\figWidth}};
\draw(8,6.75) node[anchor=south west,xshift=0pt,yshift=0pt]  {n=2};

\draw(-0.5,18.9) node[anchor=south west,xshift=0pt,yshift=0pt] {\trimfigLgd{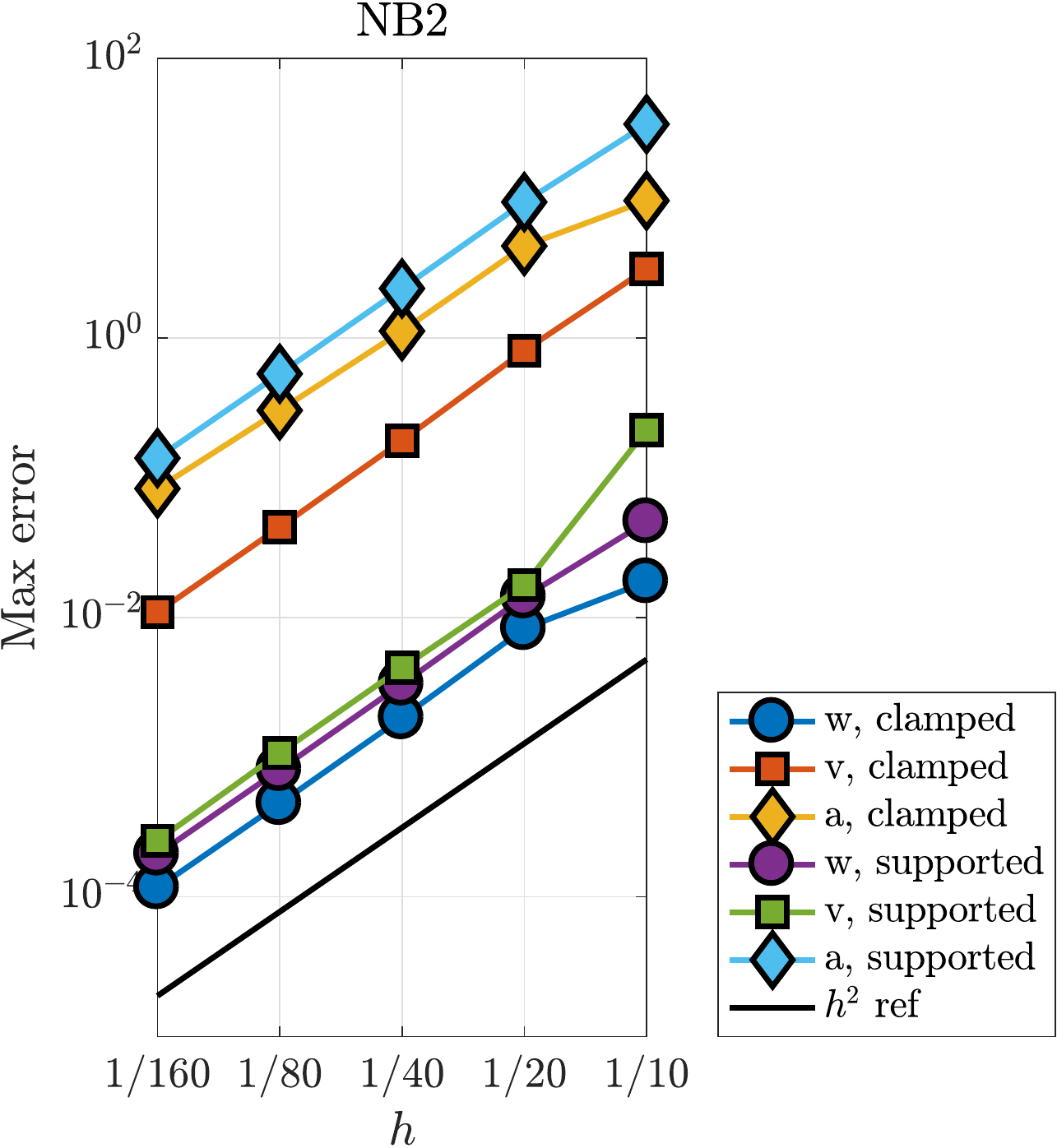}{\figWidthLgd}};


%
\end{tikzpicture}

\end{center}
\vspace{-0.3in}
\caption{Convergence rates for the cases with  the Bessel's functions of order $n=0$ (top row), $n=1$ (middle row) and $n=2$ (bottom row).} \label{fig:dpnRates}
\end{figure}
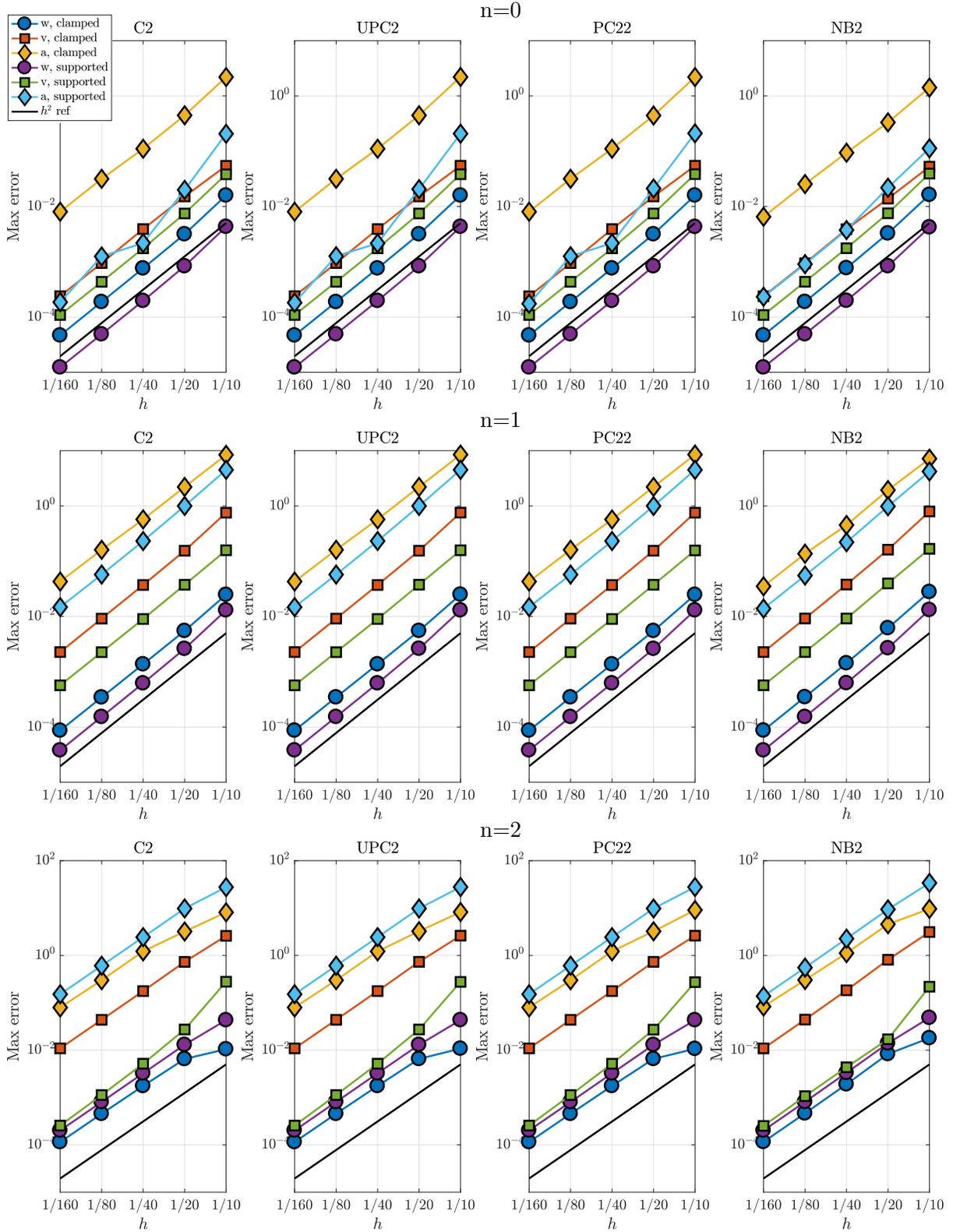
}

We note that the general solution \eqref{eq:generalSolutionDiskPlate}  consists of two sets of modal functions: one with $\cos(n\theta)$ and the other with  $\sin(n\theta)$.  To construct a standing wave solution, we specify the following initial conditions using the modal function with $\cos(n\theta)$ only;  that is, 
\begin{equation}\label{eq:circularPlateIC}
  w_0(r,\theta)= \left[J_n\left(\frac{\lambda_n r}{a}\right)-\frac{J_n(\lambda_n)}{I_n(\lambda_n)}I_n\left(\frac{\lambda_n r}{a}\right)\right] \cos(n\theta) \quad\text{and} \quad v_0(r,\theta)=0.
\end{equation}
The exact solution to the  plate equation  \eqref{eq:KLShell} with this choice of initial conditions is
\begin{equation}\label{eq:circularPlateExact}
w_e(r,\theta,t)=\left[J_n\left(\frac{\lambda_n r}{a}\right)-\frac{J_n(\lambda_n)}{I_n(\lambda_n)}I_n\left(\frac{\lambda_n r}{a}\right)\right] \cos(n\theta)\cos\left(\lambda_n^2\sqrt{\frac{D}{\rho H}}t\right).
\end{equation}

{
\newcommand{\figWidth}{4.5cm}
\def\xa{16}
\def\ya{10.5}
\newcommand{\trimfig}[2]{\trimw{#1}{#2}{0.12}{0.12}{0.12}{0.12}}
\begin{figure}[h!]
\begin{center}
\begin{tikzpicture}[scale=1]
  \useasboundingbox (0.0,0.0) rectangle (\xa,\ya);  

\def\xs{0.3}

\draw(1,5.5) node[anchor=south west,xshift=\xs cm,yshift=0pt] {\trimfig{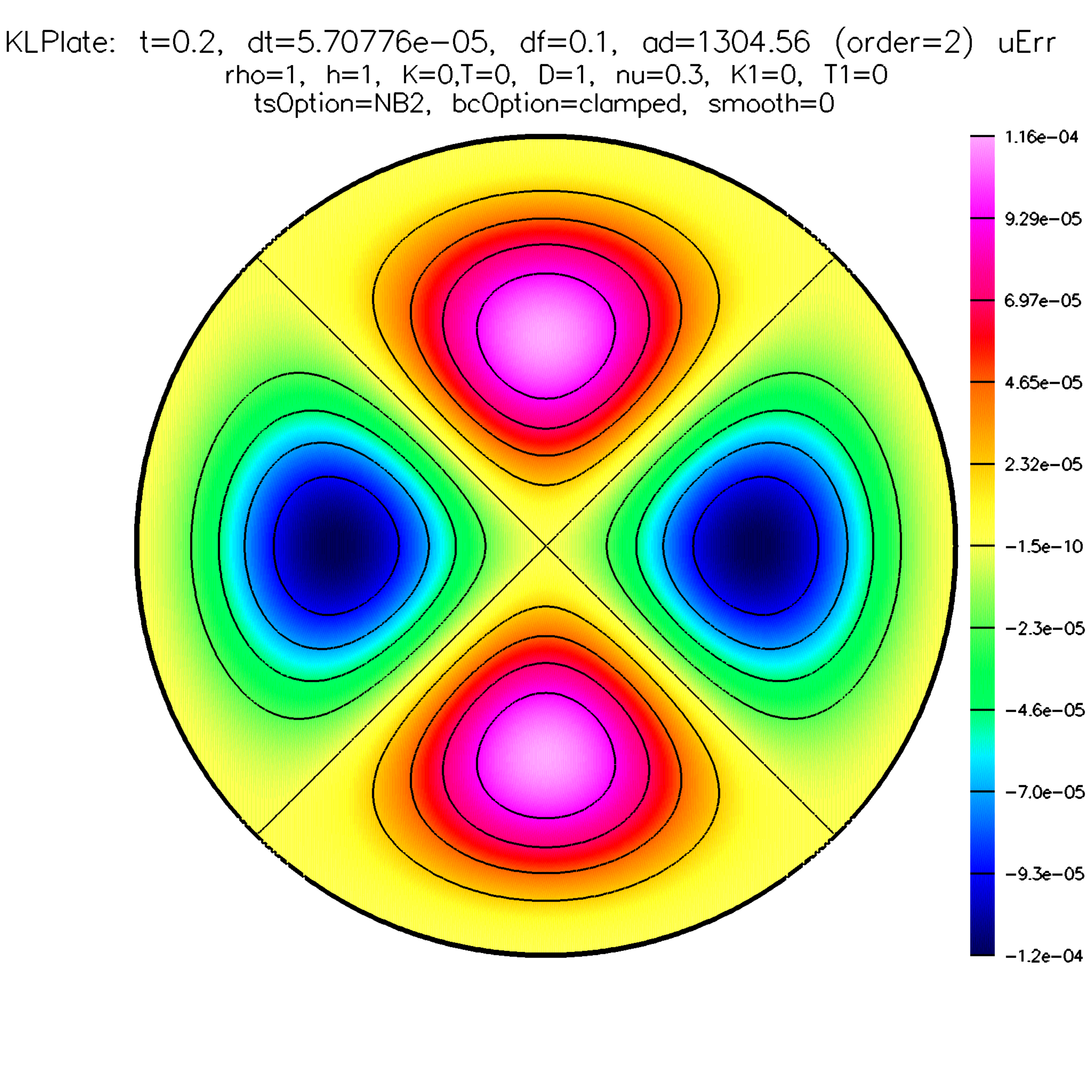}{\figWidth}};
\draw(5.75,5.5) node[anchor=south west,xshift=\xs cm,yshift=0pt] {\trimfig{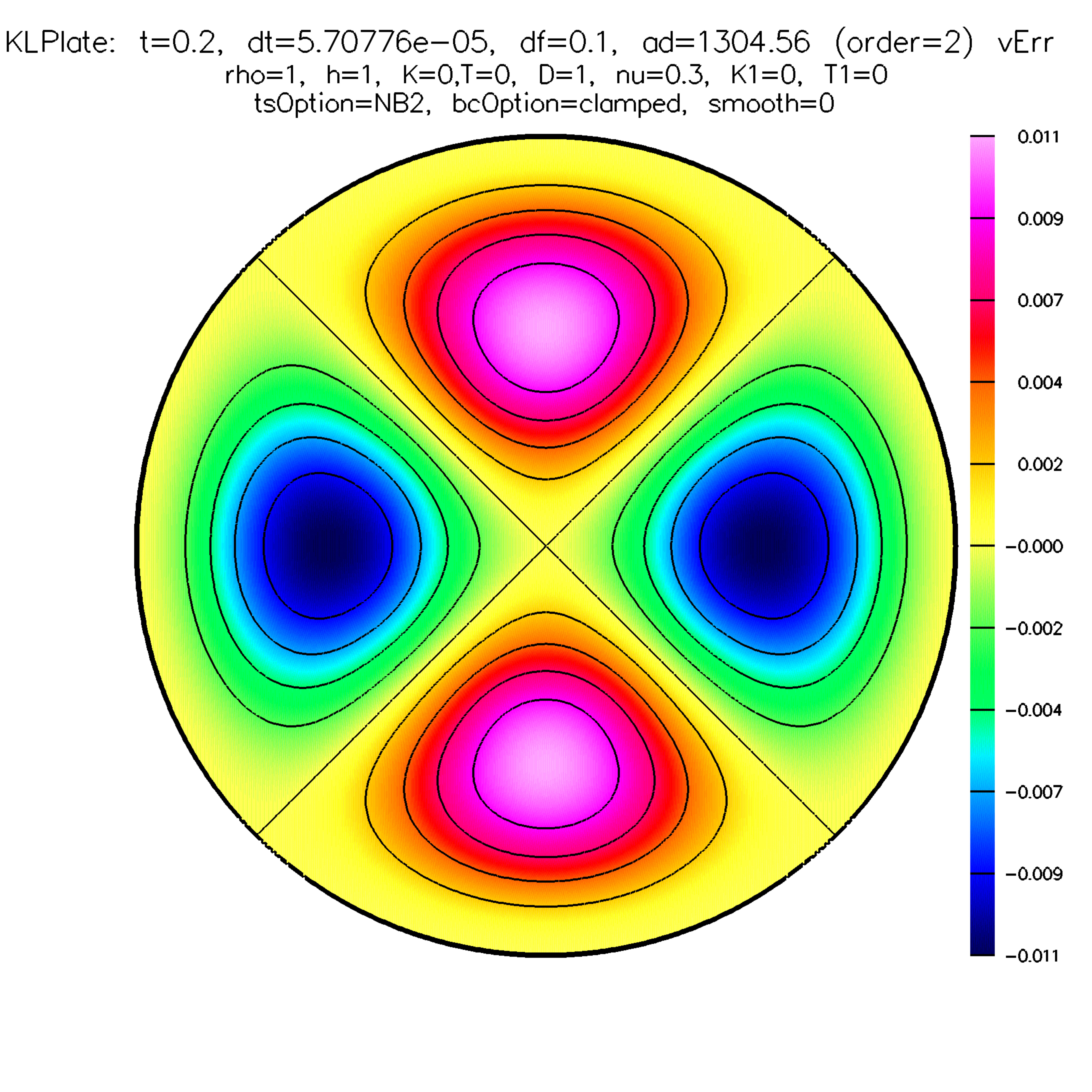}{\figWidth}};
\draw(10.5,5.5) node[anchor=south west,xshift=\xs cm,yshift=0pt] {\trimfig{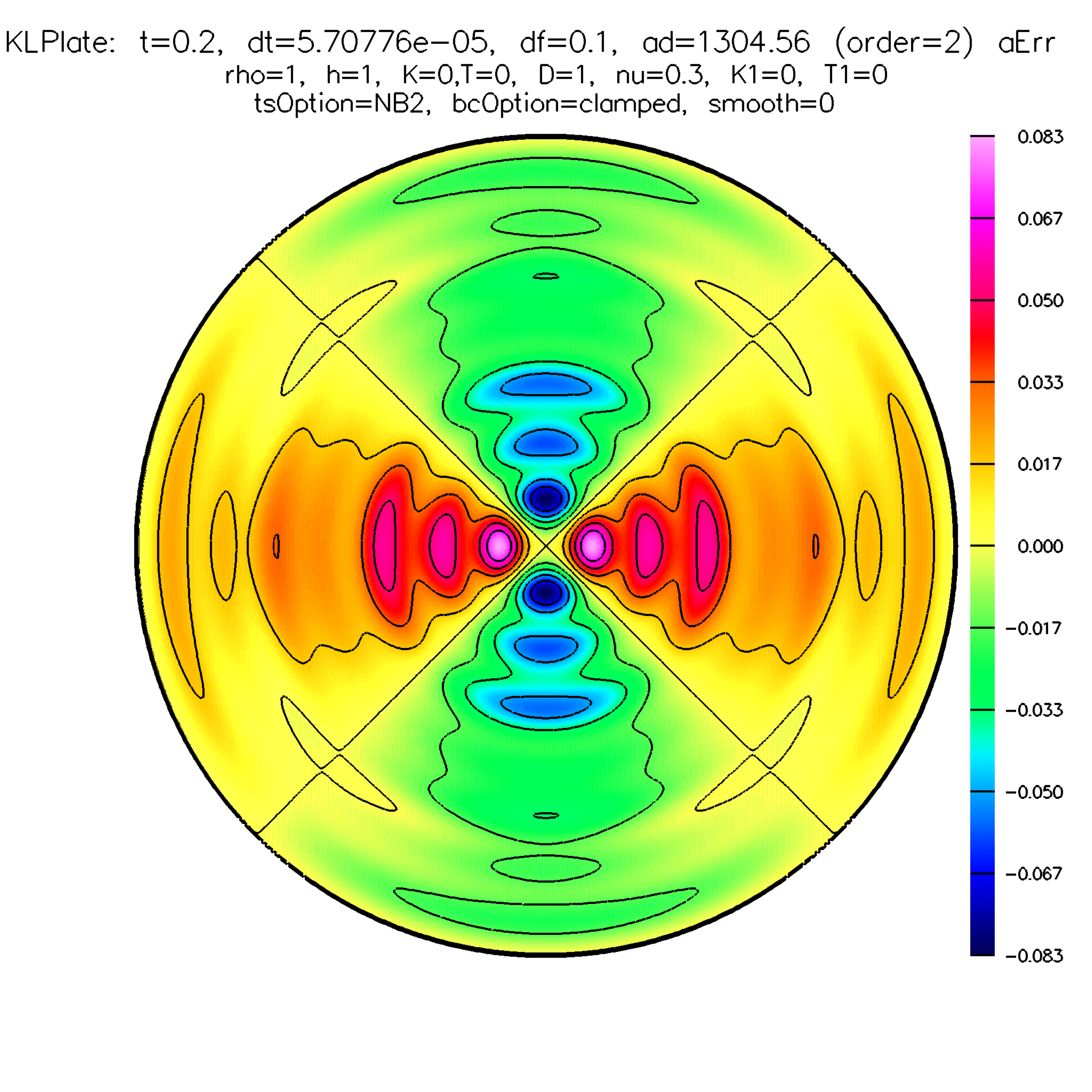}{\figWidth}};

\draw(1,0.0) node[anchor=south west,xshift=\xs cm,yshift=0pt] {\trimfig{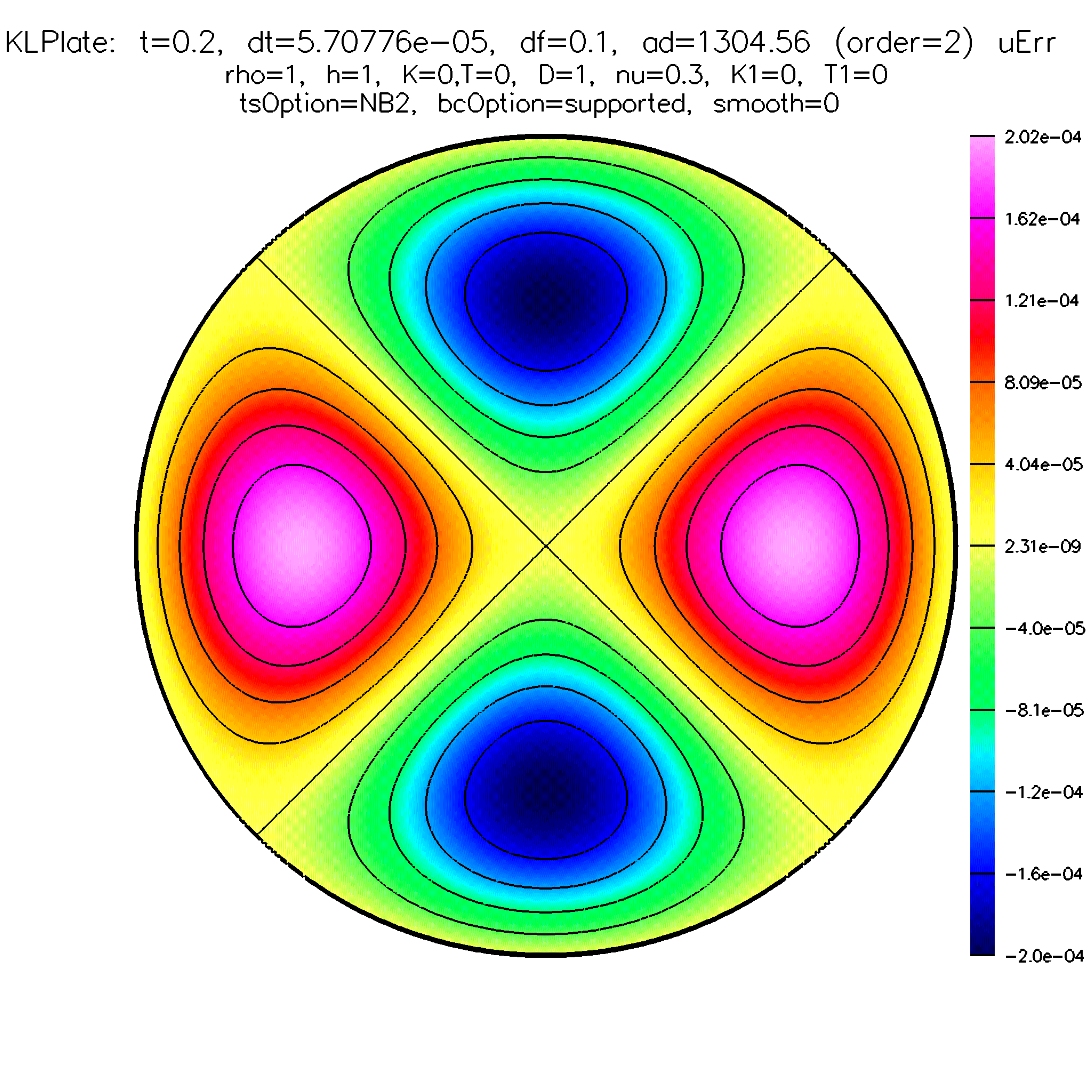}{\figWidth}};
\draw(5.75,0.0) node[anchor=south west,xshift=\xs cm,yshift=0pt] {\trimfig{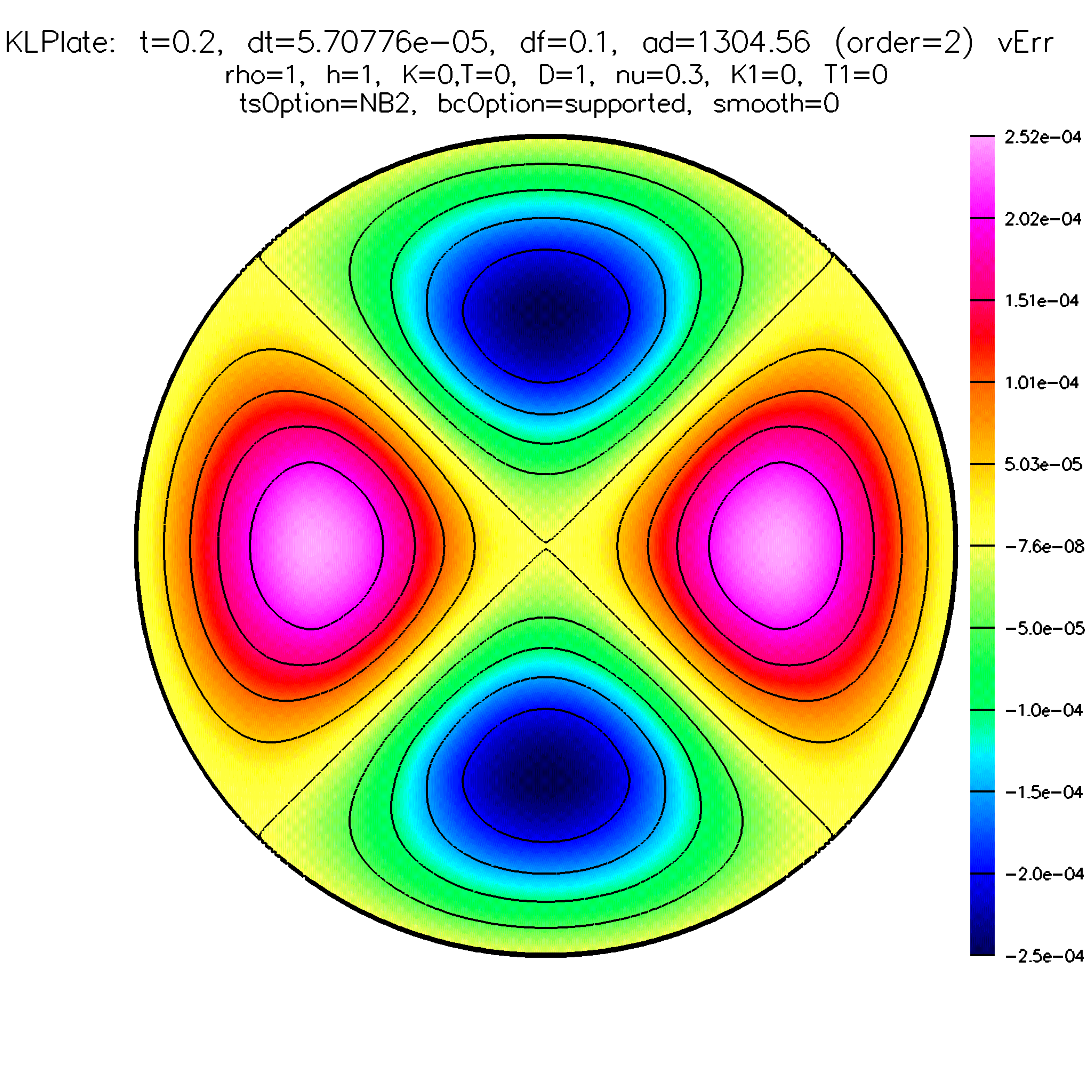}{\figWidth}};
\draw(10.5,0.0) node[anchor=south west,xshift=\xs cm,yshift=0pt] {\trimfig{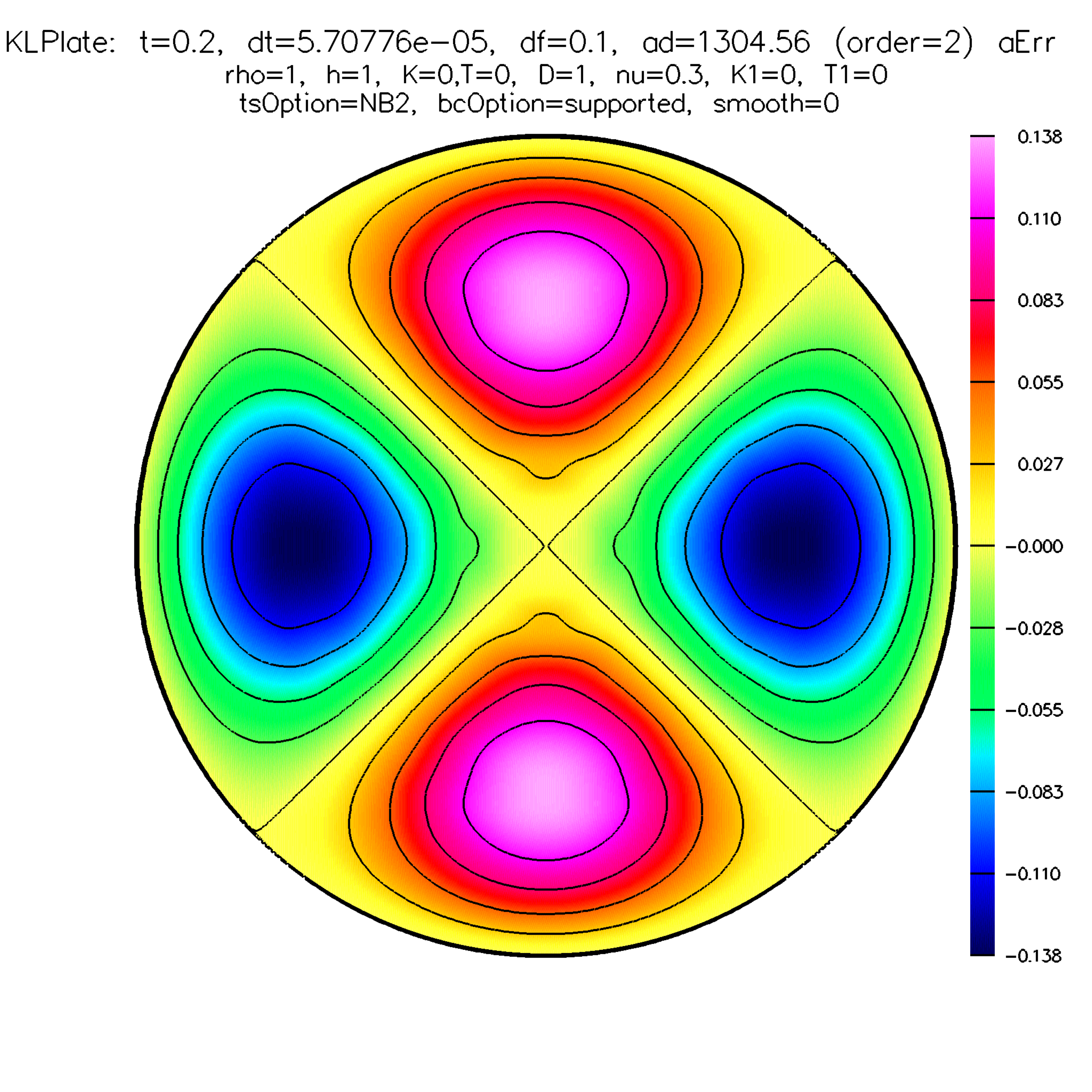}{\figWidth}};

\newcommand{\cbWidth}{.2}
\newcommand{\cbHeight}{3}
\drawColorBarH{fig/colourBarLines}{xshift=2.45cm,yshift=5.3cm}{\cbWidth}{\cbHeight}{$-1.2e-4$}{$1.16e-4$}{}
\drawColorBarH{fig/colourBarLines}{xshift=7.2cm,yshift=5.3cm}{\cbWidth}{\cbHeight}{$-0.011$}{$0.011$}{}
\drawColorBarH{fig/colourBarLines}{xshift=11.95cm,yshift=5.3cm}{\cbWidth}{\cbHeight}{$-0.083$}{$0.083$}{}

\drawColorBarH{fig/colourBarLines}{xshift=2.45cm,yshift=-0.2cm}{\cbWidth}{\cbHeight}{$-2.0e-4$}{$2.02e-4$}{}
\drawColorBarH{fig/colourBarLines}{xshift=7.2cm,yshift=-0.2cm}{\cbWidth}{\cbHeight}{$-2.5e-4$}{$2.52e-4$}{}
\drawColorBarH{fig/colourBarLines}{xshift=11.95cm,yshift=-0.2cm}{\cbWidth}{\cbHeight}{$-0.138$}{$0.138$}{}

\draw(3.7,10.8)  node[anchor=north,xshift=\xs cm,yshift=0pt] {$E(w)(x,y,0.2)$};
\draw(8.6,10.8)  node[anchor=north,xshift=\xs cm,yshift=0pt] {$E(v)(x,y,0.2)$};
\draw(13.3,10.8)  node[anchor=north,xshift=\xs cm,yshift=0pt] {$E(a)(x,y,0.2)$};

\draw(1.85,8)  node[anchor=east,xshift=0cm,yshift=0pt] {\footnotesize\bf Clamped};
\draw(1.85,2.5)  node[anchor=east,xshift=0cm,yshift=0pt] {\footnotesize\bf Supported};

%
\end{tikzpicture}
\end{center}
\caption{ Errors of the  NB2 scheme on grid $\Gc_{16}$  at $t=0.2$  for the case  $n=2$; errors of all the other schemes are similar. } \label{fig:dp2NB2Errors}
\end{figure}
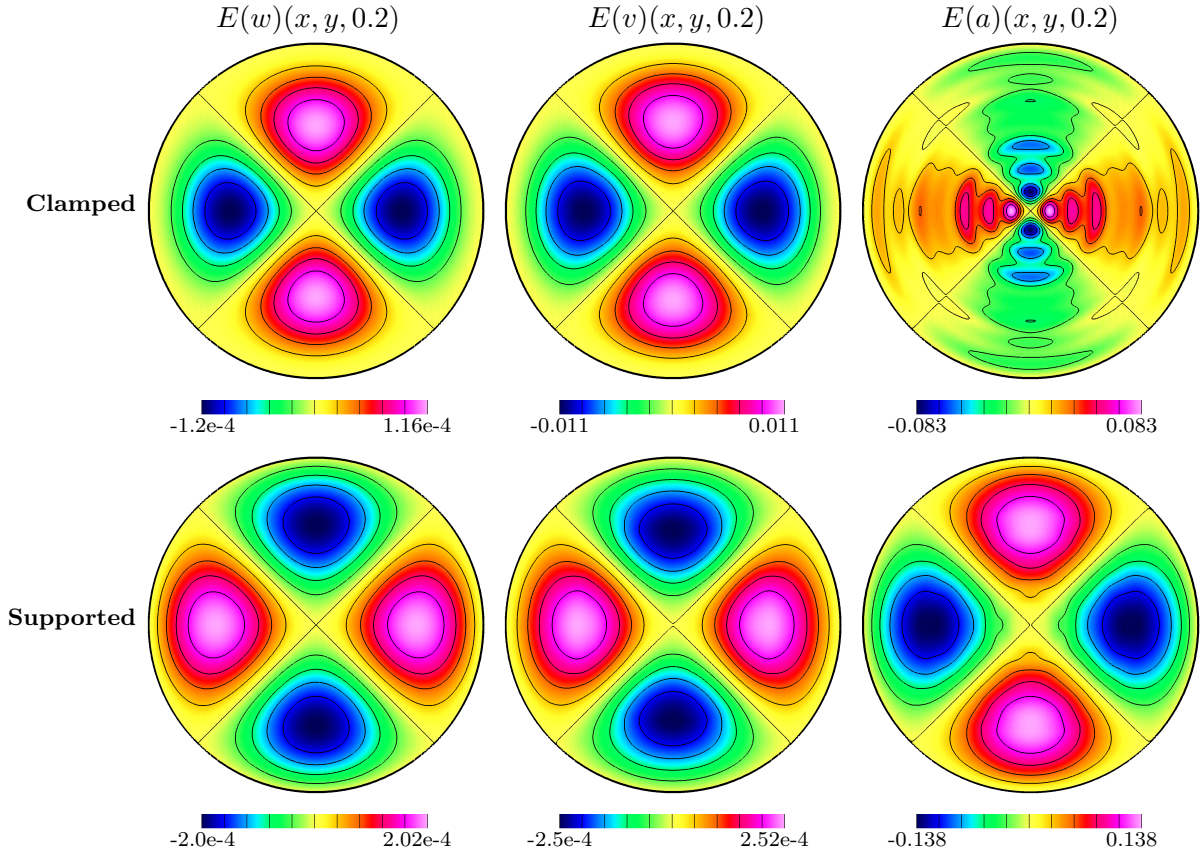
}

Prescribing the initial conditions \eqref{eq:circularPlateIC} for $n=0,1,2$, we solve the circular plate with clamped and supported boundaries numerically, and then   compare the numerical results with the corresponding  exact solution \eqref{eq:circularPlateExact} to verify the accuracy and stability of our schemes. The plate domain is assumed to be a unit circle (i.e., $a=1$); and the complete  parameter values used for this test are  $\rho=H=D=1$, $K=T=f=0$, and $\nu=0.3$.  The values of $\lambda_n$ needed by the exact solution \eqref{eq:circularPlateExact} and the initial conditions \eqref{eq:circularPlateIC}  are given in \eqref{eq:lambdac} for clamped  and in \eqref{eq:lambdas} for supported boundary conditions, respectively. In Figure~\ref{fig:dpnClampedICs}, we show the computational grids $\Gc_4$ and  the initial displacements for $n=0,1,2$  subject to the clamped boundary conditions. The initial  conditions for the supported boundary conditions are similar with slightly different scales, so their plots  are omitted here to save space.

Convergence rates for  all   the  cases are collected in Figure~\ref{fig:dpnRates}.  
It is clearly shown in this numerical experiment that the expected 2nd-order accuracy for all the solution components $(w, v, a)$ is consistently achieved by our methods proposed in this paper.  To illustrate the error behavior over the entire domain, we plot  the errors of the NB2 scheme on grid $\Gc_{16}$  at $t=0.2$  for the case  $n=2$ in Figure~\ref{fig:dp2NB2Errors}. Again, we see that  the errors of all the solution components ($w$, $v$, $a$) subject to both clamped and supported boundary conditions are well behaved; that is,  the errors  are smooth and their  magnitudes are small  throughout the domain.   The errors of all the other cases, which are  not plotted here,  behave  similarly; their  magnitudes  can be read off the corresponding log-log plots collected in Figure~\ref{fig:dpnRates}.

\subsection{Plate with  numerous holes}
Now we solve  a plate with very complicated geometrical setting to showcase the capability and efficiency of our approach for solving the \KL plate equation \eqref{eq:KLShell} using composite overlapping grids. We consider  a  circular plate with  radius $4$ and two layers of holes sitting  on  two  rings inside of the circular domain.   On the  outer ring whose radius is  $3.5$,  there are twenty-four small holes of radius $0.3$  located on equally spaced angles; the angle of the $k^\text{th}$ hole is  $\theta_k=(15k)^\circ$.  On the inner ring  whose radius is  $2.25$, there sit another twelve larger holes
of radius $0.4$ on equally spaced angles;  the angle of the $k^\text{th}$ hole is  given by $\theta_k=(30k+15)^\circ$  so that it sits in between the holes on the outer ring.

{
\newcommand{\figWidth}{6cm}
\def\xa{13.}
\def\ya{6.5}
\newcommand{\trimfig}[2]{\trimw{#1}{#2}{0.12}{0.12}{0.12}{0.12}}
\begin{figure}[h]
\begin{center}
\begin{tikzpicture}[scale=1]
  \useasboundingbox (0.0,0.0) rectangle (\xa,\ya);  

\draw(-0.5,0.0) node[anchor=south west,xshift=0pt,yshift=0pt] {\trimfig{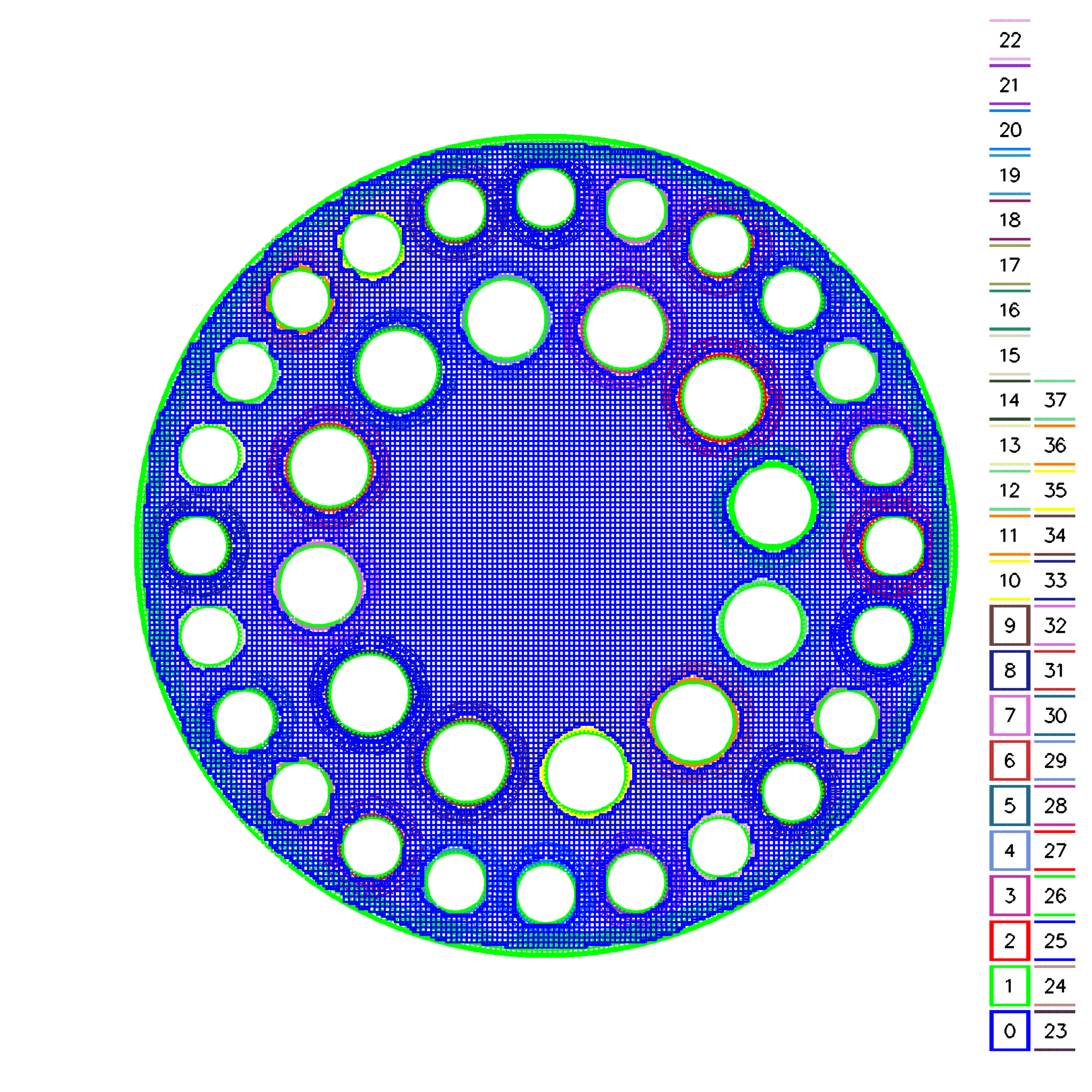}{\figWidth}};
\draw(6.5,0.0) node[anchor=south west,xshift=0pt,yshift=0pt] {\trimfig{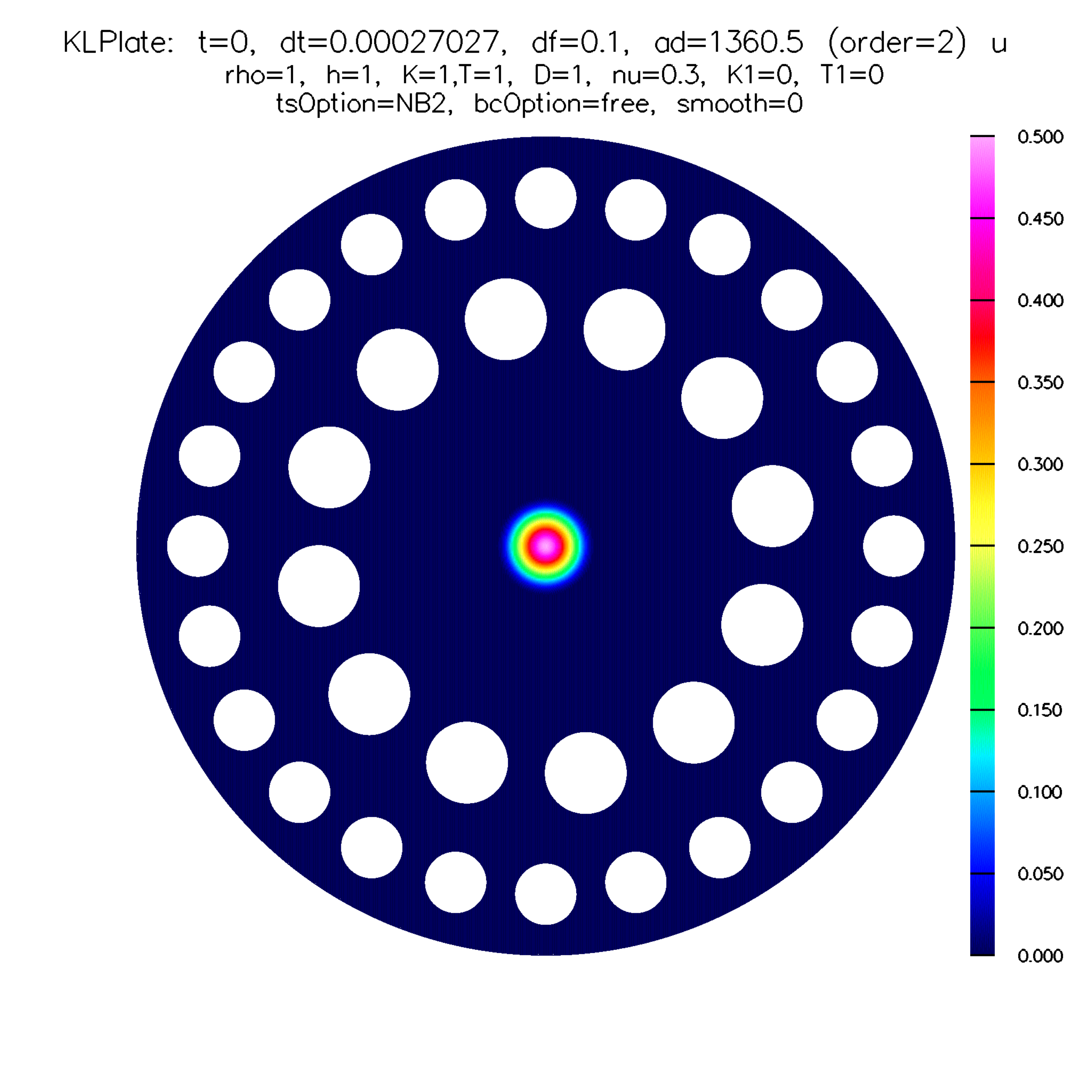}{\figWidth}};

\newcommand{\cbWidth}{.2}
\newcommand{\cbHeight}{3}
\drawColorBarH{fig/colourBarLines}{xshift=8.5cm,yshift=-.2cm}{\cbWidth}{\cbHeight}{$0$}{$0.5$}{w(x,y,0)}


%
\end{tikzpicture}

\end{center}
\caption{Computational setup for the plate with 24 outer and 12 inner holes. Left: coarsened version of the computational grid. Right: initial displacement.} \label{fig:plateWithHolesGridsAndIC}
\end{figure}
}

We discretize this complicated domain using a composite overlapping grid with 38 component grids, where a coarsened version  is  presented in Figure~\ref{fig:plateWithHolesGridsAndIC}.  The bulk of the domain is covered by a Cartesian grid, and the  outer boundary of domain is descritized by an  annulus boundary-fitted grid. The  holes in the plate are represented by  annulus boundary-fitted grids of smaller sizes that account for 
the remaining 36 component grids.  The actual computational grid is  $\Gc_{16}$ so that the target grid spacing is $1/160$. The computational  grid  has  1.8 million grid points, of which 76 thousand are  interpolation points.   We note that 
a plate with similar configuration but modeled with 3D linear elasticity  was studied in \cite{smog2012}; the computational grid with the same target grid spacing ($h=1/160$) for the 3D plate possess  42 million grid points. In this regard, the 2D \KL plate model is much cheaper to solve numerically   than a 3D plate model. Since our goal here is not to demonstrate  the  validity of the  \KL plate model for simplifying  the 3D linear elasticity model, we do not compare our results with \cite{smog2012}, which would require   tuning the parameters, initial and boundary conditions of the \KL model  to match that with   the 3D linear elasticity model.

{
\newcommand{\figWidth}{4.5cm}
\def\xa{16}
\def\ya{16.}
\newcommand{\trimfig}[2]{\trimw{#1}{#2}{0.12}{0.12}{0.12}{0.12}}
\begin{figure}[h]
\begin{center}
\begin{tikzpicture}[scale=1]
  \useasboundingbox (0.0,0.0) rectangle (\xa,\ya);  

\def\xs{0.3}
\draw(1,11.0) node[anchor=south west,xshift=\xs cm, yshift=0pt] {\trimfig{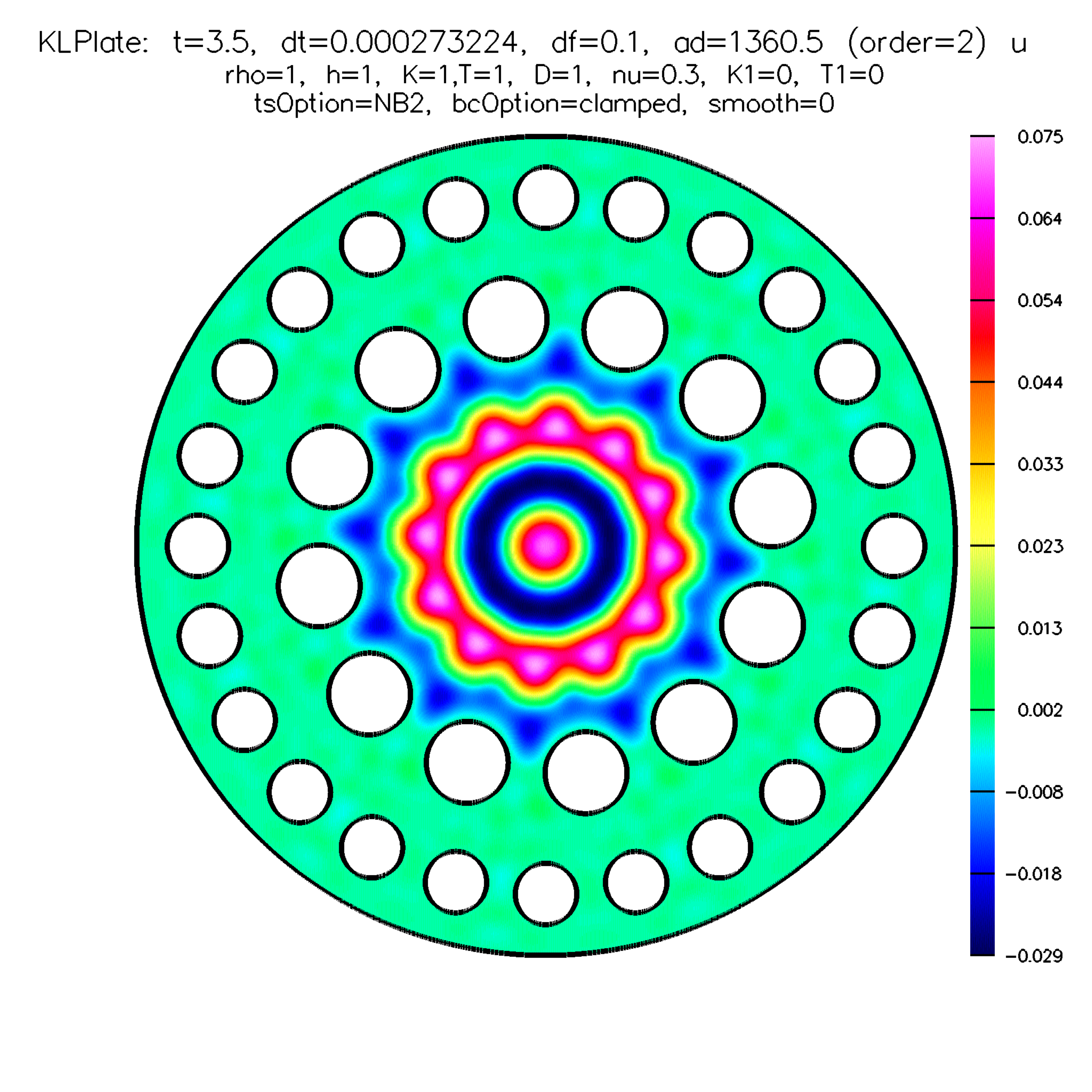}{\figWidth}};
\draw(5.75,11.0) node[anchor=south west,xshift=\xs cm,yshift=0pt] {\trimfig{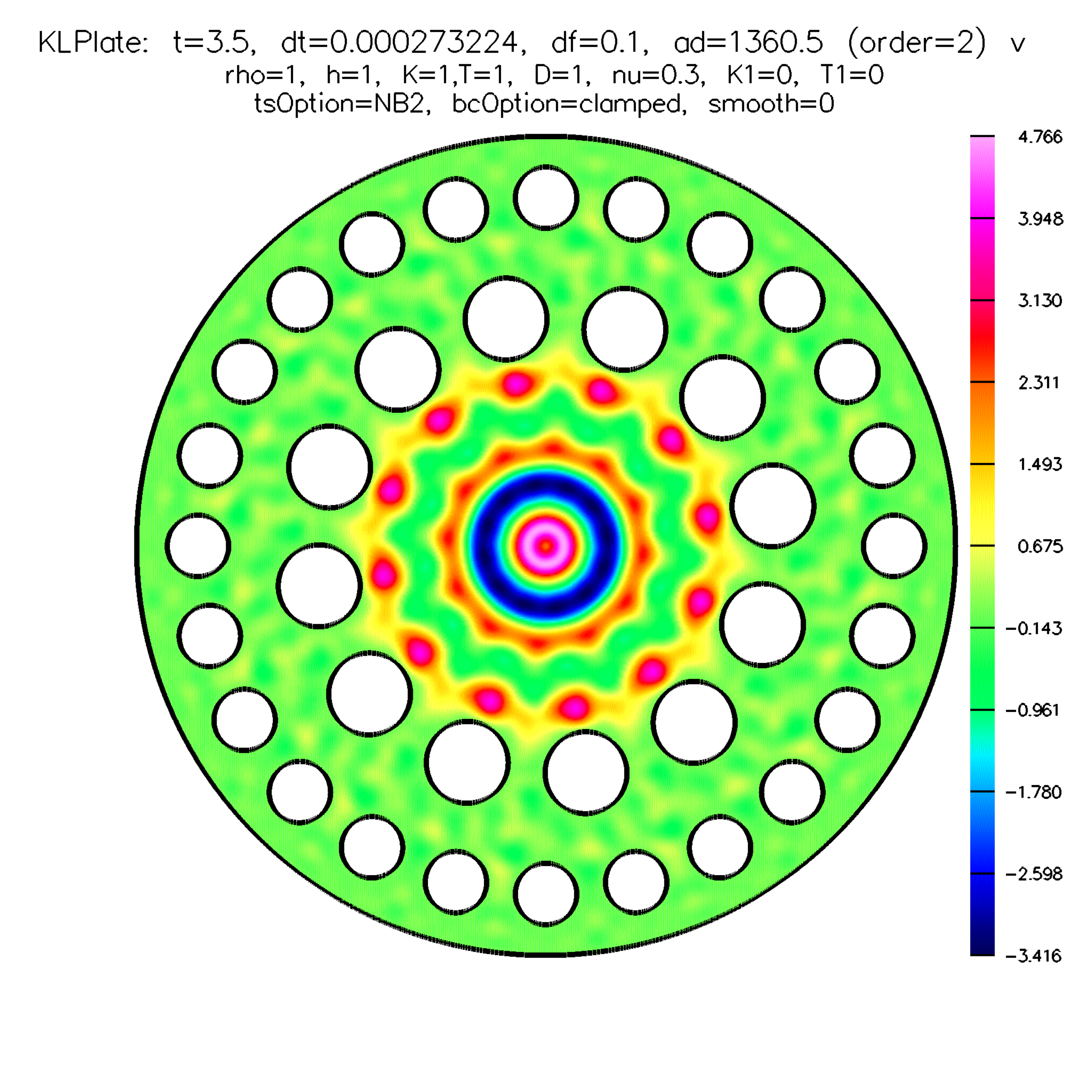}{\figWidth}};
\draw(10.5,11.0) node[anchor=south west,xshift=\xs cm,yshift=0pt] {\trimfig{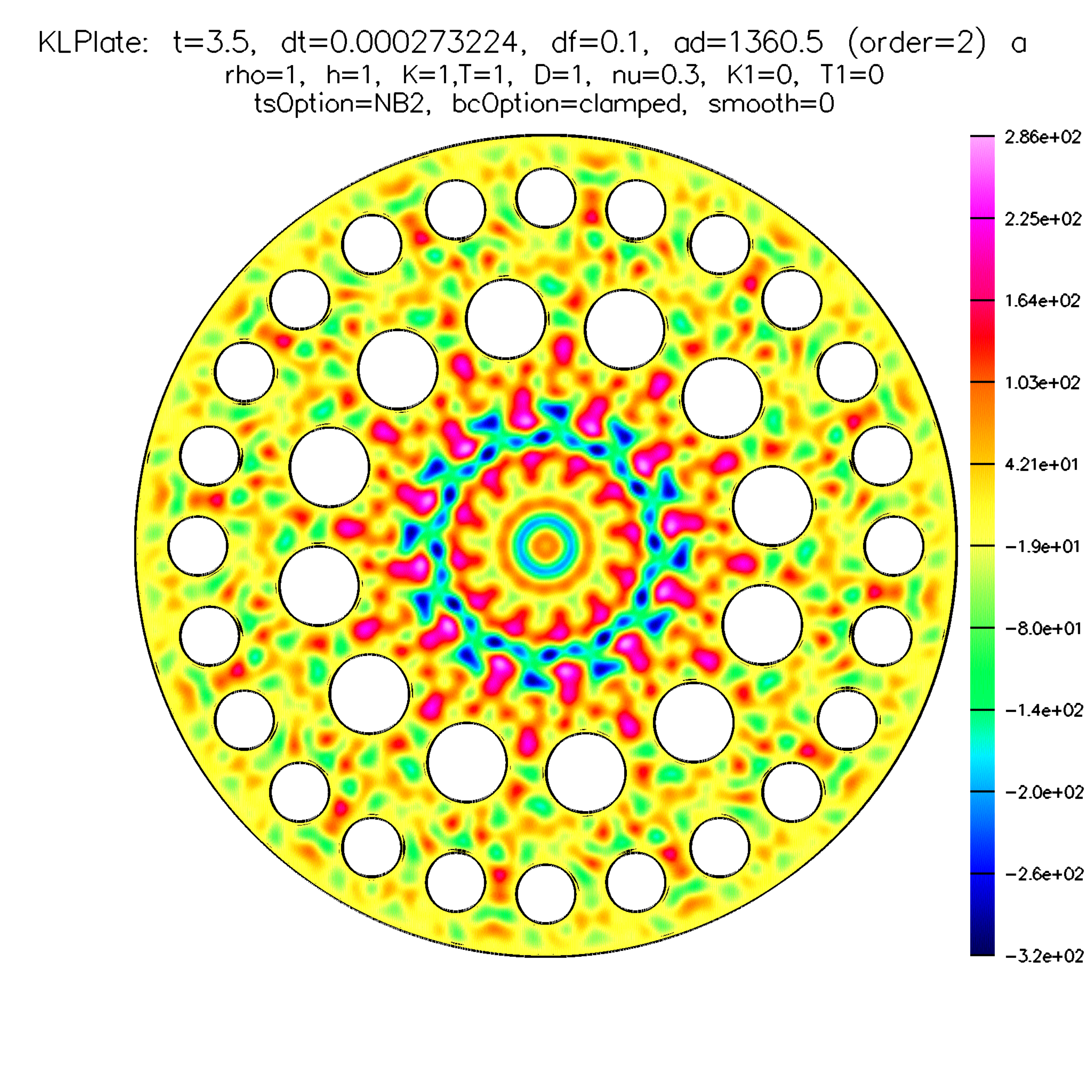}{\figWidth}};

\draw(1,5.5) node[anchor=south west,xshift=\xs cm,yshift=0pt] {\trimfig{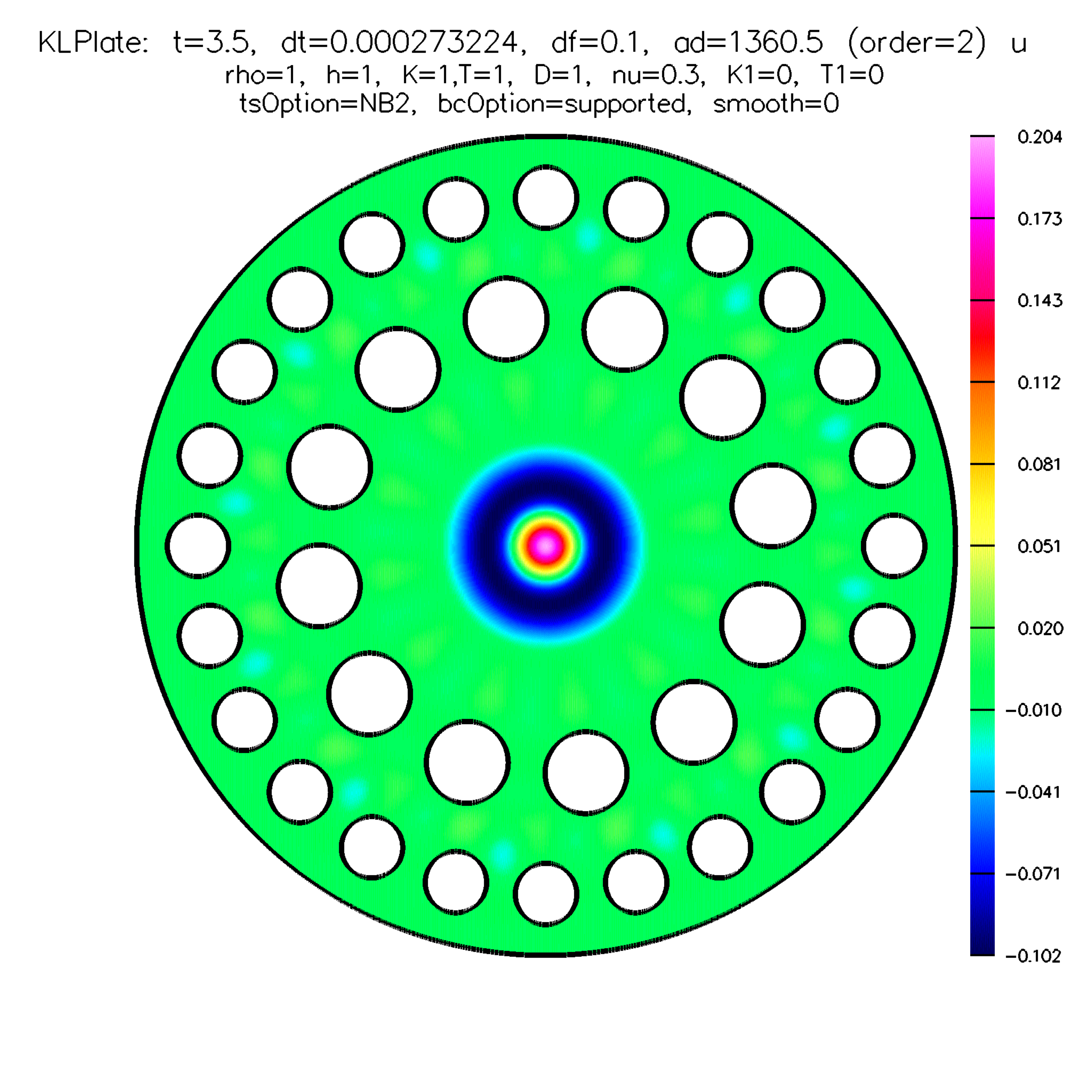}{\figWidth}};
\draw(5.75,5.5) node[anchor=south west,xshift=\xs cm,yshift=0pt] {\trimfig{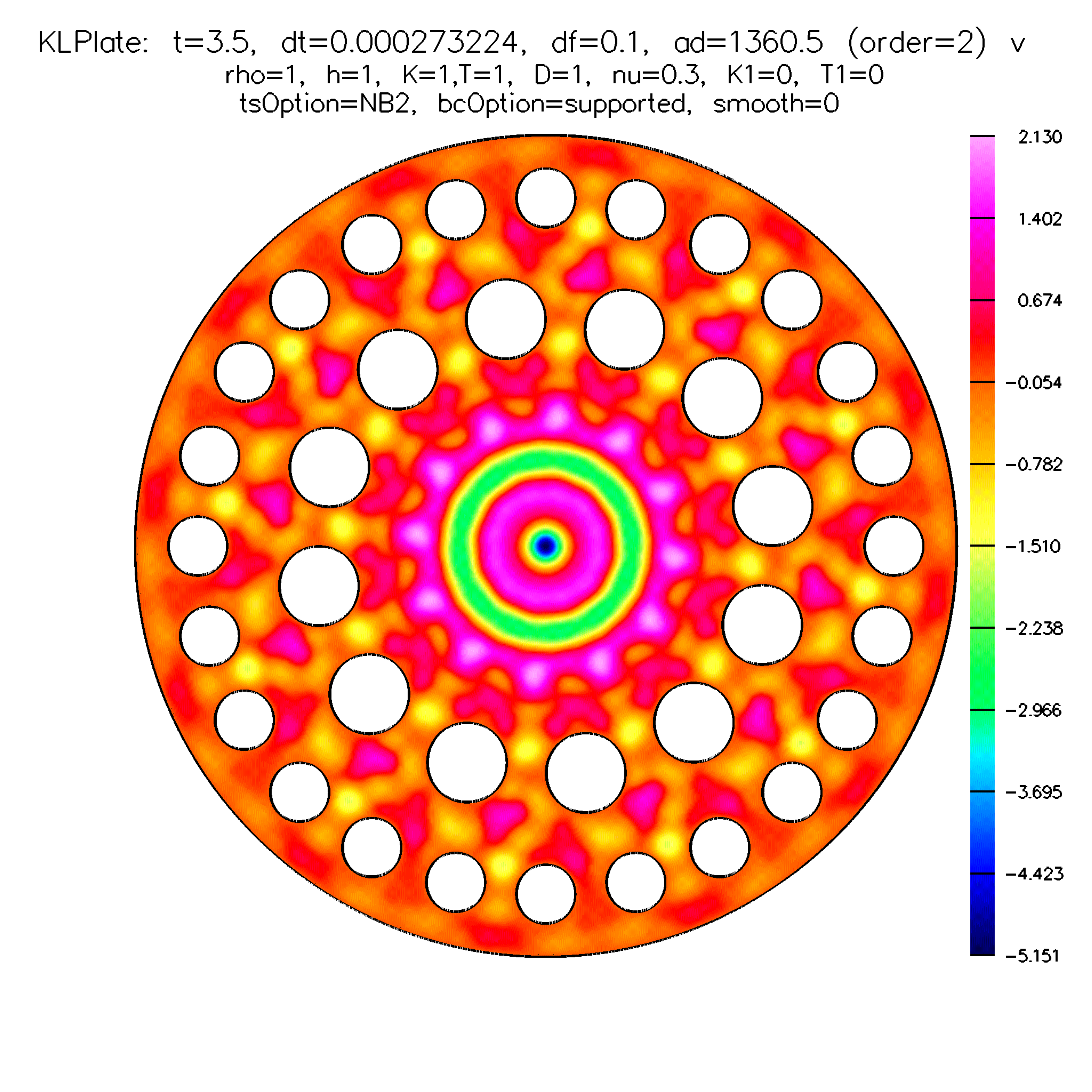}{\figWidth}};
\draw(10.5,5.5) node[anchor=south west,xshift=\xs cm,yshift=0pt] {\trimfig{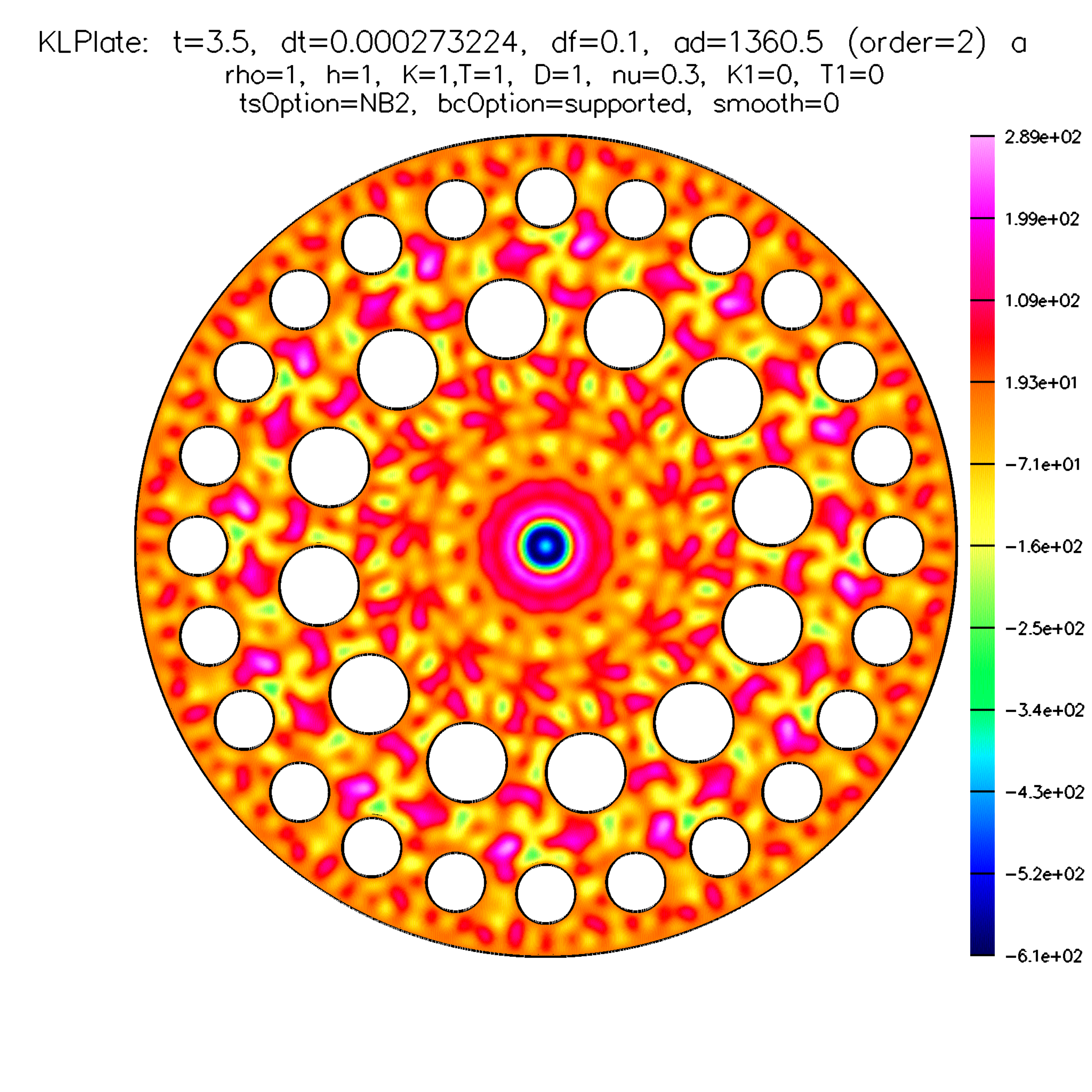}{\figWidth}};

\draw(1,0.0) node[anchor=south west,xshift=\xs cm,yshift=0pt] {\trimfig{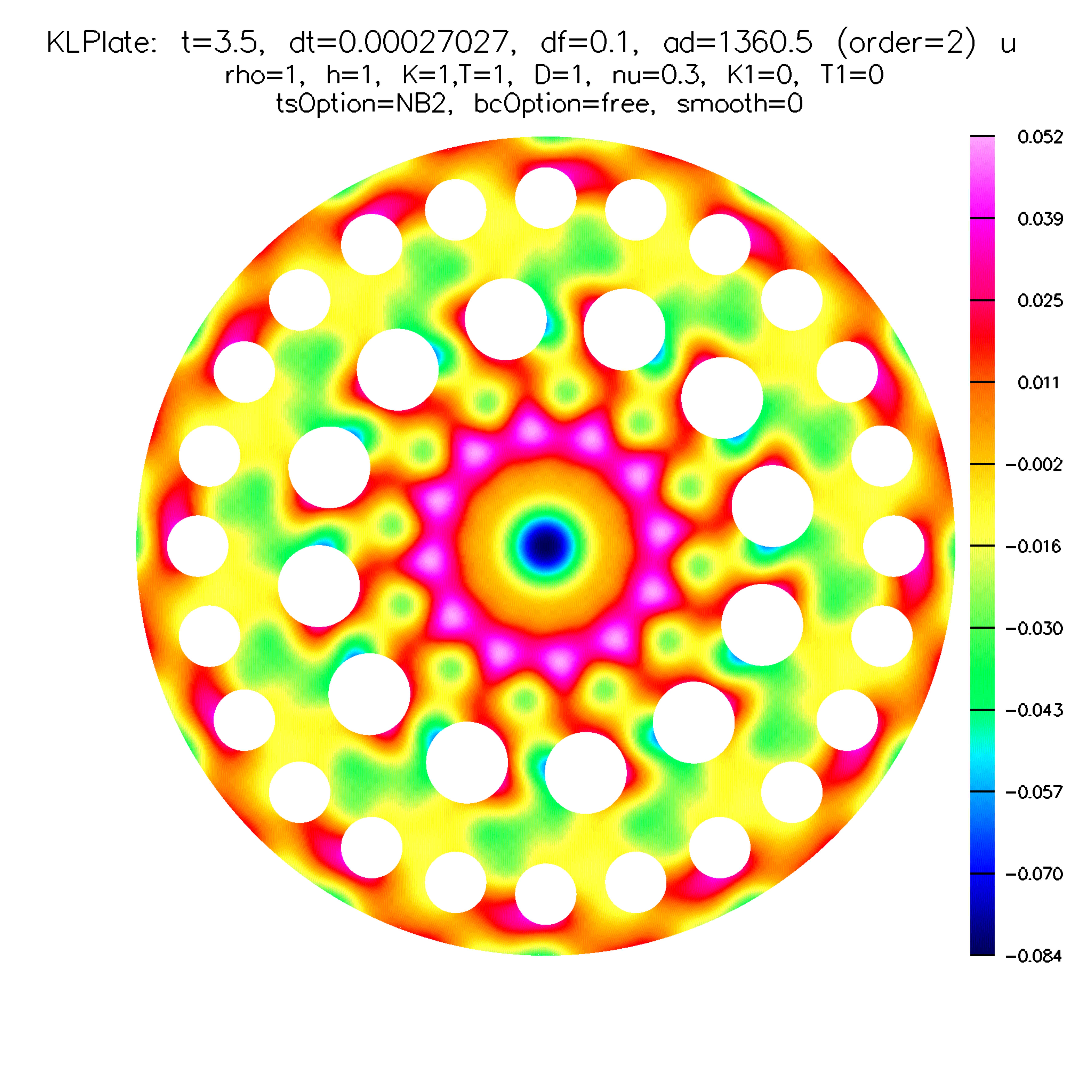}{\figWidth}};
\draw(5.75,0.0) node[anchor=south west,xshift=\xs cm,yshift=0pt] {\trimfig{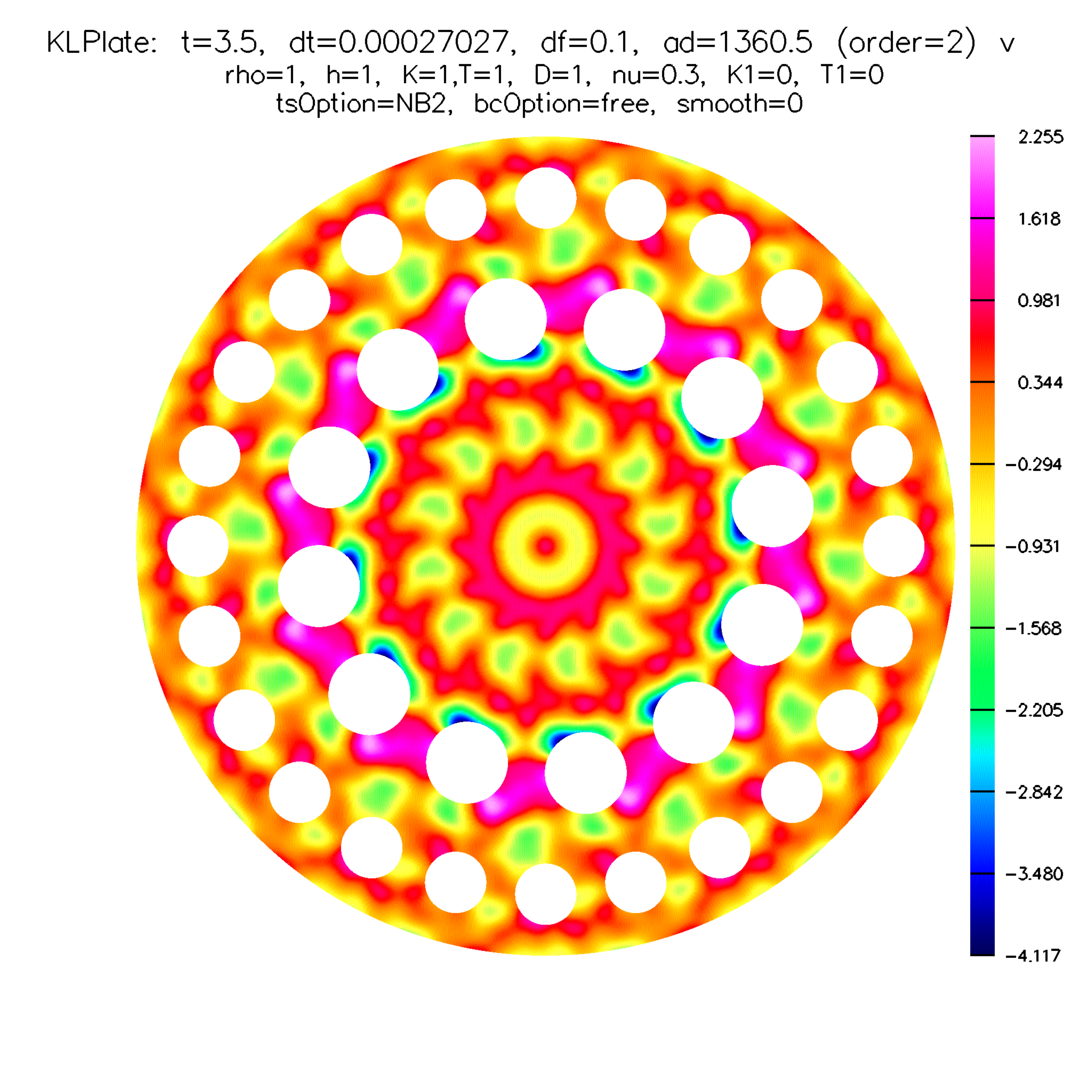}{\figWidth}};
\draw(10.5,0.0) node[anchor=south west,xshift=\xs cm,yshift=0pt] {\trimfig{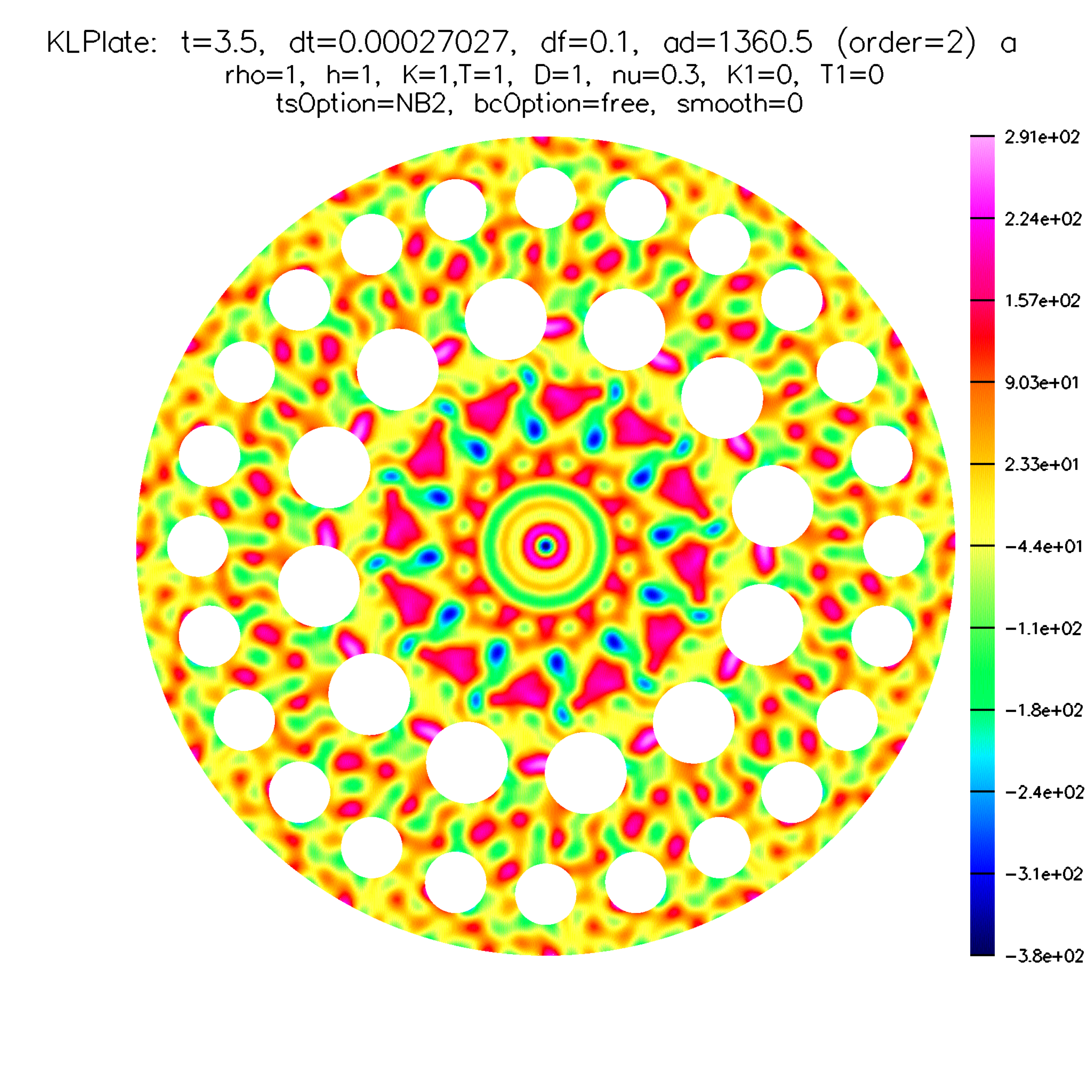}{\figWidth}};

\newcommand{\cbWidth}{.2}
\newcommand{\cbHeight}{3}

\drawColorBarH{fig/colourBarLines}{xshift=2.45cm,yshift=10.8cm}{\cbWidth}{\cbHeight}{$-0.029$}{$0.075$}{}
\drawColorBarH{fig/colourBarLines}{xshift=7.2cm,yshift=10.8cm}{\cbWidth}{\cbHeight}{$-3.416$}{$4.766$}{}
\drawColorBarH{fig/colourBarLines}{xshift=11.95cm,yshift=10.8cm}{\cbWidth}{\cbHeight}{$-3.2e2$}{$2.86e2$}{}

\drawColorBarH{fig/colourBarLines}{xshift=2.45cm,yshift=5.3cm}{\cbWidth}{\cbHeight}{$-0.102$}{$0.204$}{}
\drawColorBarH{fig/colourBarLines}{xshift=7.2cm,yshift=5.3cm}{\cbWidth}{\cbHeight}{$-5.151$}{$2.130$}{}
\drawColorBarH{fig/colourBarLines}{xshift=11.95cm,yshift=5.3cm}{\cbWidth}{\cbHeight}{$-6.1e2$}{$2.89e2$}{}

\drawColorBarH{fig/colourBarLines}{xshift=2.45cm,yshift=-.2cm}{\cbWidth}{\cbHeight}{$-0.084$}{$0.052$}{}
\drawColorBarH{fig/colourBarLines}{xshift=7.2cm,yshift=-.2cm}{\cbWidth}{\cbHeight}{$-4.117$}{$2.255$}{}
\drawColorBarH{fig/colourBarLines}{xshift=11.95cm,yshift=-.2cm}{\cbWidth}{\cbHeight}{$-3.8e2$}{$2.91e2$}{}

\draw(3.7,16.3)  node[anchor=north,xshift=\xs cm,yshift=0pt] {$w(x,y,3.5)$};
\draw(8.6,16.3)  node[anchor=north,xshift=\xs cm,yshift=0pt] {$v(x,y,3.5)$};
\draw(13.3,16.3)  node[anchor=north,xshift=\xs cm,yshift=0pt] {$a(x,y,3.5)$};

\draw(1.85,2.5)  node[anchor=east,xshift=0cm,yshift=0pt] {\footnotesize\bf Free};
\draw(1.85,8)  node[anchor=east,xshift=0cm,yshift=0pt] {\footnotesize\bf Supported};
\draw(1.85,13.5) node[anchor=east,xshift=0cm,yshift=0pt] {\footnotesize\bf Clamped};

%
\end{tikzpicture}
\end{center}
\caption{Plate with holes simulation at $t=3.5$ subject to: clamped (top row), supported (middle row), free (bottom row) boundary conditions. Results shown here are obtained using the NB2 scheme with $\csf=50$ and $\df=0.1$. } \label{fig:plateWithHolesG16NB2Feet3p5}
\end{figure}
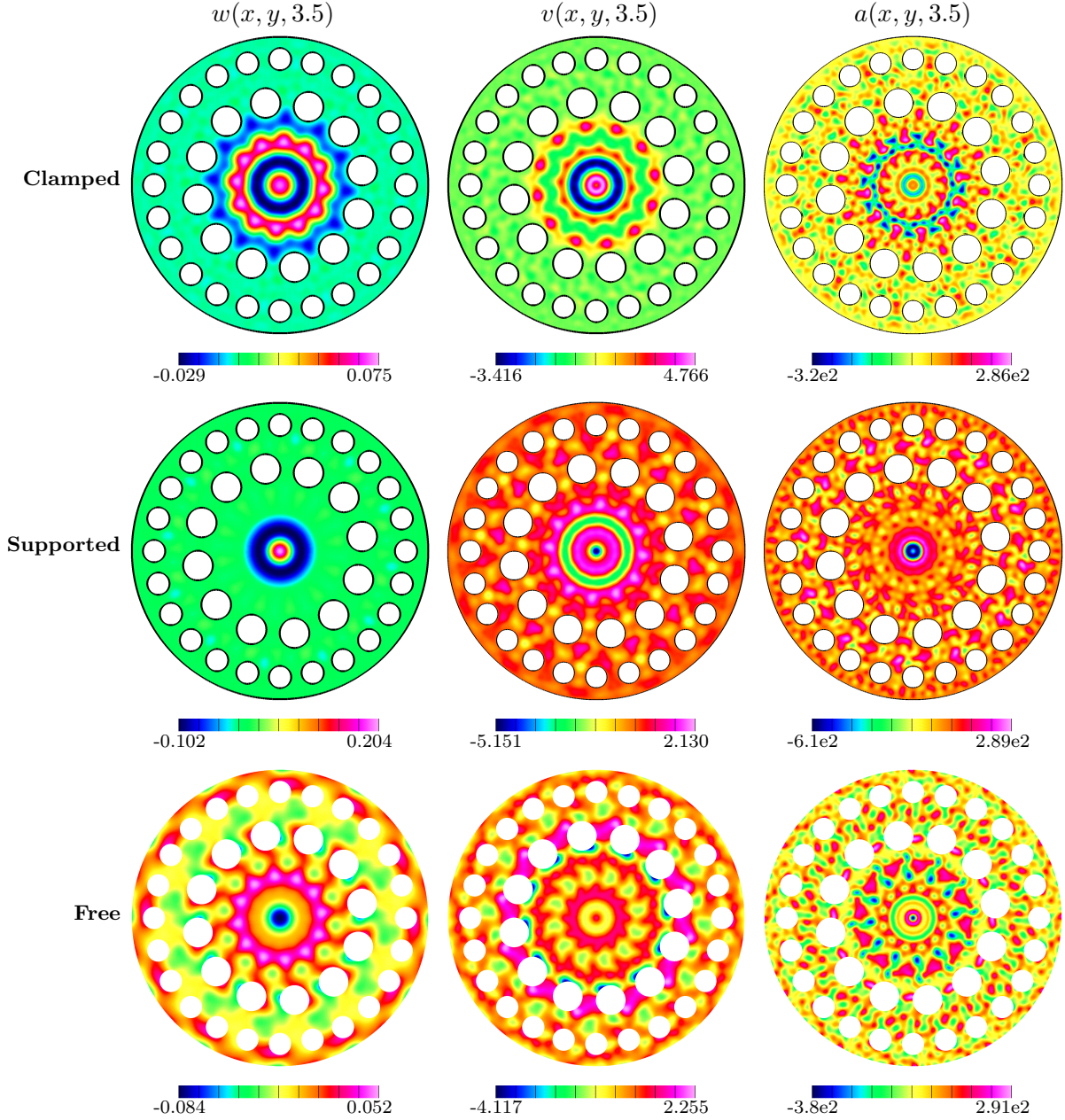
}

For this numerical test, we solve \eqref{eq:KLShell} with physical parameters being $\rho=1,H=1,K=1,T=1,D=1$ and $\nu=0.3$ to  investigate the vibrations of the plate in response to the following initial displacement  disturbance,  
\begin{align*}
w(x,y,0)&=\begin{cases}
0.25\left(\cos({\sqrt{x^2+y^2}}/{0.5})+1\right), & \sqrt{x^2+y^2} \leq 0.5,\\
0, & \text{otherwise},
\end{cases}\\
v(x,y,0)&=0.
\end{align*}
The contour of the initial displacement is plotted in Figure~\ref{fig:plateWithHolesGridsAndIC}. The governing equation is solved using all the proposed numerical methods subject to all the boundary conditions \eqref{eq:clampedBC}--\eqref{eq:freeBC}. The time step  and artificial   dissipation for each scheme is determined according to formulas summarized in Table~\ref{tab:timeStep}. The stability factor is $\csf=0.9$ for the explicit schemes, and $\csf=50$ for the NB2 scheme. And it suffices to stabilize all the schemes for this problem by setting the dissipation factor as  $\df=0.1$.

Results obtained using  all the numerical methods are similar, so only those of the NB2 scheme are shown here in Figure~\ref{fig:plateWithHolesG16NB2Feet3p5}. The figure shows the evolution of the  displacement, velocity and acceleration at time $t=3.5$ for the plate with clamped (top row), supported  (middle row), and free (bottom row) boundary conditions. We observe that the vibrations caused  by the initial disturbance      propagate  outward from the center of the plate towards  its perimeter. As the disturbance   reaches  the  holes, the behavior of   vibrations begin to differ for different boundary conditions. For the clamped and supported boundary conditions, we see that the displacement disturbance generally stops  at the   inner ring of  holes because both  boundary conditions prescribe zero displacement at the edges of the holes that serves to suppress the vibrations.  For the solutions of the  free boundary conditions, vibrations are able to propagate through both the inner and outer rings of holes since the edges of the holes are allowed to vibrate freely. All the  plots shown here  reflect the   rich and complex dynamics of the \KL plate model.

\subsection{Traveling  pulse on a  guitar soundboard }
As an illustration of the ability of our approach for solving   realistic applications, we consider the problem of smoothed  pulse  traveling through  a  guitar soundboard. See Figure~\ref{fig:guitarGrid} for the geometric configuration of the  soundboard and a coarsened version of the computational grid.
 Note that we do not intend to make a comparison for a particular guitar design or material;  the focus of this paper is to demonstrate the numerical properties of the aforementioned algorithms. Therefore, in this test, 
we  explore how an initial pulse   travel through the   guitar soundboard with hypothetical physical  parameters.

{
\newcommand{\figWidth}{6cm}
\def\xa{13.}
\def\ya{5.5}
\newcommand{\trimfig}[2]{\trimw{#1}{#2}{0.12}{0.12}{0.18}{0.18}}
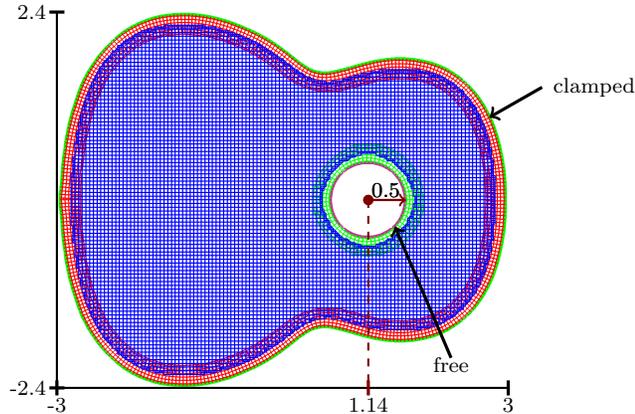
\begin{figure}[h]
\begin{center}
\begin{tikzpicture}[scale=1]
  \useasboundingbox (0.0,0.0) rectangle (\xa,\ya);  

 \begin{scope}[,xshift=3.25cm,yshift=0cm]
   \draw(0,0) node[anchor=south west,xshift=0pt,yshift=0pt] {\trimfig{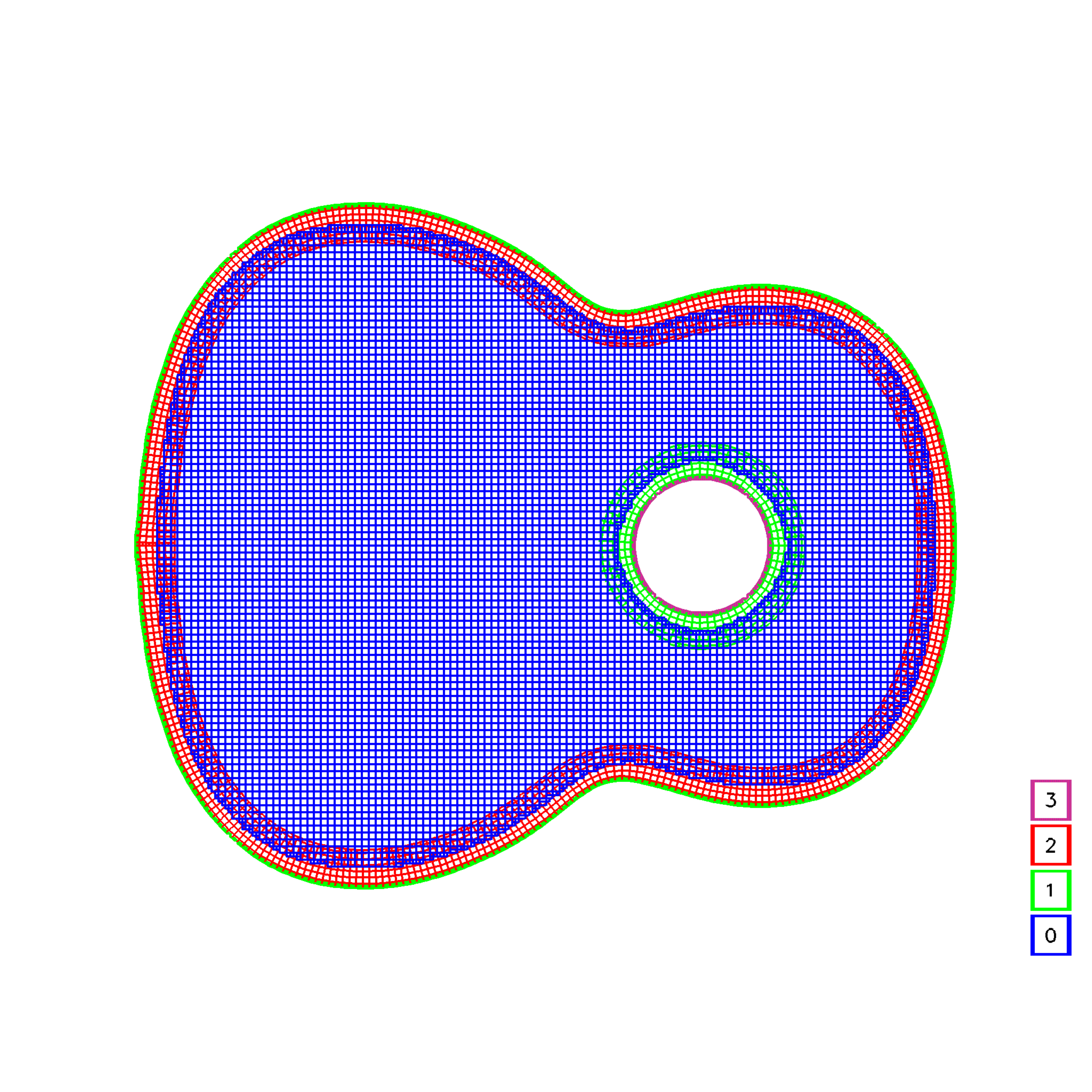}{\figWidth}};
  \def\bb{0.2}
  \def\tt{5.2}
  \def\ll{0.55}
  \def\rr{6.55}
  \def\dd{0.1}
  \def\cx{4.69} 
   \def\cy{2.7}
   \def\rad{0.5} 

  \draw[-,black,thick,xshift=0cm] (\ll,\bb) node[anchor=east,yshift=0pt] {\scriptsize $$-2.4$$} -- (\ll,\tt) node[anchor=east,yshift=0pt] {\scriptsize $$2.4$$};
  \draw[-,black,thick,xshift=0cm] (\ll,\bb)  node[anchor=north,yshift=0pt] {\scriptsize $$-3$$}  -- (\rr,\bb) node[anchor=north,yshift=0pt] {\scriptsize $$3$$};
 
\foreach \y in {\bb,\tt}
 {
  \draw[-,black,very thick,xshift=0cm] ({\ll-\dd},\y) -- ({\ll+\dd},\y);
 }
 \foreach \x in {\ll,\rr}
 {
  \draw[-,black,very thick,xshift=0cm] (\x,{\bb-\dd}) -- (\x,{\bb+\dd});
}

 \fill[red!50!black] (\cx,\cy) circle (2pt);
 \draw[dashed,red!50!black,thick,xshift=0cm] (\cx,\cy) -- (\cx,\bb) node[black,anchor=north,yshift=0pt] {\scriptsize $$1.14$$};
 \draw[-,red!50!black,very thick,xshift=0cm] (\cx,{\bb-\dd}) -- (\cx,{\bb+\dd});

 \draw[->,red!50!black,thick,xshift=0cm] (\cx,\cy) node[black,anchor=south west,xshift=-0.1cm,yshift=-0.1cm]{\scriptsize $$0.5$$} -- (\cx+\rad,\cy); 
 \draw[<-,black, very thick,xshift=0cm] ({\cx+0.705*\rad},{\cy-0.705*\rad})  -- ({\cx+1.1},{\cy-2.1})  node[black,anchor=north,yshift=0.15cm]{\scriptsize free};
 \draw[<-,black, very thick,xshift=0cm] (6.3,3.8)  -- (7,4.2)  node[black,anchor=west,yshift=0cm]{\scriptsize clamped};

 \draw[->,red!50!black,thick,xshift=0cm] (\cx,\cy) node[black,anchor=south west,xshift=-0.1cm,yshift=-0.1cm]{\scriptsize $$0.5$$} -- (\cx+\rad,\cy);

   \end{scope}


%
\end{tikzpicture}

\end{center}
\caption{Geometric configuration of the guitar soundboard and a coarsened version of the computational grids.} \label{fig:guitarGrid}
\end{figure}
}

Let us define a smoothed pulse, 
$$
w_{p}(x,y,t)=
\begin{cases}\displaystyle
0.1\left[\cos\left({\pi\sqrt{\frac{(x-x_c-ct)^2}{r^2_a}+\frac{(y-y_c)^2}{r^2_b}}}\right)+1\right], &~\text{if} ~\displaystyle{\frac{(x-x_c-ct)^2}{r^2_a}+\frac{(y-y_c)^2}{r^2_b}} \leq 1,\\
0, & \text{otherwise}.
\end{cases}\\
$$
where $x_c=-1, y_c=0, r_a=0.1, r_b=2, c=1 $. Then we  specify the initial conditions accordingly,
\begin{equation}
  w(x,y,0)=w_p(x,y,0),\quad
  v(x,y,0)=\pd{w_p}{t}(x,y,0) \approx \frac{w_p(x,y,0)-w_p(x,y,-\dt)}{\dt}.
\end{equation}
The initial conditions are plotted in Figure~\ref{fig:guitarIC}. Clamped boundary conditions are enforced for  the outer edge of the guitar domain, while the inner circular edge is allowed to move freely.

{
\newcommand{\figWidth}{6cm}
\def\xa{13.}
\def\ya{5.5}
\newcommand{\trimfig}[2]{\trimw{#1}{#2}{0.12}{0.12}{0.18}{0.18}}
\begin{figure}[h]
\begin{center}
\begin{tikzpicture}[scale=1]
  \useasboundingbox (0.0,0.0) rectangle (\xa,\ya);  

\draw(-0.5,0.0) node[anchor=south west,xshift=0pt,yshift=0pt] {\trimfig{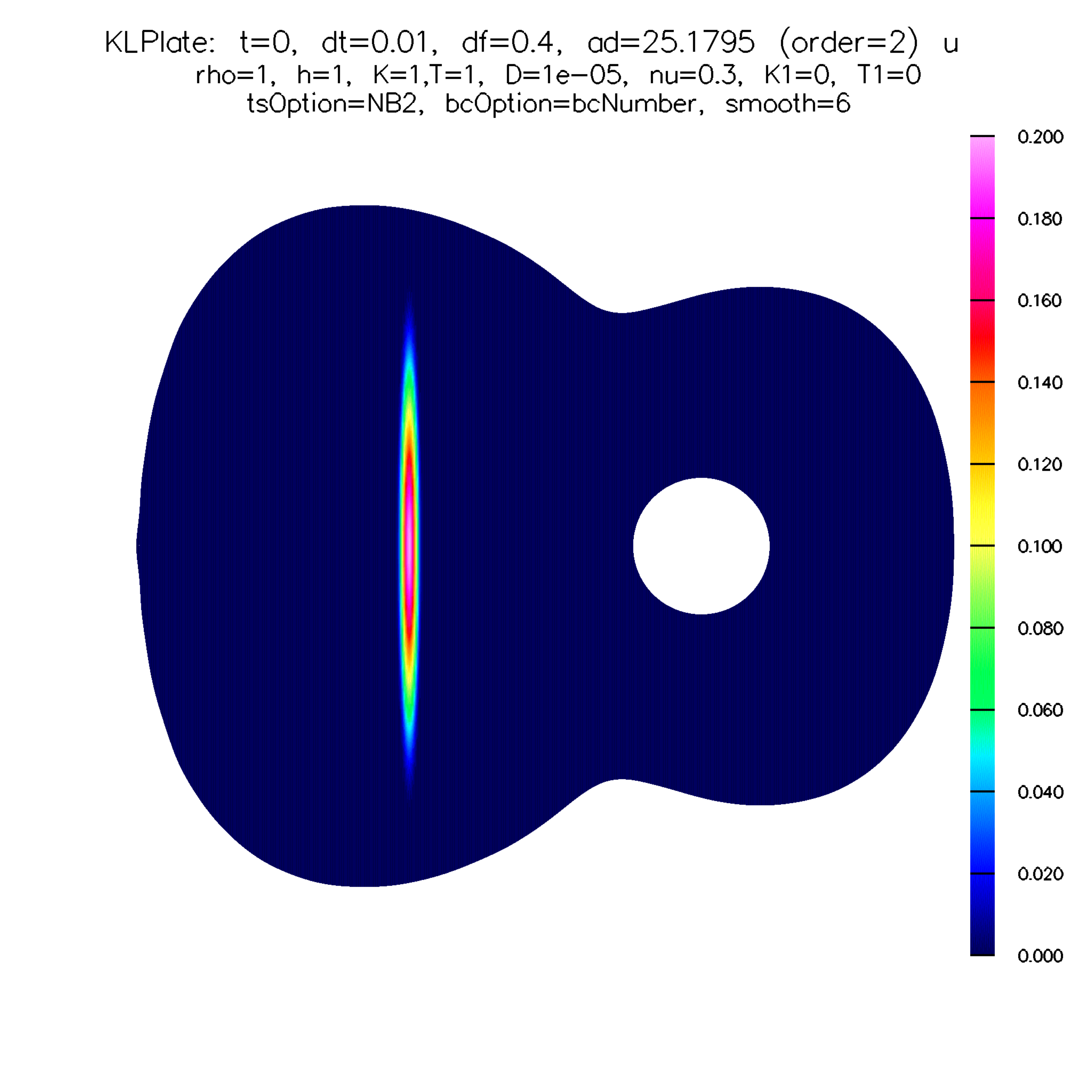}{\figWidth}};
\draw(6.5,0.0) node[anchor=south west,xshift=0pt,yshift=0pt] {\trimfig{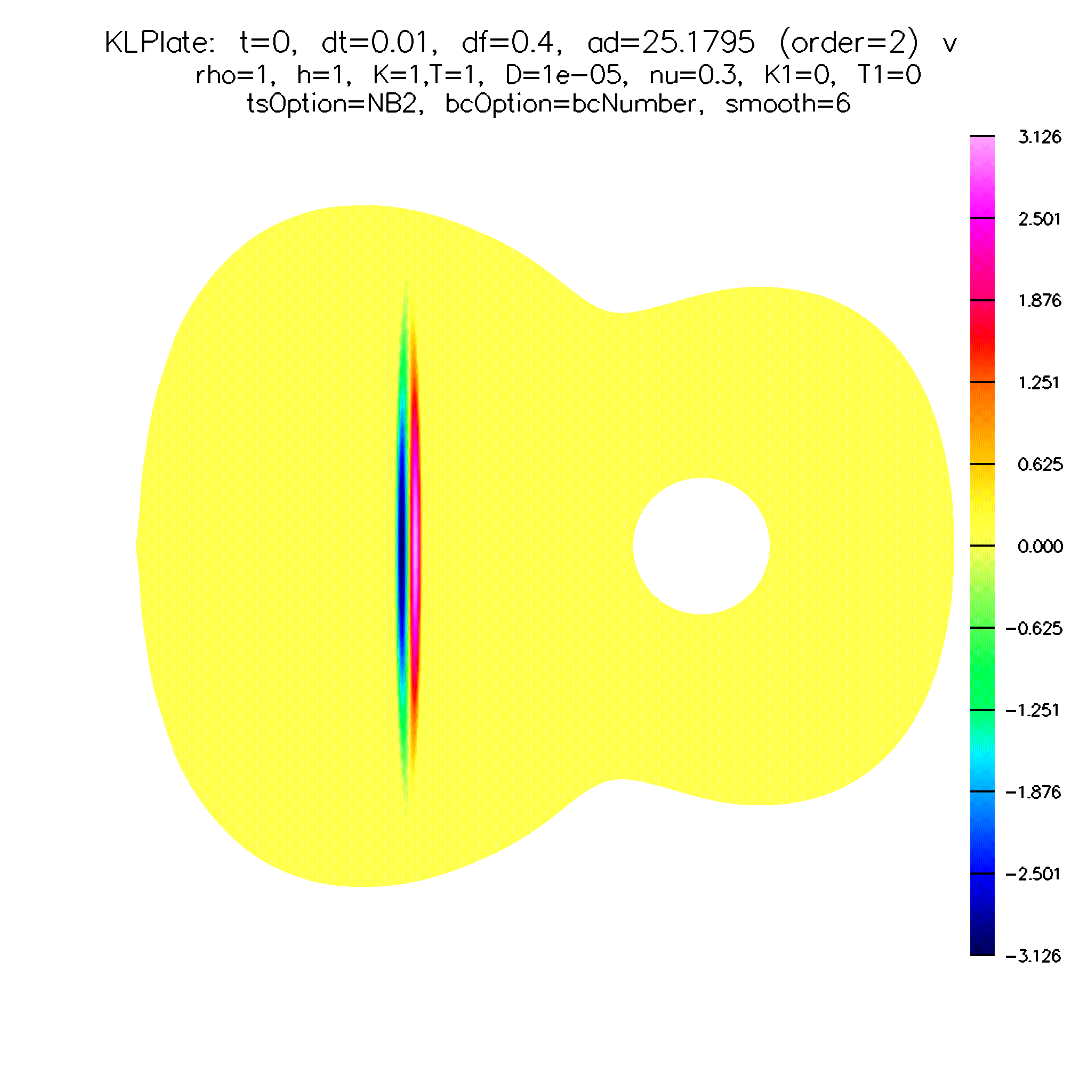}{\figWidth}};

\newcommand{\cbWidth}{.2}
\newcommand{\cbHeight}{3}

\drawColorBarH{fig/colourBarLines}{xshift=1.5cm,yshift=-0.2cm}{\cbWidth}{\cbHeight}{$0.0$}{$0.20$}{w(x,y,0)}

\drawColorBarH{fig/colourBarLines}{xshift=8.5cm,yshift=-0.2cm}{\cbWidth}{\cbHeight}{$-3.13$}{$3.13$}{v(x,y,0)}


%
\end{tikzpicture}

\end{center}
\caption{Initial displacement (left) and  velocity (right).  } \label{fig:guitarIC}
\end{figure}
}

To investigate how the initial pulse travels through the domain, we  consider two sets of physical  parameters :  (I) wave dominant case ($\rho=1, H=1,K=1,T=1,D=1\times10^{-5},\nu=0.3$) and (II) bend dominant case ($\rho=1, H=1,K=1,T=1,D=0.1,\nu=0.3$). For the wave dominant case, the bending stiffness  is set to be very small so that the initial pulse  is expected to travel like a wave, while for the bend dominant case, the bending effect dominates the vibrations.  Computations are carried out  for both cases  using  grid $\Gc_{16}$ with  UPC2 and NB2 schemes so that we can compare their respective performances. 

The simulation results of the wave dominant case are shown in Figure~\ref{fig:guitarWaveDominant}.  In the figure,  the evolution of all the solution components (i.e., $w$,$v$ and $a$)   are presented at two selected times ($t=3.5$ and $t=5$).  For this case, we observe that  the initial pulse travels  towards the hole. When the leading edge of the pulse hits the hole, more waves are generated by diffractions and reflections. By the end of the simulation ($t=5$), a complex distribution of elastic waves are observed over the entire soundboard. The velocity and acceleration evolve  similarly as the displacement.

{
\newcommand{\figWidth}{5cm}
\def\xa{16.}
\def\ya{10.5}
\newcommand{\trimfig}[2]{\trimw{#1}{#2}{0.12}{0.12}{0.18}{0.18}}
\begin{figure}[h]
\begin{center}
\begin{tikzpicture}[scale=1]
  \useasboundingbox (0.0,0.0) rectangle (\xa,\ya);  

\draw(-0.5,5.5) node[anchor=south west,xshift=0pt,yshift=0pt] {\trimfig{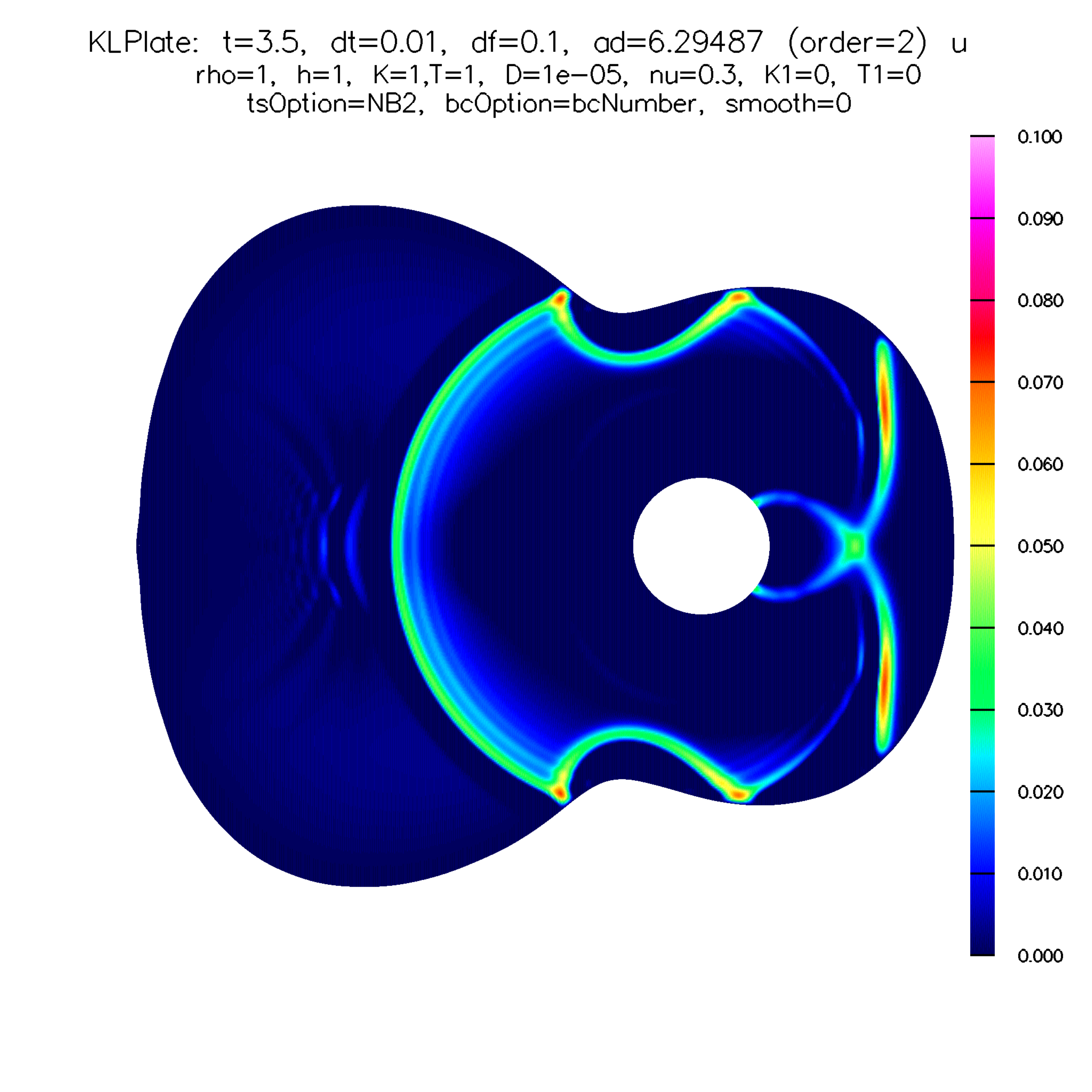}{\figWidth}};
\draw(5.,5.5) node[anchor=south west,xshift=0pt,yshift=0pt] {\trimfig{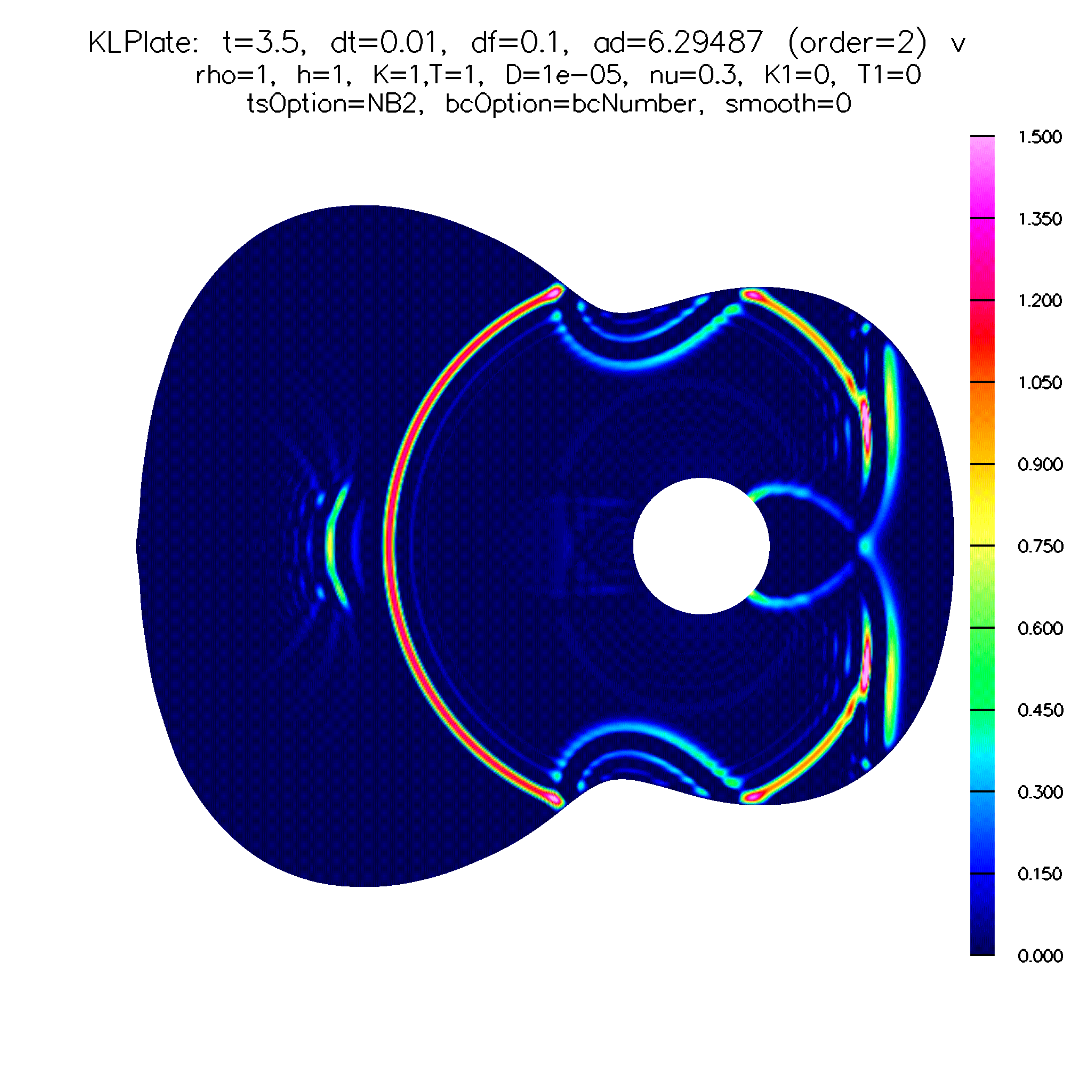}{\figWidth}};
\draw(10.5,5.5) node[anchor=south west,xshift=0pt,yshift=0pt] {\trimfig{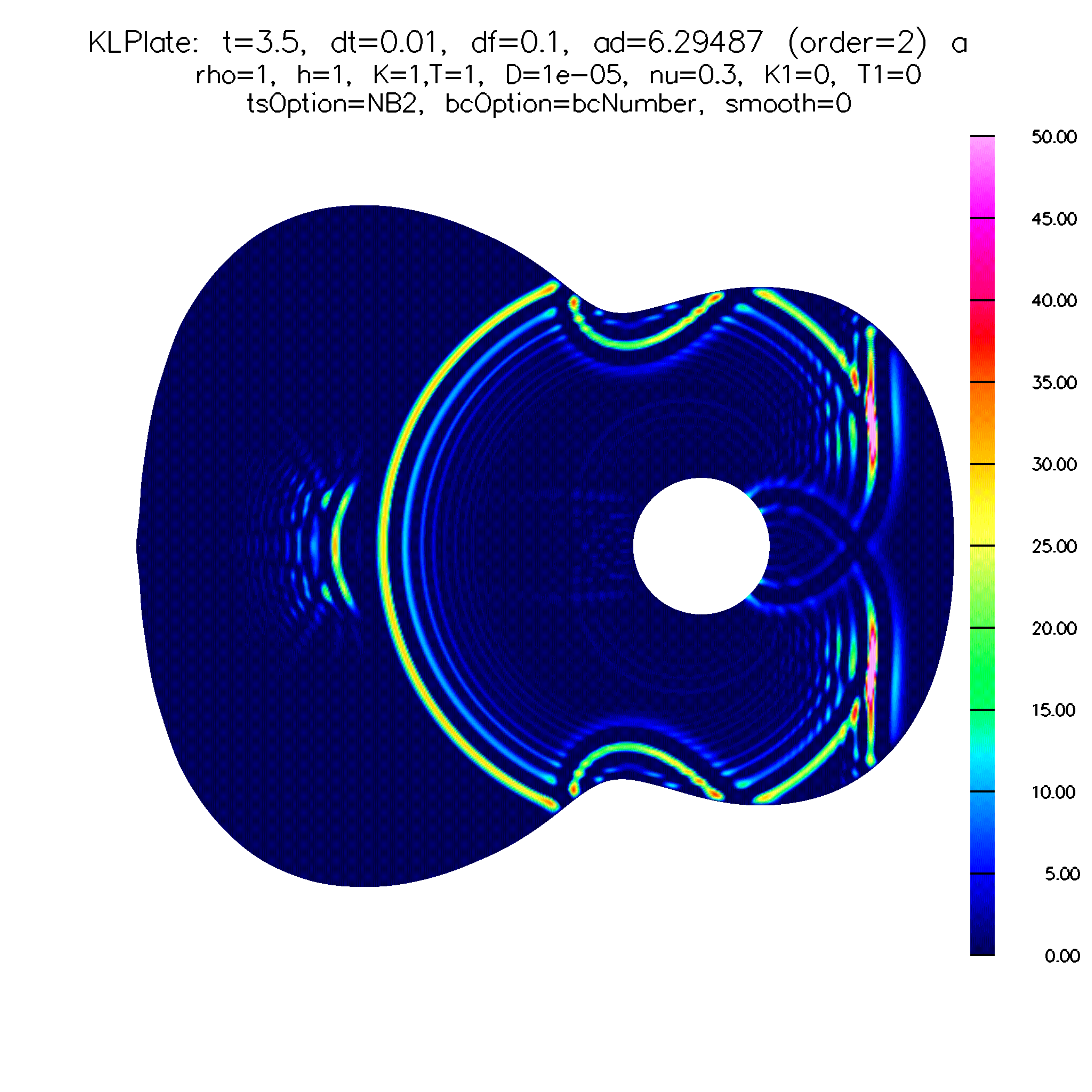}{\figWidth}};

\draw(-0.5,0.0) node[anchor=south west,xshift=0pt,yshift=0pt] {\trimfig{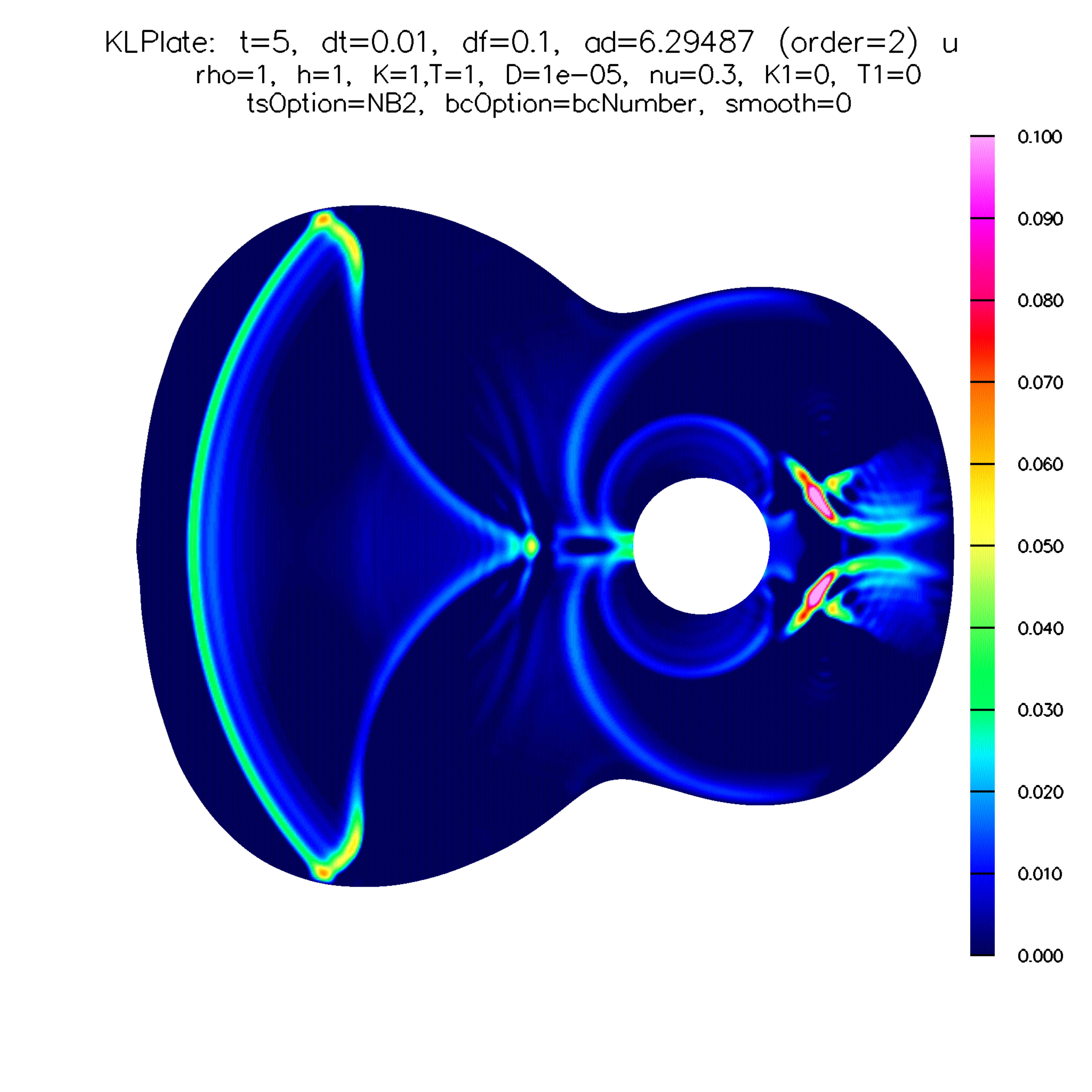}{\figWidth}};
\draw(5.,0.0) node[anchor=south west,xshift=0pt,yshift=0pt] {\trimfig{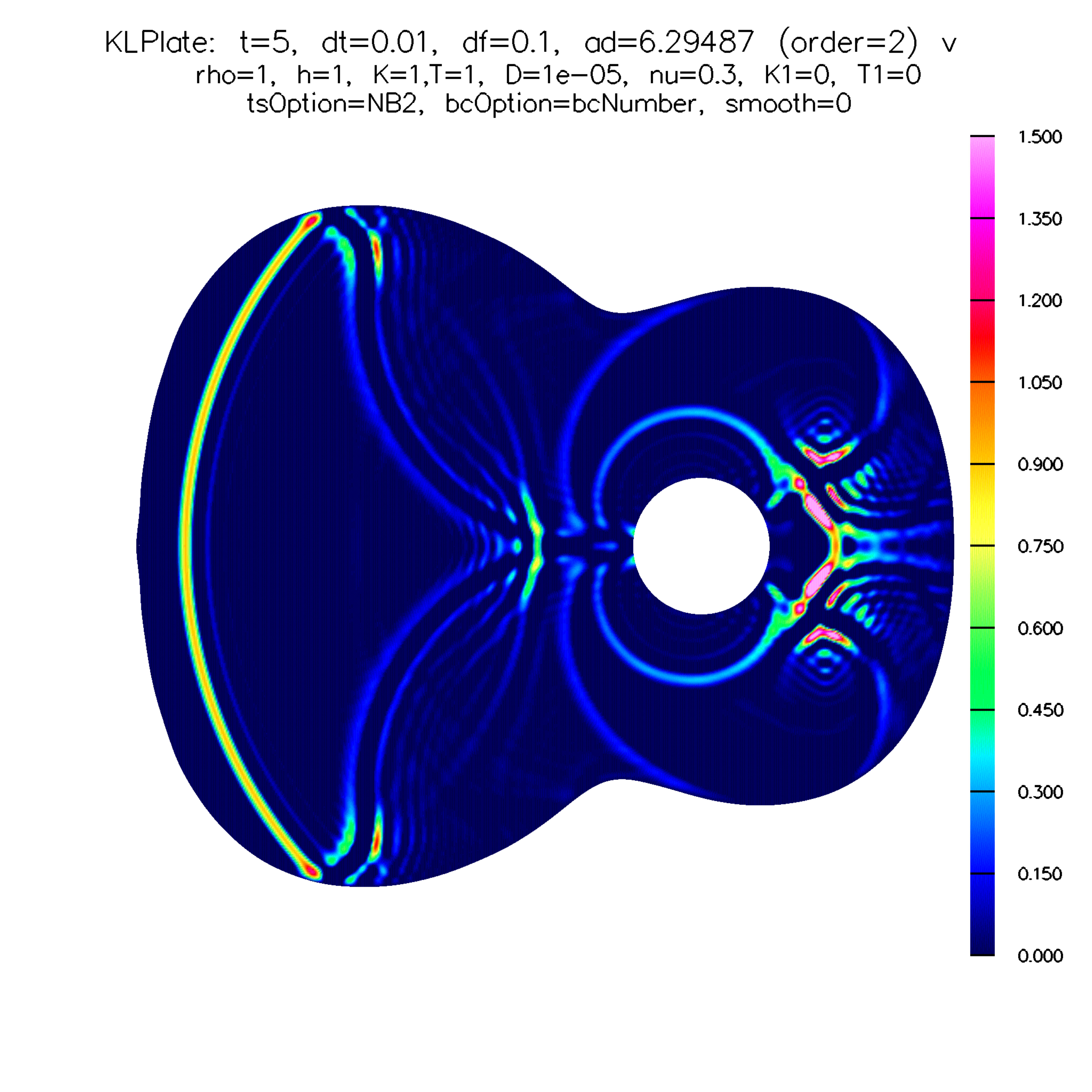}{\figWidth}};
\draw(10.5,0.0) node[anchor=south west,xshift=0pt,yshift=0pt] {\trimfig{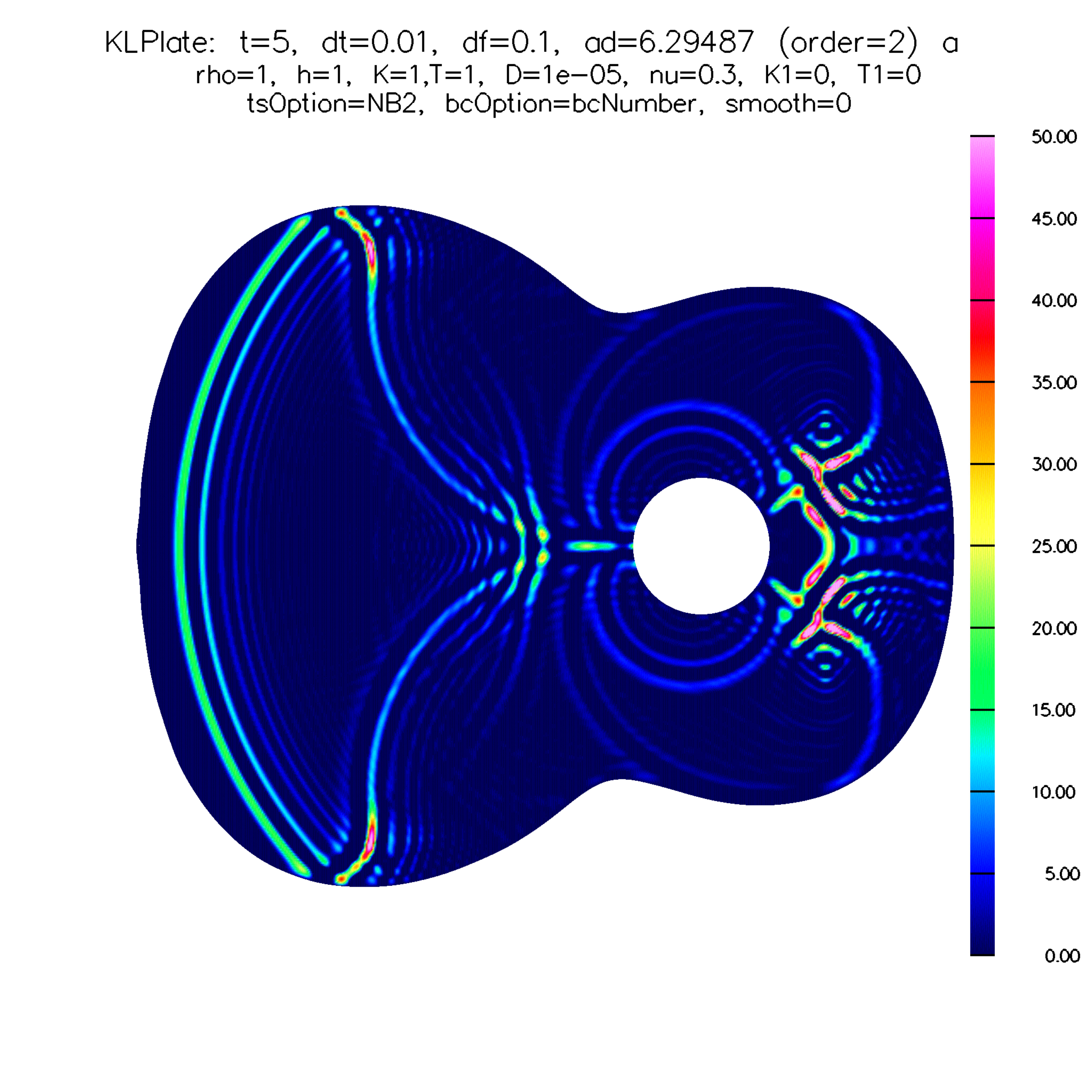}{\figWidth}};

\newcommand{\cbWidth}{.2}
\newcommand{\cbHeight}{3}

\drawColorBarH{fig/colourBarLines}{xshift=1cm,yshift=-0.2cm}{\cbWidth}{\cbHeight}{$0.0$}{$0.10$}{w}
\drawColorBarH{fig/colourBarLines}{xshift=6.6cm,yshift=-0.2cm}{\cbWidth}{\cbHeight}{$0.0$}{$1.50$}{v}
\drawColorBarH{fig/colourBarLines}{xshift=12.2cm,yshift=-0.2cm}{\cbWidth}{\cbHeight}{$0.0$}{$50.00$}{a}

\draw(2.5,10.)  node[anchor=south,xshift=0pt,yshift=0pt] {\footnotesize $w(x,y,3.5)$};
\draw(8.,10.)  node[anchor=south,xshift=0pt,yshift=0pt] {\footnotesize  $v(x,y,3.5)$};
\draw(13.5,10.)  node[anchor=south,xshift=0pt,yshift=0pt] {\footnotesize  $a(x,y,3.5)$};

\draw(2.5,4.5)  node[anchor=south,xshift=0pt,yshift=0pt] {\footnotesize $w(x,y,5.0)$};
\draw(8.,4.5)  node[anchor=south,xshift=0pt,yshift=0pt] {\footnotesize  $v(x,y,5.0)$};
\draw(13.5,4.5)  node[anchor=south,xshift=0pt,yshift=0pt] {\footnotesize  $a(x,y,5.0)$};

%
\end{tikzpicture}

\end{center}
\caption{Simulation results for the wave dominant case ($\rho=1, H=1,K=1,T=1,D=1\times10^{-5},\nu=0.3$).} \label{fig:guitarWaveDominant}

\end{figure}
}

{
\newcommand{\figWidth}{5cm}
\def\xa{16.}
\def\ya{10.5}
\newcommand{\trimfig}[2]{\trimw{#1}{#2}{0.12}{0.12}{0.18}{0.18}}
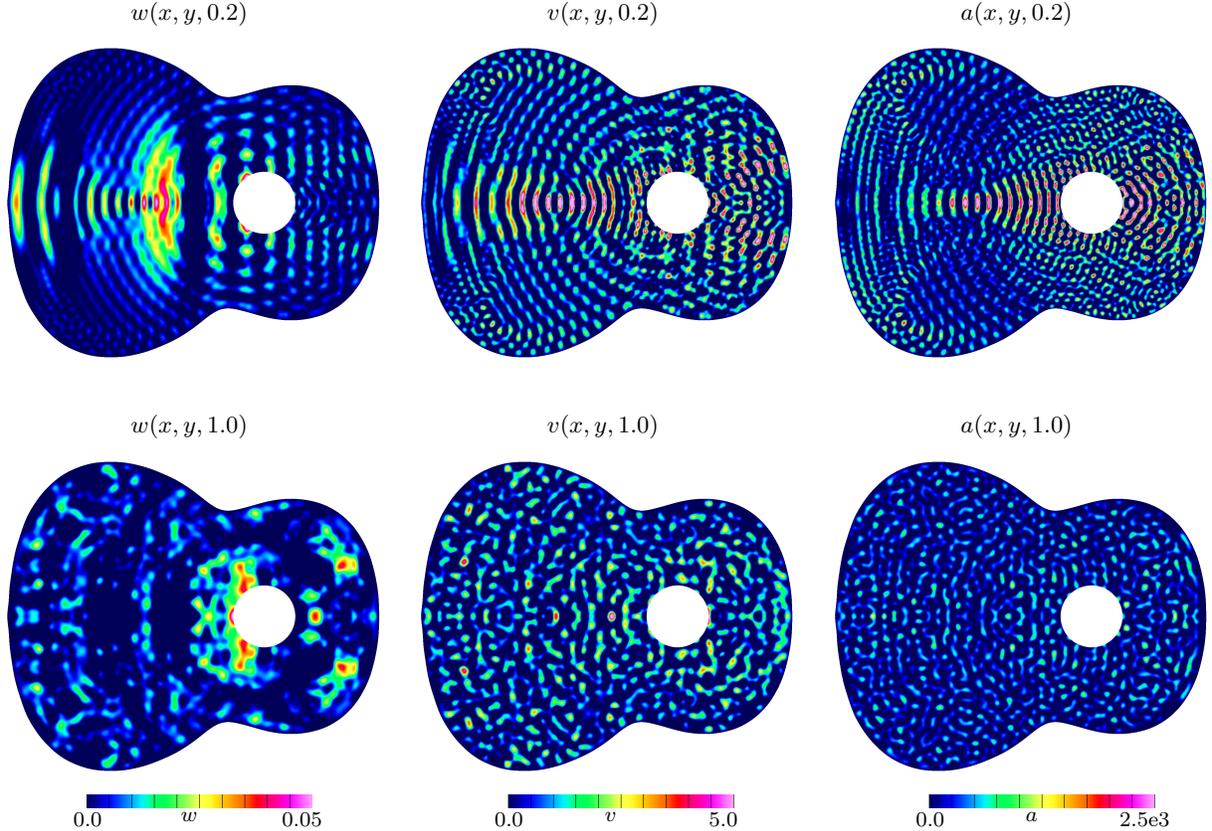
\begin{figure}[h]
\begin{center}
\begin{tikzpicture}[scale=1]
  \useasboundingbox (0.0,0.0) rectangle (\xa,\ya);  

\draw(-0.5,0.0) node[anchor=south west,xshift=0pt,yshift=0pt] {\trimfig{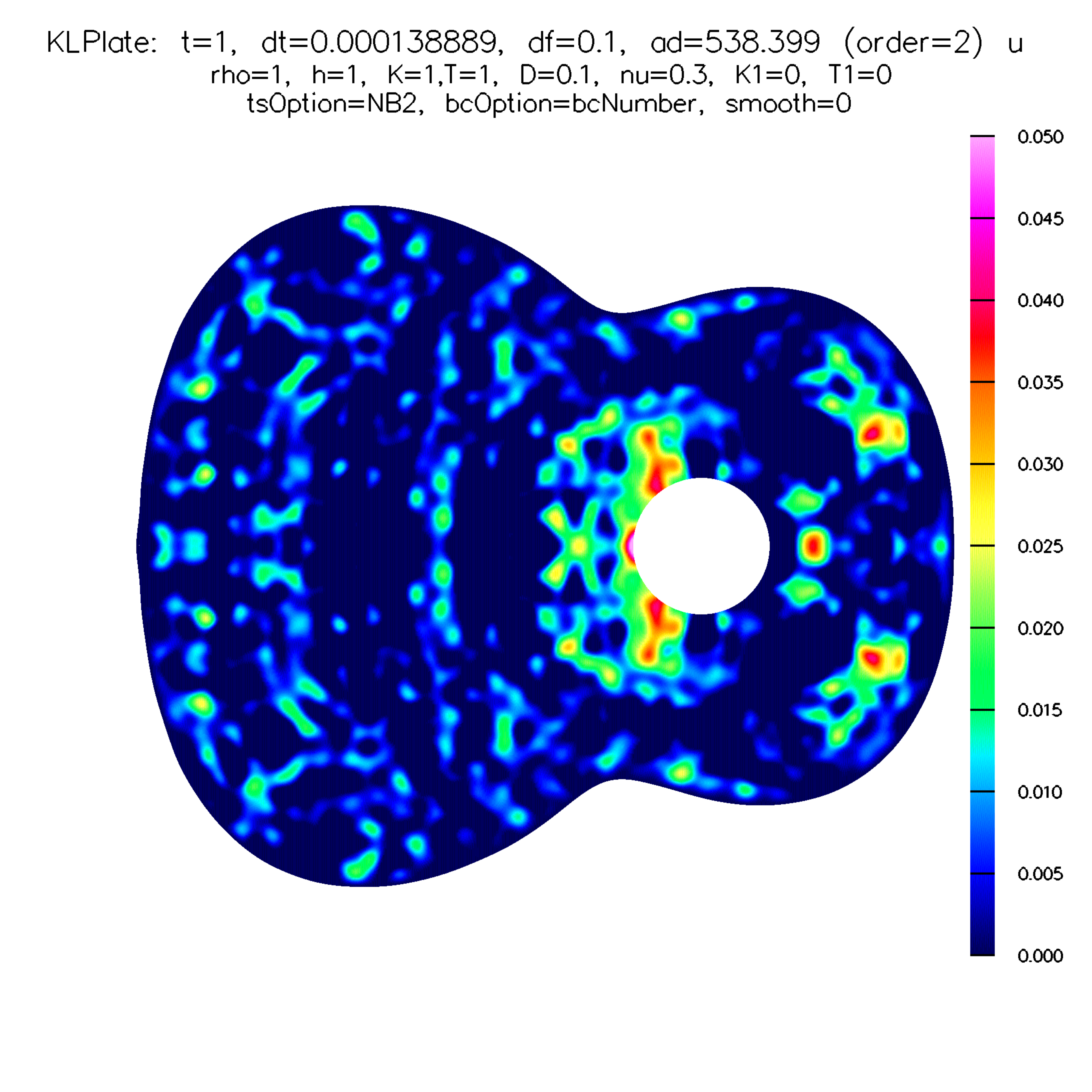}{\figWidth}};
\draw(5.,0.0) node[anchor=south west,xshift=0pt,yshift=0pt] {\trimfig{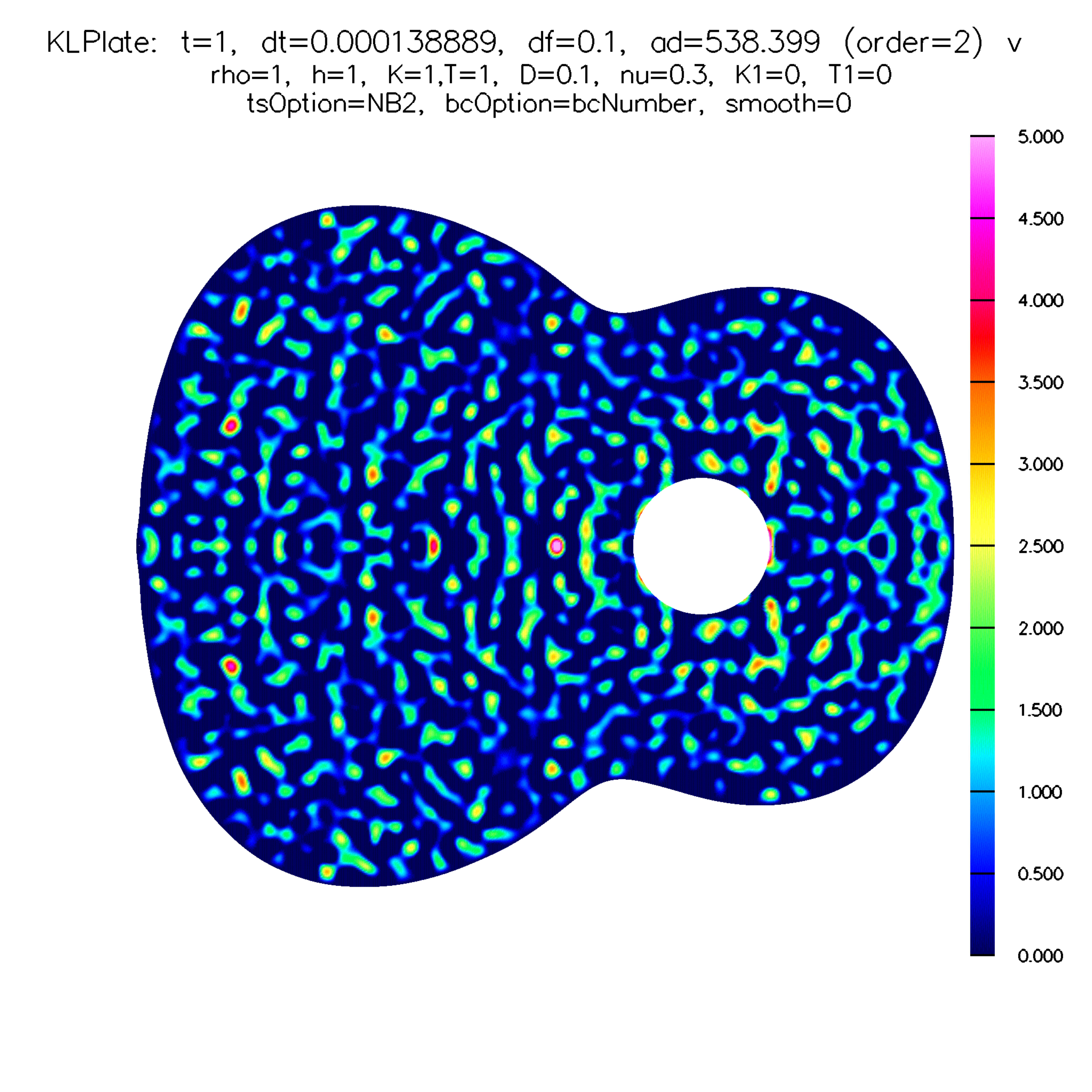}{\figWidth}};
\draw(10.5,0.0) node[anchor=south west,xshift=0pt,yshift=0pt] {\trimfig{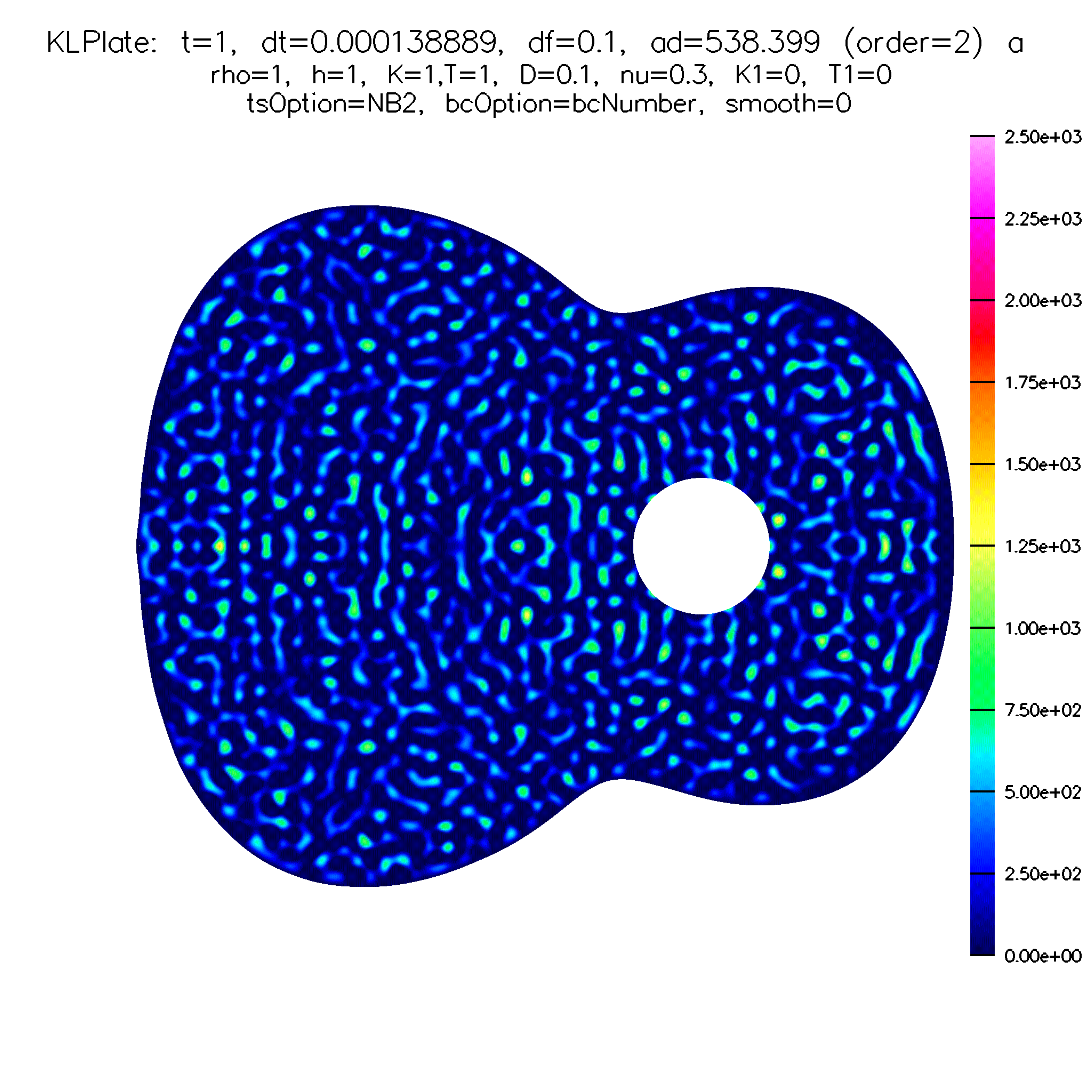}{\figWidth}};

\draw(-0.5,5.5) node[anchor=south west,xshift=0pt,yshift=0pt] {\trimfig{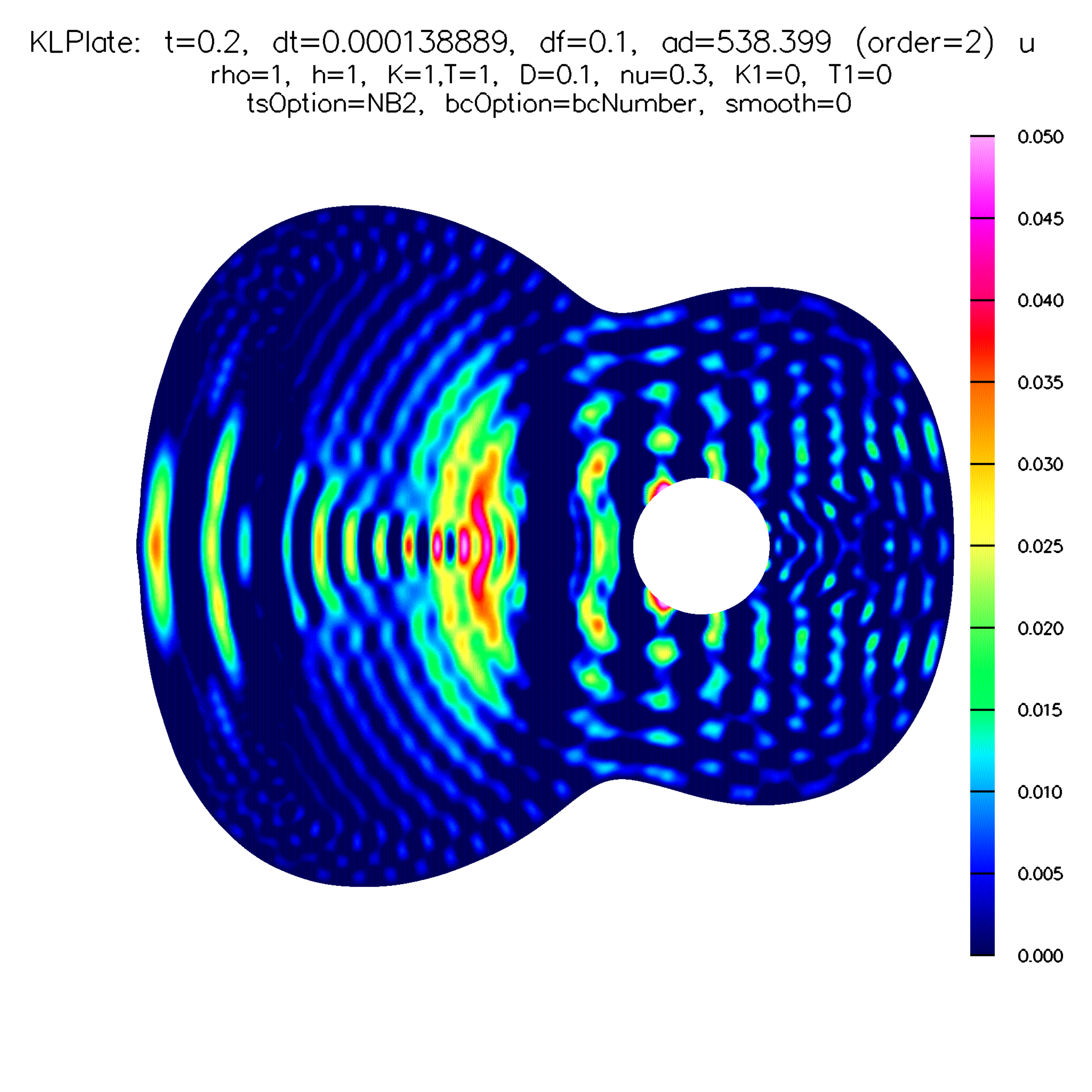}{\figWidth}};
\draw(5.,5.5) node[anchor=south west,xshift=0pt,yshift=0pt] {\trimfig{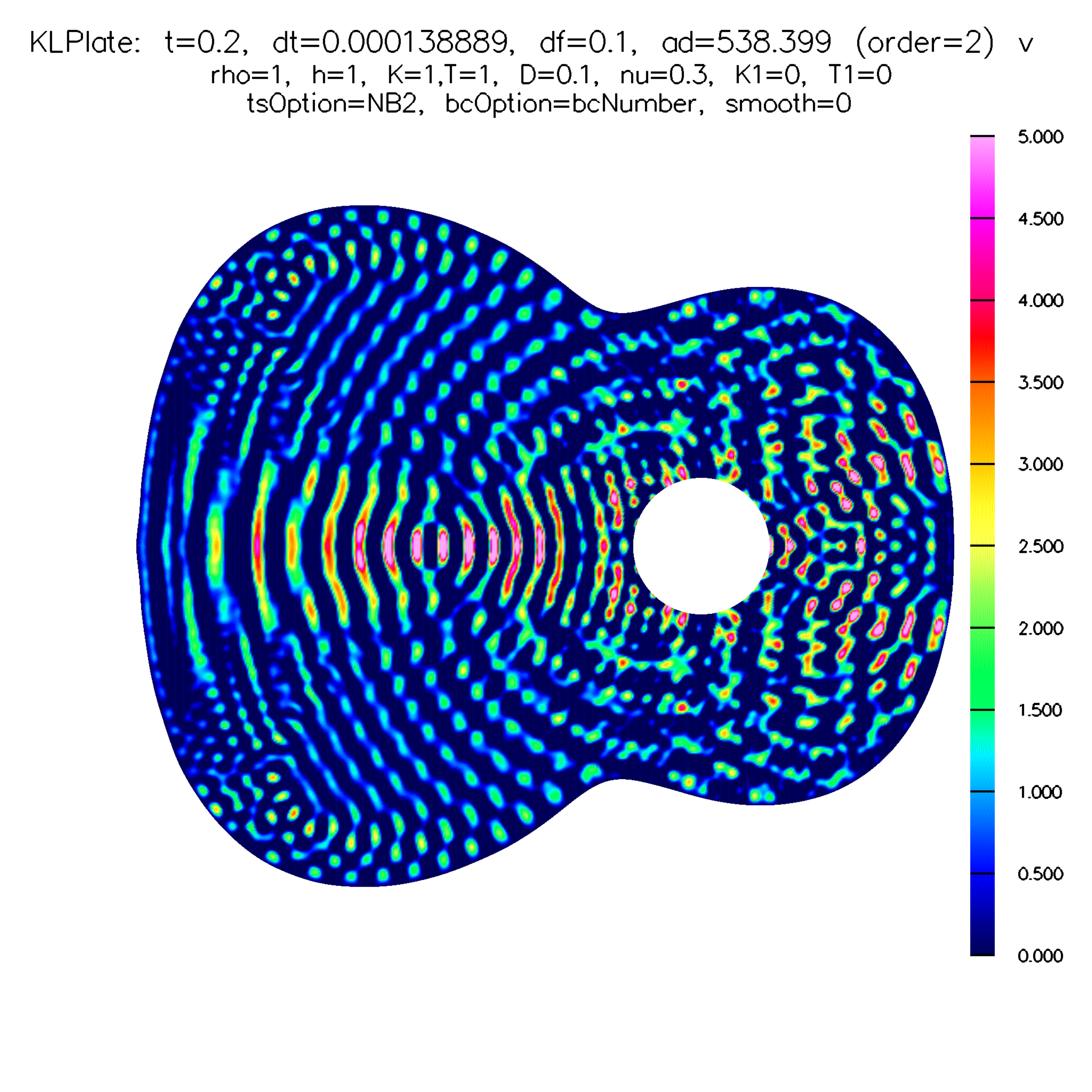}{\figWidth}};
\draw(10.5,5.5) node[anchor=south west,xshift=0pt,yshift=0pt] {\trimfig{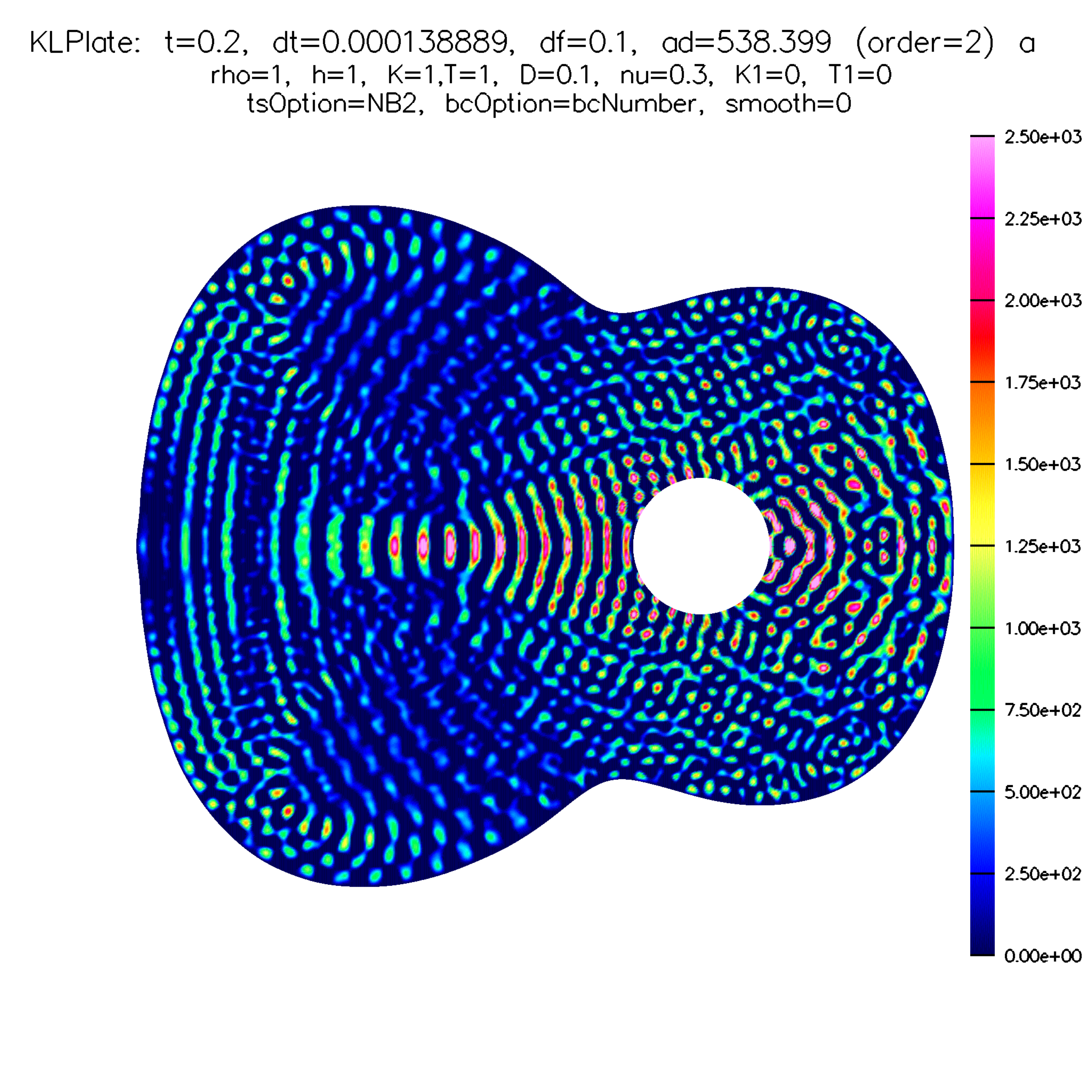}{\figWidth}};

\newcommand{\cbWidth}{.2}
\newcommand{\cbHeight}{3}

\drawColorBarH{fig/colourBarLines}{xshift=1cm,yshift=-0.2cm}{\cbWidth}{\cbHeight}{$0.0$}{$0.05$}{w}
\drawColorBarH{fig/colourBarLines}{xshift=6.6cm,yshift=-0.2cm}{\cbWidth}{\cbHeight}{$0.0$}{$5.0$}{v}
\drawColorBarH{fig/colourBarLines}{xshift=12.2cm,yshift=-0.2cm}{\cbWidth}{\cbHeight}{$0.0$}{$2.5e3$}{a}

\draw(2.5,10.)  node[anchor=south,xshift=0pt,yshift=0pt] {\footnotesize $w(x,y,0.2)$};
\draw(8.,10.)  node[anchor=south,xshift=0pt,yshift=0pt] {\footnotesize  $v(x,y,0.2)$};
\draw(13.5,10.)  node[anchor=south,xshift=0pt,yshift=0pt] {\footnotesize  $a(x,y,0.2)$};

\draw(2.5,4.5)  node[anchor=south,xshift=0pt,yshift=0pt] {\footnotesize $w(x,y,1.0)$};
\draw(8.,4.5)  node[anchor=south,xshift=0pt,yshift=0pt] {\footnotesize  $v(x,y,1.0)$};
\draw(13.5,4.5)  node[anchor=south,xshift=0pt,yshift=0pt] {\footnotesize  $a(x,y,1.0)$};

%
\end{tikzpicture}

\end{center}
\caption{Simulation results for the bend dominant case ($\rho=1, H=1,K=1,T=1,D=0.1,\nu=0.3$).} \label{fig:guitarBendDominant}
\end{figure}
}

The simulation results of the bend dominant case are shown in Figure~\ref{fig:guitarBendDominant}, where the displacement, velocity and acceleration of the plate are plotted at $t=0.2$ and $t=1$, respectively.  For this case, we observe that the initial pulse induced  vibrations  radiate rapidly  through the entire domain that  exhibits  a  much richer dynamics   than the  wave dominant case. The patterns for the velocity and acceleration are also observed to quickly travel towards the outer edge of the guitar and display a very complex dynamics by the time $t=1$. 

It is worth noting that, even though the above simulations are done for hypothetical physical  parameters, interesting guitar designers or artists could exploit our approach and code  to  improve their  guitar designs. For example, they may modify the shape of the guitar and specify physical parameters for specific materials in the above simulations  to investigate the acoustic quality of their intended designs.







{
\newcommand{\figWidth}{7cm}
\def\xa{15.}
\def\ya{5.}
\newcommand{\trimfig}[2]{\trimw{#1}{#2}{0}{0}{0}{0}}
\begin{figure}[h]
\begin{center}
\begin{tikzpicture}[scale=1]
  \useasboundingbox (0.0,0.0) rectangle (\xa,\ya);  

\draw(-0.5,0.0) node[anchor=south west,xshift=0pt,yshift=0pt] {\trimfig{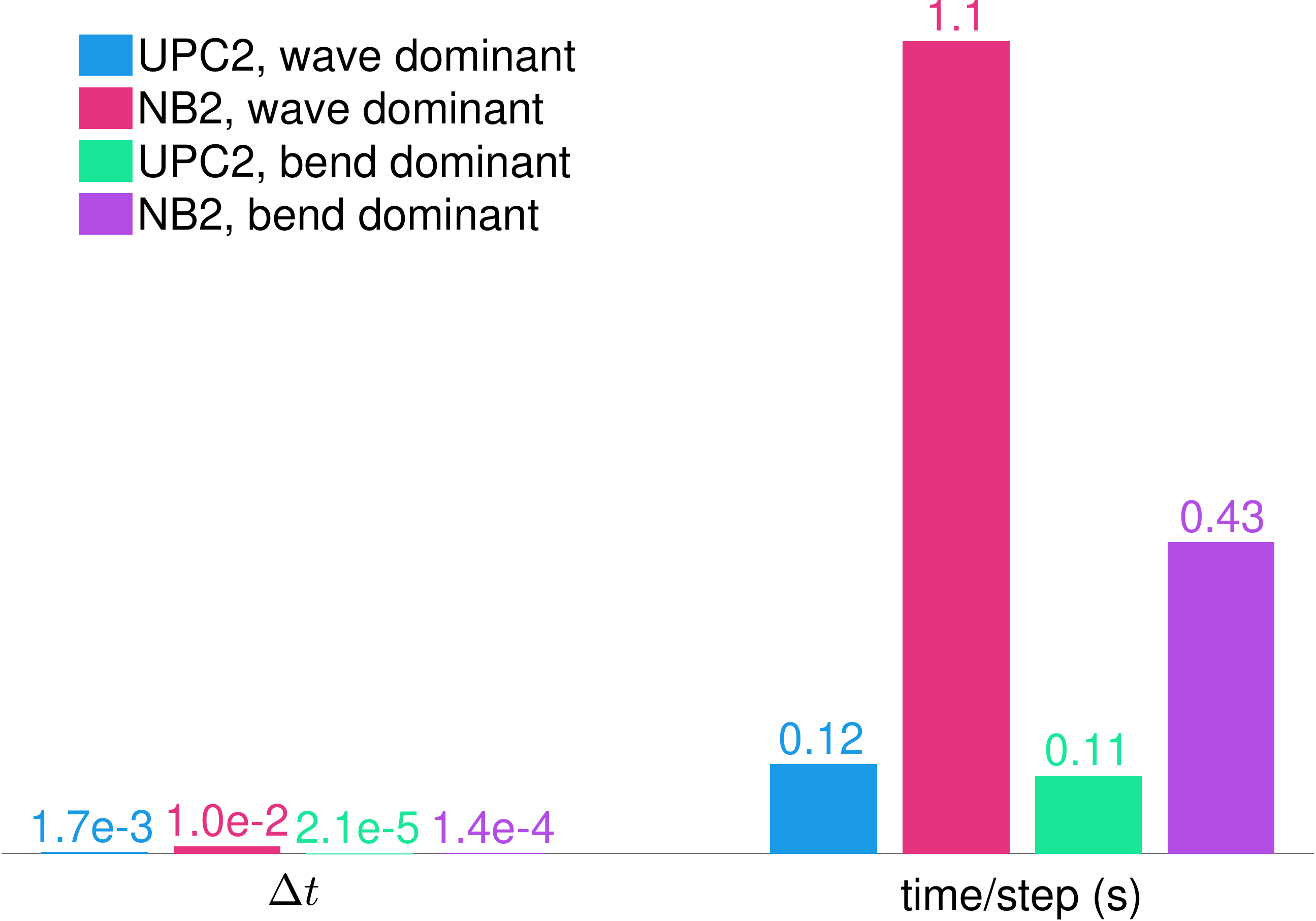}{\figWidth}};
\draw(7.5,0.0) node[anchor=south west,xshift=0pt,yshift=0pt] {\trimfig{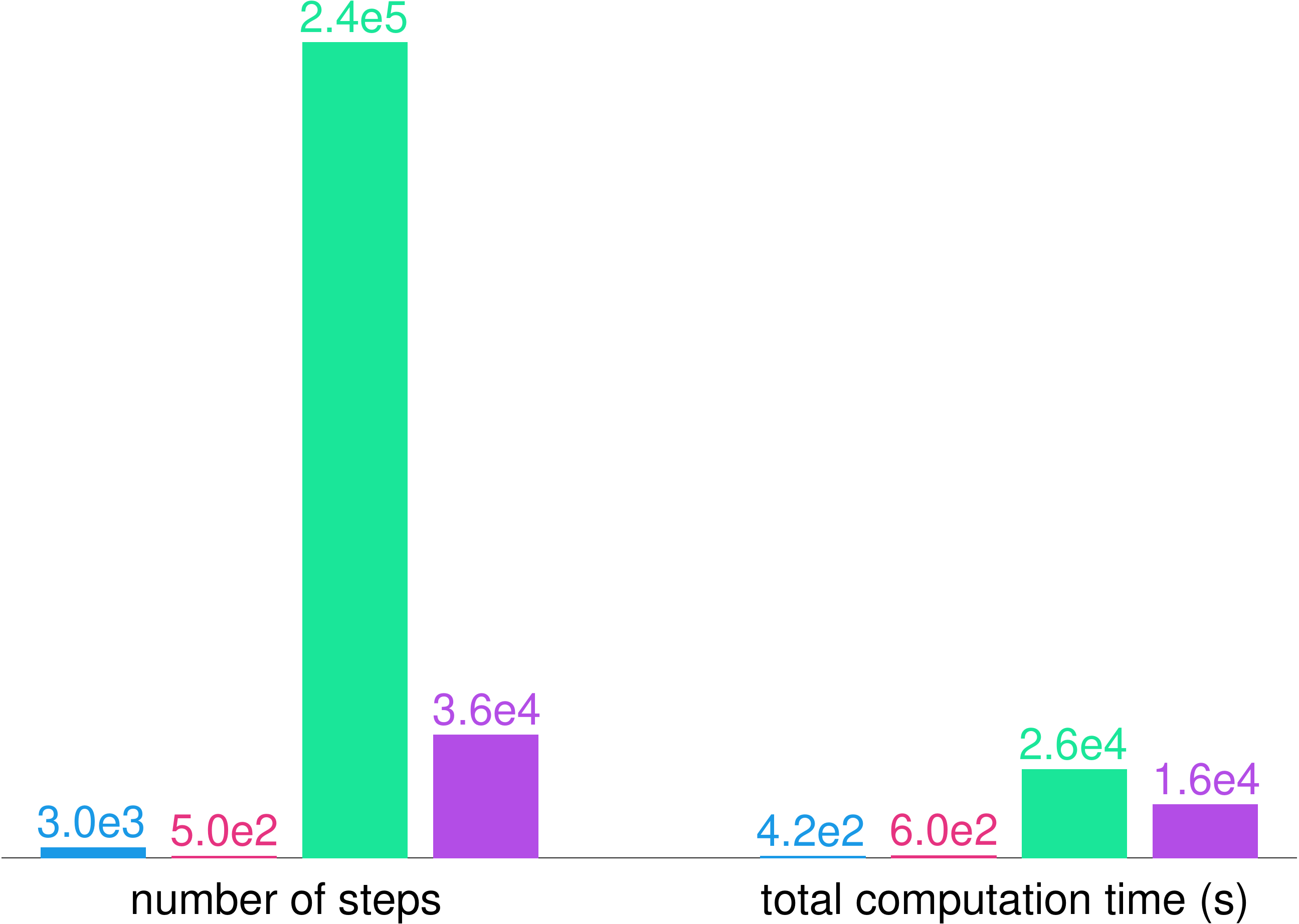}{\figWidth}};

\newcommand{\cbWidth}{.2}
\newcommand{\cbHeight}{3}


%
\end{tikzpicture}

\end{center}
\vspace{-0.2in}
\caption{performance comparison on $\Gc_{16}$} \label{fig:guitarPerformanceG16}
\end{figure}
}

We also make  a rough comparison between the run-time CPU costs for
the simulations shown in Figure~\ref{fig:guitarWaveDominant} and \ref{fig:guitarBendDominant}. The problems are solved until $t=5$
using a single processor of a MacBook Pro equipped with 16GB memory. The performance results are summarized in Figure~\ref{fig:guitarPerformanceG16}, in which we compare the step size $\dt$,  the CPU time in seconds per step (time/step), the total number of steps and the total computational times. We see  that  $\dt$ of the implicit NB2 scheme is generally 10 times bigger than that of the explicit UPC2 scheme, so  it takes 10 times more time steps for the UPC2 scheme to reach the final simulation time $t=5$. For each time step, the UPC2 scheme is more efficient  and takes much less time than the NB2 scheme  since the former scheme does not require  solving a linear system. In terms of the overall simulation, the total computation times for both schemes are more or less on the same order.  However, for  harder (stiffer) problems with larger values of  $D$, the time step of the explicit schemes can be significantly smaller. Under this circumstance, the implicit NB2 scheme can be more efficient for the  overall simulation. The performance comparison presented here can serve as a guidance  for us to choose the appropriate algorithm for a particular problem. 
 
\section{Conclusions}\label{sec:conclusions}
One novelty of this work lies in the manifestation  that finite difference methods are well-suited for solving  plate  models with  realistic and complex  geometries. In fact, our finite-difference-based approaches can be more advantageous  and efficient  for  solving  plate problems because it  solves the strong formulation of the governing equation \eqref{eq:KLShell}  directly using structured overlapping grids  without the need of any reformulation or extra complexity of non-standard finite elements.  
In this work,   four  novel numerical algorithms, referred to as the C2, UPC2, PC22 and  NB2 schemes,  are developed for the efficient and accurate solution of   the \KL plate model in complex domains. The proposed schemes are based on the standard spatial discretization on  composite overlapping grid which involves curvilinear  finite-difference   approximations for spatial  derivatives on each component grid and interpolating formulas for coupling together  solutions on the different component grids. To the best of our knowledge, the methods presented here are the first finite-difference based  algorithms for solving \KL model  with complex geometries that are suitable for realistic applications.

Our approaches resolve the numerical challenge due to 
the weak instability  excited by the presence of interpolating equations in the  discretization formulas on  overlapping grids.
The time-stepping  schemes are stabilized  by including a  novel artificial hyper-dissipation term to  the  governing equations.
Analysis on the frequency domain of the problem is performed to illustrate the stability of our methods. The analysis also   leads to the derivation of  explicit formulas (another  novelty  of the current work) for determining  stable time steps and sufficient artificial dissipations  that are applicable in real computations.
Moreover, unlike many existing algorithms which solve the plate equation for the displacement only,  we have seamlessly  incorporated the solution of velocity and acceleration into our algorithms. This  contribution can be useful for multi-physics coupling such as FSI, since the plate velocity and acceleration are important information  that needs to be accurately transferred  across the multi-physical interfaces. 

Quadratic eigenvalue problems for a simplified model plate on 1D overlapping grids are derived  for all the   proposed algorithms to  reveal the weak instability associated with the interpolation among component grids.
Using this model problem, we also investigate   how the artificial hyper-dissipation stabilizes the unstable modes for each of the algorithms.  
Carefully designed test problems are solved to demonstrate the properties and applications of our numerical approaches. In particular, the stability and accuracy of the schemes are verified by mesh refinement studies using method of manufactured solutions and using the analytical solutions of a circular plate. For the demonstration of applying our approach for  plates with complex geometries arising  from realistic applications,  two benchmark problems are considered. In the first problem, we solve a plate with numerous holes subject to  an initial disturbance of the displacement. In the second problem, we look at the traveling pulse on a guitar soundboard with two sets of  hypothetical physical parameters.   The future direction of this research  is to deploy these novel algorithms for   FSI  applications involving plates with complex geometries  in conjunction with our previously developed   FSI scheme  \cite{LiHenshaw2016}.

\section{Acknowledgment}\label{sec:Acknowledgment}
 L. Li is grateful to Professor W.D. Henshaw of Rensselaer Polytechnic Institute (RPI) for helpful conversations.

\appendix
\clearpage
\section{}

\subsection{Transformed PDE on reference domain}\label{sec:KLShellReferenceFormula}
The specific formula for Kirchoff-Love shell model \eqref{eq:KLShell} transformed into the reference domain is given by 
$$
{\rho}{H}\pdn{W}{t}{2}=\Lc W + F, 
$$
where   
\renewcommand{\cc}[1]{{\color{red} #1}}
\def\wo{W}\def\rx{r_{x}}
\def\sx{s_{x}}
\def\wr{\cc{W_{r}}}
\def\ry{r_{y}}
\def\sy{s_{y}}
\def\ws{\cc{W_{s}}}
\def\rxx{r_{xx}}
\def\sxx{s_{xx}}
\def\wrr{\cc{W_{rr}}}
\def\rxy{r_{xy}}
\def\sxy{s_{xy}}
\def\wrs{\cc{W_{rs}}}
\def\ryy{r_{yy}}
\def\syy{s_{yy}}
\def\wss{\cc{W_{ss}}}
\def\rxxx{r_{xxx}}
\def\sxxx{s_{xxx}}
\def\wrrr{\cc{W_{rrr}}}
\def\rxxy{r_{xxy}}
\def\sxxy{s_{xxy}}
\def\wrrs{\cc{W_{rrs}}}
\def\rxyy{r_{xyy}}
\def\sxyy{s_{xyy}}
\def\wrss{\cc{W_{rss}}}
\def\ryyy{r_{yyy}}
\def\syyy{s_{yyy}}
\def\wsss{\cc{W_{sss}}}
\def\rxxxx{r_{xxxx}}
\def\sxxxx{s_{xxxx}}
\def\wrrrr{\cc{W_{rrrr}}}
\def\rxxxy{r_{xxxy}}
\def\sxxxy{s_{xxxy}}
\def\wrrrs{\cc{W_{rrrs}}}
\def\rxxyy{r_{xxyy}}
\def\sxxyy{s_{xxyy}}
\def\wrrss{\cc{W_{rrss}}}
\def\rxyyy{r_{xyyy}}
\def\sxyyy{s_{xyyy}}
\def\wrsss{\cc{W_{rsss}}}
\def\ryyyy{r_{yyyy}}
\def\syyyy{s_{yyyy}}
\def\wssss{\cc{W_{ssss}}}
{\footnotesize
\begin{align*}
\mathcal{L} W &=\left(-K\right)\,\wo\\
&+\left(T\,\left(\rxx+\ryy\right)-D\,\left(\rxxxx+2\,\rxxyy+\ryyyy\right)\right)\,\wr\\
&+\left(T\,\left(\sxx+\syy\right)-D\,\left(\sxxxx+2\,\sxxyy+\syyyy\right)\right)\,\ws\\
&+\left(T\,\left(\rx^2+\ry^2\right)-D\,\left(3\,\rxx^2+2\,\rxx\,\ryy+4\,\rxy^2+3\,\ryy^2+4\,\rx\,\rxxx+4\,\rx\,\rxyy+4\,\rxxy\,\ry+4\,\ry\,\ryyy\right)\right)\,\wrr\\
&+(T\,\left(2\,\rx\,\sx+2\,\ry\,\sy\right)-D\,(4\,\rx\,\sxxx+4\,\rxxx\,\sx+4\,\rx\,\sxyy+4\,\rxyy\,\sx+6\,\rxx\,\sxx\\
&+8\,\rxy\,\sxy+4\,\rxxy\,\sy+4\,\ry\,\sxxy+2\,\rxx\,\syy+2\,\ryy\,\sxx+4\,\ry\,\syyy+4\,\ryyy\,\sy+6\,\ryy\,\syy))\,\wrs \\
&+\left(T\,\left(\sx^2+\sy^2\right)-D\,\left(3\,\sxx^2+2\,\sxx\,\syy+4\,\sxy^2+3\,\syy^2+4\,\sx\,\sxxx+4\,\sx\,\sxyy+4\,\sxxy\,\sy+4\,\sy\,\syyy\right)\right)\,\wss\\
&+\left(-D\,\left(6\,\rx^2\,\rxx+2\,\rxx\,\ry^2+6\,\ry^2\,\ryy+2\,\rx\,\left(\rx\,\ryy+2\,\rxy\,\ry\right)+4\,\rx\,\rxy\,\ry\right)\right)\,\wrrr\\
&+(-D\,(2\,\rx\,\left(\rx\,\syy+2\,\rxy\,\sy+2\,\ry\,\sxy+\ryy\,\sx\right)+3\,\rx^2\,\sxx+2\,\ry^2\,\sxx+3\,\ry^2\,\syy\\
&+2\,\sx\,\left(\rx\,\ryy+2\,\rxy\,\ry\right)+\rx\,\left(3\,\rx\,\sxx+3\,\rxx\,\sx\right)+4\,\rxy\,\left(\rx\,\sy+\ry\,\sx\right)+\ry\,\left(3\,\ry\,\syy+3\,\ryy\,\sy\right)\\
&+9\,\rx\,\rxx\,\sx+4\,\rx\,\ry\,\sxy+4\,\rxx\,\ry\,\sy+9\,\ry\,\ryy\,\sy))\,\wrrs\\
&+(-D\,(2\,\sx\,\left(\rx\,\syy+2\,\rxy\,\sy+2\,\ry\,\sxy+\ryy\,\sx\right)+3\,\rxx\,\sx^2+2\,\rxx\,\sy^2+3\,\ryy\,\sy^2\\
&+\sx\,\left(3\,\rx\,\sxx+3\,\rxx\,\sx\right)+4\,\sxy\,\left(\rx\,\sy+\ry\,\sx\right)+\sy\,\left(3\,\ry\,\syy+3\,\ryy\,\sy\right)+2\,\rx\,\left(\sx\,\syy+2\,\sxy\,\sy\right)\\
&+9\,\rx\,\sx\,\sxx+4\,\rxy\,\sx\,\sy+4\,\ry\,\sxx\,\sy+9\,\ry\,\sy\,\syy))\,\wrss\\
&+\left(-D\,\left(6\,\sx^2\,\sxx+2\,\sxx\,\sy^2+6\,\sy^2\,\syy+2\,\sx\,\left(\sx\,\syy+2\,\sxy\,\sy\right)+4\,\sx\,\sxy\,\sy\right)\right)\,\wsss\\
&+\left(-D\,\left(\rx^4+2\,\rx^2\,\ry^2+\ry^4\right)\right)\,\wrrrr\\
&+\left(-D\,\left(\rx\,\left(\sx\,\ry^2+2\,\rx\,\sy\,\ry\right)\,2+4\,\rx^3\,\sx+4\,\ry^3\,\sy+2\,\rx\,\ry^2\,\sx\right)\right)\,\wrrrs\\
&+\left(-D\,\left(6\,\rx^2\,\sx^2+6\,\ry^2\,\sy^2+\rx\,\left(\rx\,\sy^2+2\,\ry\,\sx\,\sy\right)\,2+\sx\,\left(\sx\,\ry^2+2\,\rx\,\sy\,\ry\right)\,2\right)\right)\,\wrrss\\
&+\left(-D\,\left(\sx\,\left(\rx\,\sy^2+2\,\ry\,\sx\,\sy\right)\,2+4\,\rx\,\sx^3+4\,\ry\,\sy^3+2\,\rx\,\sx\,\sy^2\right)\right)\,\wrsss\\
&+\left(-D\,\left(\sx^4+2\,\sx^2\,\sy^2+\sy^4\right)\right)\,\wssss. 
\end{align*}
}
The definitions of  coefficients $a_i(\rv), b_i(\rv),c_i(\rv),d_i(\rv)$ (coefficients in front of the red terms)  can be readily read off the above formula for $\Lc W$.

\subsection{Formula of the Discrete transformation}\label{sec:formulaDFT}
The Fourier transform (symbol) of the discrete operator $\Lc_h$ is
  {\footnotesize
\begin{align*}
\hat{Q}(\xi_r,\xi_s;\rv)=&-K-b_{11}\frac{4\sin^2(\xi_r/2)}{h_r^2}-b_{12} \frac{\sin(\xi_r)\sin(\xi_s)}{h_r h_s}-b_{22}\frac{4\sin^2(\xi_s/2)}{h_s^2}+d_{1111}\frac{16\sin^4(\xi_r/2)}{h_r^4}\\
&+d_{1112}\frac{4\sin^2(\xi_r/2)\sin(\xi_r)\sin(\xi_s)}{h_r^3h_s}+d_{1122}\frac{16\sin^2(\xi_r/2)\sin^2(\xi_s/2)}{h_r^2h_s^2}+d_{1222}\frac{4\sin(\xi_r)\sin^2(\xi_s/2)\sin(\xi_s)}{h_rh_s^3}\\
&+d_{2222}\frac{16\sin^4(\xi_s/2)}{h_s^4}
\end{align*}
  }

\clearpage
\bibliographystyle{elsart-num}
\bibliography{journal-ISI,henshaw,henshawPapers,fsi,LongfeiRef,LongfeiPapers,plate}

\end{document}